%% file: Hurst_fSDE.tex
\documentclass[a4paper,11pt]{article}
\usepackage{amsmath, amssymb}
\usepackage{amsthm}% for \proof,\newtheorem
\usepackage{latexsym}

\usepackage{subfiles}
\usepackage{comment}
\usepackage{setspace} % for gyoukan chousei
\usepackage{authblk} %\author関連

\usepackage{enumerate} %\begin{enumerate}[{(}i{)}]
\usepackage{cases}
  
\usepackage{ascmac} % for \screen,\itembox,\shadebox,highlight,\yen
\usepackage{framed}
\usepackage{fancybox}

\usepackage{cancel}
\usepackage{url}

\usepackage[left=19.91truemm,right=19.91truemm,top=26.4truemm,bottom=42truemm]{geometry}

\usepackage[colorlinks=true]{hyperref}
\hypersetup{
    colorlinks=true,
    citecolor=blue,
    linkcolor=blue,
}
\usepackage{cleveref}

\usepackage[backend=bibtex,style=alphabetic,isbn=false,eprint=false,doi=false]{biblatex}
\usepackage[normalem]{ulem}

\usepackage{xcolor}
\usepackage{layout} % for layout 
\usepackage{tikz} 
\usepackage{graphicx}
\usepackage[labelformat=parens, font=footnotesize]{caption}

\usepackage{subcaption}

\usepackage{here}

\usepackage{dsfont}
\usepackage{bbm} %\mathbbm{ABCabc\nabla}
\usepackage[bbgreekl]{mathbbol} %\mathbb{abcABC\nabla\bbalpha\bbbeta\Delta  }
\usepackage{upgreek}
\usepackage{mathrsfs}

\newtheorem{theorem}{Theorem}[section]
\newtheorem{ass}[theorem]{Assumption}

\newtheorem{lemma}[theorem]{Lemma}

\newtheorem{proposition}[theorem]{Proposition}

\newtheorem{remark}[theorem]{Remark}

\numberwithin{equation}{section}
\newtheorem{theorem*}{Theorem}
\newtheorem{ass*}[theorem*]{Assumption}
\newtheorem{note*}[theorem*]{Note}
\newtheorem{lemma*}[theorem*]{Lemma}
\newtheorem{definition*}[theorem*]{Definition}
\newtheorem{proposition*}[theorem*]{Proposition}
\newtheorem{corollary*}[theorem*]{Corollary}
\newtheorem{remark*}[theorem*]{Remark}
\newtheorem{example*}[theorem*]{Example}
\numberwithin{equation}{section}

\input{subfiles/hayamacro-240702.tex}
\input{subfiles/nakamacro291103+++-240702.tex}

\input{subfiles/fRV_symbol-240702.tex}
\input{subfiles/fSDE_Hurst_symbol-240702.tex}

\DeclareMathOperator\supp{supp}
\DeclareMathOperator*{\argmax}{arg\,max}

\DeclareMathOperator*{\Dom}{Dom}

\excludecomment{en-text}
\includecomment{scrap}
\excludecomment{scrap}

\includecomment{todelete}
\excludecomment{todelete}

\allowdisplaybreaks
\title{
Asymptotic expansion of a Hurst index estimator for a stochastic differential equation driven by fBm%
\footnote{
This work was in part supported by 
Japan Science and Technology Agency CREST JPMJCR14D7, JPMJCR2115; 
Japan Society for the Promotion of Science Grants-in-Aid for Scientific Research 
No. 17H01702 (Scientific Research); 
and by a Cooperative Research Program of the Institute of Statistical Mathematics. %%
}
}

\author{Hayate Yamagishi}
\affil{Graduate School of Mathematical Sciences, University of Tokyo
\footnote{Graduate School of Mathematical Sciences, 
University of Tokyo: 3-8-1 Komaba, Meguro-ku, Tokyo 153-8914, Japan. 
e-mail: 
yhayate@ms.u-tokyo.ac.jp
}}
\affil{CREST, Japan Science and Technology Agency%\footnote{}
}
\begin{document}
\maketitle
\begin{abstract}
We study the asymptotic properties of an estimator of Hurst parameter of 
a stochastic differential equation driven by a fractional Brownian motion with $H>1/2$.
Utilizing the theory of asymptotic expansion of Skorohod integrals introduced 
by Nualart and Yoshida \cite{nualart2019asymptotic},
we derive an asymptotic expansion formula of the distribution of the estimator.
As an corollary, we also obtain a mixed central limit theorem for the statistic,
indicating that the rate of convergence is $n^{-\half}$,
which improves the results in the previous literature.
% $n^{-\frac14}$.
To handle second-order quadratic variations appearing in the estimator,
a theory of exponent has been developed based on weighted graphs 
to estimate asymptotic orders of norms of functionals involved.
\end{abstract}

\keywords{%
  fractional Brownian motion, 
  stochastic differential equation,
  estimator of Hurst parameter,
  % ratio statistic,
  quadratic variation,
  second-order difference, 
  mixed central limit theorem,
  % mixed normal distribution,
  asymptotic expansion,
  % random symbol,
  Malliavin calculus,
  Skorohod integral,
  multiple stochastic integral,
  exponent, 
  weighted graph,
  asymptotic order,
  product formula
}

\tableofcontents

%%%%%%%%%%%%%%%%%%%%%%%%%%%%%%%%%%%%%%%%%
%%%%%%%%%%%%%%%%%%%%%%%%%%%%%%%%%%%%%%%%%

\subfile{subfiles/1-sec_introduction}
%
\subfile{subfiles/2-sec_main_result}
\subfile{subfiles/3-sec_preliminaries}

%
\subfile{subfiles/4-1-sec_exponent_def}

\subfile{subfiles/4-2-sec_exponent_main}
%
\subfile{subfiles/5-sec_asy_exp}
\subfile{subfiles/6-sec_technical_lemmas}
\subfile{subfiles/7-sec_numerical_analysis}

\printbibliography

\end{document}

%% file: subfiles/hayamacro-240702.tex
\usepackage{comment}

\makeatletter
\newsavebox{\@brx}
\newcommand{\llangle}[1][]{\savebox{\@brx}{\(\m@th{#1\langle}\)}%
  \mathopen{\copy\@brx\kern-0.5\wd\@brx\usebox{\@brx}}}
\newcommand{\rrangle}[1][]{\savebox{\@brx}{\(\m@th{#1\rangle}\)}%
  \mathclose{\copy\@brx\kern-0.5\wd\@brx\usebox{\@brx}}}
\makeatother

\providecommand{\keywords}[1]
{
  \small	
  \textbf{Keywords:} #1
}

\newcommand\vspsm{\vspace{5mm}}
\newcommand\hspsm{\hspace{5mm}}
\newcommand\vspssm{\vspace{2.5mm}}

\renewcommand{\emptyset}{\varnothing}

\newcommand\sfz{{\sf z}}

\newcommand\sfi{{\tt i}}
\newcommand\tti{{\tt i}}

\newcommand{\stout}[1]{\ifmmode\text{\sout{\ensuremath{#1}}}\else\sout{#1}\fi}

\newcommand{\myred}{\color[rgb]{0.8,0,0}}
\newcommand{\myblue}{\color[rgb]{0.01, 0.28, 1.0}}

\newcommand{\delc}[1]{}
\newcommand{\delb}[1]{[]}

\newcommand\vpr{{v'}}

\newcommand\jon{{j_{1}}}
\newcommand\jtw{{j_{2}}}
\newcommand\jth{{j_3}}

\newcommand{\nnorm}[1]{{\| #1 \|}}
\newcommand{\bnorm}[1]{{\bigl\| #1 \bigr\|}}
\newcommand{\Bnorm}[1]{{\Bigl\| #1 \Bigr\|}}

\newcommand{\snorm}[1]{{\left\| #1 \right\|}}
\newcommand{\norm}[1]{{\left\| #1 \right\|}}

\newcommand{\abs}[1]{{\left| #1 \right|}}
\newcommand{\nabs}[1]{{| #1 |}}
\newcommand{\babs}[1]{{\bigl| #1 \bigr|}}

\newcommand{\bbabs}[1]{{\biggl| #1 \biggr| }}

\newcommand{\cbr}[1]{{\left\{ #1 \right\} }}
\newcommand{\ncbr}[1]{{\{ #1 \} }}
\newcommand{\bcbr}[1]{{\bigl\{ #1 \bigr\} }}
\newcommand{\Bcbr}[1]{{\Bigl\{ #1 \Bigr\} }}

\newcommand{\rbr}[1]{{\left( #1 \right) }}
\newcommand{\nrbr}[1]{{( #1 ) }}
\newcommand{\brbr}[1]{{\bigl( #1 \bigr) }}
\newcommand{\Brbr}[1]{{\Bigl( #1 \Bigr) }}
\newcommand{\bbrbr}[1]{{\biggl( #1 \biggr) }}

\newcommand{\sbr}[1]{{\left[ #1 \right] }}
\newcommand{\bsbr}[1]{{\bigl[ #1 \bigr] }}
\newcommand{\Bsbr}[1]{{\Bigl[ #1 \Bigr] }}
\newcommand{\bbsbr}[1]{{\biggl[ #1 \biggr] }}

\newcommand{\abr}[1]{{\left\langle #1 \right\rangle }}
\newcommand{\babr}[1]{{\bigl\langle #1 \bigr\rangle }}
\newcommand{\Babr}[1]{{\Bigl\langle #1 \Bigr\rangle }}
\newcommand{\bbabr}[1]{{\biggl\langle #1 \biggr\rangle }}

\newcommand{\subotimes}[1]{{\underset{ #1 }{\otimes}}}

\newcommand\qtor{{\tt qTor}}
\newcommand\mqtor{{\tt mqTor}}
\newcommand\qtan{{\tt qTan}}

\newcommand{\mS}{\mathfrak S}
\newcommand{\mT}{\mathfrak T}

\newcommand{\scrC}{{\mathscr C}}

\newcommand{\bbf}{{\mathbb f}}

\newcommand{\bbone}{{\mathbb 1}}

\newcommand{\tfor}{\quad\text{for}\ }
\newcommand{\tforsm}{\ \text{for}\ }

\newcommand{\tifsm}{\ \text{if}\ }

\newcommand{\tand}{\quad\text{and}\quad}

\newcommand{\tandsm}{\ \text{and}\ }
\newcommand{\tor}{\quad\text{or}\quad}
\newcommand{\torsm}{\ \text{or}\ }

%% file: subfiles/nakamacro291103+++-240702.tex
\def\mH{{\mathfrak H}}

\def\bd{\begin{description}}
\def\ed{\end{description}}

\def\D2{\bbD_{2,\infty-}}

\def\tj{{t_j}}
\def\tjm{{t_{j-1}}}
\def\tjp{{t_{j+1}}}
\def\C{{\bf C}}
\def\D{{\bf D}}

\def\V{{\bf V}}

\def\cala{{\cal A}}

\def\cale{{\cal E}}
\def\calf{{\cal F}}

\def\calh{{\cal H}}
\def\cali{{\cal I}}
\def\calj{{\cal J}}

\def\caln{{\cal N}}

\def\calr{{\cal R}}

\def\yeq{\>=\>}

\def\simleq{\ \raisebox{-.7ex}{$\stackrel{{\textstyle <}}{\sim}$}\ }

\def\half{\frac{1}{2}}

\def\nn{\nonumber}
\def\be{\begin{equation}}
\def\ee{\end{equation}}
\def\bea{\begin{eqnarray}}
\def\eea{\end{eqnarray}}
\def\beas{\begin{eqnarray*}}
\def\eeas{\end{eqnarray*}}
\def\bi{\begin{itemize}}
\def\ei{\end{itemize}}

\def\bd{\begin{description}}
\def\ed{\end{description}}
\newcommand{\bbD}{{\mathbb D}}

\newcommand{\bbJ}{{\mathbb J}}

\newcommand{\bbN}{{\mathbb N}}

\newcommand{\bbR}{{\mathbb R}}
\newcommand{\bbS}{{\mathbb S}}
\newcommand{\bbT}{{\mathbb T}}

\newcommand{\bbV}{{\mathbb V}}

\newcommand{\bbX}{{\mathbb X}}
\newcommand{\bbY}{{\mathbb Y}}
\newcommand{\bbZ}{{\mathbb Z}}

%% file: subfiles/fRV_symbol-240702.tex
\newcommand{\barq}{{{\bar  q}}}
\newcommand{\bartheta}{{ {\bar  \theta}}}
\newcommand{\Comp}{\scrC}

%% file: subfiles/fSDE_Hurst_symbol-240702.tex
\usepackage{stmaryrd}
\newcommand{\secDiff}[2]{\Delta^{(2)}_{#1,#2}}
\newcommand{\secVar}{\V^{(2)}}

\newcommand{\jp}{{j+1}}

\newcommand{\ntwo}{{\bbN_2}}

\newcommand{\resint}{\hat\Xi_{n,j}}
\newcommand{\resinc}[1]{\Xi_{n,#1}}
\newcommand{\resincj}{{\resinc{j}}}
\newcommand{\convDiff}{{R}}

\newcommand{\varconst}{c_\infty}

\newcommand{\constInProd}{c_H}

\newcommand{\negTerm}{\Gamma}
\newcommand{\psiKer}{F}
\newcommand{\decayDelX}{\calr^\Delta}

\newcommand{\calii}{\cali^{(i)}}

\newcommand{\bbd}{{\mathbb d}}

\newcommand{\bbh}{{\mathbb h}}

\newcommand{\kerinc}{\nrbr{{\bbone}_+}}
\newcommand{\kerdec}{\nrbr{{\bbone}_-}}
\newcommand{\wt}[1]{\widehat{#1}}
\newcommand{\wtv}{\wt{V}}
\newcommand{\wh}[1]{\widehat{#1}}
\newcommand{\whv}{\widehat{V}}
\newcommand{\extv}{\widetilde{V}}
\newcommand{\exti}{\tilde{i}}
\newcommand{\extj}{\tilde{j}}
\newcommand{\extbbf}{\tilde{\bbf}}
\newcommand{\tp}{\hat{p}}
\newcommand{\hp}{\hat{p}}
\newcommand{\tv}{\hat{v}}
\newcommand{\hv}{\hat{v}}
\newcommand{\tvp}{\hat{v}'}

\newcommand{\sftv}{{V}_{\star}}
\newcommand{\sftwhv}{{\whv}_{\star}}
\newcommand{\sftc}{{C}_{\star}}
\newcommand{\sfthv}{{\hv}_{\star}}
\newcommand{\dblv}{\widetilde{V}}
\newcommand{\dbli}{\tilde{\imath}}
\newcommand{\dblj}{\tilde{\jmath}}
\newcommand{\dblbbf}{\tilde{\bbf}}
\newcommand{\dblvwq}{\tilde{\vwq}}
\newcommand{\dblewt}{\tilde{\ewt}}

\newcommand{\diffker}{\bbd}
\newcommand{\dker}{\diffker}
\newcommand{\vwq}{{\mathbf{q}}}
\newcommand{\bfq}{{\mathbf{q}}}
\newcommand{\ewt}{{\mathbf{\uptheta}}}
\newcommand{\cewt}{\check{\ewt}}
\newcommand{\bff}{{\mathbf{f}}}
\newcommand{\bfp}{{\mathbf{p}}}

\newcommand{\kap}{\kappa}
\newcommand{\kapr}{{\kappa'}}

\newcommand{\iv}{{\bbV_\infty}}
\newcommand{\rv}{{\bbV}}

\newcommand{\tnj}[2]{{t^{#1}_{#2}}}

\usepackage{accents}

\newcommand{\convDNeg}{\convDiff^{({\sf N})}}
\newcommand{\caliN}[1]{\cali^{({\sf N})(#1)}}
\newcommand{\AN}{A^{({\sf N})}}

\newcommand{\eq}{e_{q}}
\newcommand{\et}{e_{\theta}}

\newcommand{\totpi}{\mathring{\bfp}}

\newcommand{\eltw}{\ell_2}

\newcommand{\lambdainv}{\sbr\lambda}
\newcommand{\DLamCali}{D^{(\lambda)}\cali_n}

\newcommand{\newComp}{{\tt C}}

\newcommand{\ttj}{{t^{2n}_{2j}}}
\newcommand{\ttjm}{{t^{2n}_{2j-1}}}
\newcommand{\ttjp}{{t^{2n}_{2j+1}}}

%% file: subfiles/1-sec_introduction.tex
%
\section{Introduction}
In this paper,
we consider the following one-dimensional SDE driven by a fractional Brownian motion (fBm)
with its Hurst parameter $H>1/2$:
\begin{align}\label{eq:230925.1452}
  dX_t &= V^{[2]}(X_t) dt + V^{[1]}(X_t) dB^H_t
  \\
  X_0&=x_0\in\bbR,
  \nn
\end{align}
where $V^{[i]}$ ($i=1,2$) are functions from $\bbR$ to $\bbR$,
$t\in[0,T]$ and $T>0$ is a fixed terminal time.
To make the presentation simple, we set $T=1$.
The stochastic integral in (\ref{eq:230925.1452}) is a pathwise Riemann-Stieltjes
integral (Young integral), 
and it is known that 
there exists a unique solution under some regularity conditions on $V^{[i]}$ ($i=1,2$), 
as Nualart and Rascanu \cite{rascanu2002differential} detailed.

For $n\geq2$ and $j=1,...,n-1$,
we denote the second order difference of $X_t$ by
\begin{align*}
  \secDiff{n}{j}X = (X_{t^n_{j+1}} - X_{t^n_{j}}) - (X_{t^n_{j}} - X_{t^n_{j-1}}), 
\end{align*}
where $t^n_{j}= j/n$ for $j=1,..,n$.
We define the quadratic variation $\secVar_n$ based on the second order difference by
\begin{align}\label{eq:231002.0958}
  \secVar_n&= \secVar_n(X) = 
  \sum_{j=1}^{n-1}\rbr{\secDiff{n}{j}X}^2,
  % =O(n^{-2H+1}),
\end{align}
and an estimator $\hat H_n$ of the Hurst parameter $H$ by
% Recall the estimator defined at \eqref{??}
\begin{align}\label{eq:231001.1046}
  \hat H_n &=
  0\vee
  \rbr{\half + \frac1{2\log2} \log\frac{\secVar_n}{\secVar_{2n}}}
  \wedge1.
\end{align}
We consider the asymptotics of $\hat H_n$ as $n$ tends to $\infty$, 
that is, the length of the interval of the difference goes to zero.
This estimator is a slight modification of the one proposed in \cite{istas1997quadratic}
for the Hurst index of a fractional Brownian motion.
For a SDE driven by fBm,
the above estimator was previously studied (with a slightly different truncation) by 
Kubilius and his coauthors
% This estimator (with a slightly different truncation) 
in \cite{kubilius2010quadratic}, \cite{kubilius2012rate} and \cite{kubilius2017rate}.
% Kubilius and Mishura in \cite{kubilius2012rate}.
They established the strong consistency and the rate of convergence 
by analyzing the asymptotic behavior with probability $1$ 
of the quadratic variations of the solution to SDE \eqref{eq:230925.1452}.
The rate of convergence established in \cite{kubilius2017rate} is as follows: 
\begin{align}\label{eq:231004.1859}
  \babs{\hat H_n-H} = O((n/(\log n))^{\alpha})
  % \babs{\hat H_n-H} = O((\log n/n)^{1/4+\epsilon/2})
  % \babs{\hat H_n-H} = O(n^{-1/4}(\log n)^{1/4+\epsilon})
\end{align}
almost surely for any $-1/4>\alpha>-H/2$.
% almost surely for any $-1/4>-(1/4+\epsilon/2)>-H/2$.
% almost surely for any $0<\epsilon<H-\half$.

In the present paper, with more regularity imposed
on the coefficient functions $V^{[i]}$ ($i=1,2$),
we go further into considering the asymptotic expansion of 
the distribution of the rescaled error $\sqrt{n}(\hat H_n-H)$
associated with the above convergence.
In particular, we obtain a mixed central limit theorem 
(i.e. weak convergence to a mixed normal distribution) of this error.
It is noteworthy that the order of the normalizing factor in our result is $n^{1/2}$, 
and hence the rate of convergence of this estimator
in almost-sure sense is $n^{-1/2+\epsilon}$ for any $\epsilon>0$,
which improves the result \eqref{eq:231004.1859} obtained in \cite{kubilius2017rate}.
% {\myred at least in $L^p$ sense.}
%%
The approach we adopt in this paper is 
the asymptotic expansion technique of Skorohod integrals,
a methodology recently developed by \cite{nualart2019asymptotic},
which belongs to the active research area of investigating limit theorems 
using Malliavin calculus.
The same Hurst estimator for a fractional Brownian motion was considered in 
\cite{mishura2023asymptotic}, 
where the authors derived the asymptotic expansion formula for the error of that estimator.
The result presented in this article generalizes their findings to a more intricate model.

When applying the general theories of asymptotic expansion 
to the context of the estimator $\hat H_n$,
the major steps are finding the asymptotic (conditional) variance, 
identifying the limit random symbols, and 
verifying the conditions related to
the asymptotic orders of functionals and 
the non-degeneracy of the distribution of the variable in question.
Additionally, prior to deriving the asymptotic expansion, 
it is necessary to write a stochastic expansion of the variable 
as a sum of a Skorohod integral and terms of smaller order
that still contribute to the first order terms in the asymptotic expansion.
Through the steps described above, 
there appear many different functionals
% \sout{a functional is often decomposed into a sum of functionals, }
and it becomes essential to estimate 
the asymptotic order of each functional.
%\sout{ of each summand}. 
%%
In general, these functionals are expressed as 
randomly weighted polynomials of multiple Wiener integrals 
with respect to the fractional Brownian motion, 
and estimating their magnitudes can be challenging. 
To overcome this difficulty, we introduce a theory of exponent which provides 
an estimate of an upper bound of the order of a functional having a certain form 
through easy calculations based on a graphical representation of a structure of the functional.
In the paper \cite{yamagishi2023order}, motivated by the problem of the asymptotic expansion 
of the quadratic variation of fractional SDE,
Yamagishi and Yoshida introduced other theories of exponent 
that also rely on graph theoretical arguments.
While the theory of exponent in the current paper
shares similar arguments with the one in the previous paper,
it differs from \cite{yamagishi2023order} in that we mainly treat functionals related to 
the second-order differences of fractional Brownian motion.

% {\mygreen 
% To apply the asymptotic expansion for Skorohod integrals, 
% the essential steps are identification of the limit of the random symbols and 
% estimation of the Sobolev norms of variables obtained by 
% repeated Malliavin derivatives and their projections. 
% Generally, both steps were not easy so far 
% because of lack of a methodology for systematic assessment of 
% the magnitude of random polynomials of multiple Wiener integrals 
% with respect to the fractional Brownian motion.
% To deal with this problem, we introduce two exponents which give an upper bound of the order of functionals 
% with some easy calculation.
% These exponents are applicable to a wider class of functionals related to fBm than treated in this specific problem of 
% the asymptotic expansion of the quadratic variations. 
% Our exponents in a sense generalize the exponent defined by \cite{yoshida2020asymptotic} 
% for functionals of a Brownian motion 
% (i.e. fBm of $H=1/2$) having a structure similar to that we encounter in this paper.
% The stability and contraction properties under the Malliavin derivatives and the projections will be clarified 
% with the help of the exponents of functionals.
% }

The organization of the paper is as follows.
In Section \ref{sec:231002.2414},
we present the main results of the current paper 
and provide an overview of our approach to proving them. 
Section \ref{sec:231002.2418} offers some preliminaries
on the theories of asymptotic expansion and fractional Brownian motion.
Section \ref{sec:231002.2156} introduces the theory of exponent mentioned above.
\begin{comment}
  {\myred 具体的な定理を名指しでどの命題が何を示しているかを書いてあげると良い説．}
\end{comment}
The proof of the asymptotic expansion formula is provided in Section \ref{sec:231002.2423},
while the proofs of some technical lemmas 
% required for establishing the results 
are deferred to 
Section \ref{sec:231002.2426}.
The theory of exponent developed in Section \ref{sec:231002.2156}
plays a significant role throughout 
Sections \ref{sec:231002.2423} and \ref{sec:231002.2426}.
In Section \ref{sec:240701.2110},
the results of numerical simulations verify that the asymptotic expansion formula provides a better approximation than that by the (mixed) central limit theorem.

To conclude the introduction, 
we briefly mention some references related to
the recently developed theories of asymptotic expansion and 
their applications to stochastic models.
A general theory for the asymptotic expansion of martingales 
with a mixed normal limit was introduced in \cite{yoshida2013martingale}.
This theory has been applied to various contexts, 
including 
% For the asymptotic expansion of martingales having a mixed normal limit,
% a general theory was given in \cite{yoshida2013martingale} and 
% it is applied to 
the quadratic variation of a diffusion process in \cite{yoshida2012asymptotic},
power variations of a diffusion process in \cite{podolskij2016edgeworth},
a pre-averaging estimator in \cite{podolskij2017edgeworth}, and 
the Euler approximation of continuous diffusion processes in 
\cite{podolskij2020edgeworth}.
Although these stochastic models are within the realms of 
continuous martingale theory and It\^{o} calculus, 
the asymptotic expansion formulas obtained in these works 
involve the Malliavin derivatives of functionals.

For the expansion theory 
targeting objects within the Malliavin calculus,
general theories have been obtained in 
\cite{tudor2019asymptotic}, \cite{tudor2023high}, and 
\cite{nualart2019asymptotic}.
The first two papers utilize the so-called Gamma factor and 
are suitable to handle functionals residing in Wiener chaos of 
fixed orders.
In contrast, the last one considers functionals expressed as 
a Skorohod integral perturbed by another functional of a minor order,
allowing us to deal with more general and broader functionals,
such as the one in the present paper.
The applications of the first two papers
can be found in 
\cite{tudor2020asymptotic} and 
\cite{mishura2023asymptotic},
where the authors address
the quadratic variation of a mixed fractional Brownian motion and 
a Hurst estimator of a fractional Brownian motion, respectively.
Meanwhile, the papers 
\cite{yoshida2023asymptotic} and \cite{2024Yamagishi-AsymptoticEO}
build upon the theory presented in \cite{nualart2019asymptotic}, and 
explore functionals involving randomly-weighted variations of 
a classical Brownian motion and a fractional Brownian motion, 
respectively.

% {\myblue 
%   \begin{itemize}
%   \item martingale expansion in mixed normal limit \cite{yoshida2013martingale}
%   \item Asymptotic expansion for the quadratic form of the diffusion process 
%   \cite{yoshida2012asymptotic}
%   \item Edgeworth expansion for functionals of continuous diffusion processes
%   \cite{podolskij2016edgeworth}
%   \item Edgeworth expansion for the pre-averaging estimator
%   \cite{podolskij2017edgeworth}
%   \item Edgeworth expansion for Euler approximation of continuous
%   diffusion processes
%   \cite{podolskij2020edgeworth}

%   \item  Asymptotic expansion for vector-valued sequences of random variables with
%   focus on Wiener chaos. 
%   \cite{tudor2019asymptotic}
%   \item High order asymptotic expansion for Wiener functionals.
%   \cite{tudor2023high}
%   \item Asymptotic expansion of Skorohod integrals.
%   \cite{nualart2019asymptotic}
%   \item Asymptotic expansion of the quadratic variation of a mixed
%   fractional Brownian motion
%   \cite{tudor2020asymptotic}

%   % \item Parameter estimation in fractional diffusion models
%   % \cite{kubilius2017parameter}

%   \item Asymptotic expansion of an estimator for the Hurst coefficient
%   \cite{mishura2023asymptotic}

%   \item Asymptotic expansion and estimates of Wiener functionals
% \cite{yoshida2023asymptotic}

% \end{itemize}}

% \begin{itembox}[l]{Memo}
%   \begin{itemize}
%     \item 漸近展開について
%     \item skorohod integral 周りの漸近展開について
%   \end{itemize}
% \end{itembox}

\subsection{Notations}
The following notations are repeatedly used in the following sections.
\begin{itemize}
  \setlength{\parskip}{0cm} \setlength{\itemsep}{5pt}

  \item We denote by $\bbZ$, $\bbZ_{\geq0}$ and $\bbN$
  the set of integers, that of nonnegative integers and
  that of positive integers.
  We write $\ntwo=\cbr{j\in\bbZ \mid j\geq2}$ for notational convenience.
  For $n\in\bbN$, we write $[n]=\cbr{1,...,n}$ for short.

  \item For a (finite) set $S$, we denote 
  the cardinality of $S$ by $\abs{S}$.

  \item 
  % $C^\infty_b(\bbR)$: the set of the smooth functions 
  % whose derivatives of any order are bounded together with itself.
  For $k\in\bbN\cup\cbr{\infty}$, $C^k_b(\bbR)$ denotes
  the space of $k$ times continuously differentiable functions $\bbR\to\bbR$
  which are bounded together with its derivatives of order up to $k$.

  \item The $\bbR$-valued functions $V^{[1]}$ and $V^{[2]}$ 
  defined on $\bbR$ appears in SDE \eqref{eq:230925.1452} as 
  the diffusion coefficient and drift coefficient, respectively, and 
  $V^{[i;k]}$ is the $k$-th derivative of $V^{[i]}$.

  \item For a function $f:\bbR\to\bbR$, we write $f_t=f(X_t)$ for brevity,
  where $X_t$ is the solution to SDE \eqref{eq:230925.1452}.
  In particular, we often write $V^{[i]}_t$ for the coefficient $V^{[i]}(X_t)$ of SDE \eqref{eq:230925.1452}.

  \item 
  We write $t^n_j=j/n$ for $j=0,...,n$ and $n\in\bbN$.
  In the following argument, there appear both $t^n_j$ ($j=0,...,n$)
  and $t^{2n}_j$ ($j=0,...,2n$) from $\secVar_n$ and $\secVar_{2n}$, respectively.
  % However, when there is no fear of confusion, we write $t_j$ for $t^n_j=j/n$.
  However, when there is no risk of confusion, we use the notation $t_j$ to represent $t^n_j=j/n$.
  % We denote the indicator function of the interval $\sbr{\tjm,\tj}$
  % by $1_{n,j}$ or $1_j$.
  
  \item 
  We denote by $L^p$ the $L^p$-space of
  random variables on the probability space fixed in the following arguments.
  % When we simply write $L^p$, it means the $L^p$-space of
  % random variables on the probability space.
  We also write
  $L^{1+}=\cup_{p>1} L^p$ and $L^{\infty-}=\cap_{p>1}L^p$.
  The $L^p$-norm of a random variable is denoted by 
  $\norm{.}_p$ or $\norm{.}_{L^p}$.
  % We also write $\norm{.}_{L^p(P)}$ 
  % to distinguish it from the $L^p$-norm of functions on $[0,T]^k$, 
  % which is denoted by $\norm{.}_{L^p([0,T]^k)}$.
  We denote the Sobolev norm in Malliavin calculus by 
  $\norm{.}_{i,p}$ for $i\in\bbN$ and $p\geq1$.
  
  \item 
  Let  $\alpha\in\bbR$ and 
  $(F_n)_{n\in\bbN}$ be a sequence of random variables.
  We use the following notation to express asymptotic orders.
  \begin{itemize}
    \item [\labelitemi]
    $F_n=O_{L^{p}}(n^\alpha)$ 
    [resp. $F_n=o_{L^{p}}(n^\alpha)$],\quad
    if $\norm{F_n}_{p}=O(n^\alpha)$
    [resp. $\norm{F_n}_{p}=o(n^\alpha)$]\; for $p\geq1$.
    \item [\labelitemi]
    $F_n=O_{L^{\infty-}}(n^\alpha)$ 
    [resp. $F_n=o_{L^{\infty-}}(n^\alpha)$],\quad
    if $F_n=O_{L^{p}}(n^\alpha)$
    [resp. $F_n=o_{L^{p}}(n^\alpha)$]\; for every $p>1$.
  \end{itemize}
  For $(F_n)_{n\in\bbN}$ with $F_n\in\bbD^\infty$, we write
  \begin{itemize}
    \item [\labelitemi]
    $F_n=O_M(n^\alpha)$ [resp. $F_n=o_M(n^\alpha)]$, \quad
    if $\norm{F_n}_{i,p}=O(n^\alpha)$
    [resp. $\norm{F_n}_{i,p}=o(n^\alpha)$]
    for every $i\in\bbN$ and $p>1$.
  \end{itemize}
  When $F_n=O_{L^{p}}(n^{\alpha+\epsilon})$ for any $\epsilon>0$, we write 
  $F_n=\hat O_{L^p}(n^\alpha)$.
  For the other norms, we use the similar notations.

  % \item We use the notation $\norm{.}_\beta$
  % for the $\beta$-H\"{o}lder seminorm of 
  % $\beta$-H\"{o}lder continuous functions on $[0,T]$.
  % We use the notation $\norm{.}_\infty$ for the uniform norm.
  % For the seminorms restricted on an interval $[t,t']$,
  % we write $\norm{.}_{t,t',\beta}$ and $\norm{.}_{t,t',\infty}$.
  % See Section \ref{220405.1158} for their formal definitions.
\end{itemize}
%

%% file: subfiles/2-sec_main_result.tex
%
% \section{Main result, Outline and preliminaries}\label{sec:231002.2410}
% \section{Main results}
\section{%Asymptotic expansion of quadratic variations of SDE driven by fBm
Main results %of asymptotics of $\hat H_n$ 
and strategy}
\label{sec:231002.2414}
In this section, firstly we state the statement of the main result and
overview the strategies to prove it.
Let $(X_t)_{t\in[0,T]}$ be the solution 
to SDE \eqref{eq:230925.1452}.
Recall that the second-order quadratic variation of $X_t$ is defined at \eqref{eq:231002.0958}.
We denote the rescaled variation by
\begin{align*}
  \rv_n = n^{2H-1}\secVar_n
  = n^{2H-1}\sum_{j=1}^{n-1}\rbr{\secDiff{n}{j}X}^2.
\end{align*}
Under some assumptions, Theorem 5 of \cite{kubilius2012rate} showed 
that $\rv_n$ converges almost surely to $\iv$, where $\iv$ is written as 
\begin{align*}
  \iv=c_{2,H}\int^1_0 \rbr{V^{[1]}_t}^2 dt
\end{align*}
with the constant $c_{2,H}$ defined at \eqref{230925.1607}.
When $\iv\neq0$, 
the main part of the estimator $\hat H_n$ defined at \eqref{eq:231001.1046} 
can be written as
\begin{align}
  \half + \frac1{2\log2} 
  \log\frac{\secVar_n}{\secVar_{2n}} &= 
  \half + \frac1{2\log2} 
  \log\rbr{2^{2H-1} \frac{\rv_n}{\rv_{2n}}}
  \nn\\*&=
  H+
  \frac1{2\log2} \rbr{
    \log\rbr{\frac{\rv_n-\iv}{\iv}+1}-\log\rbr{\frac{\rv_{2n}-\iv}{\iv}+1}},
  \label{eq:231003.1930}
\end{align}
and this shows the heuristics of the definition of the estimator.

%   {\myblue
% \begin{align*}
%   \log\frac{\secVar_n}{\secVar_{2n}} &= 
%   \log\rbr{2^{2H-1} \frac{\rv_n}{\rv_{2n}}}
%   =
%   \rbr{H-\half}2\log2 +
%   \log\rbr{\frac{\convDiff_n}{\iv}+1}-\log\rbr{\frac{\convDiff_{2n}}{\iv}+1}.
% \end{align*}}
% 
% {\myred 
% \begin{align*}
%   \rv_n = n^{2H-1}\secVar_n
%   = n^{2H-1}\sum_{j=1}^{n-1}\rbr{\secDiff{n}{j}X}^2
% \end{align*}}

%
We define the functional $G_\infty$, 
which plays the role of the asymptotic variance of the rescaled estimator (up to some constant),
by
\begin{align*}
  G_\infty &=
  \varconst
  % \rbr{3\hat c - 2^{2H+1}\tilde c}
  \rbr{\iv}^{-2} \int^1_0 \rbr{V^{[1]}_t}^4 dt
  % \label{230725.1534}
  % \\
  % \iv&=
  % c_{2,H}\int^1_0 \rbr{V^{[1]}_t}^2 dt,
\end{align*}
with the constant $\varconst$ defined at \eqref{eq:231001.1716}.
To derive the asymptotic expansion of the distribution of $\hat H_n$, 
we assume several conditions on the SDE \eqref{eq:230925.1452}.
% We impose the following conditions on SDE \eqref{eq:230925.1452}.
%
\begin{ass}\label{ass:230927.1617}
  (i) $V^{[i]}\in C^\infty_b(\bbR)$ for $i=1,2$. 
  % where $C^\infty_b(\bbR)$ is the set of the smooth functions whose derivatives of any order are bounded together with itself.
  
  \item[(ii)] 
  The functional $G_\infty$ satisfies 
  $G_\infty^{-1}\in L^{\infty-}=\cap_{p>1}L^p$.

  \item[(iii)]
  The functional $\iv$ satisfies 
  $\bbV_\infty^{-1}\in L^{\infty-}$.
  % $\iv^{-1}\in L^{\infty-}$.
\end{ass}
The second condition is assumed to ensure the nondegeneracy of the distribution of 
$\sqrt{n}(\hat H_n-H)$.
% $Z_n$.
The third condition allows us to use the Taylor expansion to decompose the two $\log$'s 
in \eqref{eq:231003.1930}.
In \cite{kubilius2012rate}, the authors assumed a similar condition that 
$\iv$ is separated from zero, namely
$\iv\geq c_0$ almost surely with some constant $c_0>0$.

The main result of this article is stated as follows:
\begin{theorem*}\label{231002.1320}
  Suppose that $(X_t)_{t\in[0,T]}$ is the solution to SDE \eqref{eq:230925.1452}
  satisfying Assumption \ref{ass:230927.1617}.
  Let $\hat H_n$ be as in \eqref{eq:231001.1046}.
  For $M>0$ and $\gamma>0$,
  we denote by $\cale(M,\gamma)$ 
  the set of measurable functions $f:\bbR\to\bbR$ such that 
  $\abs{f(z)}\leq M(1+\abs{z})^\gamma$ for all $z\in\bbR$.
  Then, % for $\hat H_n$ defined at \eqref{eq:231001.1046}, 
  it holds that 
  \begin{align*}
    \sup_{f\in\cale(M,\gamma)}
    \abs{E[f(\sqrt{n} (\hat H_n - H))] - \int_{z\in\bbR} f(z) p_n(z) dz} 
    = {o}(n^{-\half})
  \end{align*}
  as $n\to\infty$,
  where the approximate density $p_n$ is given in Theorem \ref{thm:231001.1041}.
  % with $p_n(z)=2\log2\times p_n^Z(2\log2\,z)$
\end{theorem*}
The approximate density $p_n$ can be written as 
\begin{align*}
  p_n(z)&= 
  E\sbr{\phi(z;0,(2\log2)^{-2}\,G_\infty)} 
 + n^{-\half} \mu(z),
\end{align*}
where $\phi(z;0,v)$ is the density function of the normal distribution 
with mean $0$ and variance $v$, and $\mu(z)$ is some $n$-independent function on $\bbR$
with good integrability.
Thus, as a corollary of the above theorem, 
we obtain a mixed central limit theorem of $\hat H_n$:
\begin{proposition*}
  % {\myred Suppose that $(X_t)_{t\in[0,T]}$ is the solution to SDE \eqref{eq:230925.1452}
  % satisfying Assumption \ref{ass:230927.1617}, and 
  % let $\hat H_n$ be as in \eqref{eq:231001.1046}.}
  Let the assumptions of Theorem \ref{231002.1320} hold.
  Then the random variable $\sqrt{n}(\hat H_n-H)$ weakly converges to 
  $\caln(0,(2\log2)^{-2}\,G_\infty)$, that is, 
  the mixed normal distribution
  with conditional variance $(2\log2)^{-2}\,G_\infty$.
\end{proposition*}
We remark that this result is new to the best of the author's knowledge, and 
this fact indicates that the rate of convergence of $\hat H_n\to H$ is $n^{-\half}$,
which improves the previous result \eqref{eq:231004.1859} established in \cite{kubilius2017rate}.

From now on, 
we describe our approach to this result.
We define a functional $Z_n$ by
\begin{align*}
  Z_n = 2\log2\,\sqrt{n} (\hat H_n - H) \psi_n
\end{align*}
with a truncation functional $\psi_n(\in\bbD^\infty)$ 
which satisfies 
$P\sbr{\psi_n\neq1}=O(n^{-L})$ for any $L>0$.
The truncation $\psi_n$ provides the (higher-)differentiability of $\hat H_n$
in the sense of Malliavin,
while the factor $2\log2$ is only to simplify the calculation.
In Proposition \ref{230720.2014},
we obtain a stochastic expansion of $Z_n$ of the form 
\begin{align*}
  Z_n=M_n+n^{-\half} N_n + n^{-\half} Y_n,
\end{align*}
where 
$M_n$ is the Skorohod integral of the $\calh$-valued functional $u_n$,
that is $M_n=\delta(u_n)$, and
% (defined at \eqref{eq:230926.1446}), and 
$N_n$ and $Y_n$ are functionals of $O_M(1)$.
We first (in Section \ref{sec:231001.1030})
consider the asymptotic expansion of 
$Z_n^\circ=M_n+n^{-\half}N_n$
by means of the theory of asymptotic expansion of Skorohod integrals
recently developed by Nualart and Yoshida \cite{nualart2019asymptotic}.
Subsequently in Section \ref{sec:231001.1649}, by the perturbation method 
(see Section \ref{sec:231001.1036}), 
% (reviewed in \ref{sec:231001.1036}),
we calculate the asymptotic expansion of 
$Z_n=Z_n^\circ+n^{-\half}Y_n$ and then that of $\sqrt{n}(\hat H_n-H)$.
% {\myred which is the main result and is written as follows:}

% {\myred [exponentについて]}
%\subsection{About estimate of functionals of a certain form}
While %most arguments in deriving the asymptotic expansion 
%are related to 
verifying the conditions of the general theories 
in the context of the estimator $\hat H_n$, or more practically the functional $Z_n$ defined above,
to derive the asymptotic expansion,
% in the context of $Z_n^\circ$, $u_n$, $N_n$ and so on,
there often appear functionals of a certain common form and 
we need to have estimates of their norms (in $L^p$ or Malliavin sense).
An example of such functional is as follows:
% \begin{align*}
%   \rbr{\iv}^{-1} \convDiff^{(2,1)}_n &=
%   2\;\rbr{\iv}^{-1} 
%   n^{2H-1} \sum_{j=1}^{n-1} V^{[1]}_\tj V^{[(1;1),1]}_\tjm
%   I_3\rbr{\diffker^n_j \otimes (\bbone^n_j)^{\otimes2}} 
%   =:n^{2H-1}\caliN{1}_n
%   \\
%   D_{u_n}N_n &= 
%   \sum_{i_1,i_2=1}^2 (-1)^{i_1+i_2}  
%   \rbr{
%     i_1^{(2H-1)}%2^{(2H-1)(i_1-1)} 
%   n^{2H}\sum_{k=1}^3 D_{u_n^{(i_2)}}\caliN{k}_{i_1n} 
%   + i_1^{-1}%2^{-(i_1-1)}
%   D_{u_n^{(i_2)}}\caliN{4}_{i_1n}}
% \end{align*} 
\begin{align*}
  % +%\\&\qquad+
  % 2^{-1}\sum_{i=1,2}(-1)^{i-1}D_{u_n^{(i)}}\AN
  % \\&\qquad+ O_M(n^{(\half-H)\vee(-1+H)}),
  % \\
  % &(\text{a part of }D_{u_n^{(i_2)}}\caliN{1}_{i_1n})
  % \\&=
  % \abr{
  % \sum_{j=1}^{n-1} {\rbr{\iv}^{-1} }
  % V^{[1]}_\tj V^{[(1;1),1]}_\tjm
  % I_2\rbr{\diffker^n_j \otimes \bbone^n_j} \bbone^n_j,
  % n^{2H-\half} 
  % \sum_{j\in[2n-1]} (\iv)^{-1} \Brbr{V^{[1]}_{t^{2n}_j}}^2
  % I_1(\diffker^{2n}_j) \diffker^{2n}_j}
  % \\&=
  &%n^{2H-\half} 
  \sum_{j_1\in[n-1]} 
  \sum_{j_2\in[2n-1]} 
  (\iv)^{-2} 
  V^{[1]}_{t^n_{j_1}} 
  V^{[(1;1)]}_{t^n_{j_1-1}}
  V^{[1]}_{t^n_{j_1-1}}
  \Brbr{V^{[1]}_{t^{2n}_{j_2}}}^2
  I_2\rbr{\diffker^n_{j_1} \otimes \bbone^n_{j_1}} 
  I_1(\diffker^{2n}_{j_2}) 
  \abr{\bbone^n_{j_1},\diffker^{2n}_{j_2}}_\calh
  % \\
  % % n^{2H-1}
  % \tilde{\cali}^{(\sf N)(1)}_n&=
  % % n^{2H-1} 
  % \sum_{j=1}^{n-1} {\rbr{\iv}^{-1} }
  % V^{[1]}_\tj V^{[(1;1),1]}_\tjm
  % I_3\rbr{\diffker^n_j \otimes (\bbone^n_j)^{\otimes2}} 
  % \\ 
  % \tilde u_n&=
  % n^{2H-\half} 
  % \sum_{j\in[2n-1]} (\iv)^{-1} \Brbr{V^{[1]}_{t^{2n}_j}}^2
  % I_1(\diffker^{2n}_j) \diffker^{2n}_j
\end{align*}
Here,
$\bbone^m_{j}$ is the indicator function of $[(j-1)/m,j/m]$,
$\diffker^{m}_{j} = \bbone^{m}_{j+1} - \bbone^{m}_{j}$, 
and 
$\abr{\bbone^n_{j_1},\diffker^{2n}_{j_2}}_\calh$ is the inner product of 
$\bbone^n_{j_1}$ and $\diffker^{2n}_{j_2}$
in the Hilbert space $\calh$ associated with a representation of fBm; 
see Section \ref{sec:231004.1910}.
The factors
$I_2\rbr{\diffker^n_{j_1} \otimes \bbone^n_{j_1}}$ and 
$I_1(\diffker^{2n}_{j_2})$ are (multiple) stochastic integrals.
To obtain estimates of the asymptotic order of $L^p$-norms or norms in Malliavin calculus, 
it requires us to consider a great number of patterns of functionals,
which emerges from the use of the IBP formula and Leibniz rule.

In \cite{yamagishi2022order}, Yamagishi and Yoshida considered 
the asymptotic expansion of quadratic variations of SDE \eqref{eq:230925.1452},
where they encounter similar problems as above, and 
% and developed a theory of exponent to overcome similar obstacles.
developed a theory of exponent based on some graphical representation 
capturing the structure of functionals.
In this paper, we introduce another theory of exponent to obtain 
estimates of the norms of functionals in Section \ref{sec:231002.2156}, 
since the previous theory only can handle functionals having $\bbone^n_j$ 
as the kernels of stochastic integrals or in the inner products,
and it also restricts the Hurst parameter $H$ to $(\frac12,\frac34)$.

%% file: subfiles/3-sec_preliminaries.tex
\section{Preliminaries}
\label{sec:231002.2418}
\subsection{Asymptotic expansion of Skorohod integrals}
\label{230901.1140}
We review the theory of asymptotic expansion of Skorohod integrals 
described in \cite{nualart2019asymptotic} in the one-dimensional case.
Let $(\Omega, \calf, P)$ be a complete probability space equipped with an isonormal Gaussian process 
$W=\cbr{W(h)_{h\in\calh}}$ 
on a separable real Hilbert space $\calh$.
We denote the Malliavin derivative operator by $D$ and 
its adjoint operator, namely the divergence operator or the Skorohod integral, by $\delta$.
For $p\geq1$, $k\in\bbN$ and a real separable Hilbert space $V$, 
we write $\bbD^{k,p}(V)$ for the Sobolev space of $V$-valued random variables 
that have the Malliavin derivatives up to $k$-th order which have finite moments of order $p$.
We write $\bbD^{k,p}=\bbD^{k,p}(\bbR)$, 
$\bbD^{k,\infty}(V)=\cap_{p>1}\bbD^{k,p}(V)$ and
$\bbD^{\infty}(V)=\cap_{p>1,k\geq1}\bbD^{k,p}(V)$.
The $(k,p)$-Sobolev norm is denoted by $\norm{\cdot}_{k,p}$.
We refer to the monograph \cite{nualart2006malliavin} for a detailed account on this subject.

Consider a sequence of random variables $Z_n$ defined on the probability space 
$(\Omega, \calf, P)$ written as 
\begin{align}\label{eq:231004.1810}
  Z_n=M_n+r_nN_n,
\end{align}
where
$M_n=\delta(u_n)$ is the Skorohod integral of an $\calh$-valued random variable
$u_n\in\Dom(\delta)$,
$N_n$ is a random variable and 
$(r_n)_{n\in\bbN}$ is a sequence of positive numbers such that 
$\lim_{n\to\infty}r_n=0$.
%
% {\mygreen
The variable $Z_n$ of \eqref{eq:231004.1810} is a perturbation of $M_n$ when $N_n=O_p(1)$ as $n\to\infty$.
Such a perturbation always appears when one derives 
a stochastic expansion of a statistic $Z_n$ around a principal part $M_n$ that is easy to handle by a limit theorem. 
%
% In the case of the quadratic variation 
% $\bbV_n =n^{2H-1}\sum_{j=1}^n (\Delta_jX)^2$ of (\ref{220420.1130}), 
% the scaled variable $Z_n =n^{1/2}\big(\bbV_n-\bbV_\infty)$ admits the stochastic expansion (\ref{202308260454}) 
% with 
% $r_n = n^{2H-3/2}$ and 
% $M_n=\delta(u_n)$ for $u_n$ of (\ref{220404.1632}): 
% $u_n = n^{2H-1/2} \sum_{j\in[n]} a_{t_{j-1}}  I_1(1_j) 1_j$. 
% The variable $N_n$ has an expression given in (\ref{220404.1633}). 
% }
%

We consider an asymptotic expansion of the distribution of $Z_n$ in the situation where 
$M_n$ stably converges to 
a mixed normal distribution, that is, 
$M_n\overset{d_s}{\to}M_\infty=G_\infty^{1/2}\zeta$ as $n\to\infty$,
where $G_\infty$ is a positive random variable 
and $\zeta$ is a standard Gaussian random variable independent of $\calf$. 
The variable $\zeta$ is given on an extension of the probability space $(\Omega,\calf,P)$. 
Here the stable convergence means that the convergence $(M_n,Y)\to^d(M_\infty,Y)$ holds for any 
random variable $Y$ measurable with respect to $\sigma[W]$, the $\sigma$-field generated 
by the isonormal Gaussian process $W$.  
%
% We consider a stable convergence of $Z_n$ to 
% a mixed normal distribution $G_\infty^{1/2}\zeta$,
% where $G_\infty$ is a positive random variable 
% and $\zeta$ is a standard normal distribution independent of $\calf$.}

Nualart and Yoshida \cite{nualart2019asymptotic} introduced the following random symbols.
To write random symbols, we use a simplified notation  
due to the one-dimensional setting.
%We write $\mS(\xi)=\sum_{k}\chi_k\,\xi^k$ (a finite sum) for a polynomial random symbol,
% where $\chi_k$ are coefficient random variables
% %\redb{$\chi_k\,\xi^k$ is the random symbol of degree $k$ with its coefficient random variable $\chi_k$}
%  and $\xi$ stands for a dummy variable.
We denote
$D_{u_n}F=\abr{DF,u_n}_\calh$ for a random variable $F$ regular enough.
The quasi-tangent is defined by 
\begin{align*}
  \qtan_n{[\sfi\sfz]^2}=
  r_n^{-1}\brbr{\babr{DM_n{[\sfi\sfz]}, u_n{[\sfi\sfz]}}_\calh
  -G_\infty{[\sfi\sfz]^2}}
  =r_n^{-1}\brbr{D_{u_n}M_n-G_\infty}{[\sfi\sfz]^2}.
\end{align*}
The quasi-torsion and modified quasi-torsion are defined by 
\begin{align*}
  \qtor_n{[\sfi\sfz]^3}&=
  r_n^{-1}\Babr{D\babr{DM_n{[\sfi\sfz]}, u_n{[\sfi\sfz]}}_\calh,
  u_n{[\sfi\sfz]}}_\calh
  =r_n^{-1}(D_{u_n})^2M_n{[\sfi\sfz]^3}
  \\
  \mqtor_n{[\sfi\sfz]^3}&=
  r_n^{-1}\babr{DG_\infty{[\sfi\sfz]^2},u_n{[\sfi\sfz]}}_\calh
  =r_n^{-1}D_{u_n}G_\infty{[\sfi\sfz]^3},
\end{align*}
respectively.

We write
\begin{align*}
  G^{(2)}_n
  %G^{(2)}_n(\sfz) 
  &%= r_n\qtan_n
  = D_{u_n}M_n - G_\infty
  \tand% \\
  G^{(3)}_n
  %G^{(3)}_n(\sfz) 
  % &%= r_n\mqtor_n 
  = D_{u_n}G_\infty,
\end{align*}
and define the following random symbols
\begin{align*}
  \mS^{(3,0)}_n(\tti\sfz)&=
  \frac13\qtor_n[\tti\sfz]^3
  =\frac13r_n^{-1}(D_{u_n})^2 M_n[\tti\sfz]^3
  \\
  \mS^{(2,0)}_{0,n}(\tti\sfz)&
  % \frac12 r_n^{-1}G^{(2)}_n(\sfz)=
  =\frac12 \qtan_n[\tti\sfz]^2
  =\frac12 r_n^{-1}\brbr{D_{u_n} M_n - G_\infty}[\tti\sfz]^2
  % =\frac12 r_n^{-1}\brbr{D_{u_n[\tti\sfz]} M_n[\tti\sfz] - G_\infty[\tti\sfz]^2}
  \\
  \mS^{(1,0)}_n(\tti\sfz)&=N_n[\tti\sfz]\\
  \mS^{(2,0)}_{1,n}(\tti\sfz)&=
  D_{u_n} N_n[\tti\sfz]^2
  % D_{u_n[\tti\sfz]} N_n[\tti\sfz]
\end{align*}
for $\tti\sfz\in\tti\bbR$.
%
% We consider random symbols 
We denote by
$\mS^{(3,0)}$,
$\mS^{(2,0)}_{0}$,
$\mS^{(1,0)}$ and
$\mS^{(2,0)}_{1}$
% These work as 
the limits of the above random symbols
$\mS^{(3,0)}_n$,
$\mS^{(2,0)}_{0,n}$,
$\mS^{(1,0)}_n$ and
$\mS^{(2,0)}_{1,n}$, respectively.
(The meaning of the limit is explained in the condition {\bf[D]} (iii).)
Let
\begin{align*}%\label{eq:230810.1905}
  \Psi(\sfz)=
  \exp\rbr{2^{-1} G_\infty[\tti\sfz]^2}
  \yeq \exp\rbr{-2^{-1} G_\infty\>\sfz^2}
\end{align*}
for $\sfz\in\bbR$.
The following set of conditions {\bf [D]} from Nualart and Yoshida \cite{nualart2019asymptotic}
is a sufficient condition to validate the asymptotic expansion of $Z_n$ in \eqref{eq:231004.1810}.
The parameter $l$ about differentiability below is $l=9$ in this case.
For a one-dimensional functional $F$, we write 
$\Delta_F %=\sigma_F
=\det(\abr{DF,DF}_\calh)=\abr{DF,DF}_\calh=\norm{DF}_\calh^2$
%for the Malliavin covariance (matrix) of $F$,
for the Malliavin covariance of $F$.

\begin{itemize}
  \item [{\bf [D]}]
  \begin{itemize}
    \item [(i)]
    $u_n\in\bbD^{l+1,\infty}(\mH)$,
    %$u_n\in\bbD^{l+1,\infty}(\mH\otimes\bbR^\sfd)$,
    $G_\infty \in \bbD^{l+1,\infty}(\bbR_+)$,
    %$G_\infty \in \bbD^{l+1,\infty}(\bbR^\sfd \otimes_+ \bbR^\sfd)$,
    %$W_n, 
    and
    $N_n\in\bbD^{l,\infty}$.
    %$N_n\in\bbD^{l,\infty}(\bbR^\sfd)$,
    %$W_\infty\in\bbD^{l\vee\sfd_2,\infty}(\bbR^\sfd)$,
    %$X_n\in\bbD^{l,\infty}(\bbR^{\sfd_1})$,
    %$X_\infty\in\bbD^{l\vee(\sfd_2+1),\infty}(\bbR^{\sfd_1})$,
  
    \item [(ii)]
    There exists a positive constant $\kappa$ such that the following estimates hold for every $p>1$:
    \begin{align}
      &\norm{u_n}_{l,p}=O(1)
      \label{220215.1241}\\
      &\bnorm{G^{(2)}_n}_{l-2,p}=O(r_n)
      \label{220215.1242}\\
      &\bnorm{G^{(3)}_n}_{l-2,p}=O(r_n)
      \label{220215.1243}\\
      &\bnorm{D_{u_n}G^{(3)}_n}_{l-1,p}=O(r_n^{1+\kappa})
      \label{220215.1244}\\
      &\bnorm{D^2_{u_n}G^{(2)}_n}_{l-3,p}=O(r_n^{1+\kappa})
      \label{220215.1245}\\
      &\norm{N_n}_{l-1,p}=O(1)
      \label{220215.1246}\\
      &\bnorm{D^2_{u_n} N_n}_{l-2,p}=O(r_n^{\kappa})
      \label{220215.1247}
    \end{align}
  
    \item [(iii)]
    For each pair 
    $(\mT_n,\mT)= (\mS_n^{(3,0)},\mS^{(3,0)})$,
    $(\mS_{0,n}^{(2,0)},\mS_0^{(2,0)})$,
    $(\mS_n^{(1,0)},\mS^{(1,0)})$ and 
    $(\mS_{1,n}^{(2,0)},\mS_1^{(2,0)})$,
    the following conditions are satisfied.
    \begin{itemize}
      \item [(a)]
      $\mT$ is a polynomial random symbol the coefficients of which are in
      $L^{1+}=\cup_{p>1} L^p$.
  
      \item [(b)]
      For some $p>1$, there exists a polynomial random symbol $\bar\mT_n$ such that
      the coefficients of $\bar\mT_n$ belongs to $L^p$,
      the equation
      $E \sbr{\Psi (\sfz)\mT_n(\tti\sfz)} = E[\Psi (\sfz)\bar\mT_n(\tti\sfz)]$
      holds for $\sfz\in\bbR$, and
      the convergence $\bar\mT_n\to\mT$ in $L^p$ holds.
    \end{itemize}
  
    \item [(iv)]
    \begin{itemize}
      \item [(a)] $G^{-1}_{\infty}\in L^{\infty-}$
      \item [(b)]
      There exist 
      $\kappa'>0$ such that
      \begin{align*}
      P[\Delta_{M_n}<s_n]=O(r_n^{1+\kappa'})
      \end{align*}
      for some positive random variables $s_n\in\bbD^{l-2,\infty}$ satisfying 
      $\sup _{n\in\bbN}(\norm{s_n^{-1}}_p + \norm{s_n}_{l-2,p})<\infty$
      for every $p>1$.
    \end{itemize}
  \end{itemize}
\end{itemize}
Writing 
$\bar\mT_n(\sfi\sfz)=\sum_{k\geq0}c_n^k[\sfi\sfz]^k$ and 
$\mT(\sfi\sfz)=\sum_{k\geq0}c^k[\sfi\sfz]^k$,
the convergence $\bar\mT_n\to\mT$ in $L^p$ means that
there exists $k_0\in\bbZ_{\geq0}$ such that 
$c_n^k$ and $c^k$ are all zero for $k>k_0$, and
$c_n^k\to c^k$ in $L^p$ for every $k\leq k_0$.
%
\begin{comment}
  {\myred [なにか説明を入れるか？ quad-varの文章を参照せよとjournal投稿版に書くかな？]
  The functional $s_n$ in (iv) (b) is used to make
  a truncation functional to gain local asymptotic non-degeneracy.
  See Section 7 of \cite{nualart2019asymptotic} for a construction of a truncation function.}
\end{comment}

We write $\phi(z;\mu,v)$ for the density function of the normal distibution 
with mean $\mu\in\bbR$ and variance $v>0$.
For a (polynomial) random symbol 
$\varsigma(\sfi\sfz) = \sum_\alpha c_\alpha[\sfi\sfz]^\alpha$
% $\varsigma(\xi) = \sum_\alpha c_\alpha\,\xi^\alpha$
with a random variable $c_\alpha$ and $\alpha\in\bbZ_{\geq0}$,
% where $\xi$ is a dummy variable, 
the action of the adjoint $\varsigma(\partial_z)^*$ to $\phi(z;0,G_\infty)$ under the expectation
is defined by 
\begin{align}
  E\sbr{\varsigma(\partial_z)^* \phi(z;0,G_\infty)}
  =\sum_\alpha (-\partial_z)^\alpha E\sbr{c_\alpha\phi(z;0,G_\infty)}.
  \label{eq:230927.1623}
  %{220401.1952}
\end{align} 
Define the random symbol 
$\mS_n = 1 +r_n\mS$ with
\begin{align*}
  \mS(\tti\sfz) = 
  \mS^{(3,0)}(\tti\sfz) + \mS^{(2,0)}_{0}(\tti\sfz) +
  \mS^{(1,0)}(\tti\sfz) + \mS^{(2,0)}_{1}(\tti\sfz)
\end{align*}
and the approximate density 
$p_n(z)$ by
% $\hat p_n(z)$ by
\begin{align*}
  p_n(z) = 
  % \hat p_n(z) = 
  E\sbr{\mS_n(\partial_z)^* \phi(z;0,G_\infty)}.
\end{align*}
% where the action of the adjoint of a random symbol is defined at (\ref{220401.1952}).
%
For $M,\gamma>0$, we denote by
$\cale(M,\gamma)$ 
% $\hat\cale(M,\gamma)$
the set of measurable functions $f:\bbR\to\bbR$
such that 
$\abs{f(z)}\leq M(1+\abs{z})^\gamma$ for all $z\in\bbR$.
The following theorem rephrases Theorem 7.7 of \cite{nualart2019asymptotic}.
\begin{theorem}\label{thm:230927.1318}
  [Theorem 7.7 of \cite{nualart2019asymptotic}]
  Suppose that Condition {\bf [D]} is satisfied.
  Then, for each $M,\gamma\in\bbR_{>0}$,
  \begin{align*}
    \sup_{f\in\cale(M,\gamma)}
    % \sup_{f\in\hat\cale(M,\gamma)}
    \abs{E\sbr{f(Z_n)}- \int_{\bbR}f(z) p_n(z)dz}
    % \abs{E\sbr{f(Z_n)}- \int_{\bbR}f(z)\hat p_n(z)dz}
    =o(r_n)
    \quad\text{as }n\to\infty.
  \end{align*}
\end{theorem}

\subsection{Perturbation method}\label{sec:231001.1036}
Here we review the perturbation method 
used in Yoshida 
\cite{yoshida1997malliavin,
yoshida2001malliavin,
yoshida2013martingale, yoshida2023asymptotic}
and in Sakamoto and Yoshida \cite{sakamoto2003asymptotic}.
While we only need the $1$-dimensional version 
in the argument of Section \ref{sec:231001.1649}, 
we give an exposition in a general dimension.
Consider a $d$-dimensional random variable $\bbS_n$ that decomposes as 
\begin{align*}
  \bbS_n = \bbX_n + r_n\bbY_n,
\end{align*}
where $\bbX_n$ and $\bbY_n$ are $d$-dimensional random variables, % for $n>0$, 
and $r_n$ is a positive number converging to $0$ as $n\to\infty$.
Suppose that positive numbers $M,\gamma,\ell_1$ and $\ell_2$, and 
a functional $\xi_n$ satisfy the following conditions:
\begin{itemize}
  \item [(i)] 
  There exists $m\in\bbN$ such that
  $\ell_1\geq\ell_2-1$, 
  $\ell_2\geq(2m+1)\vee(d+3)$, and 
  $m>\half\gamma+1$.
  % {\myred$m$はSakamoto Yoshidaの証明の中では超関数関連のregularityのorderとして使われている}
  \item[(ii)] 
  $\sup_{n>0}\norm{\bbX_n}_{\ell_2,p}+ \sup_{n>0}\norm{\bbY_n}_{\ell_2,p}<\infty$
  for every $p>1$.
  \item[(iii)]
  $(\bbX_n,\bbY_n)\overset{d}{\to}(\bbX_\infty, \bbY_\infty)$ 
  as $n\to\infty$ for some random variables $\bbX_\infty$ and $\bbY_\infty$.
  \item[(iv)]
  $\sup_{n>0}\norm{\xi_n}_{\ell_1,p}<\infty$ for every $p>1$.
  \item[(v)]
  $P\sbr{\abs{\xi_n}>\half}=O(r_n^\alpha)$ with some $\alpha>1$.
  \item[(vi)]
  $\sup_{n>0} E\sbr{\bbone_{\cbr{\abs{\xi_n}<1}} \Delta_{\bbX_n}^{-p}}<\infty$ 
  for every $p>1$.
  \item[(vii)]
  There exists a signed measure $\Psi_n$ on $\bbR^d$ such that 
  \begin{align*}
    \sup_{f\in\cale(M,\gamma)}
    \abs{E\sbr{f(\bbX_n)}-\int_{\bbR^d}f(z)\Psi_n(dz)}
    =o(r_n)
  \end{align*}
  as $n\to\infty$, where 
  $\cale(M,\gamma)$ is the set of measurable functions 
  $f:\bbR^d\to\bbR$ satisfying $\abs{f(z)}\leq M(1+\abs{z})^\gamma$.
\end{itemize}

Under these conditions, Sakamoto and Yoshida \cite{sakamoto2003asymptotic} showed the following result:
\begin{theorem}\label{thm:230927.1632}
  It holds that 
  \begin{align*}
    \sup_{f\in\cale(M,\gamma)}
    \abs{E\sbr{f(\bbS_n)} 
    - \int_{\bbR^d}f(z)\cbr{\Psi_n(dz)+r_n g_\infty(z)dz}}
    =o(r_n)
  \end{align*}
  as $n\to\infty$, where 
  \begin{align*}
    g_\infty(z) = 
    -\partial_z\cdot\cbr{E\sbr{\bbY_\infty\mid\bbX_\infty=z}p^{\bbX_\infty}(z)}.
  \end{align*}
\end{theorem}

\subsection{Fractional Brownian motion}\label{sec:231004.1910}
We review the definitions and notations about a fractional Brownian motion and the related Malliavin calculus.
Although we only consider the case $T=1$ in the following sections,
we give an exposition for a general $T>0$.
Let $B=\cbr{B_t\mid t\in[0,T]}$ a fractional Brownian motion 
with Hurst parameter $H\in(1/2,1)$ %$H\in\rbr{1/2,3/4}$
defined on some complete probability space $(\Omega,\calf, P)$.
An inner product on the set $\cale$ of step functions on $[0,T]$ is defined by 
\begin{align*}
  \abr{1_{[0,t]},1_{[0,s]}}_\calh = E\sbr{B_s\,B_t}
  = \frac12\rbr{\abs{t}^{2H}+\abs{s}^{2H}-\abs{t-s}^{2H}}.
\end{align*}
and $\calh$ is the closure of $\cale$ with respect to
$\snorm{\cdot}_\calh=\abr{\cdot,\cdot}_\calh^{1/2}$.
The map $\cale\ni1_{[0,t]}\mapsto B_t\in L^2(\Omega, \calf, P)$
can be extended linear-isometrically to $\calh$.
We denote this map by $\phi\mapsto B(\phi)$ and 
the process $\cbr{B(\phi), \phi\in \calh}$ is an isonormal Gaussian process. 
In the following sections, tools from Malliavin calculus are based on 
this isonormal Gaussian process.
We write $I_q(h)$ for the $q$-th multiple stochastic integral of $h\in\calh^{\otimes q}$.

It is known that the Hilbert space $\calh$ contains 
not only measurable functions on $[0,T]$ but also distributions as its elements;
see \cite{pipiras2000integration} and \cite{pipiras2001classes} for detailed accounts.
As a subspace of the Hilbert space $\calh$,
we have the linear space $\abs\calh$ of measurable functions $\phi:[0,T]\to\bbR$ such that
\begin{align*}
  \int_0^T\int_0^T \abs{\phi(s)}\abs{\phi(t)}\abs{t-s}^{2H-2}dsdt<\infty,
\end{align*}
and for $\phi$ and $\psi$ in $\abs\calh$, the inner product of $\calh$ is written as 
\begin{align*}
  \abr{\phi,\psi}_\calh=
  \constInProd
  \int^T_0\int^T_0\phi(s)\psi(t)\abs{t-s}^{2H-2}dsdt,
\end{align*}
where $\constInProd=H(2H-1)$.
For a measurable function $\phi$ on $[0,T]^l$, we define
\begin{align*}
  \norm\phi_{\abs\calh^{\otimes l}}^2
  =\constInProd^l \int_{[0,T]^{l}} \int_{[0,T]^{l}} \abs{\phi(u)} \abs{\phi(v)}
  \abs{u_1-v_1}^{2H-2}...\abs{u_{l}-v_{l}}^{2H-2}dudv,
\end{align*}
with $u=(u_1,...,u_l)$ and $v=(v_1,...,v_l)$,
and the space
$\abs\calh^{\otimes l}=\bcbr{\phi:[0,T]^l\to\bbR\mid \norm\phi_{\abs\calh^{\otimes l}}<\infty}$ 
of measurable functions
forms a subspace of the $l$-fold tensor product space $\calh^{\otimes l}$ of $\calh$.
We often drop $\calh$ from the notation $\abr{\cdot,\cdot}_\calh$, and 
instead write $\abr{\cdot,\cdot}$ for brevity, when there is no risk of confusion.
We refer to \cite{nualart2006malliavin} or 
\cite{nourdin2012selected} for a detailed account on fBm. 
We collect some basic estimates related to SDE's driven by a fBm with Hurst index $H>1/2$
in Section \ref{sec:231004.1533}.

Recall that we write 
$t^n_j=j/n$ for $n\in\bbN$ and $j\in\cbr{0,...,n}$.
We denote by $\bbone^n_j$ the indicator function $\bbone_{[t^n_{j-1},t^n_j]}$
of the interval $[t^n_{j-1},t^n_j]$,
and we write $\diffker^n_j=\bbone^n_{j+1}-\bbone^n_j$.
We regard $\bbone^n_j$ and $\diffker^n_j$ as an element of $\abs{\calh}\subset\calh$.
%%%
%%%
% {\mygreen 
% We write $[n]=\cbr{1,..,n}$.
% Define $\rho_H(k)$, $c_H$ and 
% $\beta_n\rbr{j_1,j_2}$ by %=\beta_{j_1,j_2}
%  %:= \abr{1_{\jon}, 1_\jtw}_\calh$  for $j_1,j_2\in[n]$ by
% \begin{align}
%   \rho_H(k) &= \frac12 \rbr{\abs{k+1}^{2H} +\abs{k-1}^{2H} -2\abs{k}^{2H}}
%   =\abr{1_{[0,1]}, 1_{[k,k+1]}}_\calh
%   &&\tfornsp k\in\bbZ,
%   \nn\\
%   c_H^2 &%= c_{H,(\ref{210417.1805})}^2
%   =\sum_{k\in\bbZ} \rho_H(k)^2 ,
%   \label{210417.1805}
%   \\
%   \beta_n\rbr{j_1,j_2}=\beta_{j_1,j_2}%=\beta_{j_1,j_2;T}
%   &= \abr{1_{\jon}, 1_\jtw}_\calh
%   = T^{2H} n^{-2H} \rho_H(\jon-\jtw)
%   &&\tfornsp \jon,\jtw\in[n],
%   \nn%\label{220422.1911}
% \end{align}
% respectively.}
%%%
%%%
We define a constant $c_{2,H}$ by 
\begin{align}
  % \bbone^n_j&=\bbone_{[t^n_{j-1},t^n_j]}
  % \label{230925.1608}
  % \\
  % \diffker^n_j&=\bbone^n_\jp - \bbone^n_j
  % \label{230925.1609}
  % \\
  c_{2,H}&=4-2^{2H}.
  % =\abr{\diffker^1_j,\diffker^1_j}
  \label{230925.1607}
\end{align}
Note that the relation 
$c_{2,H}=n^{2H}\abr{\diffker^n_j,\diffker^n_j}_\calh$
holds.
%%

%% file: subfiles/4-1-sec_exponent_def.tex
%
% \tableofcontents
\section{Theory of exponent}\label{sec:231002.2156}
As is explained in Section \ref{sec:231002.2414},
to derive the asymptotic expansion of functionals related to some variations of 
fractional Brownian motion,
it is necessary to establish estimates of norms %(in $L^p$ or Malliavin sense)
of functionals with certain common structure.
To capture this structure, we define a {\it  weighted graph}
and an {\it exponent} for each weighted graph, denoted by $G$ and $e(G)$, respectively.
We can roughly say that a functional corresponding to a weighted graph $G$ 
has the order of $n^{e(G)}$ at most as $n\to\infty$.
Specifically, 
Propositions \ref{221208.1720} and \ref{230725.1133} show this result in
$L^p$ norms and in the Sobolev norms of Malliavin calculus respectively.
Additionally, in Proposition \ref{221223.1757}  
we consider how 
the action of the operator $D_{u_n}$ on the functionals
changes the exponent of the corresponding weighted graph,
where $u_n$ has the form such as \eqref{eq:230926.1446},
and $D_{u_n}$ means $D_{u_n}F=\abr{DF,u_n}_\calh$ for a random variable $F$ regular enough.

The organization of this section is as follows.
In Section \ref{sec:231005.1616}, we give the definitions to illustrate the theory of exponent.
We prove the result about the asymptotic order of
the expectation and norms of functionals in Section \ref{sec:231005.1622},
and study the problem of $D_{u_n}$ in Section \ref{sec:231005.1626}.
The preliminary lemmas are collected in Section \ref{230519.1730}.

\subsection{Definitions of
weighted graph, exponent and related functional}
\label{sec:231005.1616}

First we introduce some definitions on which 
the weighted graphs and the theory of exponent are built.
%
%
%%%%% set of pairs %%%%%
For a nonempty finite set $V$, 
we denote by $p(V)$
the set of subsets of $V$ having two distinct elements.
When $V$ is a singleton, i.e. $V=\cbr{v_0}$ with some $v_0$,
we define $p(V)=\emptyset$.
The element $\cbr{v,v'}$ of $p(V)$ is 
written as $[v,v']$ hereafter to emphasize that 
it is a set consisting of distinct elements ($v\neq v'$).
%
%%%%% set of `expanded' vertices $\whv$ %%%%%
% For nonempty finite set $V$, w
We denote 
$\whv=V\times\cbr{1,2}$.
We often denote elements of $\whv$ by $\hv, \hv', \hv_0, \hv_1...$.
%
%%%%% expanded pairs %%%%%
We define $\hp(V)$ to be
the set of subsets 
$\cbr{(v,\kap),(v',\kap')}$
of $\whv=V\times \cbr{1,2}$ having two (distinct) elements such that 
$v\neq v'(\in V)$.
If $V$ is a singleton, we define $\hp(V)=\emptyset$.
Again we write $[(v,\kap),(v',\kap')]$ or $[\hv,\hv']$ for elements of $\hp(V)$.
%
%%%%% `expanded' pair for two disjoint sets %%%%%
For disjoint finite sets $V_0$ and $V_1$ such that 
$\abs{V_0\sqcup V_1}\geq1$,
we define $\hp(V_0,V_1)$ by 
\begin{align*}
  \hp(V_0,V_1) = 
  \cbr{[\hv,\hv']\in\hp(V_0\sqcup V_1)\mid 
  (\hv,\hv')\text{ or }(\hv',\hv)\in \whv_0\times \whv_1}.
\end{align*}
When either $V_0$ or $V_1$ is empty, we read 
$\hp(V_0,V_1) = \emptyset$.
Notice that there is a bijection $\hp(V_0,V_1)\simeq \whv_0\times\whv_1$ and 
$\hp(V_0\sqcup V_1)$ decomposes as
$\hp(V_0\sqcup V_1)=\hp(V_0)\sqcup\hp(V_0,V_1)\sqcup\hp(V_1)$.

%
% \subsubsection*{weighted graph}
%%%%% weighted graph %%%%%
Let $V$ be a nonempty finite set.
Given functions 
$\ewt:\hp(V)\to\bbZ_{\geq0}$ and 
$\vwq:\whv\to\bbZ_{\geq0}$,
we call the triplet $(V,\ewt,\vwq)$ {\it a weighted graph}.
When $V$ is a singleton, the domain of $\ewt$ is $\hp(V)=\emptyset$, and 
$\ewt$ is the function from the empty set to $\bbZ_{\geq0}$,
which we denote by $0$.
For a weighted graph $G=(V,\ewt,\vwq)$, 
we call $V$, $\ewt$ and $\vwq$ 
{\it the set of vertices, 
the weights on edges and % the edge weight (function) and 
the weights on vertices}, % the vertex weight (function), 
respectively.
% the set of vertices, the weight function on edges and the weight function on vertices, respectively.
We denote the set of vertices of $G$ (i.e. $V$) by $V(G)$, and
write $\whv(G)=V(G)\times\cbr{1,2}$.
%
%%%%% abuse of notation for weights with wider domains %%%%%
Given $V\subsetneq V'$, 
$\ewt:\hp(V')\to\bbZ_{\geq0}$ with $\supp(\ewt)\subset \hp(V)$ and 
$\vwq:\whv'\to\bbZ_{\geq0}$ with $\supp(\vwq)\subset\whv$,
we denote the weighted graph 
$(V,\ewt|_{\hp(V)}, \vwq|_{\whv})$ by $(V,\ewt,\vwq)$
for notational conciseness.
%
%%%%% concatenation of multiple weighted graphs %%%%%
For two weighted graphs $G=(V,\ewt,\vwq)$ and $G'=(V',\ewt',\vwq')$
such that $V\cap V'=\emptyset$,
we define a weighted graph 
$\widetilde{G}=(\extv,\tilde\ewt,\tilde\vwq)$ by
$\extv=V\sqcup V'$,
\begin{align*}
  \tilde\ewt|_{\hp(V)} = \ewt, 
  \tilde\ewt|_{\hp(V')} = \ewt',
  \tilde\ewt|_{\hp(V,V')} = 0\tand
  \tilde\vwq|_{\whv} = \vwq,
  \tilde\vwq|_{\whv'} = \vwq',
\end{align*}
and denote this weighted graph $\widetilde{G}$ by $G\vee G'$.
In general, we write 
\begin{align*}
  \mathop\vee_{k'\in\cbr{1,...,k}} G_{k'} = G_1 \vee ... \vee G_{k}
\end{align*}
for $k\in\bbN$ and weighted graphs $G_1,...,G_{k}$ whose sets of vertices are disjoint.

%
% \subsubsection*{projected unweighted graph (ちょっと違う)}
%%%%% projected edge weight %%%%%
For a nonempty finite set $V$ and 
$\ewt:\hp(V)\to\bbZ_{\geq0}$,
define the function $\cewt:p(V)\to\bbZ_{\geq0}$ by 
\begin{align*}
  \cewt([v,v']) = \sum_{\kap,\kapr\in\cbr{1,2}}\ewt([(v,\kap),(\vpr,\kapr)])
\end{align*}
for $[v,v']\in p(V)$.
We call this function $\cewt$ {\it the projected edge weight (function)}.
%
%%%%% Set of edges %%%%%
For another nonempty set $V'$ with $V'\subset V$,
we denote 
\begin{align*}
  E(V',\ewt) = \cbr{[v,\vpr]\in p(V')\mid \cewt([v,\vpr])>0}.
\end{align*}
This set consists of the pairs of $V'$ which has positive projected edge weight.
%
%%%%% projected graph %%%%%
Then for a weighted graph $G=(V,\ewt,\vwq)$, we define $E(G)=E(V,\ewt)$,
and we call 
$(V(G),E(G))(=(V,E(G)))$
{\it the projected graph} of the weighted graph $G$.
\begin{comment}
  {\myred projected graphに$\check{G}$という記号を与えては？}
\end{comment}
% 
% \subsubsection*{
% Connectedness, component, tree, spanning tree}
% treeはここでは定義されず，beta-chainの和の評価のところでtree-shaped w.g.として定義されることになった．
%%%%% connectedness %%%%%
If the projected graph $(V(G), E(G))$ of $G$ is a connected graph in a usual sense,
we say a weighted graph $G$ is {\it connected}.
%
%%%%% component %%%%%
Denoting the set of connected components of the projected graph $(V(G),E(G))$
of $G$ by $\check{\Comp}(G)$,
we define 
\begin{align*}
  \Comp(G) = \cbr{(V',\ewt|_{\hp(V')}, \vwq|_{\whv'})\mid 
  %\Comp(G) = \cbr{(V',\ewt\bbone_{\hp(V')}, \vwq\bbone_{\whv'})\mid 
  (V',E')\in\check\Comp(G)},
\end{align*}
and call an element of $\Comp(G)$ a {\it component} of the weighted graph $G$.

% 
% \subsubsection*{
%   Some quantity about weights on edges and vertices towards definition of exponent}
To define the exponent for weighted graphs, we introduce some summarizing quantities.
We start with those related to weights on vertices. 
Let $G=(V,\ewt,\vwq)$ be a weighted graph.
%
%%%%% Total weight on vertices $\barq(G)$, $\barq_1(G)$, $\barq_2(G)$ %%%%%
We denote
\begin{align*}
  \barq(G)=\sum_{\hv\in\whv(G)}\vwq(\hv),\tand
    % =\sum_{v\in V(G)}\sum_{\kap=1,2}\vwq(v,\kap)
  \barq_\kap(G)=\sum_{v\in V(G)}\vwq(v,\kap)
  \tfor \kap=1,2.
\end{align*}
Since 
$\barq(G)=\sum_{v\in V(G)}\sum_{\kap=1,2}\vwq(v,\kap)$,
it holds that 
$\barq(G)=\barq_1(G)+\barq_2(G)$.
%
%%%%% Classification of components %%%%%
% $\Comp_0(G)$, $\Comp_+(G)$, $\Comp_1(G)$, $\Comp_2(G)$;
With the above quantities, we classify the components of $G$ as follows:
\begin{align*}
  \Comp_0(G)&=\cbr{C\in\Comp(G)\mid \barq(C)=0},&
  \Comp_+(G)&=\cbr{C\in\Comp(G)\mid \barq(C)>0},
  \\*
  \Comp_1(G)&=\bcbr{C\in\Comp(G)\mid \barq_1(C)>0},&
  \Comp_2(G)&=\bcbr{C\in\Comp(G)\mid \barq_1(C)=0 \tandsm \barq_2(C)>0}.
\end{align*}
Notice that 
$\Comp_+(G) = \Comp_1(G)\sqcup\Comp_2(G)$ and 
$\Comp(G) = \Comp_+(G)\sqcup\Comp_0(G)$.
%
%%%%% Total weight on edges $\bartheta(G)$ %%%%%
As for edges, we define %for a weighted graph $G=(V,\ewt,\vwq)$
\begin{align*}
  \bartheta(G)=\sum_{[\hv,\hv']\in\hp(V(G))}\ewt([\hv,\hv']).
\end{align*}
Notice that we can rewrite $\bartheta(G)$ as 
\begin{align*}
  \bartheta(G)
  =\sum_{[v,v']\in p(V(G))}\sum_{\kap,\kap'=1,2}\ewt([(v,\kap),(v',\kap')])
  =\sum_{[v,v']\in p(V(G))}\check\ewt([v,v']).
\end{align*}
% {\myred $C$(connected)に対してdefすべきか？一般のgraphでもいいが，使わない気がする．}

To obtain a sharper estimate, we classify the edges of 
the projected graph of a weighted graph.
%
%%%%% `projected' edge weight $\cewt_1, \cewt_2$ %%%%%
For a nonempty finite set $V$ and $\ewt:\hp(V)\to\bbZ_{\geq0}$,
we define 
$\cewt_1,\cewt_2:p(V)\to\bbZ_{\geq0}$ by 
\begin{align*}
  \cewt_1([v,\vpr]) = \ewt([(v,1),(\vpr,1)])
  \tand
  \cewt_2([v,\vpr]) = 
  \sum_{\substack{\kap,\kapr\in\cbr{1,2}\\(\kap,\kapr)\neq(1,1)}}
    \ewt([(v,\kap),(\vpr,\kapr)])
\end{align*}
for $[v,\vpr]\in p(V)$.
Obviously,
$\cewt([v,\vpr]) = \cewt_1([v,\vpr]) + \cewt_2([v,\vpr])$.
%
%%%%% edges of `fast decay', edges of `heavy tail' $E_1(V,\ewt),E_2(V,\ewt)$ %%%%%
For another nonempty finite set $V'$ with $V'\subset V$,
we denote 
\begin{align*}
  E_1(V',\ewt) &= \cbr{[v,\vpr]\in p(V')\mid \cewt_1([v,\vpr])>0 \tandsm 
  \cewt_2([v,\vpr])=0}
  \\
  E_2(V',\ewt) &= \cbr{[v,\vpr]\in p(V')\mid \cewt_2([v,\vpr])>0}.
\end{align*}
Notice that 
$E(V',\ewt)=E_1(V',\ewt)\sqcup E_2(V',\ewt)$.
%
%%%%% Sets of edges of a weighted graph $E_1(G), E_2(G)$ %%%%%
For a weighted graph $G=(V,\ewt,\vwq)$, we define
\begin{align*}
  % E(G)=E(V,\ewt) \tand
  E_i(G)=E_i(V,\ewt) \tforsm i=1,2,
\end{align*}
and obviously we have 
$E(G)=E_1(G)\sqcup E_2(G)$.

%
% \subsubsection*{Exponent}
%
%%%%% spanning tree of a connected weighted graph %%%%%
For a connected weighted graph $C$,
we say that $\bbT\subset E(C)$ is a {\it spanning tree} of $C$,
if the subgraph $(V(C), \bbT)$ of the $(V(C), E(C))$ 
is a tree in a usual sense, that is a connected graph without cycles.
When $\abs{V(C)}=1$ and hence $p(V(C))=\emptyset$, 
we regard $\emptyset$ as the only spanning tree of $C$.
%
%%%%% exponent for connected weighted graphs, $\ell_2(C)$ %%%%%
Then 
the exponent $e(C)$ of a connected weighted graph $C$ is defined by
\begin{align}\label{def:231005.1440}
  e(C) = \et(C) + \eq(C)
\end{align}
with
\begin{align*}
  \eq(C) &=
  \begin{cases}
    -\half -H\barq(C)
    % \rbr{-H-\half} -H(\barq(C)-1) {\myred(=-\half -H\barq(C))}
    &\tifsm \barq_2(C)>0 \tandsm \barq_1(C)=0 
    \\
    -1 -H(\barq(C)-1)               
    &\tifsm \barq_1(C)>0 
    \\
    0 
    & \tifsm \barq(C)=0
  \end{cases}
  \\
  \et(C) &=
  % 1 - (\abs{V(C)}-1) + (1-2H) \ell_2(C) 
  % - 2H (\bartheta(C)-(\abs{V(C)}-1))
  % \\&=
  1 - 2H\bartheta(C) + (2H-1) \brbr{\abs{V(C)}-1-\ell_2(C)},
\end{align*}
and 
\begin{align*}
  \ell_2(C) = \max_\bbT\abs{E_2(C)\cap \bbT}
\end{align*}
where the maximum is taken over the set of all spanning trees of $C$. 
Notice that 
when $\abs{V(C)}=1$,
% for a connected weighted graph $C$ such that $\abs{V(C)}=1$,
it holds that $e_\theta(C)=1$, since $\bartheta(C)=0$ and $\ell_2(C)=0$.
%
%%%%% Exponent for a weighted graph %%%%%
The exponent of a weighted graph $G$ is defined as 
the sum of the exponent of its components, namely 
\begin{align*}
  e(G)=\sum_{C\in\Comp(G)} e(C).
\end{align*}

The estimates of the next lemma is basic in the following arguments.
\begin{lemma}
(i) For any connected graph $C$,
  \begin{align}
    \eq(C)\leq-H\barq(C). \label{221206.1156}  
  \end{align}
(ii) For any connected graph $C$ such that $\barq(C)>0$,
  \begin{align}
    \eq(C)&\geq-\half-H\barq(C), \label{221206.1202}
    \\
    \eq(C)&\leq-1-H(\barq(C)-1). \label{eq:230818.1803}
  \end{align}
\end{lemma}

% 
% \subsubsection*{Realization of functional}
Now we introduce 
% the functionals represented by aforementioned weighted graphs 
the correspondence between the aforementioned weighted graphs and functionals.
%
%%%%% \beta %%%%%
Suppose that 
$V$, $V_\bff$ and $V_\ewt$ are nonempty finite sets
satisfying $V\subset V_\bff, V_\ewt$.
Let
$\bff=(\bff^{(\hv)})_{\hv\in\wtv_\bff}\in\calh^{\wtv_\bff}$.
% $\bff\in\abs{\calh}^{\wtv_\bff}$
For 
$\ewt:\hp(V_\ewt)\to\bbZ_{\geq0}$ such that 
$\supp(\ewt)\subset\hp(V)$,
we denote
\begin{align*}
  \beta_{V}(\bff, \ewt) = 
  \prod_{[\tv, \tv']\in\hp(V)} \babr{\bff^{(\tv)}, \bff^{(\tv')}}^{\ewt[\tv,\tv']}.
\end{align*}
If $\abs{V}=1$, and therefore $\hp(V)=\emptyset$, we define 
$\beta_{V}(\bff, \ewt)=1$.
% Notice that $\beta_{V}(\bff, \ewt)=1$
% when $\sum_{[\tv,\tv']\in\tp(V)}\ewt[\tv,\tv']=0$.
%
%%%%% \delta %%%%%
For a nonempty finite set $V_\vwq$ satisfying $V\subset V_\vwq$ and 
$\vwq:\whv_\vwq\to\bbZ_{\geq0}$ such that $\supp(\vwq)\subset\whv$,
we denote 
\begin{align*}
  \delta_V(\bff, \vwq)
  = 
  \delta^{\bar\vwq}
  \Brbr{
    \subotimes{\hv\in\whv}
    \brbr{\nrbr{\bff^{(\hv)}}^{\otimes \vwq(\hv)}}},
\end{align*}
which is a $\bar\vwq$-th Skorohod integral with
$\bar\vwq = \sum_{\hv\in\whv}\vwq(\hv)$.
When $\bar\vwq=0$ (i.e. $\vwq=0$), 
we define $\delta_V(\bff, \vwq)=1$.
%
%%%%% B %%%%%
For a connected weighted graph $C=(V, \ewt, \vwq)$,
% $\bff\in{\calh}^{\wtv_\bff}$ with some $V_\bff(\supset V)$,
we define 
\begin{align*}
  \beta(C, \bff) = \beta_{V}(\bff, \ewt),\qquad
  \delta(C,\bff) = \delta_{V}(\bff, \vwq)\tand
  B(C,\bff) = \beta(C, \bff)\;\delta(C, \bff).
\end{align*}

%
%%%%% functionals %%%%%
%
%%%%% $i:V\to\cbr{1,2},j\in\bbJ_n(V,i)$;  %%%%%
For nonempty finite sets $V$ and $V'$ with $V\subset V'$
and $i:V'\to\cbr{1,2}$,
we write 
$\bbJ_n(V,i)=\prod_{v\in V} [i(v)n-1]$.
%
%%%%% $T_{m,j}, T^{(1)}_{m,j}, T^{(2)}_{m,j}$: 関数のtranslation %%%%%
For $f:\bbR\to\bbR$, $m\in\bbN$ and $j\in\bbZ$, 
we define a function $T_{m,j}(f):\bbR\to\bbR$ by 
\begin{align*}
  T_{m,j}(f)(x) = f(mx-j),
\end{align*}
and we denote 
\begin{align*}
  T_{m,j}^{(1)}(f) = T_{m,j}(f)\tand
  T_{m,j}^{(2)}(f) = T_{m,j+1}(f) - T_{m,j}(f).
\end{align*}
% \begin{itembox}[l]{Memo}
%   If $\supp(f)\subset[0,1]$, then
%   $\supp(T_{m,j}(f))\subset[\frac{j}{m},\frac{j+1}m]$.
% \end{itembox}

%
% \item $\calf(V), \bbf, \cala(V,i), \cali_n^{(i)}$
%%%%% \calf(V) %%%%%
For a nonempty finite set $V$,
consider the following set of conditions for 
$\bbf\in\rbr{L^\infty(\bbR)}^{\whv}$:
\begin{itemize}
  \item For $v\in V$, 
  the support of the function $\bbf(v,1)\in L^\infty(\bbR)$ 
  is included in $[a,a+1]$ with some $a\in[-1,0]$. 

  \item For $v\in V$, 
  the support of the function $\bbf(v,2)\in L^\infty(\bbR)$ 
  is included in $[-1,0]$. 
\end{itemize}
Denote by $\calf(V)$ the subset of 
$\rbr{L^\infty(\bbR)}^{\whv}$ satisfying the above conditions.
%
%%%%% \bbf %%%%%
For $\bbf\in\calf(V)$,
$n\in\ntwo$, $v\in V$, $\kappa=1,2$, $i:V\to\cbr{1,2}$ and $j\in\bbJ_n(V,i)$,
we define
\begin{align*}
  \bbf^{(v,\kappa)}_{n,j} = \bbf^{(v,\kappa)}_{n,i,j}
    &= T_{i(v)n,j_v}^{(\kappa)}(\bbf^{(v,\kappa)}),
  %\bbf^{(v,2)}_{n,j} &= T_{i(v)n,j_v}^{(2)}(\bbf^{(v,2)})=
  %T_{i(v)n,j_v+1}(\bbf^{(v,2)}) - T_{i(v)n,j_v}(\bbf^{(v,2)}).
\end{align*}
where
we recall $\ntwo=\cbr{j\in\bbZ\mid j\geq2}$ and
we write $\bbf^{(v,\kap)}:=\bbf(v,\kap)$ for notational convenience.
Observe that
the support of the function $\bbf^{(v,\kap)}_{n,j}\in L^\infty(\bbR)$ is included in $[0,1]$
thanks to the conditions imposed on $\calf(V)$,
and we can naturally consider $\bbf^{(v,\kap)}_{n,j}$ as an element of $\abs\calh$.
We write
$\bbf_{n,j}=\brbr{\bbf^{(v,\kap)}_{n,j}}_{(v,\kap)\in\wtv}
\brbr{\in\abs{\calh}^{\wtv}}$.

%
%%%%% $\cala(V,i)$ %%%%%
Let $V$ be a nonempty finite set and $i:V\to\cbr{1,2}$.
Writing $\bbJ(V,i) = \bigsqcup_{n\geq2} \cbr{n}\times\bbJ_n(V,i)$,
we consider the following conditions
for $A=(A_{n,j})_{(n,j)\in\bbJ(V,i)}\in(\bbD^\infty)^{\bbJ(V,i)}$:
\begin{itemize}
  \item For any $p\geq1$,
  $\sup_{n\geq2, j\in\bbJ_n(V,i)}\norm{A_{n,j}}_{L^p(P)}<\infty$.
  % \item $\sup_{n\geq2, j\in\bbJ_n(V,i)}\norm{A_{n,j}}_{L^p(P)}<\infty$ for any $p\geq1$.
  % \item For $k\in\bbN$ and $j\in\bbJ_n(V,i)$, 
  % $D^kA_{n,j}$ is an $\abs{\calh}^{\otimes k}$-valued random variable.
  \item For $k\geq1$, 
  $D^kA_{n,j}$ is $\abs{\calh}^{\otimes k}$-valued,
  and represented by
  $(D^k_{s_1,...,s_k}A_{n,j})_{s_1,...,s_k\in[0,1]}$ 
  such that
  \begin{align*}
    \sup_{n\geq2, j\in\bbJ_n(V,i)}
    \sup_{s_1,...,s_k\in[0,1]}\norm{D^k_{s_1,...,s_k}A_{n,j}}_{L^p(P)}
    <\infty
    \tfor p\geq1.
  \end{align*}
\end{itemize}
We denote by $\cala(V,i)$ the subset of $(\bbD^\infty)^{\bbJ(V,i)}$ 
satisfying the above conditions.

%
%%%%% $\cali_n^{(i)}(G,A,\bbf)$ %%%%%
Having prepared the definitions above, we define the functional to consider in the following arguments.
For $n\in\ntwo$, a weighted graph $G=(V,\ewt, \vwq)$, $i:V\to\cbr{1,2}$, 
$A\in\cala(V,i)$ and $\bbf\in\calf(V)$,
we define the functional $\calii_n\rbr{G, A,\bbf}$ by
\begin{align*}
  \calii_n\rbr{G, A,\bbf} &=
  \sum_{j\in\bbJ_n(V,i)}
  A_{n,j} %A_n(j)
  \prod_{C\in\Comp(G)} B\rbr{C, \bbf_{n,j}}.
\end{align*}

%% file: subfiles/4-2-sec_exponent_main.tex
\subsection{Estimates of the expectation and norms of $\cali_n^{(i)}\rbr{G, A,\bbf}$}
\label{sec:231005.1622}
First, we will give an estimate of the expectation 
$E\bsbr{\calii_n\rbr{G, A,\bbf}}$ 
of the functional defined for 
a weighted graph $G=(V,\ewt, \vwq)$, 
$i:V\to[2]$,
$A\in\cala(V,i)$ and 
$\bbf\in\calf(V)$ (Proposition \ref{prop:230605.2133}).
Lemma \ref{221207.2648} shows a decomposition of 
$E\bsbr{\calii_n\rbr{G, A,\bbf}}$ 
to the sum of the expectations of functionals, 
whose corresponding weighted graphs have 
fewer weights on vertices, and exponents less than or equal to that of $G$.
\begin{lemma}\label{221207.2648}
  For $n\in\ntwo$, 
  $G=(V,\ewt, \vwq)$ with $\barq(G)>0$, $i:V\to[2]$,
  $A\in\cala(V,i)$ and $\bbf\in\calf(V)$, let
  \begin{align*}
    \cali_n =
    \calii_n\rbr{G, A,\bbf} &=
    \sum_{j\in\bbJ_n(V,i)}
    A_{n,j}
    \prod_{C\in\Comp(G)} B\rbr{C, \bbf_{n,j}}.
    %\prod_{C\in\Comp(G)} \check B_n^C\rbr{j; \bbf_{V(C)}}.
  \end{align*}
  Then, with 
  %some $(v_0,\kap_0)\in\wtv$ and 
  some finite set $\Gamma$,
  there exist 
  $\alpha^{(\gamma)}\in\bbR$, 
  a weighted graph $G^{(\gamma)}=(V,\ewt^{(\gamma)}, \vwq^{(\gamma)})$, 
  and
  $A^{(\gamma)}\in\cala(V,i)$
  %$\bbf_\gamma\in\calf(\vwq_\gamma)$
  for each $\gamma\in\Gamma$
  satisfying the following conditions:
  \begin{itemize}
    \item $\barq(G^{(\gamma)}) \leq \barq(G)-1$ for all $\gamma\in\Gamma$
    \item $\alpha^{(\gamma)} + e(G^{(\gamma)})\leq e(G)$ for all $\gamma\in\Gamma$
    \item The expectation $E\sbr{\cali_n}$ can be decomposed as follows:
    \begin{align*}
      E\sbr{\cali_n} &=
      \sum_{\gamma\in\Gamma}
      E\sbr{n^{\alpha^{(\gamma)}}\calii_n\rbr{G^{(\gamma)}, A^{(\gamma)},\bbf}}
    \end{align*}
  \end{itemize}

\end{lemma}

\begin{proof}
  Until we specify how to choose $\tv_0=(v_0,\kap_0)\in\wtv$ afterwards, 
  fix $\tv_0\in\wtv$ satisfying $\vwq(\tv_0)>0$.
  We write $C_0\in\Comp(G)$ such that $v_0\in V(C_0)$.
  For notational convenience, we denote $V_0=V(C_0)$ and we can write 
  $C_0 = (V_0, \ewt\bbone_{\tp(V_0)}, \vwq\bbone_{\wtv_0})$.
  By the duality of $\delta$ and $D$, 
  for $j\in\bbJ_n(V,i)$
  %\newpage
  \begin{align*}
    &\hspsm 
    E\Bsbr{
      A_{n,j}
      \prod_{\substack{C\in\Comp(G)}} B\rbr{C, \bbf_{n,j}}
    }
    \\&=
    E\bbsbr{
      A_{n,j}\;
      %\beta_{C_0}\rbr{\bbf_{n,j}}
      \beta\rbr{C_0, \bbf_{n,j}}
      %\beta^{C_0}\rbr{\bbf_{n,j}}
      %\delta\rbr{\bbf_{n,j},\vwq\bbone_{V(C_0)}}
      %\delta\rbr{\bbf_{n,j}^{\otimes\vwq\bbone_{\wtv(C_0)}}}
      %\delta\rbr{\bbf_{n,j}^{\otimes\vwq\bbone_{\wtv_0}}}
      \delta_{V_0}\rbr{\bbf_{n,j},\vwq\bbone_{\wtv_0}}
      \prod_{\substack{C\in\Comp(G)\\C\neq C_0}} B\rbr{C, \bbf_{n,j}}
    }
    \\&=E\bbsbr{
      \bbabr{D\bbrbr{
        A_{n,j}
        \prod_{\substack{C\in\Comp(G)\\C\neq C_0}} B\rbr{C, \bbf_{n,j}}},
        \bbf_{n,j}^{(\tv_0)}
      }\;
      \beta\rbr{C_0, \bbf_{n,j}}\;
      \delta_{V_0}\rbr{\bbf_{n,j}, \vwq\bbone_{\wtv_0}-\bbone_\cbr{\tv_0}}
    }
    \\&=
    E\bbsbr{
      \abr{DA_{n,j},\bbf_{n,j}^{(\tv_0)}}\;
      \beta\rbr{C_0, \bbf_{n,j}}\;
      \delta_{V_0}\rbr{\bbf_{n,j}, \vwq\bbone_{\wtv_0}-\bbone_\cbr{\tv_0}}\;
      %\delta\rbr{\bbf_{n,j}^{\otimes\vwq\bbone_{\wtv(C_0)}-\bbone_\cbr{\tv_0}}}
      \prod_{\substack{C\in\Comp(G)\\C\neq C_0}} B\rbr{C, \bbf_{n,j}}
    }
    \\&\hspsm+
    \sum_{\substack{C_1\in\Comp_+(G)\\C_1\neq C_0}}
    E\bbsbr{
      A_{n,j}
      \abr{DB\rbr{C_1, \bbf_{n,j}},\bbf_{n,j}^{(\tv_0)}}\;
      \beta\rbr{C_0, \bbf_{n,j}}\;
      \delta_{V_0}\rbr{\bbf_{n,j}, \vwq\bbone_{\wtv_0}-\bbone_\cbr{\tv_0}}\;
      %\delta\rbr{\bbf_{n,j}^{\otimes\vwq\bbone_{\wtv(C_0)}-\bbone_\cbr{\tv_0}}}
      \prod_{\substack{C\in\Comp(G)\\C\neq C_0, C_1}} B\rbr{C, \bbf_{n,j}}
    }
  \end{align*}
  Again for notational convenience, we denote $V_1=V(C_1)$ and we can write 
  $C_1 = (V_1, \ewt\bbone_{\tp(V_1)}, \vwq\bbone_{\wtv_1})$.
  Then the factor 
  $\babr{DB\rbr{C_1, \bbf_{n,j}},\bbf_{n,j}^{(\tv_0)}}$ 
  can be decomposed as
  \begin{align*}
    \abr{DB\rbr{C_1,\bbf_{n,j}},\bbf_{n,j}^{(\tv_0)}}
    &=
    \beta\rbr{C_1,\bbf_{n,j}}
    \abr{D\delta_{V_1}\rbr{\bbf_{n,j}, \vwq\bbone_{\wtv_1}},\bbf_{n,j}^{(\tv_0)}}
    %\abr{D\delta\rbr{\bbf_{n,j}^{\otimes\vwq\bbone_{\wtv(C_1)}}},\bbf_{n,j}^{(\tv_0)}}
    \\&=
    \beta\rbr{C_1,\bbf_{n,j}}
    \sum_{\substack{\tv_1\in \wtv_1\\\vwq(\tv_1)>0}}
    %\sum_{\substack{\tv_1\in \wtv(C_1)\\\vwq(\tv_1)>0}}
    \vwq(\tv_1)\;
    \delta_{V_1}\rbr{\bbf_{n,j}, \vwq\bbone_{\wtv_1}-\bbone_\cbr{\tv_1}}
    %\delta\rbr{\bbf_{n,j}^{\otimes\vwq\bbone_{\wtv(C_1)}-\bbone_\cbr{\tv_1}}}
    \abr{\bbf_{n,j}^{(\tv_1)}, \bbf_{n,j}^{(\tv_0)}}
    \\&=
    \beta\rbr{C_1,\bbf_{n,j}}
    \sum_{\substack{\tv_1\in \wtv_1\\\vwq(\tv_1)>0}}
    \vwq(\tv_1)\;
    \delta_{V_1}\rbr{\bbf_{n,j}, \vwq\bbone_{\wtv_1}-\bbone_\cbr{\tv_1}}
    \beta_{V_1\sqcup V_0}(\bbf_{n,j}, \bbone_{[\tv_1, \tv_0]})
    %\abr{\bbf_{n,j}^{(\tv_1)}, \bbf_{n,j}^{(\tv_0)}}
  \end{align*}

  Let us write
  $\vwq'_0:=\vwq\bbone_{\wtv_0}-\bbone_{\cbr{\tv_0}}$,
  $\vwq'_1:=\vwq\bbone_{\wtv_1}-\bbone_{\cbr{\tv_1}}$ and
  $\vwq':=\vwq'_0+\vwq'_1
  =\vwq\bbone_{\wtv_0\sqcup\wtv_1}-\bbone_{\cbr{\tv_0,\tv_1}}$.
  By Lemma \ref{lemma:230526.1230}, the product of the two Skorohod integrals is decomposed as
  \begin{align*}
    &\hspsm
    \delta_{V_1}\rbr{\bbf_{n,j},\vwq\bbone_{\wtv_1}-\bbone_\cbr{\tv_1}}
    %\delta\rbr{\bbf_{n,j}^{\otimes\vwq\bbone_{\wtv(C_1)}-\bbone_\cbr{\tv_1}}}
    \delta_{V_0}\rbr{\bbf_{n,j},\vwq\bbone_{\wtv_0}-\bbone_\cbr{\tv_0}}
    %\delta\rbr{\bbf_{n,j}^{\otimes\vwq\bbone_{\wtv(C_0)}-\bbone_\cbr{\tv_0}}}
    = 
    \delta_{V_1}\rbr{\bbf_{n,j},\vwq'_1}
    \delta_{V_0}\rbr{\bbf_{n,j},\vwq'_0}
    \\&=
    \sum_{\bfp\in\Pi_{V_0, V_1}(\vwq'_0, \vwq'_1)}
    c(\bfp)~
    \delta_{V_1\sqcup V_0}\rbr{\bbf_{n,j}, \vwq'-\bar \bfp}\;
    %\delta\rbr{\bbf_{n,j}^{\otimes\vwq'-\bar \bfp}}
    \beta_{V_1\sqcup V_0}(\bbf_{n,j}, \bfp)
  \end{align*}
  with
  \begin{align*}
    \Pi_{V_0, V_1}(\bfq'_0, \bfq'_1) =
    \Big\{\bfp:\tp(V_0\sqcup V_1)\to\bbZ_{\geq0}\mid
    \supp(\bfp)\subset\tp(V_0,V_1), %&\;
    \bar\bfp(\tv)\leq\bfq'_k(\tv)\tforsm\tv\in\wtv_k \tforsm k=0,1 \Big\},
  \end{align*}
  where we define
  $\bar\bfp:\whv_0\sqcup\whv_1\to\bbZ_{\geq0}$ by
  $\bar\bfp(\tv_0) = \sum_{\tv_1\in\wtv_1}\bfp[\tv_0,\tv_1]$ for $\tv_0\in\wtv_0$ and
  $\bar\bfp(\tv_1) = \sum_{\tv_0\in\wtv_0}\bfp[\tv_0,\tv_1]$ for $\tv_1\in\wtv_1$, and
  %$\bar\bfp(\tv) = 0$ for $\tv\not\in\wtv_0\sqcup\wtv_1$.
  $c(\bfp)\in\bbZ_{>0}$ are constants.
  Then we have
  \begin{align*}
    &\hspsm
    \abr{DB\rbr{C_1,\bbf_{n,j}},\bbf_{n,j}^{(\tv_0)}}\;
    \beta\rbr{C_0, \bbf_{n,j}}\;
    \delta_{V_0}\rbr{\bbf_{n,j},\vwq\bbone_{\wtv_0}-\bbone_\cbr{\tv_0}}
    %\delta\rbr{\bbf_{n,j}^{\otimes\vwq\bbone_{\wtv(C_0)}-\bbone_{\tv_0}}}
    \\&=
    \sum_{\substack{\tv_1\in \wtv_1\\\vwq(\tv_1)>0}}
    \sum_{\bfp\in\Pi_{V_0, V_1}(\vwq'_0, \vwq'_1)}
    %\sum_{\substack{\tv_1\in \wtv_1\text{ s.t. }\vwq(\tv_1)>0
    %\\\bfp\in\Pi_{V_0, V_1}(\vwq'_0, \vwq'_1)}}
    c(\tv_1,\bfp)\;
    \beta\rbr{C_0, \bbf_{n,j}}
    \beta\rbr{C_1,\bbf_{n,j}}
    \beta_{V_1\sqcup V_0}(\bbf_{n,j}, \bfp
    +\bbone_{\sbr{\tv_1,\tv_0}})\;
    \delta_{V_1\sqcup V_0}\rbr{\bbf_{n,j}, \vwq'-\bar \bfp}
    %\delta\rbr{\bbf_{n,j}^{\otimes\vwq'-\bar \bfp}}
    \\&=
    \sum_{\substack{\tv_1\in \wtv_1\\\vwq(\tv_1)>0}}
    \sum_{\bfp\in\Pi_{V_0, V_1}(\vwq'_0, \vwq'_1)}
    %\sum_{\substack{\tv_1\in \wtv_1\text{ s.t. }\vwq(\tv_1)>0
    %\\\bfp\in\Pi_{V_0, V_1}(\vwq'_0, \vwq'_1)}}
    c(\tv_1,\bfp)\;
    B\brbr{C^{(\tv_1,\bfp)}, \bbf_{n,j}},
  \end{align*}
  where we define the weighted graph
  $C^{(\tv_1,\bfp)}$ by 
  \begin{align*}
    C^{(\tv_1,\bfp)} &= 
    %(V^{(v,\kap,\bfp)}, \ewt^{(v,\kap,\bfp)}, \vwq^{(v,\kap,\bfp)}) 
    %\\&=
    (V_0\sqcup V_1,\;
    \ewt\bbone_{\tp(V_0)\sqcup \tp(V_1)}
    +\bfp+\bbone_{\sbr{\tv_1,\tv_0}},\;
    \vwq\bbone_{\wtv_0\sqcup \wtv_1}-\bbone_{\cbr{\tv_1,\tv_0}}
    -\bar \bfp
    )
  \end{align*}
  and $c(\tv_1,\bfp)$ is a positive integer-valued constant.
  Hence we obtain 
  \begin{align}
    E\sbr{\cali_n} &=
    \sum_{j\in\bbJ_n(V,i)} 
    E\Bsbr{A_{n,j} \prod_{\substack{C\in\Comp(G)}} B\rbr{C, \bbf_{n,j}}}
    =%\\&=
    E\bsbr{\cali_n^{(0)}} +
    \sum_{\substack{C_1\in\Comp_+(G)\\C_1\neq C_0}}
    \sum_{\substack{\hv_1\in \wtv_1\\\vwq(\hv_1)>0}}
    \sum_{\bfp\in\Pi_{V_0, V_1}(\vwq'_0, \vwq'_1)}
    E\bsbr{\cali_n^{(\hv_1,\bfp)}}
    \label{eq:230816.2300}
  \end{align}
  with 
  \begin{align*}
    \cali_n^{(0)} &=
    \sum_{j\in\bbJ_n(V,i)} 
    \babr{DA_{n,j},\bbf_{n,j}^{(\tv_0)}}\;
    \beta\brbr{C_0, \bbf_{n,j}}\;
    \delta_{V_0}\brbr{\bbf_{n,j},\vwq\bbone_{\wtv_0}-\bbone_\cbr{\tv_0}}
    %\delta\rbr{\bbf_{n,j}^{\otimes\vwq\bbone_{\wtv(C_0)}-\bbone_\cbr{\tv_0}}}
    \prod_{\substack{C\in\Comp(G)\\C\neq C_0}} B\rbr{C, \bbf_{n,j}}
    \\
    \cali_n^{(\hv_1,\bfp)} &=
    \sum_{j\in\bbJ_n(V,i)} 
    c(\hv_1,\bfp)\,A_{n,j}\;
    B\brbr{C^{(\hv_1,\bfp)}, \bbf_{n,j}}
    \prod_{\substack{C\in\Comp(G)\\C\neq C_0, C_1}} B\rbr{C, \bbf_{n,j}}
  \end{align*}
  We set
  \begin{align*}
    \Gamma_1 = 
    \bigsqcup_{\substack{C_1\in\Comp_+(G)\\C_1\neq C_0}}
    \bigsqcup_{\substack{\tv_1\in \wtv_1\\\vwq(\tv_1)>0}}
    \cbr{(\tv_1,\bfp)\mid \bfp\in\Pi_{V_0, V_1}(\vwq'_0, \vwq'_1)}
  \end{align*}
  and 
  $\Gamma = \cbr{0} \sqcup\Gamma_1$.

  \vspsm
  We let $\alpha^{(0)}=-1$ if $\kap_0=1$ and $\alpha^{(0)}=-2H$ if $\kap_0=2$, 
  where we recall that we write $\hv_0=(v_0,\kappa_0)$.
  We define 
  \begin{align*}
    A^{(0)} 
    &= \brbr{A^{(0)}_{n,j}}_{j\in\bbJ_n(V,i),n\in\ntwo}
    = \Brbr{n^{-\alpha^{(0)}}\babr{DA_{n,j},\bbf_{n,j}^{(\tv_0)}}_{\calh}}_{j\in\bbJ_n(V,i),n\in\ntwo}.
  \end{align*}
  Notice that $A^{(0)}\in\cala(V,i)$ in both cases $\kappa_0=1,2$ by Lemma \ref{lemma:230616.1754}.
  We define a connected weighted graph $C^{(0)}$ and a weighted graph $G^{(0)}$ by
  \begin{align*}
    C^{(0)} = 
    (V_0, 
    \ewt\bbone_{\tp(V_0)},
    \vwq\bbone_{\wtv_0}-\bbone_\cbr{\tv_0})\tand
    G^{(0)} = 
    C^{(0)} \vee \Brbr{\mathop\vee_{C\in \Comp(G)\setminus\cbr{C_0}} C}
  \end{align*}
  Since 
  $\beta(C^{(0)}, \bbf_{n,j}) = \beta(C_0, \bbf_{n,j})$,
  we have 
  $B(C^{(0)}, \bbf_{n,j}) = 
  \beta(C_0, \bbf_{n,j}) 
  \delta_{V_0}\brbr{\bbf_{n,j},\vwq\bbone_{\wtv_0}-\bbone_\cbr{\tv_0}}$,
  and hence 
  \begin{align*}
    \cali_n^{(0)} =
    n^{\alpha^{(0)}}\sum_{j\in\bbJ_n(V,i)} A^{(0)}_{n,j}\;
    \prod_{\substack{C\in(\Comp(G)\setminus\cbr{C_0})\sqcup\ncbr{C^{(0)}}}} B\rbr{C, \bbf_{n,j}}
    = n^{\alpha^{(0)}}\; \cali_n^{(i)}(G^{(0)}, A^{(0)}, \bbf).
  \end{align*}
  Notice that $\barq(G^{(0)})=\barq(G)-1$.
  
  \vspsm\noindent
  For $(\hv_1,\bfp)\in\Gamma_1$, we let
  $\alpha^{(\hv_1,\bfp)}=0$, 
  $A^{(\hv_1,\bfp)}=(c(\hv_1,\bfp)A_{n,j})_{j\in\bbJ_n(V,i), n\in\ntwo}$ and 
  \begin{align*}
    G^{(\hv_1,\bfp)} &=
    %(V, \ewt^{(\hv_1,\bfp)}, \vwq^{(\hv_1,\bfp)})=
    C^{(\hv_1,\bfp)} \vee \Brbr{\mathop\vee_{C\in \Comp(G)\setminus\cbr{C_0, C_1}} C}.
  \end{align*}
  % \blub{
  % We may write $ G^{(\hv_1,\bfp)} = (V, \ewt^{(\hv_1,\bfp)}, \vwq^{(\hv_1,\bfp)})$.
  % :使っていない記号のようなので，後ろで別の意味で使うために消す予定.}
  Then obviously we have 
  \begin{align*}
    \cali_n^{(\hv_1,\bfp)} = 
    n^{\alpha^{(\hv_1,\bfp)}} \cali_n^{(i)}(G^{(\hv_1,\bfp)}, A^{(\hv_1,\bfp)}, \bbf),
  \end{align*}
  % Setting $\Gamma = \cbr{0} \sqcup \Gamma_1$,
  and by \eqref{eq:230816.2300}
  we can write 
  \begin{align*}
    E\sbr{\cali_n} &=
    \sum_{\gamma\in\Gamma}
    E\sbr{n^{\alpha^{(\gamma)}}\calii_n\rbr{G^{(\gamma)}, A^{(\gamma)},\bbf}}.
  \end{align*}

  We prepare some estimates of the exponent of $C^{(\gamma)}$ for 
  $\gamma=(\hv_1,\bfp)\in\Gamma_1$. 
  Notice that 
  \begin{align}
    \bartheta(C^{(\gamma)}) &=
    \bartheta(C_0)+\bartheta(C_1)+1+\totpi
    \label{eq:230526.1607}
    \\
    \barq(C^{(\gamma)}) &= 
    \barq(C_0) + \barq(C_1) -2 - 2 \totpi,
    \label{221206.2046}
  \end{align}
  where we define
  $\totpi =\sum_{[\tv,\tv']\in\tp(V_0, V_1)} \bfp[\tv,\tv']
  = \half\sum_{\tv\in\wtv_0\sqcup \wtv_1} \bar\bfp(\tv)$.
  By $\totpi\geq0$ and \eqref{221206.2046},
  we can verify that 
  $\barq(G^{(\gamma)}) = \barq(G) -2 - 2 \totpi
  % \leq \barq(G) -2 
  < \barq(G) -1$.
  By the inequality \eqref{221206.1156},
  \begin{align*}
    e(C^{(\gamma)}) 
    % &=\et(C^{(\gamma)}) + \eq(C^{(\gamma)}) 
    &\leq%\\&\leq%
    \et(C^{(\gamma)}) -H \barq(C^{(\gamma)})
    && (\text{by }\eqref{221206.1156})
    \\&=
    1 + (2H-1) \brbr{\nabs{V(C^{(\gamma)})}-1-\ell_2(C^{(\gamma)})}
    -2H\bartheta(C^{(\gamma)})
    -H\barq(C^{(\gamma)})
    \\&=
    1 + (2H-1) \brbr{\nabs{V_0} + \nabs{V_1}-1-\ell_2(C^{(\gamma)})}
    \\&\hspsm
    -2H\rbr{\bartheta(C_0)+\bartheta(C_1)}
    -H \rbr{\barq(C_0)+\barq(C_1)}
    && (\text{by }\eqref{eq:230526.1607}\tandsm\eqref{221206.2046})
  \end{align*}
  We also have 
  $e(C_k) =
    1 + (2H-1) \brbr{\nabs{V_k}-1-\ell_2(C_k)}
    - 2H\bartheta(C_k) +\eq(C_k)$
  for $k=0,1$.
  Hence we have 
  \begin{align}
    &e(C^{(\gamma)}) - \rbr{e(C_0)+e(C_1)}
    %&e(C^{(\hv_1,\bfp)}) - \rbr{e(C_0)+e(C_1)}
    % \nn\\&\leq
    % 1 + (2H-1) \brbr{\nabs{V_0} + \nabs{V_1}-1-\ell_2(C^{(\gamma)})} 
    % %\\&\hspsm
    % -2H\rbr{\bartheta(C_0)+\bartheta(C_1)}
    % -H \rbr{\barq(C_0)+\barq(C_1)}
    % \nn\\&\hspsm 
    % -\cbr{1 + (2H-1) \brbr{\nabs{V_0}-1-\ell_2(C_0)} - 2H\bartheta(C_0) +\eq(C_0)}
    % \nn\\&\hspsm 
    % -\cbr{1 + (2H-1) \brbr{\nabs{V_1}-1-\ell_2(C_1)} - 2H\bartheta(C_1) +\eq(C_1)}
    % \nn\\&=
    \nn\\&\leq
    -1 + (2H-1) \brbr{1-\ell_2(C^{(\gamma)}) +\ell_2(C_0)+\ell_2(C_1)} 
    -H \rbr{\barq(C_0)+\barq(C_1)} - \eq(C_0) -\eq(C_1).
    \label{221206.1340}
  \end{align}

  Concerning to the estimate of $\ell_2$, we decompose $E_2(C^{(\gamma)})$ 
  for $\gamma\in\Gamma_1$.
  Recall that for $\ewt:\hp(V)\to\bbZ_{\geq0}$, we define 
  $\check\ewt_2([v,v']) = \sum_{\substack{\kap,\kap'=1,2\\(\kap,\kap')\neq(1,1)}}
  \ewt([(v,\kap), (v',\kap')])$ for $[v,v']\in p(V)$.
  We write $ C^{(\gamma)} = (V_0\sqcup V_1, \ewt^{(\gamma)}, \vwq^{(\gamma)})$.
  We have 
  \begin{align*}
    (\ewt^{(\gamma)}\check{)}_2 &= 
    (\ewt\bbone_{\hp(V_0)}\check{)}_2 + 
    (\ewt\bbone_{\hp(V_1)}\check{)}_2 +
    (\bbone_{[\hv_0,\hv_1]}\check{)}_2 +
    \check{\bfp}_2
    \\&=
    \check\ewt_2\bbone_{p(V_0)} + 
    \check\ewt_2\bbone_{p(V_1)} + 
    \bbone_{[v_0,v_1]}\bbone_{(\kappa_0,\kappa_1)\neq(1,1)} + 
    \check{\bfp}_2,
  \end{align*}
  where we write $(v_i,\kappa_i)=\hv_i$ for $i=0,1$.
  Notice that the supports of 
  $\check\ewt_2\bbone_{p(V_0)}$, $\check\ewt_2\bbone_{p(V_1)}$ and
  $\bbone_{[v_0,v_1]}\bbone_{(\kappa_0,\kappa_1)\neq(1,1)} + \check{\bfp}_2$
  are included in $p(V_0)$, $p(V_1)$ and $p(V_0,V_1)$, respectively.
  Hence,
  \begin{align}
    E_2(C^{(\gamma)}) &= 
    E_2(V_0\sqcup V_1, \ewt^{(\gamma)}) =
    \cbr{[v,v']\in p(V_0\sqcup V_1)\mid (\ewt^{(\gamma)}\check{)}_2 ([v,v'])>0}
    \nn\\&= 
    \cbr{[v,v']\in p(V_0)\mid \check{\ewt}_2 ([v,v'])>0}\sqcup
    \cbr{[v,v']\in p(V_1)\mid \check{\ewt}_2 ([v,v'])>0}
    \nn\\&\quad\sqcup
    \cbr{[v,v']\in p(V_0, V_1)\mid 
    \bbone_{[v_0,v_1]}([v,v'])\bbone_{(\kappa_0,\kappa_1)\neq(1,1)} + \check{\bfp}_2([v,v'])>0}
    %(\ewt^{(\gamma)}\check{)}_2 ([v,v'])>0}
    \nn\\&=
    E_2(C_0)\sqcup E_2(C_1)\sqcup E_2^{(\gamma)},\label{eq:230602.1932}
  \end{align}
  where we set 
  $E_2^{(\gamma)} = \cbr{[v,v']\in p(V_0, V_1)\mid 
  \bbone_{[v_0,v_1]}([v,v'])\bbone_{(\kappa_0,\kappa_1)\neq(1,1)} + \check{\bfp}_2([v,v'])>0}$.

  %\newpage
  \vspsm
  Since we have checked $\barq(G^{(\gamma)})\leq\barq(G)-1$ for any $\gamma\in\Gamma$,
  it suffices to show $\alpha^{(\gamma)} + e(G^{(\gamma)}) \leq e(G)$ for $\gamma\in\Gamma$.
  \item [(i)] First consider the case where $\Comp_2(G)\neq\emptyset$.
  Choose $C_0\in\Comp_2(G)$ and 
  $\tv_0=(v_0,\kappa_0)\in\wtv(C_0)$ such that $\vwq(\tv_0)>0$
  at the beginning of the proof.
  Here notice that $\kap_0=2$.
  We have $\eq(C_0) = -\half -H\barq(C_0)$.

  \item[(i-i)] Suppose that $\gamma=0$.
  We have
  $\et(C^{(0)}) = \et(C_0)$ and 
  $\eq(C^{(0)}) \leq -H \barq(C^{(0)}) = -H(\barq(C_0)-1)$,
  where we used the inequality \eqref{221206.1156}, and hence
  $e(C^{(0)})$ is bounded as
  \begin{align*}
    e(C^{(0)}) 
    % = \et(C^{(0)}) + \eq(C^{(0)}) 
    &\leq%\et(C_0)  -H \barq(C^{(0)}) =
    \et(C_0) -H (\barq(C_0)-1)
    % \\&
    =\et(C_0) + \eq(C_0) - (-H-\half)
    =e(C_0) +H+\half.
  \end{align*}
  Since $\kap_0=2$ and we have set $\alpha^{(0)}=-2H$, we obtain
  \begin{align*}
    \alpha^{(0)} + e(G^{(0)}) %&= 
    % -2H + e(C^{(0)}) + \sum_{C\in \Comp(G)\setminus\cbr{C_0}} e(C)
    &\leq%\\&\leq
    -2H
    +\rbr{H+\half} + e(C_0)  + \sum_{C\in \Comp(G)\setminus\cbr{C_0}} e(C)
    <%\leq
    \sum_{C\in \Comp(G)} e(C) = e(G).
  \end{align*}

  \item[(i-ii)] Suppose that $\gamma\in\Gamma_1$. (We write $\gamma=(\tv_1,\bfp)$ and $\tv_1=(v_1,\kappa_1)$.)
  We have 
  $\eq(C_0) = -\half -H\barq(C_0)$ from $C_0\in\Comp_2(G)$, 
  $\eq(C_1) \geq -\half -H\barq(C_1)$ by \eqref{221206.1202}, and hence
  \begin{align}
    &e(C^{(\gamma)}) - \rbr{e(C_0)+e(C_1)}
    % \nn\\&\leq
    % -1 + (2H-1) \brbr{1-\ell_2(C^{(\gamma)}) + \ell_2(C_0) + \ell_2(C_1)} 
    % -H \rbr{\barq(C_0)+\barq(C_1)} - \eq(C_0) -\eq(C_1)
    % \nn\\&\leq
    \leq
    (2H-1) \brbr{1-\ell_2(C^{(\gamma)}) + \ell_2(C_0) + \ell_2(C_1)}.
    \label{eq:230602.1953}
  \end{align}
  by \eqref{221206.1340}.
  We would like to show 
  $1 + \ell_2(C_0) + \ell_2(C_1)\leq\ell_2(C^{(\gamma)})$.

  For $k=0,1$, let 
  $\bbT_k\subset E(C_k)$ be a spanning tree of $C_k$
  satisfying $\eltw(C_k) = \abs{E_2(C_k)\cap\bbT_k}$.
  Then the subset $\bbT = \bbT_0\sqcup \bbT_1\sqcup\cbr{[v_0,v_1]}$ of 
  $p(V_0\sqcup V_1)=p(V(C^{(\gamma)}))$ 
  is a spanning tree of $C^{(\gamma)}$.
  Recall the decomposition \eqref{eq:230602.1932}
  $E_2(C^{(\gamma)})=E_2(C_0)\sqcup E_2(C_1)\sqcup E_2^{(\gamma)}$ and 
  notice that $E_2^{(\gamma)}\supset\cbr{[v_0,v_1]}$ since $\kap_0=2$.
  From the inclusions
  $E_2(C_k), \bbT_k\subset p(V_k)$ for $k=0,1$ and 
  $E_2^{(\gamma)}, \cbr{[v_0,v_1]} \subset p(V_0,V_1)$, we have
  \begin{align*}
    \abs{E_2(C^{(\gamma)})\cap\bbT}=
    \abs{E_2(C_0)\cap\bbT_0}+\abs{E_2(C_1)\cap\bbT_1} + \babs{E_2^{(\gamma)}\cap \cbr{[v_0,v_1]}}=
    \eltw(C_0) + \eltw(C_1) + 1
  \end{align*}
  and 
  $\eltw(C^{(\gamma)})\geq\abs{E_2(C^{(\gamma)})\cap\bbT}
  =\eltw(C_0) + \eltw(C_1) + 1$.
  By \eqref{eq:230602.1953}, we obtain
  $e(C^{(\gamma)})\leq e(C_0)+e(C_1)$
  and 
  $\alpha^{(\gamma)} + e(G^{(\gamma)}) 
  =0 + e(C^{(\gamma)}) + \sum_{C\in \Comp(G)\setminus\cbr{C_0, C_1}} e(C)
  =\sum_{C\in \Comp(G)} e(C)
  \leq e(G)$.

  \vspsm
  (ii) Consider the case where $\Comp_2(G)=\emptyset$.
  Choose $C_0\in\Comp_1(G)=\Comp_+(G)$ and 
  $\tv_0=(v_0,\kap_0)\in\wtv(C_0)$ such that $\vwq(\tv_0)>0$.

  \item[(ii-i)] The case $\gamma=0$. By a similar argument in (i-i) and $e_q(C_0)=-1-H(\barq(C_0)-1)$, we have 
  \begin{align*}
    e(C^{(0)}) 
    \leq\et(C_0) -H (\barq(C_0)-1)
    =\et(C_0) + \eq(C_0)+1
    =e(C_0)+1.
  \end{align*}
  Whether $\kap_0=1$ or $2$, we have $\alpha^{(0)}\leq-1$ and 
  \begin{align*}
    \alpha^{(0)} + e(G^{(0)}) &\leq
    -1 + e(C^{(0)}) + \sum_{C\in \Comp(G)\setminus\cbr{C_0}} e(C)
    \leq
    \sum_{C\in \Comp(G)} e(C) = e(G).
  \end{align*}

  \item[(ii-ii)] The case $\gamma\in\Gamma_1$. (We write $\gamma=(\tv_1,\bfp)$ and $\tv_1=(v_1,\kappa_1)$.)
  %(ii-ii) The case $\gamma=(v_1,\kap_1,\bfp)$.
  Again, for $k=0,1$, let $\bbT_k\subset E(C_k)$ be a spanning tree of $C_k$ satisfying 
  $\eltw(C_k)=\abs{E_2(C_k)\cap\bbT_k}$.
  Then the subset $\bbT=\bbT_0\sqcup \bbT_1\sqcup\cbr{[v_0,v_1]}$ of 
  $p(V_0\sqcup V_1)=p(V(C^{(\gamma)}))$
  is a spanning tree of $C^{(\gamma)}$.
  By the decomposition \eqref{eq:230602.1932} of $E_2(C^{(\gamma)})$,
  \begin{align*}
    \babs{E_2(C^{(\gamma)})\cap\bbT} =
    \abs{E_2(C_0)\cap\bbT_0}+\abs{E_2(C_1)\cap\bbT_1} + \babs{E_2^{(\gamma)}\cap \cbr{[v_0,v_1]}}\geq
    %\abs{E_2(C_0)\cap\bbT_0} + \abs{E_2(C_1)\cap\bbT_1}=
    \eltw(C_0) + \eltw(C_1).
  \end{align*}
  Hence we have 
  $\eltw(C^{(\gamma)})\geq\abs{E_2(C^{(\gamma)})\cap\bbT}\geq\eltw(C_0) + \eltw(C_1)$.

  Since we assumed $\Comp_2(G)=\emptyset$, we have 
  $\eq(C_i) = -1 -H(\barq(C_i)-1)$ for $i=0,1$.
  Therefore, we obtain using \eqref{221206.1340}
  \begin{align*}
    &e(C^{(\gamma)}) - \rbr{e(C_0)+e(C_1)}
    % \\&\leq
    % -1 + (2H-1) \brbr{1-\ell_2(C^{(\gamma)}) +\ell_2(C_0)+\ell_2(C_1)} 
    % -H \rbr{\barq(C_0)+\barq(C_1)} - \eq(C_0) -\eq(C_1)
    % \\&\leq 
    \leq
    -1 + 2H-1 + 2(1-H) = 0
  \end{align*}
  and 
  $\alpha^{(\gamma)} + e(G^{(\gamma)}) 
  =0 + e(C^{(\gamma)}) + \sum_{C\in \Comp(G)\setminus\cbr{C_0, C_1}} e(C)
  \leq e(G)$.
\end{proof}
By the decomposition of %the expectation 
$E\bsbr{\calii_n\rbr{G, A,\bbf}}$ 
obtained at the previous lemma, 
we can deduce an estimate of the order of the expectation.
\begin{proposition}\label{prop:230605.2133}
  For $n\in\ntwo$, %$\alpha\in\bbR$, 
  $G=(V,\ewt, \vwq)$, $i:V\to[2]$,
  $A\in\cala(V,i)$ and $\bbf\in\calf(V)$, let
  \begin{align*}
    %\sum_{j\in[n-1]^{V_1}\times [2n-1]^{V_2}}
    \cali_n =
    \cali_n^{(i)}\rbr{G, A,\bbf} &=
    \sum_{j\in\bbJ_n(V,i)} A_{n,j} \prod_{C\in\Comp(G)} B\rbr{C, \bbf_{n,j}}.
  \end{align*}
  Then, the expectation of $\cali_n$ is estimated as 
  \begin{align}
    E\sbr{\cali_n} = O(n^{e(G)}).
    \label{221207.0332}
  \end{align}
\end{proposition}

\begin{proof}
  We will prove the statement by induction with respect to $\barq(G)$.
  First we consider the case $\barq(G)=0$.
  In this case, we have 
  $B\rbr{C, \bbf_{n,j}}=\beta\rbr{C,\bbf_{n,j}}$.
  \begin{align}
    %\sum_{j\in[n-1]^{V_1}\times [2n-1]^{V_2}}
    \abs{E\sbr{\cali_n}} 
    % &\leq
    % \sum_{j\in\bbJ_n(V,i)}
    % \abs{E\sbr{A_{n,j}}}
    % \prod_{C\in\Comp(G)} \abs{\beta\rbr{C,\bbf_{n,j}}}
    % \nn\\
    &\leq
    \sup_{n\in\ntwo, j\in\bbJ_n(V,i)} 
    \norm{A_{n,j}}_{L^1} % \abs{E\sbr{A_{n,j}}}
    \sum_{j\in\bbJ_n(V,i)}
    \prod_{C\in\Comp(G)} \abs{\beta\rbr{C,\bbf_{n,j}}}
    % \nn\\
    % &=
    % \sup_{n\in\ntwo, j\in\bbJ_n(V,i)} \abs{E\sbr{A_{n,j}}}
    % \sum_{j\in\bbJ_n(V,i)}
    % \prod_{C\in\Comp(G)} \babs{\beta\brbr{C,(\bbf|_{\whv(C)})_{n,j|_{V(C)}}}}
    \nn\\
    &=
    \sup_{n\in\ntwo, j\in\bbJ_n(V,i)} 
    \norm{A_{n,j}}_{L^1} % \abs{E\sbr{A_{n,j}}}
    \prod_{C\in\Comp(G)}
    \sum_{j_C\in\bbJ_n(V(C),i)}
    \babs{\beta\brbr{C,(\bbf|_{\whv(C)})_{n,j_C}}}
     \nn\\&=
    \prod_{C\in\Comp(G)} O(n^{e(C)})
    % \\&
    =O(n^{e(G)})
    \label{221129.1820}
  \end{align}
  Since $\barq(C)=0$ for all $C\in\Comp(G)$,
  we can apply Lemma \ref{230519.1737} to
  $C$, $i|_{V(C)}$ and $\bbf|_{\whv(C)}\in\calf(V(C))$
  to obtain 
  $\sum_{j\in\bbJ_n(V(C),i)}
  \babs{\beta\brbr{C,(\bbf|_{\whv(C)})_{n,j}}} = O(n^{e(C)})$ 
  at \eqref{221129.1820}.

  Fix $d\in\bbZ_{\geq0}$, and
  assume that the estimate \eqref{221207.0332} holds for 
  any weighted graph $G$ satisfying $\barq(G)\leq d$.
  Consider a weighted graph $G=(V,\ewt,\vwq)$ with $\barq(G)=d+1$, 
  $i:V\to[2]$, $A\in\cala(V,i)$ and $\bbf\in\calf(V)$.
  By Lemma \ref{221207.2648}, 
  the expectation of $\calii_n\rbr{G, A,\bbf}$  %$E\sbr{\cali_n}$ 
  is decomposed as
  \begin{align}
    E\sbr{\cali_n^{(i)}\rbr{G, A,\bbf}} &=
    \sum_{\gamma\in\Gamma}
    E\sbr{n^{\alpha^{(\gamma)}}
    \cali_n^{(i)}\rbr{G^{(\gamma)}, A^{(\gamma)},\bbf}}
    \label{221207.0346}
  \end{align}
  with some finite set $\Gamma$,
  $\alpha^{(\gamma)}\in\bbR$, 
  $A^{(\gamma)}\in\cala(V,i)$
  and
  a weighted graph $G^{(\gamma)}$ ($\gamma\in\Gamma$) %$=(V,\ewt^{(\gamma)}, \vwq^{(\gamma)})$
  satisfying the following conditions:
  \begin{itemize}
    \item $\barq(G^{(\gamma)}) \leq \barq(G)-1$ for all $\gamma\in\Gamma$
    \item $\alpha^{(\gamma)} + e(G^{(\gamma)})\leq e(G)$ for all $\gamma\in\Gamma$.
  \end{itemize}
  Since $\barq(G^{(\gamma)})\leq d$,
  we can apply the assumption of induction to $G^{(\gamma)}$ %$i, A^{\gamma}, \bbf$ 
  to obtain 
  \begin{align*}
    E\sbr{n^{\alpha^{(\gamma)}}
    \cali_n^{(i)}\rbr{G^{(\gamma)}, A^{(\gamma)}, \bbf}}
    =
    n^{\alpha^{(\gamma)}} O(n^{e(G^{(\gamma)})})
    = O(n^{e(G)}),
  \end{align*}
  and $E\bsbr{\cali_n^{(i)}\rbr{G, A,\bbf}} = O(n^{e(G)})$ from \eqref{221207.0346}.
  Therefore, we have the estimate \eqref{221207.0332} for any weighted graph $G$ with $\barq(G)=d+1$,
  and by induction, 
  the estimate \eqref{221207.0332} is proved for any weighted graph.
\end{proof}

As an easy corollary, we obtain an estimate of the $L^p$-norm of $\calii_n\rbr{G, A,\bbf}$.
\begin{proposition} \label{221208.1720}
  For $n\in\ntwo$, a weighted graph $G=(V,\ewt,\vwq)$, 
  $i:V\to\cbr{1,2}$, $A\in\cala(V,i)$ and $\bbf\in\calf(V)$,
  let $\cali_n=\cali_n^{(i)}(G,A,\bbf)$.
  Then, for $p>1$, $\norm{\cali_n}_{L^p}=O(n^{e(G)})$ as $n\to\infty$.
\end{proposition}
\begin{proof}
  Fix $k_0\in\bbN$ such that $2k_0\geq p$.
  Without loss of generality, we can assume that $V=\cbr{1,...,m}$ with some $m\in\bbN$.
  We define weighted graphs $G^{(k)}$ ($k=1,...,2k_0$) by 
  shifting the weighted graph $G$ by $(k-1)m$, that is 
  $G^{(k)}=(V^{(k)}, \ewt^{(k)}, \vwq^{(k)})$ with 
  \begin{align*}
    V^{(k)}&=\cbr{(k-1)m+1,...,km},
    \\
    \ewt^{(k)}([(v,\kap), (v',\kap')]) &= \ewt([(v-(k-1)m,\kap), (v'-(k-1)m,\kap')])
    &&\tforsm [(v,\kap), (v',\kap')]\in\hp(V^{(k)}), 
    \\
    \vwq^{(k)}(v,\kap) &= \vwq(v-(k-1)m,\kap) 
    &&\tforsm (v,\kap)\in\whv^{(k)}=V^{(k)}\times\cbr{1,2},
  \end{align*}
  and let 
  $%\begin{align*}
    \widetilde{G} = \mathop\vee_{k\in[2k_0]} G^{(k)}.
  $ %\end{align*}
  We can observe $e(\widetilde{G})=2k_0e(G)$.
  We write $\extv=\bigsqcup_{k\in[2k_0]} V^{(k)} = \cbr{1,...,2k_0m}$.
  Similarly, we extend $\exti:\extv\to\cbr{1,2}$ by $\exti|_{V^{(k)}}(v)=i(v-(k-1)m)$ and 
  $\extbbf\in\calf(\extv)$ by 
  $\extbbf^{(v,\kap)}=\bbf^{(v-(k-1)m,\kap)}$ for $(v,\kap)\in\whv^{(k)}$.
  Define 
  $\widetilde{A}\in\cala(\extv, \exti)$ by
  \begin{align*}
    \widetilde{A}_{n,\extj} 
    = \prod_{k\in[2k_0]} A_{n,\extj|_{V^{(k)}}}
    \quad\tforsm \extj\in\bbJ_n(\extv,\exti).
  \end{align*}
  Here we read
  $A_{n,\extj|_{V^{(k)}}}$ as $A_{n,j'}$ 
  with $j'%=(j'(v))_{v\in V}
  =(\extj(v+(k-1)m))_{v\in V}(\in\bbJ_n(V,i))$.

  Then we can write
  \begin{align*}
    (\cali_n)^{2k_0} = 
    \rbr{\sum_{j\in\bbJ_n(V,i)}A_{n,j}\prod_{C\in\Comp(G)}B(C,\bbf_{n,j})}^{2k_0}
    &=%\\&=
    \prod_{k\in[2k_0]}\rbr{\sum_{j_k\in\bbJ_n(V,i)}A_{n,j_k}\prod_{C\in\Comp(G)}B(C,\bbf_{n,j_k})}
    \\&=
    \sum_{\extj\in\bbJ_n(\extv,\exti)}
    \widetilde{A}_{n,\extj}
    \prod_{C\in\Comp(\widetilde{G})} B(C,\extbbf_{n,\extj})
    =%\\&=
    \cali_n^{(\exti)}(\widetilde{G}, \widetilde{A}, \extbbf).
  \end{align*}
  By Proposition \ref{prop:230605.2133}, we have 
  $E\sbr{(\cali_n)^{2k_0}}=
  E\sbr{\cali_n^{(\exti)}(\widetilde{G}, \widetilde{A}, \extbbf)}
  =O(n^{e(\widetilde{G})})=O(n^{2k_0e(G)})$ as $n\to\infty$, and hence 
  we obtain 
  \begin{align*}
    \norm{\cali_n}_{2k_0}=O(n^{e(G)}).
  \end{align*}
\end{proof}

% \subsection{Stability of $\cali_n^{(i)}(G,A,\bbf)$ by the action of Malliavin derivative $D$}
The following estimate shows that the exponent can be used to obtain an estimate of 
the Malliavin-Sobolev norm of $\calii_n(G,A,\bbf)$ as well as its expectation and $L^p$-norm.

\begin{proposition}\label{230725.1133}
  For a weighted graph $G=(V,\ewt, \vwq)$, $i:V\to\cbr{1,2}$, $A\in\cala(V,i)$ and $\bbf\in\calf(V)$,
  let $\cali_n = \cali_n^{(i)}(G,A,\bbf)$.
  Then, for $k_0\in\bbZ_{\geq1}$ and $p>1$,
  it holds that
  \begin{align*}
    \norm{\norm{D^{k_0}\cali_n}_{\calh^{\otimes k_0}}}_{p}=O(n^{e(G)})
  \end{align*}
  as $n\to\infty$.
  As a consequence, $\cali_n=O_M(n^{e(G)})$.
\end{proposition}
\begin{proof}
  Without loss of generality, we assume $V=\cbr{1,...,m}$ with some $m\in\bbN$.
  Let
  $\Lambda(k_0,G)=
  \{\lambda:\cbr{1,...,k_0}%[k_0]
  \to\cbr{0}\sqcup\whv\mid \abs{\lambda^{-1}(\hv)}\leq\vwq(\hv)
  \tforsm \hv\in\whv\}$.
  For $\lambda\in\Lambda(k_0,G)$,
  write
  $\Comp^\lambda_0=\ncbr{C\in\Comp(G)\mid \lambda^{-1}(\whv(C))=\emptyset}$,
  $\Comp^\lambda_+=\{C\in\Comp(G)\mid \lambda^{-1}(\whv(C))\neq\emptyset\}$ and
  $\lambdainv(\hv)=\abs{\lambda^{-1}(\hv)}$ for $\hv\in\cbr{0}\sqcup \whv$. 
  The $k_0$-th derivative of 
  $\cali_n=\cali_n^{(i)}(G,A,\bbf)$ is written as 
  \begin{align*}
    D^{k_0}\cali_n = 
    \sum_{\lambda\in\Lambda(k_0,G)}
    c(\vwq,\lambda)\; \DLamCali%\cali_n^{(\lambda)}
  \end{align*}
  with some integer-valued constants $c(\vwq,\lambda)\in\bbZ_{>0}$ and
  \begin{align*}
    \DLamCali &= 
    \Bigg(
      \sum_{j\in\bbJ_n(V,i)}
      D^{\lambdainv(0)}_{s_{k_1},..,s_{k_{\lambdainv(0)}}} A_{n,j}
      %D^{\lambdainv(0)}_{s_\sbr{\lambda,0}} A_{n,j}
      \times
      \prod_{C\in\Comp^\lambda_0} B(C, \bbf_{n,j})
      \\&\hspace{70pt}\times
      \prod_{C\in\Comp^\lambda_+}\bbrbr{
      %\prod_{C\in\Comp(G)}\bbrbr{
      \beta(C,\bbf_{n,j})\:
      \delta_{V(C)}(\bbf_{n,j},(\vwq-\lambdainv)\bbone_{\whv(C)})
      %\delta(\bbf_{n,j}^{\otimes (\vwq-\lambdainv)\bbone_{\whv(C)}})
      %\prod_{\hv\in\whv(C)}
      \prod_{k\in\lambda^{-1}(\whv(C))}
      \bbf_{n,j}^{(\lambda(k))}(s_k)}
    \Bigg)_{s_1,..,s_{k_0}\in[0,1]}.
    %\prod_{k\in\lambda^{-1}(\whv(C))}
    %\bbf_{n,j}^{(\lambda(k))}(s_k)}
  \end{align*}
  Here,
  if $\lambda^{-1}({0})$ is nonempty,
  we write $\lambda^{-1}({0})=\cbr{k_1,..,k_{\lambdainv(0)}}(\subset\cbr{1,..,k_0})$,
  and if $\lambda^{-1}({0})$ is empty, that is $\lambdainv(0)=0$,
  $D^{\lambdainv(0)}_{s_{k_1},..,s_{k_{\lambdainv(0)}}} A_{n,j}$
  reads $A_{n,j}$.
  We define $\vwq^{(\lambda)}:\whv\to\bbZ_{\geq0}$ by
  $\vwq^{(\lambda)}(\hv) = \vwq(\hv)-\lambdainv(\hv)$ for $\hv\in\wtv$.
  Notice that the expression 
  $\prod_{k\in\lambda^{-1}(\whv(C))} \bbf_{n,j}^{(\lambda(k))}(s_k)$
  can be replaced by 
  $\prod_{\hv\in\whv(C)}\prod_{k\in\lambda^{-1}(\hv)} \bbf_{n,j}^{(\hv)}(s_k)$.
  % \grnb{
  % i.e.
  % $\prod_{k\in\lambda^{-1}(\whv(C))} \bbf_{n,j}^{(\lambda(k))}(s_k) = 
  % \prod_{\hv\in\whv(C)}\prod_{k\in\lambda^{-1}(\hv)} \bbf_{n,j}^{(\hv)}(s_k)$
  % for $j, C\in\Comp^\lambda_+, s_1,...,s_{k_0}$.
  % }

  Since the following inequality holds for $p>2$
  \begin{align}\label{221208.1718}
    \Bnorm{\bnorm{D^{k_0}\cali_n}_{\calh^{\otimes k_0}}}_p 
    \simleq
    \sum_{\lambda\in\Lambda(k_0,G)}
    \Bnorm{\bnorm{\DLamCali}_{\calh^{\otimes k_0}}}_p
    =
    \sum_{\lambda\in\Lambda(k_0,G)}
      \Bnorm{\bnorm{\DLamCali}_{\calh^{\otimes k_0}}^2 }_{p/2}^{1/2},
  \end{align}
  we consider the functional
  $\bnorm{\DLamCali}_{\calh^{\otimes k_0}}^2
  =\abr{\DLamCali,\DLamCali}_{\calh^{\otimes k_0}}
  =:\cali_n^{(\lambda)}$ for $\lambda\in\Lambda(k_0,G)$.
  Fix $\lambda\in\Lambda(k_0,G)$ and we can write
  \begin{align}
    \cali_n^{(\lambda)}&=
    % \bnorm{\DLamCali}_{\calh^{\otimes k_0}}^2=
    % \abr{\DLamCali,\DLamCali}_{\calh^{\otimes k_0}}
    % \\&=
    % \Bigg\langle\Bigg(
    %   \sum_{j\in\bbJ_n(V,i)}
    %   D^{\lambdainv(0)}_{s_{k_1},..,s_{k_{\lambdainv(0)}}} A_{n,j}
    %   \times
    %   \prod_{C\in\Comp^\lambda_0} B(C, \bbf_{n,j})
    %   \\&\hspace{70pt}\times
    %   \prod_{C\in\Comp^\lambda_+}\bbrbr{
    %   \beta(C,\bbf_{n,j})\:
    %   \delta_{V(C)}(\bbf_{n,j},\vwq^{(\lambda)}\bbone_{\whv(C)})
    %   %\prod_{k\in\lambda^{-1}(\whv(C))}\bbf_{n,j}^{(\lambda(k))}(s_k)
    %   \prod_{\hv\in\whv(C)}\prod_{k\in\lambda^{-1}(\hv)} \bbf_{n,j}^{(\hv)}(s_k)}
    % \Bigg)_{s_1,..,s_{k_0}\in[0,1]},
    % \\&\qquad
    % \Bigg(
    %   \sum_{j'\in\bbJ_n(V,i)}
    %   D^{\lambdainv(0)}_{s_{k_1},..,s_{k_{\lambdainv(0)}}} A_{n,j'}
    %   \times
    %   \prod_{C\in\Comp^\lambda_0} B(C, \bbf_{n,j'})
    %   \\&\hspace{70pt}\times
    %   \prod_{C\in\Comp^\lambda_+}\bbrbr{
    %   \beta(C,\bbf_{n,j'})\:
    %   \delta_{V(C)}(\bbf_{n,j'},\vwq^{(\lambda)}\bbone_{\whv(C)})
    %   %\prod_{k\in\lambda^{-1}(\whv(C))} \bbf_{n,j'}^{(\lambda(k))}(s_k)
    %   \prod_{\hv\in\whv(C)}\prod_{k\in\lambda^{-1}(\hv)} \bbf_{n,j'}^{(\hv)}(s_k)}
    % \Bigg)_{s_1,..,s_{k_0}\in[0,1]}
    % \Bigg\rangle
    % \\&=
      \sum_{\substack{j\in\bbJ_n(V,i)\\j'\in\bbJ_n(V,i)}}
      \abr{D^{\lambdainv(0)} A_{n,j}, D^{\lambdainv(0)} A_{n,j'}}_{\calh^{\otimes\lambdainv(0)}}
      \times
      \prod_{C\in\Comp^\lambda_0} \Bcbr{B(C, \bbf_{n,j})\: B(C, \bbf_{n,j'})}
    \nn\\&\hspace{70pt}\times
      \prod_{C\in\Comp^\lambda_+}\Big(%\bbrbr{
      \beta(C,\bbf_{n,j})\:\beta(C,\bbf_{n,j'})\:
      \prod_{\hv\in\whv(C)}
      \abr{\bbf_{n,j}^{(\hv)}, \bbf_{n,j'}^{(\hv)}}^{\lambdainv(\hv)}
    \nn\\&\hspace{170pt}\times
      \delta_{V(C)}(\bbf_{n,j},\vwq^{(\lambda)}\bbone_{\whv(C)})
      \delta_{V(C)}(\bbf_{n,j'},\vwq^{(\lambda)}\bbone_{\whv(C)})
      \Big)
    \label{eq:230818.1855}
  \end{align}
  Here when $\lambdainv(0)=0$, we read 
  $\abr{D^{\lambdainv(0)} A_{n,j}, D^{\lambdainv(0)} A_{n,j'}}_{\calh^{\otimes\lambdainv(0)}}$
  as $A_{n,j}A_{n,j'}$.

  Set $\sftv=\cbr{m+1,...,2m}$, $\sftwhv=\sftv\times\cbr{1,2}$
  and $\dblv=V\sqcup\sftv=\cbr{1,...,2m}$.
  Define 
  $\dbli:\dblv\to[2]$ by 
  $\dbli|_V=i$ and 
  $\dbli(v)=i(v-m)$ for $v\in \sftv$.
  Then we can identify 
  $\bbJ_n(V,i)^2$ % $\bbJ_n(V,i)\times\bbJ_n(V,i)$ 
  with $\bbJ_n(\dblv,\dbli)$.
  We extend $\vwq$ to 
  $\dblvwq:\whv\sqcup\sftwhv=\dblv\times[2]\to\bbZ_{\geq0}$ by
  $\dblvwq|_{\whv} = \vwq$ and 
  $\dblvwq(v,\kappa) = \vwq(v-m,\kappa)$ for $(v,\kappa)\in\sftwhv$, %\sftv\times\cbr{1,2}$,
  and $\ewt$ to $\dblewt:\hp(\dblv)\to\bbZ_{\geq0}$ by
  $\dblewt|_{\hp(V)} = \ewt$, $\dblewt|_{\hp(V,\sftv)} = 0$ and 
  \begin{align*}
    \dblewt([(v,\kap),(v',\kap')]) = \ewt([(v-m,\kap),(v'-m,\kap')])
    \tforsm [(v,\kap),(v',\kap')]\in\hp(\sftv).
  \end{align*}
  We set $\dblbbf\in\calf(\dblv)$ by 
  $\dblbbf|_{\whv}=\bbf$ and 
  $\dblbbf^{(v,\kap)}=\bbf^{(v-m,\kap)}$ for $v\in \sftv$.
  We also define $A^{(\lambda)}=(A^{(\lambda)}_{n,\dblj})_{\dblj\in\bbJ_n(\dblv,\dbli)}$ by 
  $A^{(\lambda)}_{n,\dblj} = A_{n,\dblj_V} A_{n,\dblj_{\sftv}}$
  if $\lambdainv(0)=0$, and 
  \begin{align*}
    A^{(\lambda)}_{n,\dblj} = 
    \babr{D^{\lambdainv(0)} A_{n,\dblj_V}, D^{\lambdainv(0)} A_{n,\dblj_{\sftv}}}_{\calh^{\otimes \lambdainv(0)}}
  \end{align*}
  otherwise.
  We can verify $A^{(\lambda)}\in\cala(\dblv,\dbli)$ 
  by Lemma \ref{lemma:230616.1805}. 

  Hereafter we denote the object shifted by $m$ in the above ways by adding ${}_\star$ to the original notation,
  {and the object obtained by joining the original one and shifted one by adding $\tilde{}~$.}
  We define $\sftv(C) = V(C)+m$, $\wh{\sftv(C)} = \sftv(C)\times[2]$ and 
  %$(V(C))_\star= {\myred\sftv(C) =} V(C)+m$ 
  $\sftc=(\sftv(C), \dblewt\bbone_{\hp(\sftv(C))}, \dblvwq\bbone_{\wh{\sftv(C)}})$ for $C\in\Comp(G)$.
  %$\sftc=((V(C))_\star, \dblewt\bbone_{\hp((V(C))_\star)}, \dblvwq\bbone_{\widehat{(V(C))_\star}})$ 
  We also denote $\sfthv=(v+m,\kap)$ for $\hv=(v,\kap)$.
  Using these notations, \eqref{eq:230818.1855} is rewritten as
  \begin{align}
    \cali_n^{(\lambda)}&=
    \sum_{\dblj\in\bbJ_n(\dblv,\dbli)}
    A^{(\lambda)}_{n,\dblj}
    \times%\\&\hspace{50pt}\times
    \prod_{C\in\Comp^\lambda_0}\Bcbr{B(C,\dblbbf_{n,\dblj})\:B(\sftc,\dblbbf_{n,\dblj})}
    \nn\\&\hspace{50pt}\times
    \prod_{C\in\Comp^\lambda_+}\Big(
    \beta(C,\dblbbf_{n,\dblj})\beta(\sftc,\dblbbf_{n,\dblj})\:
    \beta_{V(C)\sqcup V(\sftc)}(\dblbbf_{n,\dblj}, \ewt^{(\lambda)}_{C})
    \nn\\&\hspace{130pt}\times
    \delta_{V(C)}(\dblbbf_{n,\dblj},\dblvwq^{(\lambda)}\bbone_{\whv(C)})\:
    \delta_{V(\sftc)}(\dblbbf_{n,\dblj},\dblvwq^{(\lambda)}\bbone_{\whv(\sftc)})\Big),
    \label{eq:230818.1919}
  \end{align}
  where we define 
  $\ewt^{(\lambda)}_{C}=\sum_{\hv\in\whv(C)}\lambdainv(\hv)\bbone_{[\hv,\sfthv]}$ for $C\in\Comp^\lambda_+$ and 
  $\dblvwq^{(\lambda)}$ is obtained by extending $\vwq^{(\lambda)}$ as $\dblvwq$.

  For $C\in\Comp^\lambda_+$,
  we define
  $\dblv(C)=V(C)\sqcup V(\sftc)$,
  $\widehat\dblv(C)=\whv(C)\sqcup\whv(\sftc)$ and
  \begin{align*}
    \Pi(C,\lambda)=
    %\Pi_{V(C),V(\sftc)}(\dblvwq^{(\lambda)}\bbone_{\whv(C)}, \dblvwq^{(\lambda)}\bbone_{\whv(\sftc)}) 
    %\\&=
    \Big\{\bfp:\tp(\dblv(C))\to\bbZ_{\geq0}&\mid
    \supp(\bfp)\subset\hp(V(C),V(\sftc)), 
    %\\&\quad
    %\bar\bfp(\tv)\leq\dblvwq^{(\lambda)}(\hv)\tforsm\hv\in\whv(C),\;
    \bar\bfp(\tv)\leq\dblvwq^{(\lambda)}(\hv)\tforsm\hv\in\widehat\dblv(C)\Big\},
  \end{align*}
  with
  $\bar\bfp:\widehat\dblv(C)\to\bbZ_{\geq0}$ such that
  $\bar\bfp(\tv) = \sum_{\tv'\in\whv(\sftc)}\bfp[\tv,\tv']$ for $\tv\in\whv(C)$ and
  $\bar\bfp(\tv) = \sum_{\tv'\in\whv(C)}\bfp[\tv,\tv']$ for $\tv\in\whv(\sftc)$.
  Then by Lemma \ref{lemma:230526.1230}, we have
  \begin{align*}
    &\delta_{V(C)}(\dblbbf_{n,\dblj},\dblvwq^{(\lambda)}\bbone_{\whv(C)})\:
    \delta_{V(\sftc)}(\dblbbf_{n,\dblj},\dblvwq^{(\lambda)}\bbone_{\whv(\sftc)})
    \\&=
    \sum_{\bfp\in\Pi(C,\lambda)}
    %\Pi(V(C),V(\sftc),\dblvwq^{(\lambda)}\bbone_{\whv(C)}, \dblvwq^{(\lambda)}\bbone_{\whv(\sftc)})}
    c(\bfp)~
    \delta_{\dblv(C)}%V(C)\sqcup V(\sftc)}
    (\dblbbf_{n,\dblj},\dblvwq^{(\lambda)}\bbone_{{\widehat\dblv(C)}}-\bar\bfp)\;
    %(\dblbbf_{n,\dblj},\dblvwq^{(\lambda)}\bbone_{\whv(C)\sqcup\whv(\sftc)}-\bar\bfp)\;
    \beta_{\dblv(C)}(\dblbbf_{n,\dblj},\bfp),
    %\beta_{V(C)\sqcup V(\sftc)}(\dblbbf_{n,\dblj},\bfp),
  \end{align*}
  %{\myred for $\dblj\in\bbJ_n(\dblv, \dbli)$}.
  where %we write 
  % $\Pi(C,\lambda)=
  % \Pi_{V(C),V(\sftc)}(\dblvwq^{(\lambda)}\bbone_{\whv(C)}, \dblvwq^{(\lambda)}\bbone_{\whv(\sftc)})$,
  % % {\myred$\dblv(C)=V(C)\sqcup V(\sftc)$ and
  % % $\widehat\dblv(C)=\whv(C)\sqcup\whv(\sftc)$}
  % for short, and
  $c(\bfp)\in\bbN$ is a constant dependent on $\bfp\in\Pi(C,\lambda)$.
  Applying this decomposition to \eqref{eq:230818.1919}, we have
  \begin{align*}
    \cali_n^{(\lambda)}
    % \\&=
    % \sum_{\dblj\in\bbJ_n(\dblv,\dbli)}
    % A^{(\lambda)}_{n,\dblj}
    % \times%\\&\hspace{50pt}\times
    % \prod_{C\in\Comp^\lambda_0}\Bcbr{B(C,\dblbbf_{n,\dblj})\:B(\sftc,\dblbbf_{n,\dblj})}
    % \\&\hspace{50pt}\times
    % \prod_{C\in\Comp^\lambda_+}\Bigg(
    % \beta(C,\dblbbf_{n,\dblj})\beta(\sftc,\dblbbf_{n,\dblj})\:
    % \beta_{\dblv(C)}(\dblbbf_{n,\dblj}, \ewt^{(\lambda)}_{C})\:
    % %\prod_{\hv\in\whv(C)}
    % %\abr{\dblbbf_{n,\dblj}^{(\hv)}, \dblbbf_{n,\dblj}^{(\sfthv)}}^{\lambdainv(\hv)}
    % \delta_{V(C)}(\dblbbf_{n,\dblj},\dblvwq^{(\lambda)}\bbone_{\whv(C)})\:
    % \delta_{V(\sftc)}(\dblbbf_{n,\dblj},\dblvwq^{(\lambda)}\bbone_{\whv(\sftc)})\Bigg)
    &=
    \sum_{\dblj\in\bbJ_n(\dblv,\dbli)}
    A^{(\lambda)}_{n,\dblj}
    \times%\\&\hspace{50pt}\times
    \prod_{C\in\Comp^\lambda_0}\Bcbr{B(C,\dblbbf_{n,\dblj})\:B(\sftc,\dblbbf_{n,\dblj})}
    % \\&\hspace{50pt}\times
    \prod_{C\in\Comp^\lambda_+}\Bigg(
    \beta(C,\dblbbf_{n,\dblj})\beta(\sftc,\dblbbf_{n,\dblj})\:
    \beta_{\dblv(C)}(\dblbbf_{n,\dblj}, \ewt^{(\lambda)}_{C})\:
    \\&\hspace{100pt}\times
    \sum_{\bfp_C\in\Pi(C,\lambda)}
    c(\bfp_C)~
    \delta_{\dblv(C)}
    (\dblbbf_{n,\dblj},\dblvwq^{(\lambda)}\bbone_{\widehat\dblv(C)}-\bar\bfp_C)\;
    \beta_{\dblv(C)}(\dblbbf_{n,\dblj},\bfp_C)
    \Bigg)
    % \\&=
    % \sum_{\bfp\in\prod_{C\in\Comp^\lambda_+}\Pi(C,\lambda)}
    % c(\bfp)~
    % \sum_{\dblj\in\bbJ_n(\dblv,\dbli)}
    % A^{(\lambda)}_{n,\dblj}
    % \times%\\&\hspace{50pt}\times
    % \prod_{C\in\Comp^\lambda_0}\Bcbr{B(C,\dblbbf_{n,\dblj})\:B(\sftc,\dblbbf_{n,\dblj})}
    % \\&\hspace{80pt}\times
    % \prod_{C\in\Comp^\lambda_+}\Bigg(
    % \beta(C,\dblbbf_{n,\dblj})\beta(\sftc,\dblbbf_{n,\dblj})\:
    % \beta_{\dblv(C)}(\dblbbf_{n,\dblj}, \ewt^{(\lambda)}_{C})\:
    % \delta_{\dblv(C)}
    % (\dblbbf_{n,\dblj},\dblvwq^{(\lambda)}\bbone_{\widehat\dblv(C)}-\bar\bfp_C)\;
    % \beta_{\dblv(C)}(\dblbbf_{n,\dblj},\bfp_C)
    % \Bigg)
    \\&=
    \sum_{\bfp\in\prod_{C\in\Comp^\lambda_+}\Pi(C,\lambda)}
    c'(\bfp)~
    \sum_{\dblj\in\bbJ_n(\dblv,\dbli)}
    A^{(\lambda)}_{n,\dblj}
    \times%\\&\hspace{50pt}\times
    \prod_{C\in\Comp^\lambda_0}\Bcbr{B(C,\dblbbf_{n,\dblj})\:B(\sftc,\dblbbf_{n,\dblj})}
    \\&\hspace{70pt}\times
    \prod_{C\in\Comp^\lambda_+}\Bigg(
      \beta_{\dblv(C)}(\dblbbf_{n,\dblj}, 
      \dblewt\bbone_{\hp(V(C))\sqcup\hp(V(\sftc))} + \ewt^{(\lambda)}_{C} + \bfp_C)\:
    \delta_{\dblv(C)}
    (\dblbbf_{n,\dblj},\dblvwq^{(\lambda)}\bbone_{\widehat\dblv(C)}-\bar\bfp_C)
    \Bigg),
  \end{align*}
  where $c'(\bfp)$ is a constant. %dependent on \bfp\in\prod_{C\in\Comp^\lambda_+}\Pi(C,\lambda)
  If $\Comp^\lambda_+=\emptyset$, we read $\prod_{C\in\Comp^\lambda_+}\Pi(C,\lambda)$ and $c'(\bfp)$ as 
  a singleton set and $1$, respectively, and 
  \begin{align*}
    \cali_n^{(\lambda)}&=
    \sum_{\dblj\in\bbJ_n(\dblv,\dbli)}
    A^{(\lambda)}_{n,\dblj}
    \times
    \prod_{C\in\Comp(G)}\Bcbr{B(C,\dblbbf_{n,\dblj})\:B(\sftc,\dblbbf_{n,\dblj})}.
  \end{align*}
  For 
  $\bfp\in\prod_{C\in\Comp^\lambda_+}\Pi(C,\lambda)$ and $C\in\Comp^\lambda_+$, 
  we define a connected weighted graph
  \begin{align*}
    \newComp^{(\lambda,\bfp,C)}=
    (\dblv(C), 
    \dblewt\bbone_{\hp(V(C))\sqcup\hp(V(\sftc))} + \ewt^{(\lambda)}_{C} + \bfp_C,
    \dblvwq^{(\lambda)}\bbone_{\widehat\dblv(C)}-\bar\bfp_C).
  \end{align*}
  Setting
  \begin{align}
    G^{(\lambda,\bfp)} &= 
    \Brbr{\mathop\vee_{C\in \Comp^\lambda_0} (C\vee \sftc)} \vee 
    \Brbr{\mathop\vee_{C\in \Comp^\lambda_+} (\newComp^{(\lambda,\bfp,C)})},
    \label{eq:230609.2153}
  \end{align}
  we can write 
  \begin{align}
    \cali_n^{(\lambda)} %\bnorm{\DLamCali}_{\calh^{\otimes k_0}}^2
    =&
    \sum_{\bfp\in\prod_{C\in\Comp^\lambda_+}\Pi(C,\lambda)}
    %\sum_{\pi\in\prod_{C\in\Comp^\lambda_+}\Pi(\whv(C),\whv(C'),\ddot\vwq-\ddot\lambdainv)}
    c'(\bfp)~
    \cali_n^{(\dbli)}(G^{(\lambda,\bfp)}, A^{(\lambda)}, \dblbbf).
    \label{221208.1719}
  \end{align}
  %for $\lambda\in\Lambda(k_0,G)$.
  
  Using \eqref{221208.1719}, we want to estimate 
  $\cali_n^{(\lambda)}$ % $\bnorm{\DLamCali}_{\calh^{\otimes k_0}}^2$ 
  in $L^p$.
  For $\lambda\in\Lambda(k_0,G)$ such that $\Comp^\lambda_+=\emptyset$,
  % where $\prod_{C\in\Comp^\lambda_+}\Pi(C,\lambda)$ is a singleton,
  we have $e(G^{(\lambda,\bfp)}) = 2e(G)$
  by the definition \eqref{eq:230609.2153} of $G^{(\lambda,\bfp)}$, and 
  \begin{align}\label{eq:230613.2101}
    \bnorm{\cali_n^{(\lambda)}}_{L^p} =
    \bnorm{\cali_n^{(\dbli)}(G^{(\lambda,\bfp)}, A^{(\lambda)}, \dblbbf)}_{L^p}
    = O(n^{e(G^{(\lambda,\bfp)})}) = O(n^{2e(G)}).
  \end{align}
  for $p>1$ by Proposition \ref{221208.1720}.
  
  Next, 
  we consider $\lambda\in\Lambda(k_0,G)$ such that $\Comp^\lambda_+\neq\emptyset$. % and $C\in\Comp^\lambda_+$,
  We observe the difference between $2e(C)$ and $e(\newComp^{(\lambda,\bfp,C)})$
  for $\bfp\in\prod_{C\in\Comp^\lambda_+}\Pi(C,\lambda)$ and $C\in\Comp^\lambda_+$.
  %Fix $\bfp\in\prod_{C\in\Comp^\lambda_+}\Pi(C,\lambda)$ and $C\in\Comp^\lambda_+$, and 
  We write $\gamma=(\lambda,\bfp,C)$ for short.
  By a usual calculation, we have
  $\babs{V(\newComp^{(\gamma)})}=2\babs{V(C)}$, 
  %$\babs{V(\newComp^{(\lambda,\bfp,C)})}=2\babs{V(C)}$, 
  \begin{align*}
    \bartheta(\newComp^{(\gamma)})&=
    2\bartheta(C) + \nabs{\lambda^{-1}(\whv(C))} + \totpi_C,
    \\%and 
    \barq(\newComp^{(\gamma)})&=
    2(\barq(C) -\nabs{\lambda^{-1}(\whv(C))} - \totpi_C),
  \end{align*}
  where we define
  $\totpi_C =\sum_{[\tv,\tv']\in\tp(V(C), V(\sftc))} \bfp_C[\tv,\tv']
  = \half\sum_{\hv\in\whv(C)\sqcup \whv(\sftc)} \bar\bfp_C(\hv)$.
  Hence, by the inequality \eqref{221206.1156}, we have
  $e_q(\newComp^{(\gamma)})\leq-H\barq(\newComp^{(\gamma)})$ and hence
  \begin{align}
    e(\newComp^{(\gamma)}) &=
    1+(2H-1)(\babs{V(\newComp^{(\gamma)})}-1-\eltw(\newComp^{(\gamma)}))
    -2H\bartheta(\newComp^{(\gamma)})
    +e_q(\newComp^{(\gamma)})
    \nn\\&\leq
    % &\leq
    1+(2H-1)(2\babs{V(C)}-2-2\eltw(C))
    +(2H-1)(1+2\eltw(C)-\eltw(\newComp^{(\gamma)}))
    \nn\\&\hspsm
    -2H(2\bartheta(C) + \nabs{\lambda^{-1}(\whv(C))} + \totpi_C)
    -2H\brbr{\barq(C)-\nabs{\lambda^{-1}(\whv(C))} - \totpi_C}
    \nn\\&=
    % 2\rbr{1+(2H-1)(\babs{V(C)}-1-\eltw(C)) -2H\bartheta(C) -H\barq(C)
    % -1/2+(H-1/2)(1+2\eltw(C)-\eltw(\newComp^{(\gamma)}))}
    % \nn\\&=
    2\rbr{e_\theta(C) -H\barq(C)
    -1/2+(H-1/2)(1+2\eltw(C)-\eltw(\newComp^{(\gamma)}))}.
    \label{221208.1708}
  \end{align}
  Notice that $\Comp^\lambda_+\subset\Comp_+(G)=\Comp_1(G)\sqcup\Comp_2(G)$. % for any $\lambda\in\Lambda(k_0,G)$.

  \item[(i)] The case $C\in\Comp_1(G)$.
  By an argument similar to the ones in the proof of Lemma \ref{221207.2648}, 
  we can show 
  $\eltw(\newComp^{(\gamma)})\geq2\eltw(C)$.
  By \eqref{221208.1708}, it holds~
  $e(\newComp^{(\gamma)}) \leq 2\rbr{e_\theta(C) -H(\barq(C)-1) -1}=2e(C)$.

  \item[(ii)] The case $C\in\Comp_2(G)$.
  Again, by an usual argument as in the proof of Lemma \ref{221207.2648}, 
  we can show 
  $\eltw(\newComp^{(\gamma)})\geq2\eltw(C)+1$.
  By \eqref{221208.1708}, it holds~
  $e(\newComp^{(\gamma)}) \leq 2\brbr{e_\theta(C) -H\barq(C)-1/2}=2e(C)$.
  
  \noindent
  Thus, we have $e(\newComp^{(\lambda,\bfp,C)})\leq2e(C)$ 
  % in either case.
  for any $C\in\Comp^\lambda_+$.
  Therefore, 
  for $\bfp\in\prod_{C\in\Comp^\lambda_+}\Pi(C,\lambda)$
  we obtain 
  \begin{align*}
    e(G^{(\lambda,\bfp)}) = 
    \sum_{C\in \Comp^\lambda_0} \rbr{e(C)+e(\sftc)} +
    \sum_{C\in \Comp^\lambda_+} e(\newComp^{(\lambda,\bfp,C)})
    \leq2e(G)
  \end{align*}
  from the definition \eqref{eq:230609.2153} of $G^{(\lambda,\bfp)}$, and 
  \begin{align}\label{eq:230613.2102}
    \bnorm{\cali_n^{(\lambda)}}_{L^p}
    \lesssim
    \sum_{\bfp\in\prod_{C\in\Comp^\lambda_+}\Pi(C,\lambda)}
    \bnorm{\cali_n^{(\dbli)}(G^{(\lambda,\bfp)}, A^{(\lambda)}, \dblbbf)}_{L^p}
    =O(n^{2e(G)}),
  \end{align}
  for $p>1$ %and $\lambda\in\Lambda(k_0,G)$ such that $\Comp^\lambda_+\neq\emptyset$
  by \eqref{221208.1719} and  Proposition \ref{221208.1720}.
  
  By the inequality \eqref{221208.1718} and the estimates \eqref{eq:230613.2101} and \eqref{eq:230613.2102}, 
  it holds that 
  \begin{align*}
    \Bnorm{\bnorm{D^{k_0}\cali_n}_{\calh^{\otimes k_0}}}_{L^p}
    \simleq
    \sum_{\lambda\in\Lambda(k_0,G)}
      \bnorm{\cali_n^{(\lambda)} }_{L^{\frac{p}2}}^{\frac12}
      % \Bnorm{\bnorm{\DLamCali}_{\calh^{\otimes k_0}}^2 }_{p/2}^{1/2}
    = O(n^{e(G)})
  \end{align*}
  for $p>2$.
\end{proof}

\subsection{Change of exponent by the action of $D_{u_n}$}
\label{sec:231005.1626}
When calculating the asymptotic expansion of a functional such as 
$Z_n$ defined at \eqref{eq:230926.1554},
we need estimates of norms of
$D_{u_n}\cali_n=\abr{D\cali_n,u_n}$
for functionals $\cali_n$ having the structure of $\calii_n(G,A,\bbf)$ and 
an $\calh$-valued functional $u_n$ that can be written in the following manner:
\begin{align*}
  u_n = 
  n^{2H-\half} \sum_{j\in[in-1]}
  a_{n,j} I_1(\bbd^{in}_j)\bbd^{in}_j,
\end{align*}
where 
$a_{n,j}$ is some regular (in the Malliavin sense) random variable,
$i\in\cbr{1,2}$, and 
$\bbd^m_j=\bbone^m_{j+1}-\bbone^m_j$ for $j\in[m-1]$.
The following proposition gives an estimate of $D_{u_n}\cali_n$ 
in terms of exponents.

\begin{proposition}\label{221223.1757}
  For a singleton $V_0=\cbr{v_0}$ and $i_0:V_0\to[2]$,
  let $A'\in\cala(V_0,i_0)$.
  Consider $\bff\in L^\infty(\bbR)$ such that {$\supp(\bff)\subset[-1,0]$}.
  Define a $\calh$-valued functional $u_n$ by 
  \begin{align*}
    u_n = n^{2H-1/2} \sum_{j\in\bbJ_n(V_0,i_0)} A'_{n,j} I_1(\bff_{n,j})\bff_{n,j},
  \end{align*}
  where $\bff_{n,j}=T^{(2)}_{i_0(v_0)n,j}(\bff)$.
  (Notice that $\bbJ_n(V_0,i_0)=[i_0(v_0)n-1]$.)
  Let $G=(V,\ewt, \vwq)$ be a weighted graph, $i:V\to[2]$, $A\in\cala(V,i)$ and $\bbf\in\calf(V)$,
  and denote $\cali_n = \cali_n^{(i)}(G,A,\bbf)$.
  We write $\extv=V_0\sqcup V$ and define $\exti:\extv\to\cbr{1,2}$ by 
  $\exti|_V=i$ and $\exti(v_0)=i_0(v_0)$.

  Then there exist 
  a finite set $\Gamma$, $\extbbf\in\calf(\extv)$, and 
  $\alpha^{(\gamma)}\in\bbR$, 
  a weighted graph $G^{(\gamma)}=(\extv,\ewt^{(\gamma)}, \vwq^{(\gamma)})$, 
  and
  $A^{(\gamma)}\in\cala(\extv,\exti)$
  for each $\gamma\in\Gamma$ such that
  the functional $D_{u_n}\cali_n$ can be expressed as 
  \begin{align*}
    D_{u_n}\cali_n &=
    \sum_{\gamma\in\Gamma}
    n^{\alpha^{(\gamma)}}\cali_n^{(\exti)}\rbr{G^{(\gamma)}, A^{(\gamma)},\extbbf},
  \end{align*}
  and the following conditions are satisfied:
  \begin{itemize}
    \item [(a)] If $\barq(C)\neq2$ holds for all $C\in\Comp(G)$, then 
    \begin{align*}
      \max_{\gamma\in\Gamma}\rbr{\alpha^{(\gamma)} + e(G^{(\gamma)})}\leq 
      e(G)-1/2.
    \end{align*}

    \item [(b)] If $\barq(C)=2$ holds for some $C\in\Comp(G)$, then 
    \begin{align*}
      \max_{\gamma\in\Gamma}\rbr{\alpha^{(\gamma)} + e(G^{(\gamma)})}\leq
      e(G).
    \end{align*}

    %\item [(c)] \grnc{Might be added for the case of existence of $C$ such that
    %$\barq(C)=2$ and $C\in\Comp_1(G)$. 
    %Difference: $1/2-H$}
  \end{itemize}
\end{proposition}

\begin{proof}
  For $j\in\bbJ_n(\extv,\exti)=\prod_{v\in\extv}[\exti(v)n-1]$,
  we denote $j|_V$ by $j_V$.
  We define $\extbbf\in\calf(\extv)(\subset L^\infty(\bbR)^{\extv\times\cbr{1,2}})$ by 
  $\extbbf|_{\whv} = \bbf$, $\extbbf^{(v_0,1)}=0$ and $\extbbf^{(v_0,2)}=\bff$.
  (Since $\extbbf^{(v_0,1)}$ is not used in the following argument, its definition does not matter.)
  Notice that we can write 
  $\bff_{n,j_{v_0}}(=T^{(2)}_{i_0(v_0)n,j_{v_0}}(\bff))
  =T^{(2)}_{\exti(v_0)n,j_{v_0}}(\extbbf^{(v_0,2)})
  =\extbbf_{n,j}^{(v_0,2)}$
  for $j\in\bbJ_n(\extv,\exti)$.%\grnb{ここもう一度チェック\tto Done230607.2340}
  
  Since $\cali_n$ is written as
  \begin{align*}
    \cali_n = \cali_n^{(i)}(G, A,\bbf) =
    \sum_{j\in\bbJ_n(V,i)} A_{n,j} \prod_{C\in\Comp(G)} B(C,\bbf_{n,j}),
  \end{align*}
  the functional
  $D_{u_n}\cali_n$ is decomposed as
  \begin{align*}
    D_{u_n}\cali_n &= 
    % \abr{D\Brbr{\sum_{j\in\bbJ_n(V,i)} A_{n,j} \prod_{C\in\Comp(G)}  B(C,\bbf_{n,j})},
    % n^{2H-1/2} \sum_{j'\in\bbJ_n(V_0,i_0)} A'_{n,j'}\: \delta(\bff_{n,j'})\bff_{n,j'}}
    % \\&= 
    n^{2H-1/2} \sum_{j\in\bbJ_n(V,i)} \sum_{j'\in\bbJ_n(V_0,i_0)}
    \Babr{D\Brbr{A_{n,j} \prod_{C\in\Comp(G)}  B(C,\bbf_{n,j})},
    \bff_{n,j'}}A'_{n,j'}\: I_1(\bff_{n,j'})
    % \\&= 
    % n^{2H-1/2} \sum_{j\in\bbJ_n(\extv,\exti)}
    % \Babr{D\Brbr{A_{n,j_V} \prod_{C\in\Comp(G)}  B(C,\bbf_{n,j_V})},
    % \bff_{n,j_{v_0}}}A'_{n,j_{v_0}} \delta(\bff_{n,j_{v_0}})
    \\&= 
    n^{2H-1/2} \sum_{j\in\bbJ_n(\extv,\exti)}
    \Babr{D\Brbr{A_{n,j_V} \prod_{C\in\Comp(G)}  B(C,\extbbf_{n,j})},
    \extbbf_{n,j}^{(v_0,2)}}A'_{n,j_{v_0}} \delta(\extbbf_{n,j}^{(v_0,2)}).
  \end{align*}
  Writing $\hv_0 = (v_0,2)$,
  we have for $j\in\bbJ_n(\extv,\exti)$
  \begin{align*}
    &\hspsm 
    \Babr{D\Brbr{
      A_{n,j_V} 
      \prod_{C\in\Comp(G)} B(C, \extbbf_{n,j})},
      \extbbf^{(\hv_0)}_{n,j}}
      A'_{n,j_{v_0}} \delta(\extbbf^{(\hv_0)}_{n,j}) 
    \\&= 
    \babr{DA_{n,j_V}, \extbbf^{(\hv_0)}_{n,j}}\:
    A'_{n,j_{v_0}}\;\delta(\extbbf^{(\hv_0)}_{n,j})
    \prod_{C\in\Comp(G)} B(C,\extbbf_{n,j})
    \\&\hspsm+
    \sum_{C_1\in\Comp_+(G)} 
    A_{n,j_V} A'_{n,j_{v_0}} 
    \babr{DB(C_1, \extbbf_{n,j}), \extbbf^{(\hv_0)}_{n,j}}\:
    \delta(\extbbf^{(\hv_0)}_{n,j}) 
    \prod_{\substack{C\in\Comp(G)\\C\neq C_1}} B(C, \extbbf_{n,j}).
  \end{align*}
  Denote $V_1=V(C_1)$ and $\whv_1=\whv(C_1)$ for brevity.
  For $C_1\in\Comp_+(G)$, %{\myred if $\Comp_+(G)\neq\emptyset$},
  since
  $B(C_1,\extbbf_{n,j})=\beta(C_1,\extbbf_{n,j})\delta(C_1,\extbbf_{n,j})$ and 
  $\delta(C_1,\extbbf_{n,j}) = \delta_{V_1}(\extbbf_{n,j}, \vwq\bbone_{\whv_1})$ by definition, 
  we have
  \begin{align*}
    \babr{DB(C_1,\extbbf_{n,j}),\extbbf^{(\hv_0)}_{n,j}}
    %\delta(\extbbf^{(\hv_0)}_{n,j})
    =%\\&=
    \sum_{\substack{\hv_1\in\wtv_1\\\vwq(\hv_1)>0}}
    \vwq(\hv_1)\:
    \beta(C_1,\extbbf_{n,j})\:
    \delta_{V_1}\brbr{\extbbf_{n,j},\vwq\bbone_{\wtv_1}-\bbone_\cbr{\hv_1}}
    \babr{\extbbf^{(\hv_1)}_{n,j},\extbbf^{(\hv_0)}_{n,j}}.
    %\delta(\extbbf^{(\hv_0)}_{n,j})
  \end{align*}
  By the product formula,
  \begin{align*}
    &\hspsm
    \delta_{V_1}\brbr{\extbbf_{n,j},\vwq\bbone_{\wtv_1}-\bbone_\cbr{\hv_1}}
    \delta(\extbbf^{(\hv_0)}_{n,j})
    \\&=
    \delta_{V_1\sqcup\cbr{v_0}}
    \rbr{\extbbf_{n,j},\vwq\bbone_{\wtv_1}-\bbone_\cbr{\hv_1}+\bbone_\cbr{\hv_0}}
    \\&\hspsm+
    \sum_{\substack{\hv_2\in\wtv_1\\ (\vwq-\bbone_\cbr{\hv_1})(\hv_2)>0}}
    (\vwq-\bbone_\cbr{\hv_1})(\hv_2)~
    \delta_{V_1\sqcup\cbr{v_0}}
    \rbr{\extbbf_{n,j},\vwq\bbone_{\wtv_1}-\bbone_\cbr{\hv_1}-\bbone_\cbr{\hv_2}}
    \babr{\extbbf^{(\hv_2)}_{n,j},\extbbf^{(\hv_0)}_{n,j}}.
  \end{align*}
  We define the following weighted graphs:
  \begin{align}
    C^{(\tv_1)}&=
    (V_1\sqcup \cbr{v_0}, %V_0, 
    \ewt\bbone_{\tp(V_1)} + \bbone_{[\tv_1,\hv_0]},\;
    \vwq\bbone_{\wtv_1}-\bbone_\cbr{\tv_1}+\bbone_\cbr{\hv_0})
    \label{221207.1841}
    \\
    C^{(\tv_1,\tv_2)}&=
    (V_1\sqcup \cbr{v_0}, %V_0, 
    \ewt\bbone_{\tp(V_1)} 
    + \bbone_{[\tv_1,\hv_0]} + \bbone_{[\tv_2,\hv_0]},\;
    \vwq\bbone_{\wtv_1}-\bbone_\cbr{\tv_1}-\bbone_\cbr{\tv_2}),
    \label{221207.1842}
  \end{align}
  Notice that $C^{(\tv_1)}$ and $C^{(\tv_1,\tv_2)}$ are connected.
  With these weighted graphs, we can write
  \begin{align*}
    &\babr{DB(C_1, \extbbf_{n,j}), \extbbf^{(\hv_0)}_{n,j}}\:
    \delta(\extbbf^{(\hv_0)}_{n,j}) 
    % \\&=
    % \sum_{\substack{\hv_1\in\wtv_1\\\vwq(\hv_1)>0}}
    % \vwq(\hv_1)\:
    % \beta(C_1,\extbbf_{n,j})\:
    % \delta_{V_1}\brbr{\extbbf_{n,j},\vwq\bbone_{\wtv_1}-\bbone_\cbr{\hv_1}}
    % \babr{\extbbf^{(\hv_1)}_{n,j},\extbbf^{(\hv_0)}_{n,j}}
    % \delta(\extbbf^{(\hv_0)}_{n,j}) 
    % \\&=
    % \sum_{\substack{\hv_1\in\wtv_1\\\vwq(\hv_1)>0}}
    % \vwq(\hv_1)\:
    % \beta(C_1,\extbbf_{n,j})\:
    % \babr{\extbbf^{(\hv_1)}_{n,j},\extbbf^{(\hv_0)}_{n,j}}
    % \delta_{V_1\sqcup\cbr{v_0}}
    % \rbr{\extbbf_{n,j},\vwq\bbone_{\wtv_1}-\bbone_\cbr{\hv_1}+\bbone_\cbr{\hv_0}}
    % \\&\hspsm+
    % \sum_{\substack{\hv_1\in\wtv_1\\\vwq(\hv_1)>0}}
    % \vwq(\hv_1)\:
    % \beta(C_1,\extbbf_{n,j})\:
    % \babr{\extbbf^{(\hv_1)}_{n,j},\extbbf^{(\hv_0)}_{n,j}}
    % \sum_{\substack{\hv_2\in\wtv_1\\ (\vwq-\bbone_\cbr{\hv_1})(\hv_2)>0}}
    % (\vwq-\bbone_\cbr{\hv_1})(\hv_2)~
    % \delta_{V_1\sqcup\cbr{v_0}}
    % \rbr{\extbbf_{n,j},\vwq\bbone_{\wtv_1}-\bbone_\cbr{\hv_1}-\bbone_\cbr{\hv_2}}
    % \babr{\extbbf^{(\hv_2)}_{n,j},\extbbf^{(\hv_0)}_{n,j}}
    \\&=
    \sum_{\substack{\hv_1\in\wtv_1\\\vwq(\hv_1)>0}}
    \vwq(\hv_1)\:
    B(C^{(\tv_1)},\extbbf_{n,j})
    +% \\&\hspsm+
    \sum_{\substack{\hv_1\in\wtv_1\\\vwq(\hv_1)>0}}
    \sum_{\substack{\hv_2\in\wtv_1\\ (\vwq-\bbone_\cbr{\hv_1})(\hv_2)>0}}
    \vwq(\hv_1)\,
    (\vwq-\bbone_\cbr{\hv_1})(\hv_2)~
    B(C^{(\tv_1,\tv_2)},\extbbf_{n,j}).
  \end{align*}
  Then we have
  \begin{align}
    D_{u_n}\cali_n &= 
    % n^{2H-1/2} \sum_{j\in\bbJ_n(\extv,\exti)}
    % \Babr{D\Brbr{A_{n,j_V} \prod_{C\in\Comp(G)}  B(C,\extbbf_{n,j})},
    % \extbbf_{n,j}^{(\hv_0)}}A'_{n,j_{v_0}} \delta(\extbbf_{n,j}^{(\hv_0)})
    % \nn\\&=
    % n^{2H-1/2} 
    % \sum_{j\in\bbJ_n(\extv,\exti)}
    % \babr{DA_{n,j_V},\extbbf^{(\hv_0)}_{n,j}} A'_{n,j_{v_0}}\;
    % \delta(\extbbf^{(\hv_0)}_{n,j}) 
    % \prod_{C\in\Comp(G)}  B(C,\extbbf_{n,j})
    % \nn\\&\hspsm+
    % \sum_{C_1\in\Comp_+(G)} 
    % n^{2H-1/2} \sum_{j\in\bbJ_n(\extv,\exti)}
    % A_{n,j_V} A'_{n,j_{v_0}}\;
    % \abr{DB(C_1,\extbbf_{n,j}), \extbbf^{(\hv_0)}_{n,j}}\:
    % \delta(\extbbf^{(\hv_0)}_{n,j}) 
    % \prod_{\substack{C\in\Comp(G)\\C\neq C_1}}  B(C,\extbbf_{n,j})
    % \nn\\&=
    % n^{-1/2} 
    % \sum_{j\in\bbJ_n(\extv,\exti)}
    % n^{2H} \babr{DA_{n,j_V},\extbbf^{(\hv_0)}_{n,j}} A'_{n,j_{v_0}}\;
    % \delta(\extbbf^{(\hv_0)}_{n,j}) 
    % \prod_{C\in\Comp(G)}  B(C,\extbbf_{n,j})
    % \nn\\&\hspsm+
    % \sum_{C_1\in\Comp_+(G)} 
    % \sum_{\substack{\hv_1\in\wtv_1\\\vwq(\hv_1)>0}}
    % n^{2H-1/2} \sum_{j\in\bbJ_n(\extv,\exti)}
    % \vwq(\hv_1)\:
    % A_{n,j_V} A'_{n,j_{v_0}}\;
    % B(C^{(\tv_1)},\extbbf_{n,j})
    % \prod_{\substack{C\in\Comp(G)\\C\neq C_1}}  B(C,\extbbf_{n,j})
    % \nn\\&\hspsm+
    % \sum_{C_1\in\Comp_+(G)} 
    % \sum_{\substack{\hv_1\in\wtv_1\\\vwq(\hv_1)>0}}
    % \sum_{\substack{\hv_2\in\wtv_1\\ (\vwq-\bbone_\cbr{\hv_1})(\hv_2)>0}}
    % n^{2H-1/2} \sum_{j\in\bbJ_n(\extv,\exti)}
    % \vwq(\hv_1)\,(\vwq-\bbone_\cbr{\hv_1})(\hv_2)~
    % A_{n,j_V} A'_{n,j_{v_0}}\;
    % \nn\\&\hspace{250pt}\times
    % B(C^{(\tv_1,\tv_2)},\extbbf_{n,j}).
    % \prod_{\substack{C\in\Comp(G)\\C\neq C_1}}  B(C,\extbbf_{n,j})
    % \nn\\&=
    n^{-1/2} 
    \sum_{j\in\bbJ_n(\extv,\exti)} A^{(0)}_{n,j}\;
    B(C^{(0)},\extbbf_{n,j})
    \prod_{C\in\Comp(G)}  B(C,\extbbf_{n,j})
    \nn\\&\hspsm+
    \sum_{C_1\in\Comp_+(G)} 
    \sum_{\substack{\hv_1\in\wtv_1\\\vwq(\hv_1)>0}}
    n^{2H-1/2} \sum_{j\in\bbJ_n(\extv,\exti)}
    A^{(\hv_1)}_{n,j}\;
    B(C^{(\tv_1)},\extbbf_{n,j})
    \prod_{\substack{C\in\Comp(G)\\C\neq C_1}}  B(C,\extbbf_{n,j})
    \nn\\&\hspsm+
    \sum_{C_1\in\Comp_+(G)} 
    \sum_{\substack{\hv_1\in\wtv_1\\\vwq(\hv_1)>0}}
    \sum_{\substack{\hv_2\in\wtv_1\\ (\vwq-\bbone_\cbr{\hv_1})(\hv_2)>0}}
    n^{2H-1/2} \sum_{j\in\bbJ_n(\extv,\exti)}
    A^{(\hv_1,\hv_2)}_{n,j}\;
    B(C^{(\tv_1,\tv_2)},\extbbf_{n,j}).
    \prod_{\substack{C\in\Comp(G)\\C\neq C_1}}  B(C,\extbbf_{n,j}),
    \label{eq:230608.1152}
  \end{align}
  where we set
  $C^{(0)}=(V_0,0,\bbone_\cbr{\hv_0})$,
  $A^{(0)}_{n,j} = n^{2H} \babr{DA_{n,j_V}, \extbbf_{n,j}^{(\hv_0)}} A'_{n,j_{v_0}}$, %,\quad% \\
  \begin{align*}
    A^{(\hv_1)}_{n,j} &= \vwq(\hv_1) A_{n,j_V} A'_{n,j_{v_0}}\tand % \\
    A^{(\hv_1,\hv_2)}_{n,j} = \vwq(\hv_1)\,(\vwq-\bbone_\cbr{\hv_1})(\hv_2)\; A_{n,j_V} A'_{n,j_{v_0}}
  \end{align*}
  for $j\in\bbJ_n(\extv,\exti)$.
  %{\mygreen$c\rbr{\tv_1}$ and $c\rbr{\tv_1,\tv_2}$ are integer-valued constants.}

  Set $\Gamma = \cbr{0}\sqcup\Gamma_1\sqcup\Gamma_2$ with
  \begin{align*}
    \Gamma_1 &= 
    \bigsqcup_{C_1\in\Comp_+(G)}    
    \bcbr{\tv_1\in\whv_1\mid\vwq(\tv_1)>0},
    \\
    \Gamma_2 &= 
    \bigsqcup_{C_1\in\Comp_+(G)}   
    \bcbr{(\tv_1,\tv_2)\in\wtv_1^2\mid 
    \vwq(\tv_1)>0 \tandsm (\vwq-\bbone_\cbr{\tv_1})(\tv_2)>0}.
  \end{align*}
  Notice that $A^{(\gamma)}=(A^{(\gamma)}_{n,j})_{j\in\bbJ_n(\extv,\exti), n\in\ntwo}
  \in\cala(\extv,\exti)$
  for $\gamma\in\Gamma$,
  which follows from Lemma \ref{lemma:230616.1754} (i) and Lemma \ref{lemma:230616.1805}.
  Let $\alpha^{(0)}=-1/2$ and $G^{(0)}=C^{(0)}\vee G$.
  For $\gamma=\hv_1\in\Gamma_1$ and $\gamma=(\hv_1,\hv_2)\in\Gamma_2$, we define 
  $\alpha^{(\gamma)}=2H-1/2$ and 
  $G^{(\gamma)} %(= (\extv, \ewt^{(\gamma)}, \vwq^{(\gamma)}))
  :=C^{(\gamma)} \vee 
  \brbr{\mathop\vee_{C\in \Comp(G)\setminus\cbr{C_\gamma}} C}$,
  where 
  we defined $C^{(\gamma)}$ in \eqref{221207.1841} and \eqref{221207.1842}, and
  we denote by $C_\gamma$ the component $C_1\in\Comp_+(G)$ such that $\hv_1\in\whv(C_1)$.
  Hence, by \eqref{eq:230608.1152}, we obtain 
  \begin{align*}
    D_{u_n}\cali_n &=
    \sum_{\gamma\in\Gamma} n^{\alpha^{(\gamma)}} 
    \cali_n^{(\exti)}(G^{(\gamma)}, A^{(\gamma)}, \extbbf).
  \end{align*}

  Next we observe the difference between 
  $e(G)$ and $\alpha^{(\gamma)}+e(G^{(\gamma)})$.

  \item[(i)] When $\gamma=0$, since $e(C^{(0)}) = 1 + \rbr{-1/2-H} = 1/2-H$,
  we have 
  \begin{align*}
    \alpha^{(0)}+e(G^{(0)}) = 
    %\alpha^{(\gamma)}+e(G^{(\gamma)}) = 
    \alpha^{(0)}+e(C^{(0)}) + e(G)= 
    %\alpha^{(\gamma)}+e(C^{(\gamma)}) + e(G)= 
    %-1/2 + (1/2-H) + e(G) = 
    e(G)-H. 
  \end{align*}

  \item[(ii)] Consider when $\gamma(=\hv_1=(v_1,\kappa_1))\in\Gamma_1$.
  %We write $v_\gamma$ for $v\in V(C_\gamma)$ such that 
  %$\gamma\in\cbr{v}\times\cbr{1,2}$ holds.
  For a spanning tree $\bbT\subset E(C_\gamma)$ such that 
  $\eltw(C_\gamma) = \abs{E_2(C_\gamma)\cap\bbT}$,
  let $\tilde{\bbT} = \bbT \sqcup\cbr{[v_0,v_1]}$.
  Then $\tilde\bbT$ is a subset of $E(C^{(\gamma)})=E(C_\gamma)\sqcup\cbr{[v_0,v_1]}$ and 
  a spanning tree of $C^{(\gamma)}$.
  Because $E_2(C^{(\gamma)})=E_2(C_\gamma)\sqcup\cbr{[v_0,v_1]}$,
  we have 
  $\ell_2(C^{(\gamma)})\geq\babs{E_2(C^{(\gamma)})\cap\tilde\bbT}=
  \abs{E_2(C_\gamma)\cap\bbT} + 1 = \eltw(C_\gamma) + 1$.
  (In fact the equality holds.)
  Since we have $\bartheta(C^{(\gamma)})=\bartheta(C_\gamma)+1$ and 
  $\babs{V(C^{(\gamma)})}=\babs{V(C_\gamma)}+1$,
  $e_\theta(C^{(\gamma)})$ is bounded as 
  \begin{align}
    e_\theta(C^{(\gamma)}) &= 
    1 + (2H-1)(\babs{V(C^{(\gamma)})}-1 -\eltw(C^{(\gamma)}))
    -2H \bartheta(C^{(\gamma)})
    \nn\\&\leq
    % 1 + (2H-1)\brbr{(\babs{V(C_\gamma)}+1) -1 -(\eltw(C_\gamma) +1)}
    % -2H (\bartheta(C_\gamma)+1)
    % \nn\\&=
    1 + (2H-1)\brbr{\babs{V(C_\gamma)}-1 -\eltw(C_\gamma)}
    -2H\bartheta(C_\gamma)-2H 
    = e_\theta(C_\gamma)-2H
    \label{221207.2015}
  \end{align}

  If $C_\gamma\in\Comp_1(G)$,
  then by \eqref{eq:230818.1803} it holds that
  \begin{align*}
    e_q(C^{(\gamma)})\leq 
    -1 -H(\barq(C^{(\gamma)})-1) = 
    -1 -H(\barq(C_\gamma)-1) = e_q(C_\gamma)
  \end{align*}
  since $\barq(C_\gamma)=\barq(C^{(\gamma)})>0$.
  Next consider the case $C_\gamma\in\Comp_2(G)$.
  Since $\barq_2(C^{(\gamma)})=\barq_2(C_\gamma)\geq1$ and 
  $\barq_1(C^{(\gamma)})=0$,
  %{\myred$\sum_{v\in V(C^{(\gamma)})}\vwq^{(\gamma)}(v,1)=0$,}
  we have $C^{(\gamma)}\in\Comp_2(G^{(\gamma)})$ and
  \begin{align*}
    e_q(C^{(\gamma)}) =  
    -1/2 -H\barq(C^{(\gamma)}) = 
    -1/2 -H\barq(C_\gamma) = e_q(C_\gamma).
  \end{align*}
  Thus $e_q(C^{(\gamma)})\leq e_q(C_\gamma)$ holds
  in the both cases $C_\gamma\in\Comp_1(G)$ or $\Comp_2(G)$, and
  we obtain 
  $e(C^{(\gamma)})\leq e(C_\gamma) - 2H$ from \eqref{221207.2015} and hence
  \begin{align*}
    \alpha^{(\gamma)} + e(G^{(\gamma)}) 
    %= \alpha^{(\gamma)} + e(C^{(\gamma)}) 
    %+ \sum_{C\in\Comp(G)\setminus\cbr{C_\gamma}} e(C)
    = (2H - 1/2) + e(C^{(\gamma)}) 
    + \sum_{C\in\Comp(G)\setminus\cbr{C_\gamma}} e(C)
    \leq e(G) -1/2.
  \end{align*}

  \vspssm
  (iii) Consider when $\gamma\in\Gamma_2$.
  By the same argument as in (ii), we have 
  $\eltw(C^{(\gamma)})\geq\eltw(C_\gamma) +1$. 
  Since we have $\bartheta(C^{(\gamma)})=\bartheta(C_\gamma)+2$ and 
  $\babs{V(C^{(\gamma)})}=\babs{V(C_\gamma)}+1$,
  $e_\theta(C^{(\gamma)})$ is bounded as 
  \begin{align}
    e_\theta(C^{(\gamma)}) &= 
    1 + (2H-1)(\babs{V(C^{(\gamma)})}-1 -\eltw(C^{(\gamma)}))
    -2H \bartheta(C^{(\gamma)})
    \nn\\&\leq
    % 1 + (2H-1)\brbr{(\babs{V(C_\gamma)}+1) -1 -(\eltw(C_\gamma) +1)}
    % -2H (\bartheta(C_\gamma)+2)
    % \nn\\&=
    1 + (2H-1)\brbr{\babs{V(C_\gamma)}-1 -\eltw(C_\gamma)}
    -2H\bartheta(C_\gamma)-4H 
    = e_\theta(C_\gamma)-4H.
    \label{221207.2105}
  \end{align}
  Notice that $\barq(C_\gamma)\geq2$.
  We separate $\Gamma_2$ with respect to whether $\barq(C_\gamma)=2$ or not.
  We set 
  \begin{align*}
    \Gamma_{2,2} &= 
    \bigsqcup_{\substack{C_1\in\Comp_+(G)\\\barq(C_\gamma)=2}}
    \bcbr{(\tv_1,\tv_2)\in\wtv(C_1)^2\mid 
    \vwq(\tv_1)>0 \tandsm (\vwq-\bbone_\cbr{\tv_1})(\tv_2)>0}
  \end{align*}
  and $\Gamma_{2,3}=\Gamma_2\setminus\Gamma_{2,2}$.

  \item[(iii-i)] First we treat the case $\gamma\in\Gamma_{2,3}$. %$\barq(C_\gamma)\geq3$.
  Notice that $\barq(C^{(\gamma)})=\barq(C_\gamma)-2\geq1$.
  If $C_\gamma\in\Comp_1(G)$, then \eqref{eq:230818.1803} we have
  \begin{align*}
    e_q(C^{(\gamma)})\leq 
    -1 -H (\barq(C^{(\gamma)})-1) = 
    -1 -H ( \barq(C_\gamma)-2-1) = 
    e_q(C_\gamma) + 2H.
  \end{align*}
  If $C_\gamma\in\Comp_2(G)$, then we can observe 
  $C^{(\gamma)}\in\Comp_2(G^{(\gamma)})$ and
  \begin{align*}
    e_q(C^{(\gamma)})=
    -1/2 -H \barq(C^{(\gamma)}) = 
    -1/2 -H (\barq(C_\gamma)-2) = 
    e_q(C_\gamma) + 2H.
  \end{align*}
  Hence $e_q(C^{(\gamma)})\leq e_q(C_\gamma)+2H$ holds
  in both the cases,
  using \eqref{221207.2105} we obtain 
  $e(C^{(\gamma)})\leq e(C_\gamma) - 2H$ and hence
  \begin{align*}
    \alpha^{(\gamma)} + e(G^{(\gamma)})
    \leq e(G) -1/2.
  \end{align*}

  \item[(iii-ii)]
  If $\gamma\in\Gamma_{2,2}$, %$\barq(C_\gamma)=2$.
  then we have $\barq(C^{(\gamma)})=\barq(C_\gamma)-2=0$,
  $e_q(C^{(\gamma)})=0$,
  $e_q(C_{\gamma}) \geq -1/2-2H$, and hence
  \begin{align*}
    e(C^{(\gamma)})=
    e_\theta(C^{(\gamma)}) &\leq
    e_\theta(C_\gamma)-4H = 
    e(C_\gamma)-e_q(C_\gamma)-4H
    \leq%\\&\leq
    % e(C_\gamma)+1/2+2H-4H=
    e(C_\gamma)+1/2-2H.
  \end{align*}
  Thus, we obtain
  \begin{align*}
    \alpha^{(\gamma)} + e(G^{(\gamma)}) &=
    (2H-1/2) + e(C^{(\gamma)})
    + \sum_{C\in\Comp(G)\setminus\cbr{C_\gamma}} e(C)
    \leq e(G).
  \end{align*}

  Summing up the above argument, we obtain the following conclusion:
  \item[(a)] If $\barq(C)\neq2$ for all $C\in\Comp(G)$, then 
  we have $\Gamma_{2,2}=\emptyset$ and $\Gamma=\cbr{0}\sqcup\Gamma_1\sqcup\Gamma_{2,3}$,
  % $\Gamma$ consists of $\gamma$ corresponding to one of the cases (i), (ii) and (iii-i),
   and hence
  \begin{align*}
    \max_{\gamma\in\Gamma}\rbr{\alpha^{(\gamma)} + e(G^{(\gamma)})}\leq 
    e(G)-1/2.
  \end{align*}
  \item[(b)] If $\barq(C)=2$ holds for some $C\in\Comp(G)$, then
  there exists $\gamma\in\Gamma_{2,2}$.
  Therefore, 
  \begin{align*}
    \max_{\gamma\in\Gamma}\rbr{\alpha^{(\gamma)} + e(G^{(\gamma)})}\leq
    e(G).
  \end{align*}
\end{proof}

\begin{scrap}
  If $C_\gamma\in\Comp_1(G)$, thus 
  $e_q(C_{\gamma}) = -1-H$, then we have
  \begin{align*}
    \alpha^{(\gamma)} + e(G^{(\gamma)}) \leq
    (2H-1/2) - (-1-H) -4H + e(G)=
    e(G) + 1/2 -H.
  \end{align*}
  If $C_\gamma\in\Comp_2(G)$, thus 
  $e_q(C_{\gamma}) = -1/2-2H$, then we have
  \begin{align*}
    \alpha^{(\gamma)} + e(G^{(\gamma)}) \leq
    (2H-1/2) - (-1/2-2H) -4H + e(G)=
    e(G).
  \end{align*}
\end{scrap}

%\newpage
\subsection{Preliminary lemmas}\label{230519.1730}
\subsubsection{Order estimate of sums of $\beta$}
\label{sec:231005.1445}
Here we want to obtain an estimate of the sum of 
$\beta_V(\bbf_{n,j}, \ewt)$ (or $\beta(C,\bbf_{n,j})$)
for a connected weighted graph $C=(V,\ewt,\vwq)$.
% {\myred subsubsection冒頭}
Since the estimate of the sum is complicated, in Lemma \ref{221222.1901},
we first consider the sum of $\beta_V(\bbf_{n,j}, \ewt)$
related to a rather simple weighted graph:
``tree-shaped''. %, whose definition is given below. 
Then in Lemma \ref{230519.1737}, 
we deal with the sum of $\beta_V(\bbf_{n,j}, \ewt)$
for a general connected one.

We define a weighted graph $C$ to be tree-shaped, %かつてはtreeと呼んでいた． 
if $(V(C), E(C))$ is a tree, that is a connected graph without cycles.
Note that 
any weighted graph $C=(V,\ewt,\vwq)$ with $\abs{V}=1$ is 
tree-shaped.
For notational convenience, we introduce another definition.
For a nonempty finite set $V$,
a map $\pi:\hp(V)\to p(V)$ is defined by
\begin{align}\label{230621.1526}
  \pi([(v,\kap), (v',\kap')]) = [v,v']
\end{align}
for $[(v,\kap), (v',\kap')]\in\hp(V)$(, hence $v\neq v'\in V$ and $[v,v']\in p(V)$).
For a subset $E\subset p(V)$, we write
$\widehat{E}=\pi^{-1}(E)(\subset\hp(V))$.
We start with the decomposition of
the summand $\beta_V(\bbf_{n,j}, \ewt)$
associated with a tree-shaped weighted graph $C=(V,\ewt,\vwq)$,
which can be written as a product of 
$\beta$ associated to another tree-shaped weighted graph 
having vertices one less than that of $C$
% $C'$ with $\abs{V(C')}=\abs{V}-1$
and the other factor whose estimate is easily obtained.
% Since $\vwq$ of $C=(V,\ewt,\vwq)$ does not concern the sum,
% we write $\vwq=0$ and $C=(V,\ewt,0)$ in this subsubsection.
The main result of Section \ref{sec:231005.1445} is Lemma \ref{230519.1737},
where we deal with a connected weighted graph with $\vwq=0$.
Therefore we assume the weights of vertices to be zero.

Let 
$C=(V,\ewt,0)$ be a tree-shaped weighted graph 
with $\abs{V}\geq2$.
Then there are vertices 
$v_0\neq v_1\in V$ such that 
$\cewt([v_0,v_1])>0$ and 
$\cewt([v_0,v])=0$ for $v\in V\setminus\cbr{v_0,v_1}$.
We write $V'=V\setminus\cbr{v_0}$ and 
$\whv'=V'\times\cbr{1,2}$.
For $j\in\bbJ_n(V,i)$, we have
\begin{align}
  \beta_V(\bbf_{n,j}, \ewt)
  &=% \nn\\&=
  \prod_{[\hv,\hv']\in\hp(V)}
  \babr{\bbf_{n,j}^{(\hv)},\bbf_{n,j}^{(\hv')}}
  ^{\ewt[\hv,\hv']}
  \nn\\&=
  {
  \prod_{[\hv,\hv']\in\hp(V')}
  \babr{\bbf_{n,j}^{(\hv)},\bbf_{n,j}^{(\hv')}}
  ^{\ewt[\hv,\hv']}}
  \times{
  \prod_{[\hv,\hv']\in\hp(V',\cbr{v_0})}
  \babr{\bbf_{n,j}^{(\hv)},\bbf_{n,j}^{(\hv')}}
  ^{\ewt[\hv,\hv']}}
  \label{eq:230928.1633}
\end{align}
We denote
$\bbf'=\bbf|_{\whv'}\in\calf(V')$
and
$j'=j|_{V'}$ for $j\in\bbJ_n(V,i)$.

For $\hv=(v,\kap)\in\whv'$ and $j\in\bbJ_n(V,i)$, 
we have
$\bbf_{n,j}^{(\hv)}
=T_{i(v)n,j_v}^{(\kap)}(\bbf^{(\hv)})%(\bbf^{(v,\kappa)})
=T_{i(v)n,j'_v}^{(\kap)}({\bbf'}^{(\hv)})%(\bbf^{(v,\kappa)})
={\bbf'}_{n,j'}^{(\hv)}$,
and the first factor of \eqref{eq:230928.1633} is written as
\begin{align*}
  &\prod_{[\hv,\hv']\in\hp(V')}
  \babr{\bbf_{n,j}^{(\hv)},\bbf_{n,j}^{(\hv')}}
  ^{\ewt[\hv,\hv']}
  =% \\&=
  \prod_{[\hv,\hv']\in\hp(V')}
  \babr{{\bbf'}_{n,j'}^{(\hv)},{\bbf'}_{n,j'}^{(\hv')}}
  ^{\ewt[\hv,\hv']}
  =% \\&=
  \beta_{V'}({\bbf'}_{n,j'}, \ewt\bbone_{\hp(V')}).
\end{align*}

Notice that there is one-to-one correspondence between 
$\hp(V',\cbr{v_0})=
\whv'\times\cbr{(v_0,1),(v_0,2)}$.
By the way we chose $v_0$ and $v_1$, we have 
$\ewt[\hv_0,\hv]=0$
for any $\hv_0\in\cbr{(v_0,1),(v_0,2)}$ and 
$\hv\in(V'\setminus\cbr{v_1})\times\cbr{1,2}$.
Hence the second factor of \eqref{eq:230928.1633}
is written as
\begin{align*}
  \prod_{[\hv,\hv']\in\hp(V',\cbr{v_0})}
  \babr{\bbf_{n,j}^{(\hv)},\bbf_{n,j}^{(\hv')}}
  ^{\ewt[\hv,\hv']}
  &=% \\&=
  \prod_{\kap_0,\kap_1\in\cbr{1,2}}
  \babr{\bbf_{n,j}^{(v_0,\kap_0)},\bbf_{n,j}^{(v_1,\kap_1)}}
  ^{\ewt[(v_0,\kap_0),(v_1,\kap_1)]}
  \\&=
  \prod_{\kap_0,\kap_1\in\cbr{1,2}}
  \babr{T_{i(v_0)n,j_{v_0}}^{(\kap_0)}(\bbf^{(v_0,\kap_0)}),
  T_{i(v_1)n,j_{v_1}}^{(\kap_1)}(\bbf^{(v_1,\kap_1)})}
  ^{\ewt[(v_0,\kap_0),(v_1,\kap_1)]}
\end{align*}
Thus we can rewrite the sum of $\beta_V(\bbf_{n,j}, \ewt)$ as
\begin{align}
  &\sum_{j\in\bbJ_n(V,i)}\beta_V(\bbf_{n,j}, \ewt)
  % \\&=
  % \sum_{j\in\bbJ_n(V,i)}
  % \beta_{V'}({\bbf'}_{n,j'}, \ewt\bbone_{\hp(V')})
  % \\&\hspace{50pt}\times
  % \prod_{\kap_0,\kap_1\in\cbr{1,2}}
  % \babr{T_{i(v_0)n,j_{v_0}}^{(\kap_0)}(\bbf^{(v_0,\kap_0)}),
  % T_{i(v_1)n,j_{v_1}}^{(\kap_1)}(\bbf^{(v_1,\kap_1)})}
  % ^{\ewt[(v_0,\kap_0),(v_1,\kap_1)]}
  \nn\\&=
  \sum_{j'\in\bbJ_n(V',i)}
  \sum_{j_0\in [i(v_0)n-1]}
  \beta_{V'}({\bbf'}_{n,j'}, \ewt\bbone_{\hp(V')})
  \nn\\&\hspace{50pt}\times
  \prod_{\kap_0,\kap_1\in\cbr{1,2}}
  \babr{T_{i(v_0)n,j_0}^{(\kap_0)}(\bbf^{(v_0,\kap_0)}),
  T_{i(v_1)n,j'_{v_1}}^{(\kap_1)}(\bbf^{(v_1,\kap_1)})}
  ^{\ewt[(v_0,\kap_0),(v_1,\kap_1)]}
  \nn\\&=
  \sum_{j'\in\bbJ_n(V',i)}
  \beta_{V'}({\bbf'}_{n,j'}, \ewt\bbone_{\hp(V')})
  \nn\\&\hspace{50pt}\times
  \sum_{j_0\in [i(v_0)n-1]}
  \prod_{\kap_0,\kap_1\in\cbr{1,2}}
  \babr{T_{i(v_0)n,j_0}^{(\kap_0)}(\bbf^{(v_0,\kap_0)}),
  T_{i(v_1)n,j'_{v_1}}^{(\kap_1)}(\bbf^{(v_1,\kap_1)})}
  ^{\ewt[(v_0,\kap_0),(v_1,\kap_1)]}.
  \label{eq:230928.1814}
\end{align}

Based on this expression, the following lemma gives an estimate of 
the sum of $\beta_V(\bbf_{n,j}, \ewt)$ 
for a tree-shaped weighted graph $C=(V,\ewt,0)$.
\begin{lemma}\label{221222.1901}
  Consider a tree-shaped weighted graph $C=(V,\ewt,0)$,
  $i:V\to[2]$ and $\bbf\in\calf(V)$.
  The following estimate holds as $n\to\infty$:
  \begin{align}
    \sum_{j\in\bbJ_n(V,i)}
    \abs{\beta_V(\bbf_{n,j}, \ewt)}
    = 
    O(n^{1-2H\:\bartheta(C) +(2H-1)\abs{E_1(C)}}),
    \label{221222.1602}
  \end{align}
  where 
  $\bartheta(C) = \sum_{[\tv,\tv']\in\tp(V)} \ewt([\tv,\tv'])
  =\sum_{[v,v']\in p(V)} \check\ewt([v,v'])$.
\end{lemma}

\begin{proof}
  We prove the estimate by induction with respect to $\abs{V}$.
  First consider the case $\abs{V}=1$, say, $V=\cbr{v_0}$.
  By definition,
  $\beta_V(\bbf_{n,j}, \ewt)=1$, and
  % $p(V)=\emptyset$ and hence $\ewt=0$, %we read $\beta_V(\bbf_{n,j}, \ewt)=1$, %
  the sum is estimated as
  \begin{align*}
    \sum_{j\in\bbJ_n(V,i)}
    \abs{\beta_V(\bbf_{n,j}, \ewt)}
    = 
    \sum_{j\in[i(v_0)n-1]} 1
    =O(n).
  \end{align*} 
  Since 
  $\bartheta(C)=0$ and $E_1(\ewt)=\emptyset$, 
  we have
  $1-2H\bartheta(C) +(2H-1)\abs{E_1(\ewt)}=1$
  and the estimate \eqref{221222.1602} holds.

  \vspssm 
  For given $d\geq1$,
  assume that the estimate \eqref{221222.1602} holds % for any $C=(V,\ewt,0)$ 
  in any case of $\abs{V}\leq d$.
  Consider a tree-shaped weighted graph $C=(V,\ewt,0)$ with $\abs{V}=d+1$.
  Since $\abs{V}\geq2$, 
  there are vertices 
  $v_0\neq v_1\in V$ such that 
  $\cewt([v_0,v_1])>0$ and 
  $\cewt([v_0,v])=0$ for $v\in V\setminus\cbr{v_0,v_1}$.
  We set
  $V' = V\setminus\cbr{v_0}$ and
  $\whv'=V'\times\cbr{1,2}$.
  % {\myred 
  % by an easy observation, 
  % there exist 
  % $v_0,v_1\in V$ such that 
  % $[v_0,v_1]\in E(C)(=E(V,\ewt) = \cbr{[v,\vpr]\in p(V)\mid \cewt([v,\vpr])>0})$ and 
  % $[v_0,v]\notin E(C)$
  % for any $v\in V\setminus\cbr{v_0,v_1}$.}
  %for any $v\neq v_0, v_1(\in V)$.
  % \sout{We set
  % $V' = V\setminus\cbr{v_0}$ and
  % $E' = E(C)\setminus\cbr{[v_0,v_1]}$.
  % We can see that 
  % $E'\subset p(V')$.}
  % \sout{We write 
  % $\bbf' = \bbf|_{\whv'}\in\calf(V')$.}
  Then by \eqref{eq:230928.1814}, we have
  \begin{align}
    &\sum_{j\in\bbJ_n(V,i)}\abs{\beta_V(\bbf_{n,j}, \ewt)}
    \nn\\&\leq
    \sup_{j_1\in [i(v_1)n-1]}
    \sum_{j_0\in [i(v_0)n-1]}
    \prod_{\kap_0,\kap_1\in\cbr{1,2}}
    \abs{\babr{T_{i(v_0)n,j_0}^{(\kap_0)}(\bbf^{(v_0,\kap_0)}),
    T_{i(v_1)n,j_1}^{(\kap_1)}(\bbf^{(v_1,\kap_1)})}}
    ^{\ewt[(v_0,\kap_0),(v_1,\kap_1)]}
    \nn\\&\hspace{40pt}\times
    \sum_{j'\in\bbJ_n(V',i)}
    \abs{\beta_{V'}({\bbf'}_{n,j'}, \ewt\bbone_{\hp(V')})}
    \label{221222.1730}
  \end{align}
  Let $C' = (V',\ewt\bbone_{\hp(V')},0)$ 
  and
  observe that $C'$ is also tree-shaped and $\abs{V'}=d$.
  We can apply the assumption of induction to
  $C'=(V', \ewt\bbone_{\hp(V')}, 0)$, $i|_{V'}$ and $\bbf'$
  to obtain the following estimate of the second factor of \eqref{221222.1730}
  as $n\to\infty$:
  \begin{align*}
    \sum_{j'\in\bbJ_n(V',i)}\abs{\beta_{V'}(\bbf'_{n,j'}, \ewt\bbone_{\hp(V')})}
    =O(n^{1-2H\bartheta(C')
    +(2H-1)\abs{E_1(C')}}).
    %=O(n^{1-2H(\bartheta(C)-\cewt([v_0,v_1]))
    %+(2H-1)\abs{E_1(\ewt\bbone_{\hp(V')})}}).
  \end{align*}
  Notice that $\bartheta(C')=\bartheta(C)-\cewt([v_0,v_1])$.

  We next estimate the order of the first factor of \eqref{221222.1730} depending on whether 
  $[v_0,v_1]\in E_1(C)$.

  \item[(i)] Consider the case $[v_0,v_1]\in E_1(C)
  =\cbr{[v,\vpr]\in p(V)\mid \cewt_1([v,\vpr])>0 \tandsm 
  \cewt_2([v,\vpr])=0}$.
  We have
  \begin{align*}
    &
    %\sup_{j_1\in[i(v_1)n-1]}
    \sup_{j_1\in\bbJ_n(\cbr{v_1},i)}
    %\rbr{\sum_{j_0\in[i(v_0)n-1]}
    \rbr{\sum_{j_0\in\bbJ_n(\cbr{v_0},i)}
    \prod_{\kap_0,\kap_1\in\cbr{1,2}}
    \abs{\babr{T_{i(v_0)n,j_{0}}^{(\kappa_0)}(\bbf^{(v_0,\kappa_0)}),
    T_{i(v_1)n,j_{1}}^{(\kappa_1)}(\bbf^{(v_1,\kappa_1)})}}
    ^{\ewt([(v_0,\kap_0),(v_1,\kap_1)])}}
    \\&=
    \sup_{j_1\in[i(v_1)n-1]}
    %\sup_{j_1\in\bbJ_n(\cbr{v_1},i)}
    \rbr{\sum_{j_0\in[i(v_0)n-1]}
    %\rbr{\sum_{j_0\in\bbJ_n(\cbr{v_0},i)}
    \abs{\babr{T_{i(v_0)n,j_{0}}^{(1)}(\bbf^{(v_0,1)}),
    T_{i(v_1)n,j_{1}}^{(1)}(\bbf^{(v_1,1)})}}
    ^{\ewt([(v_0,1),(v_1,1)])}}
    \\&\leq
    \sup_{\substack{j_0\in[i(v_0)n-1]\\j_1\in[i(v_1)n-1]}}
    %\sup_{\substack{j_0\in\bbJ_n(\cbr{v_0},i)\\j_1\in\bbJ_n(\cbr{v_1},i)}}
    \abs{\babr{T_{i(v_0)n,j_{0}}^{(1)}(\bbf^{(v_0,1)}),
    T_{i(v_1)n,j_{1}}^{(1)}(\bbf^{(v_1,1)})}}
    ^{\ewt([(v_0,1),(v_1,1)])-1}
    \\&\hspace{30pt}\times
    \sup_{j_1\in[i(v_1)n-1]}
    %\sup_{j_1\in\bbJ_n(\cbr{v_1},i)}
    \sum_{j_0\in[i(v_0)n-1]}
    %\sum_{j_0\in\bbJ_n(\cbr{v_0},i)}
    \abs{\babr{T_{i(v_0)n,j_{0}}^{(1)}(\bbf^{(v_0,1)}),
    T_{i(v_1)n,j_{1}}^{(1)}(\bbf^{(v_1,1)})}}
    \\&=
    O(n^{-2H\cbr{\ewt([(v_0,1),(v_1,1)])-1}-1})
  \end{align*}
  Here we used Lemma \ref{221222.1739} (i) and (ii).
  Since 
  $\cewt([v_0,v_1]) = \ewt([(v_0,1),(v_1,1)])$ and
  $\abs{E_1(C')} = \abs{E_1(C)}-1$,
  we obtain 
  \begin{align*}
    \sum_{j\in\bbJ_n(V,i)}\abs{\beta_V(\bbf_{n,j}, \ewt)}
    %\\
    &=
    O(n^{-2H\ncbr{\ewt([(v_0,1),(v_1,1)])-1}-1})
    \times
    O(n^{1-2H\ncbr{\bartheta(C)-\ewt([(v_0,1),(v_1,1)])}
    +(2H-1)(\abs{E_1(C)}-1)})
    \\&=
    O(n^{1-2H\bartheta(C)+(2H-1)\abs{E_1(C)}})
  \end{align*}

  \item[(ii)] Consider the case $[v_0,v_1]\in E_2(C)$.
  Then we can find $\kap'_0,\kap'_1\in\cbr{1,2}$ such that 
  $(\kap'_0,\kap'_1)\neq(1,1)$ and 
  $\ewt([(v_0,\kap'_0),(v_1,\kap'_1)])>0$.
  By a similar argument using Lemma \ref{221222.1739}(i) and (iii),
  \begin{align*}
    &
    \sup_{j_1\in\bbJ_n(\cbr{v_1},i)}
    \rbr{\sum_{j_0\in\bbJ_n(\cbr{v_0},i)}
    \prod_{\kap_0,\kap_1\in\cbr{1,2}}
    \abs{\babr{T_{i(v_0)n,j_{0}}^{(\kappa_0)}(\bbf^{(v_0,\kappa_0)}),
    T_{i(v_1)n,j_{1}}^{(\kappa_1)}(\bbf^{(v_1,\kappa_1)})}}
    %\abs{\babr{\bff_{n,j_0}^{(v_0,\kap_0)},\bff_{n,j_1}^{(v_1,\kap_1)}}}
    ^{\ewt([(v_0,\kap_0),(v_1,\kap_1)])}}
    \\&\leq
    \sup_{\substack{j_0\in[i(v_0)n-1]\\j_1\in[i(v_1)n-1]}}
    %\sup_{\substack{j_0\in\bbJ_n(\cbr{v_0},i)\\j_1\in\bbJ_n(\cbr{v_1},i)}}
    \prod_{\kap_0,\kap_1\in\cbr{1,2}}
    \abs{\babr{T_{i(v_0)n,j_{0}}^{(\kappa_0)}(\bbf^{(v_0,\kappa_0)}),
    T_{i(v_1)n,j_{1}}^{(\kappa_1)}(\bbf^{(v_1,\kappa_1)})}}
    %\abs{\babr{\bff_{n,j_0}^{(v_0,\kap_0)},\bff_{n,j_1}^{(v_1,\kap_1)}}}
    ^{(\ewt-\bbone_{[(v_0,\kap'_0),(v_1,\kap'_1)]})([(v_0,\kap_0),(v_1,\kap_1)])}
    %\abr{\bbf_{n,j_0}^{(v_0,1)},\bbf_{n,j_1}^{(v_1,1)}}
    %^{\ewt([(v_0,1),(v_1,1)])-1}
    \\&\hspace{50pt}\times
    \sup_{j_1\in[i(v_1)n-1]}
    %\sup_{j_1\in\bbJ_n(\cbr{v_1},i)}
    \sum_{j_0\in[i(v_0)n-1]}
    %\sum_{j_0\in\bbJ_n(\cbr{v_0},i)}
    \abs{\babr{T_{i(v_0)n,j_{0}}^{(\kappa'_0)}(\bbf^{(v_0,\kappa'_0)}),
    T_{i(v_1)n,j_{1}}^{(\kappa'_1)}(\bbf^{(v_1,\kappa'_1)})}}
    %\abs{\babr{\bff_{n,j_0}^{(v_0,\kap'_0)},\bff_{n,j_1}^{(v_1,\kap'_1)}}}
    \\&=
    O(n^{-2H(\cewt([v_0,v_1])-1)-2H})=
    O(n^{-2H\cewt([v_0,v_1])}).
  \end{align*}
  Since $\abs{E_1(C')} = \abs{E_1(C)}$, 
  we obtain 
  \begin{align*}
    \sum_{j\in\bbJ_n(V,i)}\abs{\beta_V(\bbf_{n,j}, \ewt)}
    %\\
    &=
    O(n^{-2H\cewt([v_0,v_1])})
    \times
    O(n^{1-2H(\bartheta(C)-\cewt([v_0,v_1]))
    +(2H-1)\abs{E_1(C)}})
    \\&=
    O(n^{1-2H\bartheta(C)
    +(2H-1)\abs{E_1(C)}}).
  \end{align*}

  Hence, in both of the cases (i) and (ii),
  the estimate \eqref{221222.1602} holds for $V$ with $\abs{V}=d+1$,
  and by induction for a general $V$.
\end{proof}

The next lemma is essential to prove Proposition \ref{prop:230605.2133}, 
since the estimate of the expectation of the functional $\calii_n\rbr{G, A,\bbf}$ 
can be reduced by induction to the following estimate.
\begin{lemma}\label{230519.1737}
  For a connected weighted graph $C=(V,\ewt,0)$, $i:V\to[2]$ and $\bbf\in\calf(V)$,
  the following estimate holds 
  \begin{align*}
    \sum_{j\in\bbJ_n(V,i)} \abs{\beta(C, \bbf_{n,j})} %\abs{\beta_V(\bbf_{n,j}, \ewt)}
    \bbrbr{=\sum_{j\in\bbJ_n(V,i)} \abs{\beta_V(\bbf_{n,j}, \ewt)}}
    =O(n^{e(C)})
  \end{align*}
  as $n\to\infty$,
  where $e(C)$ is defined at \eqref{def:231005.1440}.
\end{lemma}

\begin{proof}
  By the connectedness of $C$,
  let $\bbT\subset E(C)$ be a spanning tree of $C$ such that $\ell_2(C)=\abs{E_2(C)\cap\bbT}$, 
  that is $\bbT\in\argmax_{\bbT'}\abs{E_2(C)\cap\bbT'}$, 
  where $\bbT'$ varies in the set of all the spanning tree of $C$.
  Then 
  \begin{align*}
    \sum_{j\in\bbJ_n(V,i)}\babs{\beta(C, \bbf_{n,j})}
    =\sum_{j\in\bbJ_n(V,i)}\babs{\beta_V(\bbf_{n,j}, \ewt)}
    &=
    \sum_{j\in\bbJ_n(V,i)}\babs{\beta_V(\bbf_{n,j}, \ewt\bbone_{\wh{\bbT}})}\:
    \babs{\beta_V(\bbf_{n,j}, \ewt\bbone_{\tp(V)\setminus\wh{\bbT}})}
    \\&\leq
    \sup_{j\in\bbJ_n(V,i)} 
    \babs{\beta_V(\bbf_{n,j}, \ewt\bbone_{\tp(V)\setminus\wh{\bbT}})}
    \sum_{j\in\bbJ_n(V,i)} \babs{\beta_V(\bbf_{n,j}, \ewt\bbone_{\wh{\bbT}})}.
  \end{align*}
  Here, $\wh\bbT=\pi^{-1}(\bbT)$ with the projection $\pi:\hp(V)\to p(V)$ defined at \eqref{230621.1526}. 
  By Lemma \ref{221222.1739} (i), we have
  \begin{align*}
    \sup_{j\in\bbJ_n(V,i)} 
    \abs{\beta_V(\bbf_{n,j}, \ewt\bbone_{\tp(V)\setminus\wh{\bbT}})}
    &= 
    \sup_{j\in\bbJ_n(V,i)} 
    \prod_{[\tv,\tvp]\in\tp(V)}
    \abs{\babr{\bbf_{n,j}^{(\tv)}, \bbf_{n,j}^{(\tvp)}}}
    ^{(\ewt\bbone_{\tp(V)\setminus\wh{\bbT}})[\tv,\tvp]}
    \\&= 
    \sup_{j\in\bbJ_n(V,i)} 
    \prod_{[\tv,\tvp]\in\tp(V)\setminus\wh{\bbT}}
    \abs{\babr{\bbf_{n,j}^{(\tv)}, \bbf_{n,j}^{(\tvp)}}}
    ^{\ewt[\tv,\tvp]}
    =O(n^{-2H\sum_{[v,\vpr]\in p(V)\setminus \bbT}\cewt([v,\vpr])}).
  \end{align*}  
  The weighted graph defined by $C' = (V,\ewt\bbone_{\wh{\bbT}},0)$
  is obviously tree-shaped
  ($\because E(C')=E(V,\ewt\bbone_{\wh{\bbT}})=\bbT$),
  and 
  we have from Lemma \ref{221222.1901}
  \begin{align*}
    \sum_{j\in\bbJ_n(V,i)}\abs{\beta_V(\bbf_{n,j},  \ewt\bbone_{\wh{\bbT}})}
    =O(n^{1-2H\bartheta(C')+(2H-1)\abs{E_1(C')}}).
  \end{align*}
  Hence, the following estimate holds
  \begin{align}
    \sum_{j\in\bbJ_n(V,i)}\abs{\beta_V(\bbf_{n,j}, \ewt)}
    =O(n^{1-2H\bartheta(C)+(2H-1)\abs{E_1(C')}}),
    \label{221222.1913}
  \end{align}
  since 
  $\bartheta(C') + \sum_{[v,\vpr]\in p(V)\setminus \bbT}\cewt([v,\vpr])
  = \sum_{[v,\vpr]\in p(V)}\cewt([v,\vpr]) = \bartheta(C)$.

  We can show
  $E_2(C') = E_2(V,\ewt\bbone_{\wh{\bbT}}) = E_2(V,\ewt)\cap\bbT = E_2(C)\cap\bbT$,
  $\abs{E(C')}=\abs{\bbT}=\abs{V}-1$, and hence
  \begin{align*}
    \abs{E_1(C')} = \abs{E(C')} - \abs{E_2(C')}
    =\abs{V}-1 - \abs{E_2(C)\cap\bbT}
    =\abs{V}-1 - \ell_2(C).
  \end{align*}
  Thus, we obtain 
  \begin{align*}
    \sum_{j\in\bbJ_n(V,i)}\abs{\beta_V(\bbf_{n,j}, \ewt)}
    =O(n^{1-2H\bartheta(C)+(2H-1)(\abs{V}-1 - \ell_2(C))})
    =O(n^{e(C)}),
  \end{align*}
  since $e_q(C)=0$.
\end{proof}

Recall that for a function $f:\bbR\to\bbR$, $n\in\bbN$ and $j\in\bbZ$,
we have defined
$T_{n,j}(f)(x)=f(nx-j)$ for $x\in\bbR$,
$T_{n,j}^{(1)}(f) = T_{n,j}(f)$ and 
$T_{n,j}^{(2)}(f) = T_{n,j+1}(f) - T_{n,j}(f)$. 
We define
$\phi_s(x)=x-s$ for $s\in\bbR$ and 
$\psi_n(x)=nx$ for $n\in\bbN$.
We have
\begin{align*}
  T_{n,j}(f) = f\circ\phi_j\circ\psi_n = f\circ\psi_n\circ\phi_{j/n}\tand
  T^{(\kappa)}_{m,j}(f) = T_{1,0}^{(\kappa)}(f)\circ\phi_j\circ\psi_{m} \tfor \kappa=1,2.
\end{align*} 
%
%
% Though the Hilbert space $\calh$ contains 
% functions defined on $[0,1]$, 
For notational convenience we extend the inner product of $\calh$ for 
compactly supported bounded functions $\varphi$ and $\psi$ on $\bbR$, namely 
\begin{align}\label{def:231005.2143}
  \abr{\varphi,\psi}=
  c_H\int_{\bbR^2}\varphi(t)\psi(t')\abs{t-t'}^{2H-2}dtdt',
\end{align}
where $c_H=H(2H-1)$.
Then we have 
$\abr{f\circ\psi_{m}, g\circ\psi_{m}} = m^{-2H}\abr{f, g}$, and 
$\abr{f\circ\phi_{s_0}, g\circ\phi_{s_1}}=\abr{f, g\circ\phi_{s_1-s_0}}$.

\begin{lemma}\label{221222.1739}
  Suppose that $f^{(\ell)}\in L^{\infty}(\bbR)$ $(\ell=0,1)$ such that 
  $\supp (f^{(\ell)})\subset[\alpha_\ell, \alpha_\ell+1]$ with some $\alpha_\ell\in\bbR$.
  Then for any $i_0,i_1\in\cbr{1,2}$ the following estimates hold as $n\to\infty$:
  \item [(i)] For any $\kappa_0, \kappa_1\in\cbr{1,2}$,
  $\sup_{j_0,j_1\in\bbZ}
  \abs{\babr{T^{(\kappa_0)}_{i_0n,j_0}(f^{(0)}), T^{(\kappa_1)}_{i_1n,j_1}(f^{(1)})}}
  =O(n^{-2H})$.
  \item [(ii)] %For any $i_0,i_1\in\cbr{1,2}$,
  %\sup_{j_1\in[i_1 n] \redc{\bbZ??}}\sum_{j_0\in[i_0 n]}
  $\sup_{j_1\in\bbZ}\sum_{j_0\in[i_0 n]}
  \abs{\babr{T^{(1)}_{i_0n,j_0}(f^{(0)}), T^{(1)}_{i_1n,j_1}(f^{(1)})}}
  =O(n^{-1})$. 
  \item [(iii)] For any %$i_0,i_1\in\cbr{1,2}$ and 
  $\kappa_0, \kappa_1\in\cbr{1,2}$ such that $(\kappa_0,\kappa_1)\neq(1,1)$,
  $\sup_{j_1\in\bbZ}\sum_{j_0\in\bbZ}
  %$\sup_{j_1\in[i_1 n] \redc{\bbZ??}}\sum_{j_0\in[i_0 n]\redc{\bbZ??}}
  \abs{\babr{T^{(\kappa_0)}_{i_0n,j_0}(f^{(0)}), T^{(\kappa_1)}_{i_1n,j_1}(f^{(1)})}}
  =O(n^{-2H})$.
\end{lemma}

\begin{proof}
  In general, 
  for $\kappa_0,\kappa_1\in\cbr{1,2}$, $j_0,j_1\in\bbZ$ and $n\in\bbN$,
  we have the following relation:
  \begin{align}
    \babr{T^{(\kappa_0)}_{i_0n,j_0}(f^{(0)}), T^{(\kappa_1)}_{i_1n,j_1}(f^{(1)})}
    &=
    \babr{T_{1,0}^{(\kappa_0)}(f^{(0)})\circ\psi_{i_0}\circ\phi_{j_0/i_0}\circ\psi_{n},\;
    T_{1,0}^{(\kappa_1)}(f^{(1)})\circ\psi_{i_1}\circ\phi_{j_1/i_1}\circ\psi_{n}}
    \nn\\&=
    n^{-2H}
    \babr{T_{1,0}^{(\kappa_0)}(f^{(0)})\circ\psi_{i_0}\circ\phi_{j_0/i_0},\;
    T_{1,0}^{(\kappa_1)}(f^{(1)})\circ\psi_{i_1}\circ\phi_{j_1/i_1}}
    %\\&=
    %n^{-2H}
    %\babr{T_{1,0}^{(\kappa_0)}(f^{(0)})\circ\psi_{i_0},\;
    %T_{1,0}^{(\kappa_1)}(f^{(1)})\circ\psi_{i_1}\circ\phi_{j_1/i_1-j_0/i_0}}
    \nn\\&=
    n^{-2H}
    \babr{T_{i_0,0}^{(\kappa_0)}(f^{(0)}),\;
    T_{i_1,0}^{(\kappa_1)}(f^{(1)})\circ\phi_{j_1/i_1-j_0/i_0}}.
    \label{eq:230524.1818}
  \end{align}

  \item [(i)] By \eqref{eq:230524.1818}, it suffices to show that 
  \begin{align}
    \sup_{j_0,j_1\in\bbZ}
    \abs{\babr{T_{i_0,0}^{(\kappa_0)}(f^{(0)}),\;
    T_{i_1,0}^{(\kappa_1)}(f^{(1)})\circ\phi_{j_1/i_1-j_0/i_0}}}
    <\infty
    \label{eq:230524.1300}
  \end{align}
  for any $\kappa_0,\kappa_1\in\cbr{1,2}$ and $i_0, i_1\in\cbr{1,2}$.
  Since $\supp(T_{i_\ell,0}^{(\kappa_\ell)}(f^{(\ell)})) 
  \subset i_\ell^{-1} [\alpha_\ell, \alpha_\ell+\kappa_\ell]$ for $\ell=0,1$,
  we have
  \begin{align*}
    &\abs{\babr{T_{i_0,0}^{(\kappa_0)}(f^{(0)}),\;
    T_{i_1,0}^{(\kappa_1)}(f^{(1)})\circ\phi_{j_1/i_1-j_0/i_0}}}
    \\&\leq
    c_H\int_{\bbR^2}
    \abs{T_{i_0,0}^{(\kappa_0)}(f^{(0)})(t)}
    \abs{T_{i_1,0}^{(\kappa_1)}(f^{(1)})\circ\phi_{j_1/i_1-j_0/i_0}(t')}
    \abs{t-t'}^{2H-2} dtdt'
    \\&\leq
    c_H\int_{\bbR^2}
    \bnorm{f^{(0)}}_\infty 
    \bnorm{f^{(1)}}_\infty 
    \bbone_{i_0^{-1} [\alpha_0, \alpha_0+\kappa_0]}(t)\;
    \bbone_{i_1^{-1} [\alpha_1, \alpha_1+\kappa_1]}(t'- (j_1/i_1-j_0/i_0))\;
    \abs{t-t'}^{2H-2} dtdt'
    \\&\leq
    c_H \bnorm{f^{(0)}}_\infty \bnorm{f^{(1)}}_\infty 
    \int_{-1}^{1}\abs{t}^{2H-2}dt
    \times\frac{\kappa_1}{i_1}
    <\infty
  \end{align*}
  for any $j_0,j_1\in\bbZ$.

  \item [(ii)]
  Since $\supp (T^{(1)}_{i_\ell n,j_\ell}(f^{(\ell)})) = (i_\ell n)^{-1}(\supp (f^{(\ell)})+j_\ell)
  \subset(i_\ell n)^{-1}[\alpha_\ell + j_\ell, \alpha_\ell + j_\ell + 1]$, 
  we have
  $\abs{T^{(1)}_{i_\ell n,j_\ell}(f^{(\ell)})(t)}\leq 
  \norm{f^{(\ell)}}_\infty
  \bbone_{[\frac{\alpha_\ell + j_\ell}{i_\ell n}, \frac{\alpha_\ell + j_\ell + 1}{i_\ell n}]}(t)$.
  %\bbone_{[(i_0n)^{-1}(\alpha_0 + j_0), (i_0n)^{-1}(\alpha_0 + j_0 + 1)]}$
  Therefore, the following estimate holds
  \begin{align*}
    &%\sup_{j_1\in[i_1 n] \redc{\bbZ??}}
    \sum_{j_0\in[i_0 n]}
    \abs{\babr{T^{(1)}_{i_0n,j_0}(f^{(0)}), T^{(1)}_{i_1n,j_1}(f^{(1)})}}
    % \\&\leq
    % \sum_{j_0\in[i_0 n]}
    % c_H\int_{\bbR^2\myred} 
    % \abs{T^{(1)}_{i_0n,j_0}(f^{(0)})(t)}\;\abs{T^{(1)}_{i_1n,j_1}(f^{(1)})(t')}\;
    % \abs{t-t'}^{2H-2} dtdt'
    \\&\leq
    \sum_{j_0\in[i_0 n]}
    c_H \bnorm{f^{(0)}}_\infty \bnorm{f^{(1)}}_\infty
    \int_{\bbR^2\myred} 
    \bbone_{[\frac{\alpha_0 + j_0}{i_0 n}, \frac{\alpha_0 + j_0 + 1}{i_0 n}]}(t)
    \bbone_{[\frac{\alpha_1 + j_1}{i_1 n}, \frac{\alpha_1 + j_1 + 1}{i_1 n}]}(t')
    %\abs{T^{(1)}_{i_0n,j_0}(f^{(0)})(t)}\;
    %\abs{T^{(1)}_{i_1n,j_1}(f^{(1)})(t')}\;
    \abs{t-t'}^{2H-2} dtdt'
    \\&=
    c_H\bnorm{f^{(0)}}_\infty \bnorm{f^{(1)}}_\infty
    \int_{\bbR^2\myred} 
    \bbone_{[\frac{\alpha_0 + 1}{i_0 n}, \frac{\alpha_0 + 1}{i_0 n} + 1]}(t)
    \bbone_{[\frac{\alpha_1 + j_1}{i_1 n}, \frac{\alpha_1 + j_1 + 1}{i_1 n}]}(t')
    \abs{t-t'}^{2H-2} dtdt'
    \\&<
    %c_H\bnorm{f^{(0)}}_\infty \bnorm{f^{(1)}}_\infty
    %\int_{-1}^1\abs{t}^{2H-2} dt
    %\int_{\bbR} 
    %\bbone_{[\frac{\alpha_1 + j_1}{i_1 n}, \frac{\alpha_1 + j_1 + 1}{i_1 n}]}(t')dt'
    %\\&=
    c_H\bnorm{f^{(0)}}_\infty \bnorm{f^{(1)}}_\infty
    \int_{-1}^1\abs{t}^{2H-2} dt
    \times\rbr{i_1 n}^{-1},
  \end{align*}
  for any $j_1\in\bbZ$,
  and we obtain (ii).

  \item [(iii)] 
  By the relation \eqref{eq:230524.1818},
  it suffices to show that 
  \begin{align}\label{eq:230524.1925}
    \sup_{j_1\in\bbZ}\sum_{j_0\in\bbZ}
    \abs{\babr{T_{i_0,0}^{(\kappa_0)}(f^{(0)}),\;
    T_{i_1,0}^{(\kappa_1)}(f^{(1)})\circ\phi_{j_1/i_1-j_0/i_0}}}
    <\infty.
  \end{align}

  First assume that $\kappa_0=2$.
  By change of variables, we have
  \begin{align*}
    &\babr{T_{i_0,0}^{(2)}(f^{(0)}),\;
    T_{i_1,0}^{(\kappa_1)}(f^{(1)})\circ\phi_{j_1/i_1-j_0/i_0}}
    % \\&=
    % c_H\int_{\bbR^2}
    % T_{i_0,0}^{(2)}(f^{(0)})(t)\;
    % T_{i_1,0}^{(\kappa_1)}(f^{(1)})\circ\phi_{j_1/i_1-j_0/i_0}(t')
    % \abs{t-t'}^{2H-2} dt dt'
    % \\&=
    % {c_H}\int_{\bbR^2}
    % T_{i_0,0}^{(2)}(f^{(0)})(t)\;
    % T_{i_1,0}^{(\kappa_1)}(f^{(1)})(t' - \rbr{j_1/i_1-j_0/i_0})
    % \abs{t-t'}^{2H-2} dt dt'
    \\&=
    {c_H}\int_{\bbR^2}
    T_{i_0,0}^{(2)}(f^{(0)})(t)\;
    T_{i_1,0}^{(\kappa_1)}(f^{(1)})(s)\;
    \abs{t-(s + \rbr{j_1/i_1-j_0/i_0})}^{2H-2} dt ds
    % \\&{\myred=
    % {c_H}\int_{\bbR^2}
    % \rbr{T_{i_0,0}(f^{(0)})(t-1/i_0) - T_{i_0,0}(f^{(0)})(t)}\;
    % T_{i_1,0}^{(\kappa_1)}(f^{(1)})(s)\;
    % \abs{t-(s + \rbr{j_1/i_1-j_0/i_0})}^{2H-2} dt ds}
    % \\&=
    % {c_H}\int_{\bbR^2}
    % T_{i_0,0}(f^{(0)})(t-1/i_0)\;
    % T_{i_1,0}^{(\kappa_1)}(f^{(1)})(s)\;
    % \abs{t-(s + \rbr{j_1/i_1-j_0/i_0})}^{2H-2} dt ds
    % \\&\quad-
    % { c_H}\int_{\bbR^2}
    % T_{i_0,0}(f^{(0)})(t)\;
    % T_{i_1,0}^{(\kappa_1)}(f^{(1)})(s)\;
    % \abs{t-(s + \rbr{j_1/i_1-j_0/i_0})}^{2H-2} dt ds
    % \\&=
    % { c_H}\int_{\bbR^2}
    % T_{i_0,0}(f^{(0)})(t)\;
    % T_{i_1,0}^{(\kappa_1)}(f^{(1)})(s)\;
    % \abs{t+1/i_0-(s + \rbr{j_1/i_1-j_0/i_0})}^{2H-2} dt ds
    % \\&\quad-
    % { c_H}\int_{\bbR^2}
    % T_{i_0,0}(f^{(0)})(t)\;
    % T_{i_1,0}^{(\kappa_1)}(f^{(1)})(s)\;
    % \abs{t-(s + \rbr{j_1/i_1-j_0/i_0})}^{2H-2} dt ds
    \\&=
    { c_H}\int_{\bbR^2}
    T_{i_0,0}(f^{(0)})(t)\;
    T_{i_1,0}^{(\kappa_1)}(f^{(1)})(s)\;
    \\&\hspace{60pt}\times
    \rbr{\abs{t-s - j_1/i_1+(j_0+1)/i_0}^{2H-2} -
    \abs{t-s - j_1/i_1+j_0/i_0}^{2H-2}}
    dt ds
  \end{align*}
  Since $\supp (f^{(l)})\subset[\alpha_l, \alpha_l+1]$ with some $\alpha_l\in\bbR$, 
  we have\;
  \begin{align*}
    &\supp (T_{i_0,0}(f^{(0)}))\subset[\alpha_0/i_0, (\alpha_0+1)/i_0],\qquad
    \supp (T_{i_1,0}^{(\kappa_1)}(f^{(1)}))\subset[\alpha_1/i_1, (\alpha_1+\kappa_1)/i_1],
    %{ $\supp (T_{i_0,0}(f^{(0)})(\cdot)\;T_{i_1,0}^{(\kappa_1)}(f^{(1)})(\cdot))
    %\subset[\alpha_0/i_0, (\alpha_0+1)/i_0]\times[\alpha_1/i_1, (\alpha_1+\kappa_1)/i_1]
    %=:A_{i_0,i_1}$.}
  \end{align*}
  and 
  \begin{align*}
    &\supp (T_{i_0,0}(f^{(0)})(\cdot)\;T_{i_1,0}^{(\kappa_1)}(f^{(1)})(\cdot))
    \subset[\alpha_0/i_0, (\alpha_0+2)/i_0]\times[\alpha_1/i_1, (\alpha_1+2)/i_1]
    =:A_{i_0,i_1}.
    %{ Then $t-s\in[\alpha_0/i_0-(\alpha_1+\kappa_1)/i_1, (\alpha_0+1)/i_0-\alpha_1/i_1]
    %=:[\alpha_{(-)}, \alpha_{(+)}]$
    %for $(t,s)\in A_{i_0,i_1}$.}
  \end{align*}
  Then 
  for $(t,s)\in A_{i_0,i_1}$, we have
  $t-s\in[\alpha_0/i_0-(\alpha_1+2)/i_1, (\alpha_0+2)/i_0-\alpha_1/i_1]
  =:[\alpha_{(-)}, \alpha_{(+)}]$.
  Hence for $j_0,j_1\in\bbZ$, we obtain
  \begin{align*}
    &\abs{\babr{T_{i_0,0}^{(2)}(f^{(0)}),\;
    T_{i_1,0}^{(\kappa_1)}(f^{(1)})\circ\phi_{j_1/i_1-j_0/i_0}}}
    \\&\leq
    { c_H}\bnorm{f^{(0)}}_\infty\bnorm{f^{(1)}}_\infty
    \int_{A_{i_0,i_1}}
    \abs{\abs{t-s - j_1/i_1+(j_0+1)/i_0}^{2H-2} -
    \abs{t-s - j_1/i_1+j_0/i_0}^{2H-2}}
    dt ds.
  \end{align*}

  We have the following estimate:
  \begin{align*}
    \abs{\abs{x_0+\Delta x}^{2H-2} - \abs{x_0}^{2H-2}}
    &\leq
    (2-2H) {x_0}^{2H-3} \times \Delta x 
    &&\text{for }x_0>0 \tandsm \Delta x>0
    \\
    \abs{\abs{x_0}^{2H-2} - \abs{x_0-\Delta x}^{2H-2}}
    &\leq
    (2-2H) \abs{x_0}^{2H-3} \times \Delta x
    &&\text{for }x_0<0 \tandsm \Delta x>0
  \end{align*}
  Hence if $(j_0,j_1)$ satisfies $\alpha_{(-)} - j_1/i_1+j_0/i_0>1$, 
  then for $(t,s)\in A_{i_0,i_1}$, 
  which implies $t-s\geq\alpha_{(-)}$,
  we have
  % $t-s - j_1/i_1+j_0/i_0>1$ and 
  \begin{align*}
    &\quad\abs{\abs{t-s - j_1/i_1+(j_0+1)/i_0}^{2H-2} - \abs{t-s - j_1/i_1+j_0/i_0}^{2H-2}}
    \\&\leq 
    (2-2H)\rbr{t-s - j_1/i_1 + j_0/i_0}^{2H-3} \times \frac1{i_0}
    \\&\leq 
    (2-2H)\rbr{\alpha_{(-)} - j_1/i_1 + j_0/i_0}^{2H-3} \times \frac1{i_0}.
  \end{align*}
  Similarly, 
  if $(j_0,j_1)$ satisfies $\alpha_{(+)} - j_1/i_1+(j_0+1)/i_0<-1$, 
  then for $(t,s)\in A_{i_0,i_1}$, 
  which implies $t-s\leq\alpha_{(+)}$,
  we have
  % $t-s - j_1/i_1+(j_0+1)/i_0<-1$ and 
  \begin{align*}
    &\quad\abs{\abs{t-s - j_1/i_1+(j_0+1)/i_0}^{2H-2} - \abs{t-s - j_1/i_1+j_0/i_0}^{2H-2}}
    % \\&\leq 
    % (2-2H)\abs{t-s - \frac{j_1}{i_1}+\frac{j_0+1}{i_0}}^{2H-3} \times \frac1{i_0}
    % \\&\leq 
    \leq 
    (2-2H)\abs{\alpha_{(+)} - \frac{j_1}{i_1}+\frac{j_0+1}{i_0}}^{2H-3} \times \frac1{i_0}.
  \end{align*}
  As a result, 
  \begin{align*}
    &\sup_{j_1\in\bbZ}\sum_{j_0\in\bbZ}
    \abs{\babr{T_{i_0,0}^{(2)}(f^{(0)}),\;
    T_{i_1,0}^{(\kappa_1)}(f^{(1)})\circ\phi_{j_1/i_1-j_0/i_0}}}
    \\&\leq
    {c_H}\bnorm{f^{(0)}}_\infty\bnorm{f^{(1)}}_\infty \abs{A_{i_0,i_1}}\;
    (2-2H){i_0}^{-1}\;
    \sup_{j_1\in\bbZ}
    \sum_{j_0\in\bbZ: \alpha_{(-)} - \frac{j_1}{i_1} + \frac{j_0}{i_0}>1}
    \rbr{\alpha_{(-)} - \frac{j_1}{i_1}+\frac{j_0}{i_0}}^{2H-3} 
    \\&\quad+
    {c_H}\bnorm{f^{(0)}}_\infty\bnorm{f^{(1)}}_\infty \abs{A_{i_0,i_1}}\;
    (2-2H){i_0}^{-1}\;
    \sup_{j_1\in\bbZ}
    \sum_{j_0\in\bbZ: \alpha_{(+)} - \frac{j_1}{i_1}+\frac{j_0+1}{i_0}<-1}
    \abs{\alpha_{(+)} - \frac{j_1}{i_1}+\frac{j_0+1}{i_0}}^{2H-3}
    \\&\quad+
    \sup_{j_1\in\bbZ}
    \sum_{\substack{j_0\in\bbZ:
    \alpha_{(-)} - \frac{j_1}{i_1}+\frac{j_0}{i_0}\leq1 \tandsm\\\quad\quad
    \alpha_{(+)} - \frac{j_1}{i_1}+\frac{j_0+1}{i_0}\geq-1}}
    \abs{\babr{T_{i_0,0}^{(2)}(f^{(0)}),\;
    T_{i_1,0}^{(\kappa_1)}(f^{(1)})\circ\phi_{\frac{j_1}{i_1}-\frac{j_0}{i_0}}}}
    \\&\leq 
    2{c_H}\bnorm{f^{(0)}}_\infty\bnorm{f^{(1)}}_\infty \abs{A_{i_0,i_1}}\;
    (2-2H){i_0}^{-1}\;
    \sum_{j\in\bbZ: \frac{j}{i_0}\geq1} \rbr{\frac{j}{i_0}}^{2H-3} 
    \\&\quad+
    N_0
    \sup_{j_0,j_1\in\bbZ}
    \abs{\babr{T_{i_0,0}^{(2)}(f^{(0)}),\;
    T_{i_1,0}^{(\kappa_1)}(f^{(1)})\circ\phi_{\frac{j_1}{i_1}-\frac{j_0}{i_0}}}}
    %\abs{\cbr{j_0\in\bbZ\;\Big|\;
    %\alpha_{(-)} - \frac{j_1}{i_1}+\frac{j_0}{i_0}\leq1 \tandsm
    %\alpha_{(+)} - \frac{j_1}{i_1}+\frac{j_0+1}{i_0}\geq-1}}
    \\&<\infty
  \end{align*}
  We can see that the number of $j_0\in\bbZ$ satisfying
  $\alpha_{(-)} - j_1/i_1+j_0/i_0\leq1$ and 
  $\alpha_{(+)} - j_1/i_1+(j_0+1)/i_0\geq-1$ can be bounded by some constant $N_0$ independent of $j_1$.
  % ($\because$ The conditions above is equivalent to 
  % $i_0(-1-\alpha_{(+)})-1 + j_1i_0/i_1\leq
  % j_0\leq i_0(1 -\alpha_{(-)})   + j_1i_0/i_1$.)
  We also used the estimate \eqref{eq:230524.1300}.

  \vspssm
  We will prove \eqref{eq:230524.1925} for $\kappa_0=1$ and $\kappa_1=2$.
  By exchanging the index,
  we will prove 
  \begin{align*}
    \sup_{j_0\in\bbZ}\sum_{j_1\in\bbZ}
    \abs{\babr{T_{i_0,0}^{(2)}(f^{(0)}),\;
    T_{i_1,0}^{(1)}(f^{(1)})\circ\phi_{j_1/i_1-j_0/i_0}}}
    <\infty.
  \end{align*}
  By reusing the estimates about 
  $\abs{\babr{T_{i_0,0}^{(2)}(f^{(0)}),\;
  T_{i_1,0}^{(\kappa_1)}(f^{(1)})\circ\phi_{j_1/i_1-j_0/i_0}}}$,
  \begin{align*}
    &\sup_{j_0\in\bbZ}\sum_{j_1\in\bbZ}
    \abs{\babr{T_{i_0,0}^{(2)}(f^{(0)}),\;
    T_{i_1,0}^{(1)}(f^{(1)})\circ\phi_{j_1/i_1-j_0/i_0}}}
    \\&\leq
    {c_H}\bnorm{f^{(0)}}_\infty\bnorm{f^{(1)}}_\infty \abs{A_{i_0,i_1}}\;
    (2-2H){i_0}^{-1}\;
    \sup_{j_0\in\bbZ}
    \sum_{j_1\in\bbZ: \alpha_{(-)} - \frac{j_1}{i_1} + \frac{j_0}{i_0}>1}
    \rbr{\alpha_{(-)} - \frac{j_1}{i_1}+\frac{j_0}{i_0}}^{2H-3} 
    \\&\quad+
    {c_H}\bnorm{f^{(0)}}_\infty\bnorm{f^{(1)}}_\infty \abs{A_{i_0,i_1}}\;
    (2-2H){i_0}^{-1}\;
    \sup_{j_0\in\bbZ}
    \sum_{j_1\in\bbZ: \alpha_{(+)} - \frac{j_1}{i_1}+\frac{j_0+1}{i_0}<-1}
    \abs{\alpha_{(+)} - \frac{j_1}{i_1}+\frac{j_0+1}{i_0}}^{2H-3}
    \\&\quad+
    \sup_{j_0\in\bbZ}
    \sum_{\substack{j_1\in\bbZ:
    \alpha_{(-)} - \frac{j_1}{i_1}+\frac{j_0}{i_0}\leq1 \tandsm\\\quad\quad
    \alpha_{(+)} - \frac{j_1}{i_1}+\frac{j_0+1}{i_0}\geq-1}}
    \abs{\babr{T_{i_0,0}^{(2)}(f^{(0)}),\;
    T_{i_1,0}^{(\kappa_1)}(f^{(1)})\circ\phi_{\frac{j_1}{i_1}-\frac{j_0}{i_0}}}}
    \\&\leq 
    2{c_H}\bnorm{f^{(0)}}_\infty\bnorm{f^{(1)}}_\infty \abs{A_{i_0,i_1}}\;
    (2-2H){i_0}^{-1}\;
    \sum_{j\in\bbZ: \frac{j}{i_1}\geq1} \rbr{\frac{j}{i_1}}^{2H-3} 
    \\&\quad+
    N_1
    \sup_{j_0,j_1\in\bbZ}
    \abs{\babr{T_{i_0,0}^{(2)}(f^{(0)}),\;
    T_{i_1,0}^{(\kappa_1)}(f^{(1)})\circ\phi_{\frac{j_1}{i_1}-\frac{j_0}{i_0}}}}
    \\&<\infty
  \end{align*}
  Again we can see that the number of $j_1\in\bbZ$ satisfying
  $\alpha_{(-)} - j_1/i_1+j_0/i_0\leq1$ and 
  $\alpha_{(+)} - j_1/i_1+(j_0+1)/i_0\geq-1$ can be bounded by some constant $N_1$ independent of $j_0$.
\end{proof}

% \newpage
\subsubsection{The stability of the class $\cala(V,i)$}
In the proofs of Lemma \ref{221207.2648} and Proposition \ref{221223.1757},
we encountered factors having the form of 
$\babr{DA_{n,j},\bbf_{n,j}^{(\hv_0)}}$,
% $\babr{DA_{n,j_V}, \extbbf^{(\hv_0)}_{n,j}}$
which results from the use of the IBP formula or 
the operator $D_{u_n}$.
Lemma \ref{lemma:230616.1754} ensures that 
this factor also belongs to $\cala(V,i)$ when it is properly rescaled.
First we prove the following estimate 
which is essential to the proof of Lemma \ref{lemma:230616.1754}.
Recall that 
$\bbf^{(v,\kappa)}_{n,j}$ %($\hv=(v,\kap)$) 
is defined by
$\bbf^{(v,\kappa)}_{n,j} 
= T_{i(v)n,j_v}^{(\kappa)}(\bbf^{(v,\kappa)})$
with
\begin{align*}
  T_{m,j}^{(1)}(f) = T_{m,j}(f)\tand
  T_{m,j}^{(2)}(f) = T_{m,j+1}(f) - T_{m,j}(f),
\end{align*}
where 
$T_{m,j}(f)(x) = f(mx-j)$.
\begin{lemma}\label{230623.2016}
  For $H\in(1/2,1)$, any $i_0=1,2$ and $f\in L^\infty(\bbR)$ such that $\supp f\subset[-1,0]$,
  the following estimate holds as $n\to\infty$:
  \begin{align}
    \sup_{j_0\in[i_0n-1]}
    \int_{s\in[0,1]}ds
    \abs{\int_{t\in[0,1]}
    (T_{i_0n, j_0+1}(f)(t) - T_{i_0n, j_0}(f)(t))
    \abs{t-s}^{2H-2} dt}
    =O(n^{-2H})
    \label{230907.2434}
  \end{align}
\end{lemma}
\begin{proof}
  Fix $j_0\in[i_0n-1]$.
  By change of variables, we can write 
  \begin{align*}
    &\int_{t\in[0,1]}
    (T_{i_0n, j_0+1}(f)(t) - T_{i_0n, j_0}(f)(t))
    \abs{t-s}^{2H-2} dt  
    % \\&=
    % \int_{t\in[0,1]}
    % \Brbr{f(i_0n t - (j_0+1)) - f(i_0n t - j_0)}
    % \abs{t-s}^{2H-2} dt
    % \\&=
    % \int_{t\in[0,1]} f(i_0n t - (j_0+1)) \abs{t-s}^{2H-2} dt
    % -
    % \int_{t\in[0,1]} f(i_0n t - j_0) \abs{t-s}^{2H-2} dt
    \\&=
    \int_{t\in[\frac{j_0}{i_0n},\frac{j_0+1}{i_0n}]} f(i_0n t - (j_0+1)) \abs{t-s}^{2H-2} dt
    -
    \int_{t\in[\frac{j_0-1}{i_0n},\frac{j_0}{i_0n}]} f(i_0n t - j_0) \abs{t-s}^{2H-2} dt
    % \\&=
    % \int_{t'\in[\frac{j_0-1}{i_0n},\frac{j_0}{i_0n}]} f(i_0n t' - j_0) \abs{(t' + 1/(i_0n))-s}^{2H-2} dt'
    % -
    % \int_{t\in[\frac{j_0-1}{i_0n},\frac{j_0}{i_0n}]} f(i_0n t - j_0) \abs{t-s}^{2H-2} dt
    \\&=
    \int_{t\in[\frac{j_0-1}{i_0n},\frac{j_0}{i_0n}]} 
    f(i_0n t - j_0) \rbr{\abs{t+1/(i_0n)-s}^{2H-2}-\abs{t-s}^{2H-2}} dt.
  \end{align*}
  Hence,
  \begin{align}
    &
    \abs{\int_{t\in[0,1]}
    (T_{i_0n, j_0+1}(f)(t) - T_{i_0n, j_0}(f)(t))
    \abs{t-s}^{2H-2} dt}
    % \nn\\&=
    % \abs{\int_{t\in[\frac{j_0-1}{i_0n},\frac{j_0}{i_0n}]} 
    % f(i_0n t - j_0) \rbr{\abs{(t+1/(i_0n))-s}^{2H-2}-\abs{t-s}^{2H-2}} dt}
    \nn\\&\leq
    \norm{f}_{L^{\infty}}
    \int_{t\in[\frac{j_0-1}{i_0n},\frac{j_0}{i_0n}]} 
    %f(i_0n t - j_0) 
    %\bbone_{[\frac{j_0-1}{i_0n},\frac{j_0}{i_0n}]}(t)
    \abs{\abs{t+1/(i_0n)-s}^{2H-2}-\abs{t-s}^{2H-2}} dt
    \label{eq:230623.2006}
  \end{align}
  for $s\in[0,1]$ and $j_0\in[i_0n-1]$.

  For $2\leq j_0\leq i_0n-1$ and $s\in[0,\frac{j_0-2}{i_0n}]$,
  \begin{align*}
    &\int_{t\in[\frac{j_0-1}{i_0n},\frac{j_0}{i_0n}]} 
    \abs{\abs{t+1/(i_0n)-s}^{2H-2}-\abs{t-s}^{2H-2}} dt
    \\&=
    \int_{t\in[\frac{j_0-1}{i_0n},\frac{j_0}{i_0n}]} 
    \rbr{\rbr{t-s}^{2H-2}-\rbr{(t+1/(i_0n))-s}^{2H-2}} dt
    % \\&\leq
    % \int_{t\in[\frac{j_0-1}{i_0n},\frac{j_0}{i_0n}]} 
    % (2-2H)\rbr{t-s}^{2H-3} (1/(i_0n)) dt
    \\&\leq
    \rbr{\frac{1}{i_0n}}^2 (2-2H)
    \rbr{\frac{j_0-1}{i_0n}-s}^{2H-3}.
  \end{align*}
  Hence, for $j_0\in[i_0n-1]$,
  \begin{align*}
    &\int_{s\in[0,\frac{j_0-2}{i_0n}\vee0]}ds
    \int_{t\in[\frac{j_0-1}{i_0n},\frac{j_0}{i_0n}]} 
    \abs{\abs{t+1/(i_0n)-s}^{2H-2}-\abs{t-s}^{2H-2}} dt
    % \\&\leq 
    % % \int_{s\in[0,\frac{j_0-2}{i_0n}]}ds\;
    % % \frac{1}{i_0n} (2-2H)
    % % \rbr{\frac{j_0-1}{i_0n}-s}^{2H-3} (1/(i_0n))
    % % \\&=
    % \frac{1}{i_0n} 
    % (1/(i_0n))
    % \int_{s\in[0,\frac{j_0-2}{i_0n}]}ds\;
    % (2-2H)\rbr{\frac{j_0-1}{i_0n}-s}^{2H-3}
    % \\&=
    % \frac{1}{i_0n} 
    % (1/(i_0n))
    % \sbr{\rbr{\frac{j_0-1}{i_0n}-s}^{2H-2}}_{s=0}^{\frac{j_0-2}{i_0n}}
    % % \\&=
    % % \frac{1}{i_0n} 
    % % (1/(i_0n))
    % % \rbr{\rbr{\frac{j_0-1}{i_0n}-\frac{j_0-2}{i_0n}}^{2H-2}
    % % -\rbr{\frac{j_0-1}{i_0n}-0}^{2H-2}}
    % \\&=
    % \frac{1}{i_0n} 
    % (1/(i_0n))
    % \rbr{\rbr{\frac{1}{i_0n}}^{2H-2}
    % -\rbr{\frac{j_0-1}{i_0n}}^{2H-2}}
    % \leq%\\&\leq
    \leq
    % (\frac{1}{i_0n})^2
    % \rbr{\frac{1}{i_0n}}^{2H-2}
    % =
    \rbr{\frac{1}{i_0n}}^{2H}.
  \end{align*}
  For $j_0\in[i_0n-1]$,
  a similar calculation shows
  % Assume $j_0\leq i_0n-2$.
  % For $s\in[\frac{j_0+2}{i_0n},1]$, we have
  % \begin{align*}
  %   &\int_{t\in[\frac{j_0-1}{i_0n},\frac{j_0}{i_0n}]} 
  %   \abs{\abs{(t+1/(i_0n))-s}^{2H-2}-\abs{t-s}^{2H-2}} dt
  %   \\&=
  %   \int_{t\in[\frac{j_0-1}{i_0n},\frac{j_0}{i_0n}]} 
  %   \rbr{\rbr{s-(t+1/(i_0n))}^{2H-2}-\rbr{s-t}^{2H-2}} dt
  %   \\&\leq
  %   \int_{t\in[\frac{j_0-1}{i_0n},\frac{j_0}{i_0n}]} 
  %   \rbr{2-2H}\rbr{s-(t+1/(i_0n))}^{2H-3} \times (1/(i_0n)) dt
  %   \\&\leq
  %   \frac{1}{i_0n} 
  %   \rbr{2-2H}\rbr{s-(\frac{j_0+1}{i_0n})}^{2H-3} \times (1/(i_0n)).
  % \end{align*}
  % Hence, 
  \begin{align*}
    &\int_{s\in[\frac{j_0+2}{i_0n}\wedge1,1]}ds
    \int_{t\in[\frac{j_0-1}{i_0n},\frac{j_0}{i_0n}]} 
    \abs{\abs{t+1/(i_0n)-s}^{2H-2}-\abs{t-s}^{2H-2}} dt
    % \\&\leq
    % % \int_{s\in[\frac{j_0+2}{i_0n},1]}ds
    % % \frac{1}{i_0n} 
    % % \rbr{2-2H}\rbr{s-(\frac{j_0+1}{i_0n})}^{2H-3} \times (1/(i_0n))
    % % \\&=
    % \frac{1}{i_0n} \times (1/(i_0n))
    % \int_{s\in[\frac{j_0+2}{i_0n},1]}ds
    % \rbr{2-2H}\rbr{s-(\frac{j_0+1}{i_0n})}^{2H-3}
    % \\&=
    % \frac{1}{i_0n} \times (1/(i_0n))
    % \sbr{-\rbr{s-(\frac{j_0+1}{i_0n})}^{2H-2}}_{s=\frac{j_0+2}{i_0n}}^1
    % % \\&=
    % % \frac{1}{i_0n} \times (1/(i_0n))
    % % \rbr{-\rbr{1-(\frac{j_0+1}{i_0n})}^{2H-2}
    % % +\rbr{\frac{j_0+2}{i_0n}-(\frac{j_0+1}{i_0n})}^{2H-2}}
    % \\&=
    % \frac{1}{i_0n} \times (1/(i_0n))
    % \rbr{-\rbr{1-(\frac{j_0+1}{i_0n})}^{2H-2}
    % +\rbr{\frac{1}{i_0n}}^{2H-2}}
    % \leq%\\&\leq
    % \frac{1}{i_0n} \times (1/(i_0n))
    % \rbr{\frac{1}{i_0n}}^{2H-2}
    % =
    \leq
    \rbr{\frac{1}{i_0n}}^{2H}.
  \end{align*}

  For $j_0\in[i_0n-1]$, we have
  \begin{align*}
    &
    \int_{s\in[0\vee\rbr{\frac{j_0-2}{i_0n}},\rbr{\frac{j_0+2}{i_0n}}\wedge1]}ds
    \int_{t\in[\frac{j_0-1}{i_0n},\frac{j_0}{i_0n}]} 
    \abs{\abs{(t+1/(i_0n))-s}^{2H-2}-\abs{t-s}^{2H-2}} dt
    % \\&\leq
    % \int_{s\in[\frac{j_0-2}{i_0n},\frac{j_0+2}{i_0n}]}ds
    % \int_{t\in[\frac{j_0-1}{i_0n},\frac{j_0}{i_0n}]} 
    % \abs{\abs{(t+1/(i_0n))-s}^{2H-2}-\abs{t-s}^{2H-2}} dt
    \\&\leq
    \int_{s\in[\frac{j_0-2}{i_0n},\frac{j_0+2}{i_0n}]}ds
    \int_{t\in[\frac{j_0-1}{i_0n},\frac{j_0}{i_0n}]} 
    \abs{(t+1/(i_0n))-s}^{2H-2} dt
    +
    \int_{s\in[\frac{j_0-2}{i_0n},\frac{j_0+2}{i_0n}]}ds
    \int_{t\in[\frac{j_0-1}{i_0n},\frac{j_0}{i_0n}]} 
    \abs{t-s}^{2H-2} dt
    \\&\leq2
    \int_{s\in[\frac{-2}{i_0n},\frac{2}{i_0n}]}ds
    \int_{t\in[\frac{-1/2}{i_0n},\frac{1/2}{i_0n}]} 
    \abs{t}^{2H-2} dt
    \\&=
    % 2\int_{s\in[\frac{-2}{i_0n},\frac{2}{i_0n}]}ds
    % 2\sbr{\rbr{2H-1}^{-1}{t}^{2H-1}}_{t=0}^{\frac{1/2}{i_0n}}
    % =%\\&=
    2\frac{4}{i_0n}
    \times2
    \rbr{2H-1}^{-1}\rbr{\frac{1/2}{i_0n}}^{2H-1}
    % \\&=
    % \frac{16}{i_0n}\rbr{2H-1}^{-1}\rbr{\frac{1/2}{i_0n}}^{2H-1}
    =
    2^{5-2H}\rbr{2H-1}^{-1}
    \rbr{i_0n}^{-2H}.
    % \rbr{\frac{1}{i_0n}}^{2H}
    % \leq
    % 2^{4}\rbr{2H-1}^{-1}
    % \rbr{\frac{1}{i_0n}}^{2H}.
  \end{align*}

  Therefore, we obtain 
  \begin{align*}
    &\int_{s\in[0,1]}ds
    \int_{t\in[\frac{j_0-1}{i_0n},\frac{j_0}{i_0n}]} 
    \abs{\abs{t+1/(i_0n)-s}^{2H-2}-\abs{t-s}^{2H-2}} dt
    % \\&\leq
    % \rbr{\frac{1}{i_0n}}^{2H}+
    % 2^{4}\rbr{2H-1}^{-1}
    % \rbr{\frac{1}{i_0n}}^{2H}+
    % \rbr{\frac{1}{i_0n}}^{2H}
    % =
    % (2+2^{4}\rbr{2H-1}^{-1})
    % \rbr{\frac{1}{i_0n}}^{2H}
    % =
    % c'_H
    % \rbr{\frac{1}{i_0n}}^{2H},
    \leq c\; n^{-2H}
  \end{align*}
  with some constant $c$ independent of $j_0$ and $n$,
  and hence the estimate \eqref{230907.2434} by \eqref{eq:230623.2006}.
\end{proof}

\begin{lemma}\label{lemma:230616.1754}
  Consider a nonempty finite set $V$, $i:V\to\cbr{1,2}$ and $A\in\cala(V,i)$.
  Let $v_0\in V$.

  \item[(i)]
  For $f\in L^\infty(\bbR)$ such that $\supp f\subset[-1,0]$,
  define $A'=(A'_{n,j})_{n\in\ntwo, j\in\bbJ_n(V,i)}$ by 
  \begin{align*}
    A'_{n,j} = n^{2H}\abr{DA_{n,j}, T^{(2)}_{i(v_0)n, j_{v_0}}(f)}_{\calh}.
  \end{align*}
  Then, $A'\in\cala(V,i)$.

  \item[(ii)]
  For $f\in L^\infty(\bbR)$ such that $\supp f\subset[\alpha,\alpha+1]$ with some $\alpha\in[-1,0]$,
  define $A'=(A'_{n,j})_{n\in\ntwo, j\in\bbJ_n(V,i)}$ by 
  \begin{align*}
    A'_{n,j} = n\abr{DA_{n,j}, T^{(1)}_{i(v_0)n, j_{v_0}}(f)}_{\calh}.
  \end{align*}
  Then, $A'\in\cala(V,i)$.
\end{lemma}

\begin{proof}
  For notational convenience, we write $i_0=i(v_0)$ and $j_0=j_{v_0}$.

  \noindent(i)
  % Fix $j\in\bbJ_n(V,i)$.
  Since we have 
  \begin{align*}
    \abr{DA_{n,j}, T^{(2)}_{i_0n, j_0}(f)}&=
    % \abr{DA_{n,j}, T_{i_0n, j_0+1}(f) - T_{i_0n, j_0}(f)}
    % % \\&=
    % % c_H\int_{[0,1]^2}
    % % D_sA_{n,j} (T_{i_0n, j_0+1}(f)(t) - T_{i_0n, j_0}(f)(t))
    % % \abs{t-s}^{2H-2} dsdt 
    % \\&=
    c_H\int_{s\in[0,1]}
    D_sA_{n,j}\,ds
    \int_{t\in[0,1]}
    (T_{i_0n, j_0+1}(f)(t) - T_{i_0n, j_0}(f)(t))
    \abs{t-s}^{2H-2} dt,
  \end{align*}
  the $L^p$-norm of $A'_{n,j}$ is bounded as 
  \begin{align*}
    \norm{A'_{n,j}}_{L^p}%&=
    % n^{2H}
    % \norm{\abr{DA_{n,j}, T^{(2)}_{i_0n, j_0}(f)}}_{L^p}
    % \\&\leq
    % n^{2H}
    % c_H\int_{s\in[0,1]}
    % \norm{D_sA_{n,j}}_{L^p} ds
    % \abs{\int_{t\in[0,1]}
    % (T_{i_0n, j_0+1}(f)(t) - T_{i_0n, j_0}(f)(t))
    % \abs{t-s}^{2H-2} dt}
    % \\&\leq
    &\leq
    n^{2H}
    c_H
    \sup_{s\in[0,1]}\norm{D_sA_{n,j}}_{L^p}
    \int_{s\in[0,1]}ds
    \abs{\int_{t\in[0,1]}
    (T_{i_0n, j_0+1}(f)(t) - T_{i_0n, j_0}(f)(t))
    \abs{t-s}^{2H-2} dt},
  \end{align*}
  and by Lemma \ref{230623.2016} we obtain 
  $\sup_{j\in\bbJ_n(V,i)}\norm{A'_{n,j}}_{L^p}=O(1)$ as $n\to\infty$.
  Here we used the assumption on $A=(A_{n,j})_{j\in\bbJ_n(V,i),n\in\ntwo}$ 
  (i.e. 
  $\sup_{n\in\ntwo,j\in\bbJ_n(V,i)}\sup_{s\in[0,1]}\norm{D_sA_{n,j}}_{L^p}<\infty$
  ).

  For higher derivatives, similar arguments work as below.
  The Malliavin derivative of $A'_{n,j}$ is written as 
  \begin{align*}
    D^k_{s_1,...,s_k}A'_{n,j} &= 
    % n^{2H}D^k_{s_1,...,s_k}
    % \abr{DA_{n,j}, T^{(2)}_{i_0n, j_0}(f)}
    % \\&=
    % n^{2H} c_H
    % \int_{s_{k+1},t\in[0,1]}
    % D^{k+1}_{s_1,...,s_{k+1}}A_{n,j}\; T^{(2)}_{i_0n, j_0}(f)(t)
    % \abs{s_{k+1}-t}^{2H-2} ds_{k+1} dt
    % \\&=
    n^{2H}
    c_H\int_{s_{k+1}\in[0,1]}
    D^{k+1}_{s_1,...,s_{k+1}}A_{n,j}\, ds_{k+1}
    \int_{t\in[0,1]}
    (T_{i_0n, j_0+1}(f)(t) - T_{i_0n, j_0}(f)(t))
    \abs{t-s_{k+1}}^{2H-2} dt
  \end{align*}
  and we have the following estimate:
  \begin{align*}
    &\norm{D^k_{s_1,...,s_k}A'_{n,j}}_{L^p(P)} 
    % \\&\leq
    % n^{2H}
    % c_H\int_{s_{k+1}\in[0,1]}
    % \norm{D^{k+1}_{s_1,...,s_{k+1}}A_{n,j}}_{L^p(P)} ds_{k+1}
    % \abs{\int_{t\in[0,1]}
    % (T_{i_0n, j_0+1}(f)(t) - T_{i_0n, j_0}(f)(t))
    % \abs{t-s_{k+1}}^{2H-2} dt}
    \\&\leq
    n^{2H} c_H
    \sup_{s_{k+1}\in[0,1]}
    \norm{D^{k+1}_{s_1,...,s_{k+1}}A_{n,j}}_{L^p(P)} 
    \\&\hspace{20pt}\times
    \int_{s_{k+1}\in[0,1]}
    ds_{k+1}
    \abs{\int_{t\in[0,1]}
    (T_{i_0n, j_0+1}(f)(t) - T_{i_0n, j_0}(f)(t))
    \abs{t-s_{k+1}}^{2H-2} dt}
  \end{align*}
  Thus, 
  by the assumption on $A=(A_{n,j})_{j\in\bbJ_n(V,i),n\in\ntwo}$,
  that is 
  \begin{align*}
    \sup_{n\in\ntwo, j\in\bbJ_n(V,i)}\sup_{s_1,...,s_{k+1}\in[0,1]}
    \norm{D^{k+1}_{s_1,...,s_{k+1}}A_{n,j}}_{L^p(P)} 
    <\infty
  \end{align*}
  and Lemma \ref{230623.2016},
  we obtain 
  \begin{align*}
    &\sup_{n\in\ntwo, j\in\bbJ_n(V,i)}\sup_{s_1,...,s_k\in[0,1]}
    \norm{D^k_{s_1,...,s_k}A'_{n,j}}_{L^p(P)} 
    % \\&\leq
    % \sup_{n\in\ntwo, j\in\bbJ_n(V,i)}\sup_{s_1,...,s_k\in[0,1]}
    % \sup_{s_{k+1}\in[0,1]}
    % \norm{D^{k+1}_{s_1,...,s_{k+1}}A_{n,j}}_{L^p(P)} 
    % \\&\qquad\times
    % c_H
    % \sup_{n\in\ntwo, j\in\bbJ_n(V,i)}
    % n^{2H}
    % \int_{s_{k+1}\in[0,1]}
    % ds_{k+1}
    % \abs{\int_{t\in[0,1]}
    % (T_{i_0n, j_0+1}(f)(t) - T_{i_0n, j_0}(f)(t))
    % \abs{t-s_{k+1}}^{2H-2} dt}
    % \\&<\infty.
    <\infty.
  \end{align*} 

  \item[(ii)]
  The proof is straight forward,
  since we can prove 
  \begin{align*}
    \sup_{j_0\in[i_0n-1]}
    \int_{s\in[0,1]}ds 
    \abs{\int_{t\in[0,1]} T^{(1)}_{i_0n,j_0}(f)(t)
    \abs{t-s}^{2H-2} dt}
    =O(n^{-1}).
  \end{align*}
  % See for \ref{quad. var. paper}.\koko
\end{proof}

The following lemma is used in the proofs of 
Propositions \ref{221223.1757} and \ref{230725.1133}.
\begin{lemma}\label{lemma:230616.1805}
  (i)
  Let $V, \extv$ be nonempty finite sets such that $V\subset\extv$,
  $\exti:\extv\to\cbr{1,2}$ and
  $A\in\cala(V,\exti|_V)$.
  Define $A'=(A'_{n,\extj})_{n\in\ntwo, \extj\in\bbJ_n(\extv,\exti)}$ by 
  $%\begin{align*}
    A'_{n,\extj} = A_{n,\extj_V},
  $ %\end{align*}
  where $\extj_V$ is the projection of $\extj$ to $V$.
  Then, $A'\in\cala(\extv,\exti)$.

  \item[(ii)]
  Let $V$ be a nonempty finite set, $i:V\to\cbr{1,2}$, and $A,A'\in\cala(V,i)$.
  \begin{itemize}
    \item [(a)]
    Define $A''=(A''_{n,j})_{n\in\ntwo, j\in\bbJ_n(V,i)}$ by 
    $%\begin{align*}
      A''_{n,j} = A_{n,j}\: A'_{n,j}.
    $ %\end{align*}
    Then, $A''\in\cala(V,i)$.

    \item[(b)]
    For $k\geq1$,
    define $A''=(A''_{n,j})_{n\in\ntwo, j\in\bbJ_n(V,i)}$ by 
    $%\begin{align*}
      A''_{n,j} = \abr{D^k A_{n,j}, D^k A'_{n,j}}_{\calh^{\otimes k}}.
    $ %\end{align*}
    Then, $A''\in\cala(V,i)$.
  \end{itemize}
\end{lemma}
We omit the proof, since it is straightforward.

\subsubsection{Product formula}
Consider nonempty sets $V$, $V_0$, $V_1$ with
$V_i\subset V$ ($i=0,1$) and $V_0\cap V_1=\emptyset$.
Recall that we have defined
\begin{align*}
  \tp(V_0,V_1) = 
  \cbr{[\tv,\tv']\in\tp(V_0\sqcup V_1)\mid
  (\tv,\tv')\text{ or }(\tv',\tv)\in\wtv_0\times\wtv_1}.
\end{align*}
For $\vwq_i:\wtv\to\bbZ_{\geq0}$ such that $\supp{\vwq_i}\subset\wtv_i$ for $i=0,1$, 
we define
\begin{align*}
  \Pi_{V_0, V_1}(\bfq_0, \bfq_1) =
  \Big\{\bfp:\hp(V_0\sqcup V_1)\to\bbZ_{\geq0}\mid
  \supp(\bfp)\subset\tp(V_0,V_1),& 
  %\bar\bfp(\tv_0)\leq\bfq_0(\tv_0)\tforsm\tv_0\in\wtv_0,\quad
  \sum_{\tv_1\in\wtv_1}\bfp[\tv_0,\tv_1]\leq\bfq_0(\tv_0)\tforsm\tv_0\in\wtv_0,
  \\&\;
  %\bar\bfp(\tv_1)\leq\bfq_1(\tv_1)\tforsm\tv_1\in\wtv_1\}
  \sum_{\tv_0\in\wtv_0}\bfp[\tv_0,\tv_1]\leq\bfq_1(\tv_1)\tforsm\tv_1\in\wtv_1\Big\}.
\end{align*}
For $\bfp\in\Pi_{V_0,V_1}(\vwq_0,\vwq_1)$, we define 
$\bar\bfp:\whv_0\sqcup\whv_1\to\bbZ_{\geq0}$ by
$\bar\bfp(\tv_0) = \sum_{\tv_1\in\wtv_1}\bfp[\tv_0,\tv_1]$ for $\tv_0\in\wtv_0$ and
$\bar\bfp(\tv_1) = \sum_{\tv_0\in\wtv_0}\bfp[\tv_0,\tv_1]$ for $\tv_1\in\wtv_1$.

\begin{lemma}\label{lemma:230526.1230}
  Let $\vwq_i:\wtv\to\bbZ_{\geq0}$ such that $\supp{\vwq_i}\subset\wtv_i$ for $i=0,1$ and 
  $\bff\in\abs{\calh}^{\wtv}$.
  Then the following equality holds:
  \begin{align}\label{eq:230615.1050}
    \delta_{V_0}(\bff,\vwq_0)\: \delta_{V_1}(\bff,\vwq_1)
    =
    \sum_{\bfp\in\Pi_{V_0,V_1}(\vwq_0,\vwq_1)}
    c(\bfp) \;
    \delta_{V_0\sqcup V_1}(\bff, \bfq_0+\bfq_1-\bar\bfp)\;
    \beta_{V_0\sqcup V_1}(\bff, \bfp),
  \end{align}
  where $c(\bfp)$ is a positive integer-valued constant.
\end{lemma}
\begin{remark}
  In the following proof, the constant $c(\bfp)$ might be different for the same $\bfp$ line by line.
\end{remark}

\begin{proof}
  Let us write 
  $\barq_i=\sum_{\hv\in\whv_i}\vwq_i(\hv)$ for $i=0,1$.
  We prove the equality \eqref{eq:230615.1050} by induction with respect to $\min_{i=0,1}\barq_i$.
  First consider the case where $\min_{i=0,1}\barq_i=0$, that is
  $\vwq_0\equiv0$ or $\vwq_1\equiv0$.
  By symmetry, we assume $\vwq_0\equiv0$.
  By the definition, $\Pi_{V_0,V_1}(\vwq_0,\vwq_1)=\cbr{0}$, and 
  by setting the constant $c(0)=1$, we have 
  \begin{align*}
    &\sum_{\bfp\in\Pi_{V_0,V_1}(\vwq_0,\vwq_1)}
    c(\bfp)\:
    \delta_{V_0\sqcup V_1}(\bff, \bfq_0+\bfq_1-\bar\bfp)\:
    \beta_{V_0\sqcup V_1}(\bff, \bfp)
    =%\\&=
    c(0)\:\delta_{V_0\sqcup V_1}(\bff,\vwq_1)\:\beta_{V_0\sqcup V_1}(\bff,0)=
    \delta_{V_1}(\bff,\vwq_1),
  \end{align*}
  which proves \eqref{eq:230615.1050} since $\delta_{V_0}(\bff,\vwq_0)=1$.

  Next fix $d\in\bbZ_{\geq0}$ and assume that the equality \eqref{eq:230615.1050} holds for 
  any $\vwq_i$ ($i=0,1$) satisfying $\min_{i=0,1}\barq_i\leq d$.
  Consider $\vwq_i$ ($i=0,1$) with $\min_{i=0,1}\barq_i=d+1$.
  By symmetry we assume $\barq_0=d+1$.
  %assume that $\vwq_0\not\equiv0$ and $\vwq_1\not\equiv0$.
  Fix $\hv_0=(v_0,\kap_0)\in\whv_0$ such that $\vwq_0(\hv_0)>0$.
  Then the left hand side of \eqref{eq:230615.1050} decomposes as 
  \begin{align*}
    \delta_{V_0}(\bff, \vwq_0) \delta_{V_1}(\bff, \vwq_1)
    &=
    \delta(\delta_{V_0}(\bff, \vwq_0-\bbone_\cbr{\hv_0})\bff^{(\hv_0)}) \delta_{V_1}(\bff, \vwq_1)
    \\&=
    \delta(\delta_{V_0}(\bff, \vwq_0-\bbone_\cbr{\hv_0})\:\delta_{V_1}(\bff, \vwq_1)\:\bff^{(\hv_0)}) 
    +
    \delta_{V_0}(\bff, \vwq_0-\bbone_\cbr{\hv_0})\babr{\bff^{(\hv_0)}, D\delta_{V_1}(\bff, \vwq_1)}
  \end{align*}
  Let us write $\vwq'_0=\vwq_0-\bbone_\cbr{\hv_0}$.
  Since $\sum_{\hv\in\whv_0}\vwq'_0(\hv)=\barq_0-1=d$,
  we can apply the assumption of induction to $\vwq'_0$ and $\vwq_1$, and 
  \begin{align}
    \delta(\delta_{V_0}(\bff, \vwq'_0)\:\delta_{V_1}(\bff, \vwq_1)\:\bff^{(\hv_0)}) 
    &=
    \sum_{\bfp\in\Pi_{V_0,V_1}(\vwq'_0,\vwq_1)}c(\bfp)\:
    \delta(\delta_{V_0\sqcup V_1}(\bff, \vwq'_0+\vwq_1-\bar\bfp)\:
    \beta_{V_0\sqcup V_1}(\bff, \bfp)\:\bff^{(\hv_0)}) 
    \nn\\&=
    \sum_{\bfp\in\Pi_{V_0,V_1}(\vwq'_0,\vwq_1)}c(\bfp)\:
    \delta_{V_0\sqcup V_1}(\bff, \vwq_0+\vwq_1-\bar\bfp)\:
    \beta_{V_0\sqcup V_1}(\bff, \bfp).
    \label{eq:230615.1351}
  \end{align}
  
  We can write
  \begin{align*}
    D\delta_{V_1}(\bff, \vwq_1)
    =
    \sum_{\substack{\hv_1\in\whv_1\\\vwq_1(\hv_1)>0}}\vwq_1(\hv_1)
    \delta_{V_1}(\bff, \vwq_1-\bbone_\cbr{\hv_1})\bff^{(\hv_1)}.
  \end{align*}
  Let us denote $\vwq_1^{(\hv_1)}=\vwq_1-\bbone_\cbr{\hv_1}$ for $\hv_1\in\whv_1$ with $\vwq_1(\hv_1)>0$.
  Again, we can apply the assumption of induction to $\vwq'_0$ and $\vwq^{(\hv_1)}_1$ to obtain
  \begin{align*}
    \delta_{V_0}(\bff, \vwq'_0)\delta_{V_1}(\bff, \vwq_1^{(\hv_1)})
    &=
    \sum_{\bfp\in\Pi_{V_0,V_1}(\vwq'_0,\vwq_1^{(\hv_1)})}c(\bfp)\:
    \delta_{V_0\sqcup V_1}(\bff, \vwq'_0+\vwq_1^{(\hv_1)}-\bar\bfp)
    \beta_{V_0\sqcup V_1}(\bff, \bfp)
  \end{align*}
  and 
  \begin{align*}
    &\delta_{V_0}(\bff, \vwq'_0)\babr{\bff^{(\hv_0)}, D\delta_{V_1}(\bff, \vwq_1)}
    \\&=
    \delta_{V_0}(\bff, \vwq'_0)\babr{\bff^{(\hv_0)},
    \sum_{\substack{\hv_1\in\whv_1\\\vwq_1(\hv_1)>0}}\vwq_1(\hv_1)\:
    \delta_{V_1}(\bff, \vwq_1^{(\hv_1)})\:\bff^{(\hv_1)}}
    \\&=
    \sum_{\substack{\hv_1\in\whv_1\\\vwq_1(\hv_1)>0}}\vwq_1(\hv_1)\:
    \delta_{V_0}(\bff, \vwq'_0)\:
    \delta_{V_1}(\bff, \vwq_1^{(\hv_1)})\:
    \babr{\bff^{(\hv_0)},\bff^{(\hv_1)}}
    \\&=
    \sum_{\substack{\hv_1\in\whv_1\\\vwq_1(\hv_1)>0}}\vwq_1(\hv_1)
    \sum_{\bfp\in\Pi_{V_0,V_1}(\vwq'_0,\vwq_1^{(\hv_1)})}c(\bfp)\:
    \delta_{V_0\sqcup V_1}(\bff, \vwq'_0+\vwq_1^{(\hv_1)}-\bar\bfp)
    \beta_{V_0\sqcup V_1}(\bff, \bfp+\bbone_{[\hv_0,\hv_1]})
    %\babr{\bff^{(\hv_0)},\bff^{(\hv_1)}}
  \end{align*}
  Write $\bfp^{(\hv_1)}=\bfp+\bbone_{[\hv_0,\hv_1]}$ for $\bfp\in\Pi_{V_0,V_1}(\vwq'_0,\vwq^{(\hv_1)}_1)$, and 
  we have
  $\bar\bfp^{(\hv_1)}=\bar\bfp+\bbone_\cbr{\hv_0}+\bbone_\cbr{\hv_1}$.
  Setting
  $\Gamma=
  \bcbr{\bfp^{(\hv_1)}\mid 
  \hv_1\in\whv_1 \text{ s.t. }\vwq_1(\hv_1)>0 \tandsm \bfp\in\Pi_{V_0,V_1}(\vwq'_0,\vwq^{(\hv_1)}_1)}$,
  we can write 
  \begin{align}
    &\delta_{V_0}(\bff, \vwq'_0)\babr{\bff^{(\hv_0)}, D\delta_{V_1}(\bff, \vwq_1)}
    \nn\\&=
    \sum_{\substack{\hv_1\in\whv_1\\\vwq_1(\hv_1)>0}}
    \sum_{\bfp\in\Pi_{V_0,V_1}(\vwq'_0,\vwq_1^{(\hv_1)})}\vwq_1(\hv_1)c(\bfp)\:
    \delta_{V_0\sqcup V_1}(\bff, \vwq_0+\vwq_1-\bar\bfp^{(\hv_1)})
    \beta_{V_0\sqcup V_1}(\bff, \bfp^{(\hv_1)}),
    \nn\\&=
    \sum_{\bfp\in\Gamma}c_\Gamma(\bfp)\:
    \delta_{V_0\sqcup V_1}(\bff, \vwq_0+\vwq_1-\bar\bfp)
    \beta_{V_0\sqcup V_1}(\bff, \bfp),
    \label{eq:230615.1423}
  \end{align}
  where $c_\Gamma(\bfp)$ in \eqref{eq:230615.1423} is written as 
  \begin{align*}
    c_\Gamma(\bfp)=\sum\{\vwq_1(\hv_1)c(\bfp')\mid 
    \hv_1\in\whv_1 \text{ with }\vwq_1(\hv_1)>0,\;
    \bfp'\in\Pi_{V_0,V_1}(\vwq'_0,\vwq^{(\hv_1)}_1) \text{ s.t. }
    \bfp=\bfp'+\bbone_{[\hv_0,\hv_1]}\}.
  \end{align*}
  Notice that $c_\Gamma(\bfp)>0$ for any $\bfp\in\Gamma$.
  %Notice that the constant $c(\bfp)$ may be different from those in \eqref{eq:230615.1050} and \eqref{eq:230615.1351}
  %for the same $\bfp$.

  Hence, by \eqref{eq:230615.1351} and \eqref{eq:230615.1423}, we can write
  \begin{align*}
    \delta_{V_0}(\bff, \vwq_0) \delta_{V_1}(\bff, \vwq_1)
    &=
    \sum_{\bfp\in\Pi_{V_0,V_1}(\vwq'_0,\vwq_1)\cup\Gamma}
    (c(\bfp)\bbone_{\Pi_{V_0,V_1}(\vwq'_0,\vwq_1)}(\bfp)+c_\Gamma(\bfp)\bbone_{\Gamma}(\bfp))\:
    \delta_{V_0\sqcup V_1}(\bff, \vwq_0+\vwq_1-\bar\bfp)\:
    \beta_{V_0\sqcup V_1}(\bff, \bfp).
  \end{align*}
  By definition, we can see $\Pi_{V_0,V_1}(\vwq'_0,\vwq_1)\cup\Gamma\subset\Pi_{V_0,V_1}(\vwq_0,\vwq_1)$.
  On the other hand, we can show 
  $\Pi_{V_0,V_1}(\vwq'_0,\vwq_1)\cup\Gamma\supset\Pi_{V_0,V_1}(\vwq_0,\vwq_1)$ and 
  we have proved \eqref{eq:230615.1050} and $c(\bfp)>0$ for $\bfp\in\Pi_{V_0,V_1}(\vwq_0,\vwq_1)$.
\end{proof}
%

%% file: subfiles/5-sec_asy_exp.tex
% \maketitle
% \tableofcontents
%
\section{Proof of asymptotic expansion of the Hurst estimator}
\label{sec:231002.2423}
In this section, we prove the asymptotic expansion of the distribution of
the Hurst estimator with the help of the theory of exponent introduced in the last section.
In Section \ref{sec:231005.2242}, 
we exhibit a stochastic expansion of $Z_n$, which is
the functional obtained by truncating the rescaled error $\sqrt{n}(\hat H_n-H)$.
In Section \ref{230721.1156},
we estimate the order of functionals related to a mixed CLT of $Z_n$
and in particular identify the asymptotic conditional variance, namely $G_\infty$.
In Section \ref{sec:231001.1030}, we verify Condition {\bf [D]} of 
the general theory by Nualart and Yoshida 
to derive the asymptotic expansion of $Z_n^\circ$, a principal term of $Z_n$,
while we calculate the contribution from the rest $Z_n-Z_n^\circ$ and 
establish the expansion formula for $\sqrt{n}(\hat H_n-H)$ in Section \ref{sec:231001.1649}.

We say the functional $\cali_n$ corresponds to a weighted graph $G$,
if there exist $i:V\to\cbr{1,2}$, $A\in\cala(V,i)$ and $\bbf\in\calf(V)$ such that 
$\cali_n = \cali_n^{(i)}(G,A,\bbf)$ holds, 
where we write $V$ for the set of vertices $V(G)$ of $G$.
In Section \ref{230925.1250}, the weighted graphs related to the functionals 
appearing in the following arguments are collected.

\subsection{Stochastic expansion of $Z_n$}
\label{sec:231005.2242}
In Section \ref{sec:231002.2414}, 
we have defined $\rv_n$ and its limit $\iv$ by 
\begin{align*}
  \rv_n = n^{2H-1}\secVar_n
  = n^{2H-1}\sum_{j=1}^{n-1}\rbr{\secDiff{n}{j}X}^2
  \tand
  \iv=c_{2,H}\int^1_0 \rbr{V^{[1]}_t}^2 dt
\end{align*}
with the constant $c_{2,H}$ defined at \eqref{230925.1607},
and we write 
%$\convDiff_n = \rv_n-\iv$.
\begin{align}\label{230804.1403}
  \convDiff_n = \rv_n-\iv.
\end{align}
Let $\psi:\bbR\to[0,1]$ be a smooth function satisfying 
$\psi(x)=1$ if $\abs{x}\leq\half$ and 
% on $\cbr{x\mid\abs{x}\leq\half}$ and 
$\psi(x)=0$ if $\abs{x}\geq1$.
% on $\cbr{x\mid\abs{x}\geq1}$.
We define a sequence $(\psi_n)_{n\in\ntwo}$ of functionals by
\begin{align}\label{eq:230926.1352}
  \psi_n = 
  \psi\rbr{ \frac{R_n}{\eta_0\iv}}\psi\rbr{\frac{R_{2n}}{\eta_0\iv}}.
\end{align}
Here the constant $\eta_0\in\bbR_{>0}$ is taken small enough so that $\hat H_n\in(0,1)$ when $\psi_n>0$.
With this $\psi_n$, we define $Z_n$ by 
\begin{align}\label{eq:230926.1554}
  Z_n = 2\log2\,\sqrt{n} (\hat H_n - H) \psi_n.
\end{align}
% {\mygreen
% The main part of 
% $\hat H_n =
% 0\vee\rbr{\half + \frac1{2\log2} \log\frac{\secVar_n}{\secVar_{2n}}}\wedge1$
% can be decomposed as 
% \begin{align*}
%   \log\frac{\secVar_n}{\secVar_{2n}} &= 
%   %\log\rbr{2^{2H-1} \frac{n^{2H-1}\secVar_n}{{(2n)}^{2H-1}\secVar_{2n}}} = 
%   \log\rbr{2^{2H-1} \frac{\rv_n}{\rv_{2n}}}
%   % \\&=
%   % \log\rbr{\frac{\rv_n}{\iv}}-\log\rbr{\frac{\rv_{2n}}{\iv}}
%   % + (H-\half)2\log2
%   =%\\&=
%   (H-\half)2\log2 +
%   \log\rbr{\frac{\convDiff_n}{\iv}+1}-\log\rbr{\frac{\convDiff_{2n}}{\iv}+1}.
% \end{align*}
% Then, by the Taylor expansion, we have the next lemma,
% whose proof is given in {\myred Section \ref{230804.1348}}.
% \begin{lemma}\label{lem:230804.1452}
%   The functional $Z_n$ decomposes as
%   \begin{align*}%\label{230804.1353}
%     Z_n &=
%     \sqrt{n} 
%     \cbr{\frac{\convDiff_n}\iv - \frac{\convDiff_{2n}}\iv}
%     -  \sqrt{n} 
%     \cbr{\half\rbr{\frac{\convDiff_n}\iv}^2 -\half\rbr{\frac{\convDiff_{2n}}\iv}^2}
%     + \negTerm_n^{(1)}
%   \end{align*}
%   with $\negTerm_n^{(1)}=O_M\rbr{n^{-1}}$.
% \end{lemma}
% }
% With the above definitions of $\psi_n$ and $Z_n$,
% we will prove the stochastic expansion of $Z_n$ (i.e. Proposition \ref{230720.2014}).

Before we proceed to calculating the stochastic expansion of 
$Z_n$,
we decompose the functional $\convDiff_n$.
Recall that we have defined 
$\diffker^n_j = \bbone^n_\jp - \bbone^n_j$.
The second order difference 
$\secDiff{n}{j}X = (X^n_{j+1} - X^n_{j}) - (X^n_{j} - X^n_{j-1})$
% ($j=1,...,n-1$).
in the summand of $\rv_n$ can be written as
\begin{align*}
  \secDiff{n}{j} X &=
  % (X^n_{j+1}-X^n_{j}) - (X^n_{j}-X^n_{j-1})
  % \\&=
  \int^{t^n_{j+1}}_{t^n_j} \rbr{V^{[2]}(X_t) dt + V^{[1]}(X_t) dB_t} -
  \int^{t^n_j}_{t^n_{j-1}} \rbr{V^{[2]}(X_t) dt + V^{[1]}(X_t) dB_t}
  \\&=
  V^{[1]}_{\tnj{n}{j}} \Delta^n_{\jp}B + V^{[2]}_\tnj{n}{j} n^{-1} + 
  \int^{\tnj{n}{j+1}}_{\tnj{n}{j}} \rbr{V^{[1]}_t - V^{[1]}_{t^n_j}} dB_t + 
  \int^{\tnj{n}{j+1}}_{\tnj{n}{j}} \rbr{V^{[2]}_t - V^{[2]}_{t^n_j}} dt
  \\&\hspsm-\rbr{
    V^{[1]}_\tjm \Delta^n_{j}B + V^{[2]}_\tjm n^{-1} + 
  \int^\tj_\tjm \rbr{V^{[1]}_t - V^{[1]}_\tjm} dB_t + 
  \int^\tj_\tjm \rbr{V^{[2]}_t - V^{[2]}_\tjm} dt
  }
  % \\&=
  % V^{[1]}_\tj \rbr{\Delta^n_{\jp}B - \Delta^n_{j}B}
  % + \rbr{V^{[1]}_\tj - V^{[1]}_\tjm} \Delta^n_{j}B
  % + \rbr{V^{[2]}_\tj - V^{[2]}_\tjm} n^{-1}
  % +\resint
  \\&=
  V^{[1]}_\tj B(\dker^n_j) + \resincj,
\end{align*}
where we define 
\begin{align}
  \resincj &=
  \rbr{V^{[1]}_\tj - V^{[1]}_\tjm} \Delta^n_{j}B
  + \rbr{V^{[2]}_\tj - V^{[2]}_\tjm} n^{-1}
  +\resint
  \label{230925.1518}
  \\
  \resint &= 
  \int^\tjp_\tj \rbr{V^{[1]}_t - V^{[1]}_\tj} dB_t + 
  \int^\tjp_\tj \rbr{V^{[2]}_t - V^{[2]}_\tj} dt
  \nn\\&\hspsm-\rbr{
  \int^\tj_\tjm \rbr{V^{[1]}_t - V^{[1]}_\tjm} dB_t + 
  \int^\tj_\tjm \rbr{V^{[2]}_t - V^{[2]}_\tjm} dt
  }.\nn
\end{align}
By the product formula, $\convDiff_n$ is written as 
\begin{align*}
  \convDiff_n &= 
  % n^{2H-1}\sum_{j=1}^{n-1}\rbr{\secDiff{n}{j}X}^2
  % -\iv
  % \\&=
  n^{2H-1} \sum_{j=1}^{n-1} \rbr{
  \rbr{V^{[1]}_\tj B\rbr{\diffker^n_j} }^2
  +2 V^{[1]}_\tj B\rbr{\diffker^n_j} \resinc{j}
  + \resinc{j}^2}
  -\iv
  % \nn\\&=
  % c_{2,H}n^{-1} \sum_{j=1}^{n-1} \rbr{V^{[1]}_\tj }^2
  % - c_{2,H}\int^1_0 \rbr{V^{[1]}_t}^2 dt
  % \nn\\&\hspsm
  % +n^{2H-1} \sum_{j=1}^{n-1} \rbr{V^{[1]}_\tj }^2 I_2\rbr{(\diffker^n_j)^{\otimes2}}
  % \nn\\&\hspsm
  % +2 n^{2H-1} \sum_{j=1}^{n-1} V^{[1]}_\tj B\rbr{\diffker^n_j} \resinc{j}
  % + n^{2H-1} \sum_{j=1}^{n-1} \resinc{j}^2
  % \label{221121.2305}
  =%\\&=
  \sum_{k=0}^3 \convDiff^{(k)}_n,
  % \convDiff^{(1)}_n + \convDiff^{(0)}_n
  % + \convDiff^{(2)}_n + \convDiff^{(3)}_n,
  \nn
\end{align*}
where
\begin{align}
  \convDiff^{(0)}_n &=
  n^{2H-1} \sum_{j=1}^{n-1} 
  \rbr{V^{[1]}_\tj }^2 I_2\rbr{(\diffker^n_j)^{\otimes2}},
  \label{230925.1600}\\% &%\\
  \convDiff^{(1)}_n &=
  c_{2,H}\;n^{-1} \sum_{j=1}^{n-1} \rbr{V^{[1]}_\tj }^2
  - c_{2,H}\int^1_0 \rbr{V^{[1]}_t}^2 dt,
  \label{230925.1601}\\
  \convDiff^{(2)}_n &=
  2 n^{2H-1} \sum_{j=1}^{n-1} V^{[1]}_\tj 
  B\rbr{\diffker^n_j} \resinc{j},
  \label{230925.1602}\\% &%\\
  \convDiff^{(3)}_n &=
  n^{2H-1} \sum_{j=1}^{n-1} \resinc{j}^2.
  \label{230925.1603}
\end{align}
Recall that we have defined 
$c_{2,H} = 4-2^{2H}
% =\abr{\diffker^1_j,\diffker^1_j}
(=n^{2H}\abr{\diffker^n_j,\diffker^n_j})$.
Decomposing $\convDiff^{(k)}_n$ ($k=1,2,3$) further 
using the lemmas in Section \ref{230926.1330}, 
we have the following lemma:
\begin{lemma}\label{230720.2006}
  The random variable $\convDiff_n$ defined at \eqref{230804.1403} is written as follows: 
  \begin{align}\label{eq:230926.1136}
    \convDiff_n &=
    \convDiff_n^{(0)} 
    + \convDNeg_n
    + \convDiff^{(\negTerm)}_n,
  \end{align}
  where 
  $\convDiff^{(0)}_n$ is defined at \eqref{230925.1600}, and
  $\convDiff^{(\negTerm)}_n$ is some functional satisfying
  $\convDiff^{(\negTerm)}_n=O_M(n^{(-\half-H)\vee(-2+H)})$.
  % $\convDiff^{(\negTerm)}_n=
  % \convDiff^{(1,\negTerm)}_n+\convDiff^{(2,\negTerm)}_n+\convDiff^{(3)}_n$.
  The term $\convDNeg_n$ is written as
  $\convDNeg_n=\convDiff_n^{(1,1)}+\convDiff_n^{(1,2)}+\convDiff_n^{(2,1)}+\convDiff_n^{(2,2)}$ 
  with
  % and each summand is defined at 
  % \eqref{230804.1411}, \eqref{230804.1412}, \eqref{230804.1421} and 
  % \eqref{230804.1422}, respectively,
  \begin{align}
    \convDiff^{(1,1)}_n
    &=
    c_{2,H}\;
    n^{-1}\sum_{j=1}^{n-1} a'_\tj V^{[1]}_\tj 
    B(\bbh^n_j),
    \label{230804.1411}
    \\%\quad
    \convDiff^{(1,2)}_n&=
    - c_{2,H}\;
    \frac1{2n}\rbr{\rbr{V^{[1]}_0 }^2 + \rbr{V^{[1]}_1 }^2},
    \label{230804.1412}
    %\\%\quad
    \\
    \convDiff^{(2,1)}_n &=
    2\;n^{2H-1} \sum_{j=1}^{n-1} V^{[1]}_\tj V^{[(1;1),1]}_\tjm
    I_3\rbr{\diffker^n_j \otimes (\bbone^n_j)^{\otimes2}}
    \label{230804.1421}
    \\
    \convDiff^{(2,2)}_n &=
    n^{2H-1} \sum_{j=1}^{n-1} V^{[1]}_\tj {V^{[(1;1),1]}_\tj} 
    I_3\brbr{\diffker^n_j\otimes{\bbone^n_\jp}^{\otimes2}}
    -
    n^{2H-1} \sum_{j=1}^{n-1} V^{[1]}_\tj {V^{[(1;1),1]}_\tj} 
    I_3\brbr{\diffker^n_j\otimes{\bbone^n_j}^{\otimes2}},
    \label{230804.1422}
  \end{align}
  where
  $\bbh^{n}_{j} = \bbh^{(1) n}_{j} - \bbh^{(2) n}_{j+1}$ with
  $\bbh^{(1) n}_{j}(\cdot) = \bbone^{2n}_{2j}(\cdot)\; n(\cdot-\ttjm)$ and 
  $\bbh^{(2) n}_{j}(\cdot) = \bbone^{2n}_{2j-1}(\cdot)\; n(\ttjm-\cdot).$

  The asymptotic orders of their Malliavin norms are
  $\convDiff_n^{(0)}=O_M(n^{-\half})$ and
  $\convDNeg_n=O_M(n^{-1})$.
  \qed
\end{lemma}
\begin{proof}
  By Lemma \ref{lem:230720.1911}, we have 
$%\begin{align*}
  \convDiff^{(1)}_n =
  \convDiff^{(1,1)}_n + \convDiff^{(1,2)}_n + \convDiff^{(1,\negTerm)}_n
  %{\myred(=O_M(n^{-1})),}
$ %\end{align*}
with
$\convDiff^{(1,\negTerm)}_n = \hat O_M(n^{-2H})$.
Lemma \ref{lem:230925.1301} shows
$\convDiff^{(2)}_n =
\convDiff^{(2,1)}_n+\convDiff^{(2,2)}_n
+ \convDiff^{(2,\negTerm)}_n$ with
$\convDiff^{(2,\negTerm)}_n = O_M(n^{(-\half-H)\vee(H-2)})$
and
$\convDiff^{(3)}_n=O_M(n^{-2H})$.
Setting
$\convDiff^{(\negTerm)}_n=
\convDiff^{(1,\negTerm)}_n+\convDiff^{(2,\negTerm)}_n+\convDiff^{(3)}_n$,
we obtain the decomposition \eqref{eq:230926.1136}.

The order estimates 
$\convDiff_n^{(0)}=O_M(n^{-\half})$ and
$\convDNeg_n=O_M(n^{-1})$ 
are verified also by Lemmas
\ref{lem:230720.1911} and \ref{lem:230925.1301}.
\end{proof}

With the definition of $\psi_n$ at \eqref{eq:230926.1352} and
the above expansion of $\convDiff_n$,
we prove the stochastic expansion of $Z_n$.

\begin{proposition}\label{230720.2014}
  For the functionals $\psi_n$ and $Z_n$ defined at 
  \eqref{eq:230926.1352} and \eqref{eq:230926.1554},
  the following conditions hold:
  \item[(i)]
  $P\sbr{\psi_n<1}=O(n^{-L})$
  for any $L>0$.
  % $\psi_n:\Omega\to[0,1]$, $\psi_n\in\bbD_\infty$ for $n\in\ntwo$, and 
  % $\norm{\psi_n-1}_{k,p}=O(n^{-L})$ as $n\to\infty$ for any 
  % $k\geq1, p>1$ and $L>0$.

  \item[(ii)]
  The sequence $(Z_n)_{n\in\ntwo}$ of functionals
   defined at \eqref{eq:230926.1554}
  % \begin{align}\label{eq:230925.1831}
  %   Z_n = 2\log2\,\sqrt{n} (\hat H_n - H) \psi_n %\tforsm n\in\ntwo
  % \end{align}
  admits a stochastic expansion 
  \begin{align*}
    Z_n = M_n + n^{-1/2}N_n + n^{-1/2}Y_n
  \end{align*}
  with $M_n=\delta(u_n)$,
  $u_n=u_n^{(1)}-u_n^{(2)}$ and 
  \begin{align}
    % u_n^{(1)}&=
    % n^{2H-\half} 
    % \sum_{j\in[n-1]}
    % (\iv)^{-1}
    % \rbr{V^{[1]}_{t^n_j}}^2
    % I_1(\diffker^n_j) \diffker^n_j
    % \\ 
    % u_n^{(2)}&=
    % 2^{2H-1} n^{2H-\half} 
    % \sum_{j\in[2n-1]} (\iv)^{-1} \rbr{V^{[1]}_{t^{2n}_j}}^2
    % I_1(\diffker^{2n}_j) \diffker^{2n}_j
    % \\ 
    u_n^{(i)}&=
    i^{2H-1} n^{2H-\half} 
    \sum_{j\in[in-1]} (\iv)^{-1} \Brbr{V^{[1]}_{t^{in}_j}}^2
    I_1(\diffker^{in}_j) \diffker^{in}_j
    \qquad\tforsm i=1,2
    \label{eq:230926.1446}\\
    N_n&= 
    n (\iv)^{-1} \rbr{\convDNeg_n - \convDNeg_{2n}}
    +O_M(n^{(\half-H)\vee(-1+H)})
    \label{230725.1532}\\
    Y_n&= 
    -n
    \cbr{\half\rbr{\frac{\convDiff_n}\iv}^2 -\half\rbr{\frac{\convDiff_{2n}}\iv}^2},
    \label{230725.1533}
  \end{align}
  where the definition of $\convDNeg_n$ is given at Lemma \ref{230720.2006},
  and $N_n, Y_n=O_M(1)$. 
\end{proposition}

\begin{proof}%[Proof of Proposition \ref{230720.2014}]
  (i)
  The property (i) of $\psi_n$ is verified by Lemma \ref{230510.1711} (i).
  
  \item[(ii)]
  Let $r_n=n^{-\half}$. 
  Then, by Lemma \ref{lem:230926.1026} and Lemma \ref{230720.2006},
  \begin{align}
    Z_n &=
    \sqrt{n} 
    \cbr{\frac{\convDiff_n}\iv - \frac{\convDiff_{2n}}\iv}
    -  \sqrt{n} 
    \cbr{\half\rbr{\frac{\convDiff_n}\iv}^2 -\half\rbr{\frac{\convDiff_{2n}}\iv}^2}
    + \negTerm_n^{(1)}
    \nn
    % \\&=
    % \sqrt{n} \rbr{\iv}^{-1} \rbr{
    % \convDiff^{(0)}_n + \convDNeg_n + \convDiff^{(\negTerm)}_n
    % %\nn\\&\hspsm 
    % -\rbr{
    % \convDiff^{(0)}_{2n} + \convDNeg_{2n} + \convDiff^{(\negTerm)}_{2n}
    % }}
    % %\nn\\&\hspsm
    % -  \sqrt{n} 
    % \cbr{\half\rbr{\frac{\convDiff_n}\iv}^2 -\half\rbr{\frac{\convDiff_{2n}}\iv}^2}
    % + \negTerm_n^{(1)}
    % \nn
    \nn\\&=
    \sqrt{n} \rbr{\iv}^{-1} \rbr{
    \convDiff^{(0)}_n - \convDiff^{(0)}_{2n}}
    +
    \sqrt{n} \rbr{\iv}^{-1} \rbr{\convDNeg_n - \convDNeg_{2n}}
    +
    % \underbrace{
      \sqrt{n} \rbr{\iv}^{-1} \rbr{
      \convDiff^{(\negTerm)}_n - \convDiff^{(\negTerm)}_{2n}}
    % }_{=:\negTerm_n^{(2)} = O_M(n^{(-H)\vee(-\frac32+H)})}
    \nn\\&\quad
    - \sqrt{n} 
    \cbr{\half\rbr{\frac{\convDiff_n}\iv}^2 -\half\rbr{\frac{\convDiff_{2n}}\iv}^2}
    + \Gamma_n^{(1)}
    \nn\\&=
    \bar M_n + r_n\bar N_n + r_nY_n,
    \label{230728.2136}
  \end{align}
  where we set 
  $\bar M_n = \bar M_n^{(1)} - \bar M_n^{(2)}$,
  \begin{align*}
    \bar M_n^{(i)} &= 
    \sqrt{n}\rbr{\iv}^{-1} \convDiff^{(0)}_{in}
    =
    i^{2H-1} n^{2H-1/2} 
    \sum_{j=1}^{in-1} \rbr{\iv}^{-1} \Brbr{V^{[1]}_{t^{in}_j} }^2 I_2\rbr{(\diffker^{in}_j)^{\otimes2}}
    \tfor i=1,2,
    \\
    r_n\bar N_n &=
    \sqrt{n}\rbr{\iv}^{-1} \rbr{\convDNeg_n - \convDNeg_{2n}}
    + 
    % \underbrace{
      \Gamma_n^{(1)}+\Gamma_n^{(2)}
    % }_{=O_M(n^{(-H)\vee(-\frac32+H)})}
    % \\
    % Y_n &=
    % - n 
    % \cbr{\half\rbr{\frac{\convDiff_n}\iv}^2 -\half\rbr{\frac{\convDiff_{2n}}\iv}^2}.
  \end{align*}
  with 
  $\negTerm_n^{(2)} = 
  \sqrt{n} \rbr{\iv}^{-1} \rbr{
    \convDiff^{(\negTerm)}_n - \convDiff^{(\negTerm)}_{2n}}$.
  Notice that 
  $\Gamma_n^{(1)}+\Gamma_n^{(2)}=O_M(n^{(-H)\vee(-\frac32+H)})$.
  % For notational convenience, we denote 
  % \begin{align*}
  %   \bar M_n^{(1)} &= 
  %   n^{\half} \rbr{\iv}^{-1} \convDiff^{(0)}_n
  %   =
  %   n^{2H-1/2} 
  %   \sum_{j=1}^{n-1} \rbr{\iv}^{-1}\rbr{V^{[1]}_\tj }^2 I_2\rbr{(\diffker^n_j)^{\otimes2}}
  %   \\ 
  %   \bar M_n^{(2)} &= 
  %   n^{\half} \rbr{\iv}^{-1} \convDiff^{(0)}_{2n}
  %   =
  %   2^{2H-1} n^{2H-1/2} 
  %   \sum_{j=1}^{2n-1} \rbr{\iv}^{-1} \rbr{V^{[1]}_{t^{2n}_j} }^2 I_2\rbr{(\diffker^{2n}_j)^{\otimes2}}
  %   \\ 
  %   \end{align*}

  Let us denote 
  $M_n^{(i)} = \delta(u_n^{(i)})$ for $i=1,2$.
  Since $M_n^{(i)}$ is written as 
  \begin{align*}
    M_n^{(i)} &= %\delta(u_n^{(i)}) = 
    i^{2H-1} n^{2H-1/2} 
      \sum_{j=1}^{in-1} \rbr{\iv}^{-1} \brbr{V^{[1]}_{t^{in}_j} }^2 
      I_2\rbr{(\diffker^{in}_j)^{\otimes2}}
    %\\&\hspace{70pt}
    - i^{2H-1} n^{2H-1/2} 
    \sum_{j=1}^{in-1} \Babr{D\brbr{\rbr{\iv}^{-1} \brbr{V^{[1]}_{t^{in}_j}}^2 }, \diffker^{in}_j}
    I_1\brbr{\diffker^{in}_j},
  \end{align*}
  we set 
  \begin{align}\label{eq:231005.1844}
    \negTerm_n^{(3,i)} &= \bar M_n^{(i)} - M_n^{(i)} = 
    i^{2H-1} n^{2H-1/2} 
    \sum_{j=1}^{in-1} \Babr{D\brbr{\rbr{\iv}^{-1} \brbr{V^{[1]}_{t^{in}_j}}^2 }, \diffker^{in}_j}
    I_1\brbr{\diffker^{in}_j}
  \end{align}
  and 
  $\negTerm_n^{(3)}:=\negTerm_n^{(3,1)}-\negTerm_n^{(3,2)}$.
  % =\bar M_n^{(1)} - M_n^{(1)} -(\bar M_n^{(2)} - M_n^{(2)})
  The rescaled functional $n^{1/2}\negTerm_n^{(3,i)}$ corresponds to the weighted graph 
  \eqref{fig:230721.1714} in Section \ref{230925.1250},
  whose exponent is $\half-H$.
  Hence $\negTerm_n^{(3,i)}=O_M(n^{-H})$.
  We also have 
  $\bar M_n = \bar M_n^{(1)} - \bar M_n^{(2)}
  = (M_n^{(1)}+\negTerm_n^{(3,1)})
  - (M_n^{(2)}+\negTerm_n^{(3,2)})
  = M_n + \negTerm_n^{(3)}$.

  Setting
  $r_nN_n = r_n\bar N_n + \Gamma_n^{(3)}$,
  we have
  $N_n =
  n\rbr{\iv}^{-1} \rbr{\convDNeg_n - \convDNeg_{2n}}
  + n^\half\sum_{k=1}^3\Gamma_n^{(k)}$
  while the residual term 
  $n^\half\sum_{k=1}^3\Gamma_n^{(k)}$ 
  has the order of 
  $O_M(n^{(\half-H)\vee(-1+H)})$.
  Thus $N_n$ satisfies \eqref{230725.1532},
  and from \eqref{230728.2136}, we obtain 
  \begin{align*}
    Z_n&=M_n+r_nN_n+r_nY_n.
  \end{align*}

  Lemmas \ref{lem:230720.1911} and \ref{lem:230925.1301} show 
  $N_n,Y_n=O_M(1)$.
\end{proof}

\subsection{Order estimates of functionals and related mixed CLT}
\label{230721.1156}
% Here, we collect functionals appearing in the arguments in the following sections, 
% and give estimates of the order of norms of those functionals.
\subsubsection{Estimates of $u_n$}
The Malliavin norm of $u_n$ defined at \eqref{eq:230926.1446}
has the following estimate.
\begin{lemma}\label{230721.1101}
  Let $u_n^{(i)}$ $(i=1,2)$ be as in \eqref{eq:230926.1446}.
  Then the Malliavin norm of $u_n^{(i)}$ is estimated as
  $\nnorm{u_n^{(i)}}_{k,p}=O(1)$ for any $i=1,2$, $k\geq1$ and $p>1$.
  In particular, 
  $\nnorm{u_n}_{k,p}=O(1)$ holds for $k\geq1$ and $p>1$. 
\end{lemma}
\begin{proof}
  For $k_0\in\bbZ_{\geq0}$, the $k_0$-th Malliavin derivative of $u_n^{(i)}$ is expressed as 
the following $\calh^{\otimes k_0+1}$-valued random variable:
\begin{align*}
  D^{k_0} u_n^{(i)} &:= 
  \sum_{k=0}^{k_0}D^{(k_0,k)} u_n^{(i)}
\end{align*}
with
\begin{align*}
  D^{(k_0,0)} u_n^{(i)} &:= 
  i^{2H-1}\, n^{2H-\half} 
  \rbr{\sum_{j\in[in-1]} 
  D^{k_0}_{s_1,...,s_{k_0}}\rbr{\rbr{\iv}^{-1} \brbr{V^{[1]}_{t^{in}_j}}^2} 
  I_1\rbr{\diffker^{in}_j} \diffker^{in}_j(s_0)}_{s_0,s_1,...,s_{k_0}\in[0,1]}
  \\
  D^{(k_0,k)} u_n^{(i)} &:= 
  i^{2H-1}\, n^{2H-\half} 
  \rbr{\sum_{j\in[in-1]} 
  D^{k_0-1}_{s_1,...,\hat s_k,...,s_{k_0}}\rbr{\rbr{\iv}^{-1} \brbr{V^{[1]}_{t^{in}_j}}^2}
  \diffker^{in}_j(s_k)\diffker^{in}_j(s_0)}_{s_0,s_1,...,s_{k_0}\in[0,1]}
\end{align*}
for $k=1,...,k_0$.

We decompose $\nnorm{D^{(k_0,0)} u_n^{(i)}}_{\calh^{\otimes k_0+1}}^2$ as follows:
\begin{align}
  &\norm{D^{(k_0,0)} u_n^{(i)}}_{\calh^{\otimes k_0+1}}^2=
  \abr{D^{(k_0,0)} u_n^{(i)},D^{(k_0,0)} u_n^{(i)}}_{\calh^{\otimes k_0+1}}
  \nn\\&=
  i^{4H-2}\, n^{4H-1}
  \sum_{j\in[in-1]^2} 
  \abr{D^{k_0}\rbr{\rbr{\iv}^{-1} \brbr{V^{[1]}_{t^{in}_{j_1}}}^2},
  D^{k_0}\rbr{\rbr{\iv}^{-1} \brbr{V^{[1]}_{t^{in}_{j_2}}}^2}}_{\calh^{\otimes k_0}}
  I_1\rbr{\diffker^{in}_{j_1}} I_1\rbr{\diffker^{in}_{j_2}} 
  \abr{\diffker^{in}_{j_1},\diffker^{in}_{j_2}}
  \nn\\&=
  i^{4H-2}\, n^{4H-1}
  \sum_{j\in[in-1]^2} 
  \abr{D^{k_0}\rbr{\rbr{\iv}^{-1} \brbr{V^{[1]}_{t^{in}_{j_1}}}^2},
  D^{k_0}\rbr{\rbr{\iv}^{-1} \brbr{V^{[1]}_{t^{in}_{j_2}}}^2}}_{\calh^{\otimes k_0}}
  \abr{\diffker^{in}_{j_1},\diffker^{in}_{j_2}}^2
  \nn\\&\qquad+
  i^{4H-2}\, n^{4H-1}
  \sum_{j\in[in-1]^2} 
  \abr{D^{k_0}\rbr{\rbr{\iv}^{-1} \brbr{V^{[1]}_{t^{in}_{j_1}}}^2},
  D^{k_0}\rbr{\rbr{\iv}^{-1} \brbr{V^{[1]}_{t^{in}_{j_2}}}^2}}_{\calh^{\otimes k_0}}
  I_2\rbr{\diffker^{in}_{j_1}\otimes\diffker^{in}_{j_2}} 
  \abr{\diffker^{in}_{j_1},\diffker^{in}_{j_2}}
  \nn\\&=:
  n^{4H-1}\cali^{(k_0;i;0)}_n+n^{4H-1}\cali^{(k_0;i;1)}_n
  \label{230703.1713}
\end{align}
The functional $\cali^{(k_0;i;0)}_n$ and $\cali^{(k_0;i;1)}_n$
correspond to the weighted graphs 
\eqref{fig:230721.1717} and \eqref{fig:230721.1718},
whose exponents are $1-4H$ and $\half-4H$, respectively.
Hence by Proposition \ref{230725.1133}, we obtain
$\bnorm{D^{(k_0,0)} u_n^{(i)}}_{\calh^{\otimes k_0+1}}^2=
O_M(1)$.

For $k=1,...,k_0$,
\begin{align}
  &\norm{D^{(k_0,k)} u_n^{(i)}}_{\calh^{\otimes k_0+1}}^2=
  \abr{D^{(k_0,k)} u_n^{(i)} ,D^{(k_0,k)} u_n^{(i)} }_{\calh^{\otimes k_0+1}}
  \nn\\&=
  i^{4H-2}\, n^{4H-1}
  \sum_{j\in[in-1]^2} 
  \abr{D^{k_0-1}\rbr{\rbr{\iv}^{-1} \brbr{V^{[1]}_{t^{in}_{j_1}}}^2},
  D^{k_0-1}\rbr{\rbr{\iv}^{-1} \brbr{V^{[1]}_{t^{in}_{j_2}}}^2}}_{\calh^{\otimes k_0-1}}
  \abr{\diffker^{in}_{j_1}, \diffker^{in}_{j_2}}^2
  \nn\\&=
  n^{4H-1}\cali^{(k_0;i;2)}_n
  \label{230703.1728}
\end{align}
The weighted graph corresponding to the functional $\cali^{(k_0;i;2)}_n$
is \eqref{fig:230721.1717}, and its exponent is $1-4H$.
We have the estimate
$\bnorm{D^{(k_0,k)} u_n^{(i)}}_{\calh^{\otimes k_0+1}}^2
=O_M(1)$
by Proposition \ref{230725.1133}.

For $p>2$, by Proposition \ref{221208.1720}, we have
\begin{align*}
  \norm{\bnorm{D^{k_0} u_n^{(i)}}_{\calh^{\otimes k_0+1}}}_{L^p} &\leq
  \sum_{k=0}^{k_0}
  \norm{\bnorm{D^{(k_0,k)} u_n^{(i)}}_{\calh^{\otimes k_0+1}}}_{L^p}
  =
  \sum_{k=0}^{k_0}
  \norm{\bnorm{D^{(k_0,k)} u_n^{(i)}}_{\calh^{\otimes k_0+1}}^2}_{L^{\frac{p}2}}^{\half}
  =O(1)
\end{align*}
for $k_0\geq\bbZ_{\geq0}$.
Hence, for $k\geq1$ and $p>1$, we obtain 
$\nnorm{u_n^{(i)}}_{k,p}=O(1)$
for $i=1,2$, and 
$\nnorm{u_n}_{k,p}=O(1)$.
\end{proof}

\subsubsection{Estimates and limits of $\babr{u_n^{(i_1)}, u_n^{(i_2)}}$}
Since 
$\babr{u_n^{(i_1)}, u_n^{(i_2)}}$
repeatedly appears in the following discussion,
we give some estimates about the functional first.
For $i_1,i_2\in\cbr{1,2}$, we write
$\cali_n^{(i_1,i_2)} = \babr{u_n^{(i_1)}, u_n^{(i_2)}}$,
which is written as
\begin{align*}
  \cali_n^{(i_1,i_2)} &= 
  % \babr{u_n^{(i_1)}, u_n^{(i_2)}}
  % \\&= 
  (i_1\,i_2)^{(2H-1)} n^{4H-1}
  \sum_{j\in[i_1 n-1]\times[i_2 n-1]} 
  \rbr{\iv}^{-2} \brbr{V^{[1]}_{t^{i_1 n}_{j_1}}}^2 \brbr{V^{[1]}_{t^{i_2 n}_{j_2}}}^2 
  I_1\rbr{\diffker^{i_1 n}_{j_1}} I_1\rbr{\diffker^{i_2 n}_{j_2}} 
  \abr{\diffker^{i_1 n}_{j_1},\diffker^{i_2n}_{j_2}}.
\end{align*}
By the product formula, we have 
$I_1\rbr{\diffker^{i_1 n}_{j_1}} I_1\rbr{\diffker^{i_2 n}_{j_2}} =
\abr{\diffker^{i_1 n}_{j_1},\diffker^{i_2n}_{j_2}}+
I_2\rbr{\diffker^{i_1 n}_{j_1}\otimes\diffker^{i_2 n}_{j_2}}$, and hence
\begin{align}\label{eq:230926.2225}
  \cali_n^{(i_1,i_2)} = \cali_n^{(i_1,i_2)(1)} + \cali_n^{(i_1,i_2)(2)}
\end{align}
with 
\begin{align}
  \cali_n^{(i_1,i_2)(1)} &= 
  (i_1\,i_2)^{(2H-1)} n^{4H-1}
  \sum_{j\in[i_1 n-1]\times[i_2 n-1]} 
  \rbr{\iv}^{-2} \brbr{V^{[1]}_{t^{i_1 n}_{j_1}}}^2 \brbr{V^{[1]}_{t^{i_2 n}_{j_2}}}^2 
  \abr{\diffker^{i_1 n}_{j_1},\diffker^{i_2n}_{j_2}}^2
  \label{230703.1731}\\
  \cali_n^{(i_1,i_2)(2)} &= 
  (i_1\,i_2)^{(2H-1)} n^{4H-1}
  \sum_{j\in[i_1 n-1]\times[i_2 n-1]} 
  \rbr{\iv}^{-2} \brbr{V^{[1]}_{t^{i_1 n}_{j_1}}}^2 \brbr{V^{[1]}_{t^{i_2 n}_{j_2}}}^2 
  I_2\rbr{\diffker^{i_1 n}_{j_1}\otimes\diffker^{i_2 n}_{j_2}} 
  \abr{\diffker^{i_1 n}_{j_1},\diffker^{i_2n}_{j_2}}
  \label{230703.1732}
\end{align}
Then, by the theory of exponent, we have the following estimate: 
\begin{lemma}\label{230721.1102}
  For any $i_1,i_2\in\cbr{1,2}$, 
  the Malliavin norms of $\cali_n^{(i_1,i_2)(1)}, \cali_n^{(i_1,i_2)(2)}$ satisfy 
  $\cali_n^{(i_1,i_2)(1)}=O_M(1)$ and 
  $\cali_n^{(i_1,i_2)(2)}=O_M(n^{-\half})$
  as $n\to\infty$.

  In, particular, 
  as $n\to\infty$, the functional $\abr{u_n,u_n}$
  is written as
  \begin{align*}
    \abr{u_n,u_n}&=
    \sum_{i_1,i_2=1}^2 (-1)^{i_1+i_2} \cali_n^{(i_1,i_2)(1)}
    +O_M(n^{-\half}).
  \end{align*}
\end{lemma}
\begin{proof}
  The rescaled functional $n^{-4H+1}\cali_n^{(i_1,i_2)(1)}$ 
  corresponds to the weighted graph \eqref{fig:230721.1717},
  whose exponent is $1-4H$.
  Hence by Proposition \ref{230725.1133},
  $n^{-4H+1}\cali_n^{(i_1,i_2)(1)}=O_M(n^{1-4H})$ and 
  $\cali_n^{(i_1,i_2)(1)}=O_M(1)$.
  Similarly,  the rescaled functional $n^{-4H+1}\cali_n^{(i_1,i_2)(2)}$ 
  corresponds to the weighted graph \eqref{fig:230721.1718},
  whose exponent is $\half-4H$, 
  %Hence $n^{-4H+1}\cali_n^{(i_1,i_2)(2)}=O_M(n^{\half-4H})$ 
  and 
  $\cali_n^{(i_1,i_2)(2)}=O_M(n^{-\half})$.

  Since we have defined $u_n=u_n^{(1)}-u_n^{(2)}$,
  \begin{align*}
    \abr{u_n,u_n}&=
    \Babr{\sum_{i=1,2}(-1)^{i-1}u_n^{(i)},\sum_{i=1,2}(-1)^{i-1}u_n^{(i)}}
    =%\\&=
    % \sum_{i_1,i_2=1}^2 (-1)^{i_1+i_2-2} \babr{u_n^{(i_1)}, u_n^{(i_2)}} 
    % \\&=
    \sum_{i_1,i_2=1}^2 (-1)^{i_1+i_2} \cali_n^{(i_1,i_2)}
    % \\&=
    % \sum_{i_1,i_2=1}^2 (-1)^{i_1+i_2} 
    % \rbr{\cali_n^{(i_1,i_2)(1)} + \cali_n^{(i_1,i_2)(2)}}
    \\&=
    \sum_{i_1,i_2=1}^2 (-1)^{i_1+i_2} \cali_n^{(i_1,i_2)(1)}
    +O_M(n^{-\half})
  \end{align*}
\end{proof}

For notational convenience, we set 
\begin{align}
  \widehat{G}_\infty = \rbr{\iv}^{-2} \int^1_0 \rbr{V^{[1]}_t}^4 dt.
  \label{230728.2030}
\end{align}
\begin{lemma}\label{230724.1350}
  For any $\beta\in(\half,H)$,
  the functionals $\cali_n^{(1,1)(1)}$ and $\cali_n^{(1,2)(1)}$
  are written as
  \begin{align*}
    \cali_n^{(1,1)(1)} 
    &=
    \widehat c \; \widehat{G}_\infty
    +O_M(n^{(-\beta)\vee(-\frac34)})
    \tand%\\
    \cali_n^{(1,2)(1)}%=\cali_n^{(2,1)(1)} 
    =
    2^{2H-1} \tilde c\; \widehat{G}_\infty
    +O_M(n^{(-\beta)\vee(-\frac34)}),
  \end{align*}
  as $n$ tends to $\infty$,
  where the definitions of the constants $\widehat c$ and $\tilde c$ are 
  given at \eqref{eq:230925.1914}.
\end{lemma}
Notice that 
$\cali_n^{(1,2)(1)}=\cali_n^{(2,1)(1)}$ and 
$\cali_n^{(2,2)(1)}=\half\cali_{2n}^{(1,1)(1)}$.
%$\cali_n^{(2,2)(k)}=\half\cali_{2n}^{(1,1)(k)}$ for $k=1,2$.
The proof of this lemma will be given 
as Lemma \ref{lem:230926.1700} in Section \ref{230804.1348}.
We define
\begin{align}\label{eq:230926.1711}
  G_\infty =
  \varconst
  % \rbr{3\hat c - 2^{2H+1}\tilde c}
  \rbr{\iv}^{-2} \int^1_0 \rbr{V^{[1]}_t}^4 dt
  = 
  \varconst
  % \rbr{3\hat c - 2^{2H+1}\tilde c} 
  \widehat{G}_\infty.
\end{align}
with 
\begin{align}\label{eq:231001.1716}
  \varconst={3\hat c - 2^{2H+1}\tilde c}.
\end{align}
The estimate below follows readily from the previous arguments. 
\begin{lemma}\label{230503.1530}
  % {\myblue $\babr{u_n^{(i_1)}, u_n^{(i_2)}}$とlimitの差についての評価}
  As $n\to\infty$, the following estimate holds:
  \begin{align*}
    \abr{u_n, u_n} &= 
    \half G_\infty
    + O_M(n^{-\half}).
  \end{align*}
\end{lemma}
\begin{proof}
  \begin{align*}
    \abr{u_n, u_n} &= 
    % \cali_n^{(1,1)(1)} - 2\:\cali_n^{(1,2)(1)} + \cali_n^{(2,2)(1)} + O_M(n^{-\half})
    % % =
    \cali_n^{(1,1)(1)} - 2\:\cali_n^{(1,2)(1)} + \half\cali_{2n}^{(1,1)(1)} + O_M(n^{-\half})
    \\&=
    \rbr{\frac32\widehat c  - 2^{2H} \tilde c } \widehat{G}_\infty
    %\hat c \; \widehat{G}_\infty - 2^{2H} \tilde c\; 
    % \widehat{G}_\infty + \half \hat c \; \widehat{G}_\infty
    + O_M(n^{-\half})
    =
    \half G_\infty
    + O_M(n^{-\half}).
  \end{align*}
\end{proof}

\subsubsection{Estimates of $\abr{DM_n,u_n}$ and $\babr{DG_\infty, u_n}$, and 
related central limit theorem 
% {\myred $\abr{DM_n^{(i_1)}, u_n^{(i_2)}}$ とか 
% $\abr{D\bar M_n^{(i_1)}, u_n^{(i_2)}}$ でもよい．}
}
\label{230724.1104}

Recall that in the proof of Proposition \ref{230720.2014}
we have defined 
$\bar M_n=\bar M_n^{(1)}- \bar M_n^{(2)}$ with 
\begin{align*}
  \bar M_n^{(i)} &= 
  i^{2H-1} n^{2H-\half} 
  \sum_{j=1}^{in-1}  \rbr{\iv}^{-1} \brbr{V^{[1]}_\tnj{in}{j} }^2 
  I_2\rbr{(\diffker^{in}_j)^{\otimes2}}.
\end{align*}
Recall that  
$M_n^{(i)}=\bar M_n^{(i)}+O_M(n^{-H})$ and 
$M_n=\bar M_n+O_M(n^{-H})$,
which are shown in the proof of Proposition \ref{230720.2014}.

Since $D\bar M_n^{(i)}$ is written as 
\begin{align*}
  D\bar M_n^{(i)} &= 
  2\times
  i^{2H-1} n^{2H-\half} 
  \sum_{j=1}^{in-1}  \rbr{\iv}^{-1} \brbr{V^{[1]}_\tnj{in}{j} }^2 
  I_1(\diffker^{in}_j)\diffker^{in}_j
  +
  i^{2H-1} n^{2H-\half} 
  \sum_{j=1}^{in-1}  D\rbr{\rbr{\iv}^{-1} \brbr{V^{[1]}_\tnj{in}{j} }^2} 
  I_2\rbr{(\diffker^{in}_j)^{\otimes2}},
\end{align*}
we define 
\begin{align*}
  v_n^{(i)} &= 
  i^{2H-1} n^{2H-\half} 
  \sum_{j=1}^{in-1}  D\rbr{\rbr{\iv}^{-1} \brbr{V^{[1]}_\tnj{in}{j} }^2} 
  I_2\rbr{(\diffker^{in}_j)^{\otimes2}},
\end{align*}
for $i=1,2$, and 
$v_n = v_n^{(1)} - v_n^{(2)}$.
We have 
$D\bar M_n^{(i)} = 2 u_n^{(i)} + v_n^{(i)}$ for $i=1,2$ and 
$D\bar M_n = 2 u_n + v_n$.

% {\myred[todelete]
% \begin{align*}
%   u_n^{(i)} &= 
%   i^{2H-1} n^{2H-\half} 
%   \sum_{j\in[in-1]} \rbr{\iv}^{-1} \brbr{V^{[1]}_{t^{in}_j}}^2 
%   I_1\rbr{\diffker^{in}_j} \diffker^{in}_j.
% \end{align*}
% }

For $i_1,i_2\in\cbr{1,2}$,
\begin{align}\label{eq:230926.2226}
  \hspsm &\babr{D\bar M^{(i_1)}_n, u_n^{(i_2)}} =
  \babr{2 u_n^{(i_1)} + v_n^{(i_1)}, u_n^{(i_2)}}
  =
  2\babr{u_n^{(i_1)}, u_n^{(i_2)}} + \babr{v_n^{(i_1)}, u_n^{(i_2)}}
  =
  2\,\cali_n^{(i_1,i_2)} + \calj_n^{(i_1;i_2)},
\end{align}
where we set 
\begin{align}
  % \cali_n^{(i_1,i_2)} &= \babr{u_n^{(i_1)}, u_n^{(i_2)}}
  % \\&= 
  % (i_1\,i_2)^{(2H-1)} n^{4H-1}
  % \sum_{j\in[i_1 n-1]\times[i_2 n-1]} 
  % \rbr{\iv}^{-2} \brbr{V^{[1]}_{t^{i_1 n}_{j_1}}}^2 \brbr{V^{[1]}_{t^{i_2 n}_{j_2}}}^2 
  % I_1\rbr{\diffker^{i_1 n}_{j_1}} I_1\rbr{\diffker^{i_2 n}_{j_2}} 
  % \abr{\diffker^{i_1 n}_{j_1},\diffker^{i_2n}_{j_2}}
  % \\
  \calj_n^{(i_1;i_2)} &= \babr{v_n^{(i_1)}, u_n^{(i_2)}}
  \nn\\&=
  (i_1\,i_2)^{(2H-1)} n^{4H-1} 
  \sum_{\substack{j\in[i_1n-1]\\\times[i_2n-1]}}
  \rbr{\iv}^{-1} \brbr{V^{[1]}_{t^{i_2n}_{j_2}}}^2 
  \Babr{D\Brbr{\brbr{\iv}^{-1} \brbr{V^{[1]}_{t^{i_1n}_{j_1}}}^2},\diffker^{i_2n}_{j_2}}
  I_2\rbr{(\diffker^{i_1n}_{j_1})^{\otimes2}}
  I_1\rbr{\diffker^{i_2n}_{j_2}}.
  \label{eq:231005.1823}
\end{align}
Then, by the theory of exponent, we have the following estimate: 
\begin{lemma}\label{230721.1104}
  % [\sout{{\myblue $\cali_n^{(i_1,i_2)(1)}, \cali_n^{(i_1,i_2)(2)}, 
  % \calj_n^{(i_1;i_2)}$についての評価}}\tto
  % {\mygreen $\abr{u_n,v_n}, \calj_n^{(i_1;i_2)}$についての評価} 
  For any $i_1,i_2\in\cbr{1,2}$, 
  the Malliavin norms of
  $\calj_n^{(i_1;i_2)}$ satisfy 
  $\calj_n^{(i_1;i_2)}=O_M(n^{-H})$,
  as $n\to\infty$.

  As a corollary,
  $\babr{v_n, u_n}=O_M(n^{-H})$.
\end{lemma}
\begin{proof}
  By Lemma \ref{lemma:230616.1754}, 
  $n^{-2H+1}\calj_n^{(i_1;i_2)}$ 
  corresponds to the weighted graph \eqref{fig:230721.1719},
  whose exponent is $1-3H$, 
  %hence $n^{-2H+1}\calj_n^{(i_1;i_2)}=O_M(n^{1-3H})$ 
  and 
  $\calj_n^{(i_1;i_2)}=O_M(n^{-H})$.

  Since 
  $v_n=v_n^{(1)}-v_n^{(2)}$, we have
  $% \begin{align*}
    \abr{v_n,u_n}=
    % \Babr{\sum_{i=1,2}(-1)^{i-1}v_n^{(i)},\sum_{i=1,2}(-1)^{i-1}u_n^{(i)}}
    % =%\\&=
    % \sum_{i_1,i_2=1}^2 (-1)^{i_1+i_2-2} \abr{v_n^{(i_1)},u_n^{(i_2)}}
    % =%\\&=
    \sum_{i_1,i_2=1}^2 (-1)^{i_1+i_2} \calj_n^{(i_1;i_2)}
    %\\&
    =O_M(n^{-H}).
  $%\end{align*}
\end{proof}

Let 
\begin{align*}
  \begin{pmatrix}
    c^{(1,1)}& c^{(1,2)}\\
    c^{(2,1)}& c^{(2,2)}
  \end{pmatrix} =
  \begin{pmatrix}
    2\hat c& 2^{2H} \tilde c\\
    2^{2H} \tilde c& \hat c
  \end{pmatrix}.
\end{align*}
We obtain the following estimates about 
$\babr{DM_n^{(i_1)},u_n^{(i_2)}}$ and
$\babr{DM_n,u_n}$.
\begin{lemma}\label{230721.1235}%{230721.1103}
  (i)
  As $n\to\infty$, it holds that
  \begin{align*}
    \begin{pmatrix}
      \babr{DM_n^{(1)},u_n^{(1)}}& \babr{DM_n^{(1)},u_n^{(2)}}\\
      \babr{DM_n^{(2)},u_n^{(1)}}& \babr{DM_n^{(2)},u_n^{(2)}}
    \end{pmatrix} =
    % \begin{pmatrix}
    %   2\hat c& 2^{2H} \tilde c\\
    %   2^{2H} \tilde c& \hat c
    % \end{pmatrix}
    \begin{pmatrix}
      c^{(1,1)}& c^{(1,2)}\\
      c^{(2,1)}& c^{(2,2)}
    \end{pmatrix}
    \widehat{G}_\infty
    + O_M(n^{-\half}),
  \end{align*}
  where the estimate holds elementwise.

  \item[(ii)]
  As $n\to\infty$, the functionals $\abr{DM_n, u_n}$ 
  is written as
  \begin{align*}
    \abr{DM_n, u_n} &=
    % 2\abr{u_n, u_n} + \abr{v_n, u_n} + O_M(n^{-H})
    % =%\\&=
    G_\infty + O_M(n^{-\half}).
  \end{align*}
\end{lemma}
\begin{proof} 
  (i)
  By Lemmas \ref{230721.1102}, \ref{230724.1350} and \ref{230721.1104},
  $\babr{DM_n^{(i_1)}, u_n^{(i_2)}}$ decomposes as 
  \begin{align*}
    \babr{DM_n^{(i_1)}, u_n^{(i_2)}}
    &=\babr{D\bar M_n^{(i_1)}, u_n^{(i_2)}} + O_M(n^{-H})
    %=2\,\cali_n^{(i_1,i_2)} + \calj_n^{(i_1;i_2)} + O_M(n^{-H})
    \\&=
    2\,(\cali_n^{(i_1,i_2)(1)} + \cali_n^{(i_1,i_2)(2)}) + \calj_n^{(i_1;i_2)} + O_M(n^{-H})
    \\&=
    c^{(i_1,i_2)}\widehat G_\infty + O_M(n^{-\half}).
  \end{align*}

  \item[(ii)]
  The estimate follows from (i) and the following equation:
  \begin{align*}
    \abr{DM_n, u_n}&=
    (1\;-1)
    \begin{pmatrix}
      \babr{DM_n^{(1)},u_n^{(1)}}& \babr{DM_n^{(1)},u_n^{(2)}}\\
      \babr{DM_n^{(2)},u_n^{(1)}}& \babr{DM_n^{(2)},u_n^{(2)}}
    \end{pmatrix} 
    \begin{pmatrix}1 \\-1 \end{pmatrix}.
    % = c_\infty \widehat{G}_\infty + O_M(n^{-\half})
    % = G_\infty + O_M(n^{-\half}),
  \end{align*}
\end{proof}

\begin{lemma}\label{230721.1236}
  As $n\to\infty$, the following estimates holds:
  \item[(i)]
  $D_{u_n}G_\infty=\babr{DG_\infty,u_n}=O_M(n^{-H})$.
  \item[(ii)]
  $D_{u_n}^2G_\infty=O_M(n^{-H})$
\end{lemma}
\begin{proof}
The functional $D_{u_n}G_\infty$ decomposes as 
\begin{align*}
  D_{u_n}G_\infty&=
  \abr{DG_\infty,u_n}
  % =%\\&=
  % \abr{DG_\infty, 
  % \sum_{i=1,2}(-1)^{i-1}u_n^{(i)}}
  =%\\&=
  \sum_{i=1,2}(-1)^{i-1}
  \abr{DG_\infty,u_n^{(i)}},
\end{align*}
and for $i=1,2$, we have
\begin{align}
  \babr{DG_\infty,u_n^{(i)}}&=
  \abr{DG_\infty,
  i^{(2H-1)}%2^{(2H-1)(i-1)}
  n^{2H-\half} 
  \sum_{j\in[in-1]} \rbr{\iv}^{-1} \brbr{V^{[1]}_{t^{in}_j}}^2 
  I_1\rbr{\diffker^{in}_j} \diffker^{in}_j}
  % \\&=
  % i^{(2H-1)}%2^{(2H-1)(i-1)}
  % n^{2H-\half} 
  % \sum_{j\in[in-1]} 
  % \rbr{\iv}^{-1} \brbr{V^{[1]}_{t^{in}_j}}^2 
  % \abr{DG_\infty,
  % \diffker^{in}_j}
  % I_1\rbr{\diffker^{in}_j} 
  \nn\\*&=
  i^{(2H-1)}%2^{(2H-1)(i-1)}
  n^{-\half} 
  \sum_{j\in[in-1]} 
  \rbr{\iv}^{-1} \brbr{V^{[1]}_{t^{in}_j}}^2 
  n^{2H} \abr{DG_\infty,\diffker^{in}_j}
  I_1\rbr{\diffker^{in}_j}.
  \label{eq:231005.1842}
\end{align}
The rescaled functional $n^{\half}\babr{DG_\infty,u_n^{(i)}}$ 
corresponds to the weighted graph \eqref{fig:230721.1714},
whose exponent is $(1+(-\half-H))=\half-H$.
Hence $n^{\half}\babr{DG_\infty,u_n^{(i)}}=O_M(n^{\half-H})$ and 
$\babr{DG_\infty,u_n^{(i)}}=O_M(n^{-H})$.
We obtain 
$\babr{DG_\infty,u_n}=O_M(n^{-H})$, 
and by Lemma \ref{230721.1101},
$D_{u_n}^2G_\infty=O_M(n^{-H})$.
\end{proof}
\begin{remark}
  Although we only prove 
  $D_{u_n}^2G_\infty=O_M(n^{-H})$ here,
  for it is enough for the following argument,
  we can prove 
  $D_{u_n}^2G_\infty=O_M(n^{-H-\half})$ 
  by Proposition \ref{221223.1757}.
\end{remark}

% \subsubsection{CLT}
Define a $2\times 2$ random matrix $A_\infty$ by
\begin{align}
  A_\infty = 
  \begin{pmatrix}
    \sqrt{2\hat c} & \frac{2^{2H}\tilde c}{\sqrt{2\hat c}}\\
    0 & \sqrt{\hat c -\frac{\rbr{2^{2H}\tilde c}^2}{2\hat c}}
  \end{pmatrix}
  \brbr{\widehat{G}_\infty}^{\half}.
  % = 
  % \rbr{2\hat c}^{-\half}\;\rbr{\widehat{G}_\infty}^{\half}
  % \begin{pmatrix}
  %   2\hat c & {2^{2H}\tilde c}\\
  %   0 & \sqrt{2{\hat c}^2 -{\rbr{2^{2H}\tilde c}^2}}
  % \end{pmatrix}
  \label{230726.1711}
\end{align}
Let $\zeta=(\zeta_1,\zeta_2)$ be 
a $2$-dimensional standard normal random variable independent of $\calf$.
Then the random variable $\zeta A_\infty$
is the $2$-dimensional mixed normal random variable whose 
conditional covariance matrix is 
\begin{align}
  % \begin{pmatrix}
  %   \babr{DM_n^{(1)},u_n^{(1)}}& \babr{DM_n^{(1)},u_n^{(2)}}\\
  %   \babr{DM_n^{(2)},u_n^{(1)}}& \babr{DM_n^{(2)},u_n^{(2)}}
  % \end{pmatrix} =
  A_\infty^T A_\infty=
  \begin{pmatrix}
    2\widehat c& 2^{2H} \tilde c\\
    2^{2H} \tilde c& \widehat c
  \end{pmatrix}
  \widehat{G}_\infty
  =
  \begin{pmatrix}
    c^{(1,1)}& c^{(1,2)}\\
    c^{(2,1)}& c^{(2,2)}
  \end{pmatrix} 
  \widehat{G}_\infty
  =:G_\infty^{[2]},
  \label{230726.1712}
\end{align}
which is the limit of 
$\brbr{\babr{DM_n^{(i_1)},u_n^{(i_2)}}}_{i_1,i_2\in\cbr{1,2}}$.
(See Lemma \ref{230721.1235} (i).)

\begin{proposition}\label{230721.1408}
  As $n\to\infty$, the following weak convergences hold:
  \item[(i)]
  $%\begin{align*}
    \brbr{M_n^{(1)},M_n^{(2)}}=
    \brbr{\delta(u_n^{(1)}),\delta(u_n^{(2)})}
    \overset{d}{\to} \zeta A_\infty 
    \sim\caln_2(0,G_\infty^{[2]})
  $%\end{align*}
  \item[(ii)]
  $Z_n\overset{d}{\to} \caln_1(0,G_\infty)$.
  % \begin{comment}
  %   {\myred 話のつながり的には$Z_n^\circ$についてのほうがいいのでは？}
  % \end{comment}
\end{proposition}
\begin{proof}
  By an argument similar to the proof of Lemma \ref{230721.1236} (i), we have 
  \begin{align}\label{eq:230926.2252}
    D_{u_n^{(i_3)}}(G_\infty^{[2]})_{i_1,i_2}=O_M(n^{-H}).
  \end{align}
  %$D_{u_n^{(i)}}\widehat{G}_\infty=O_M(n^{-H})$
  By this estimate and  Lemma \ref{230721.1235} (i), 
  we have the weak convergence (i).

  Since $Z_n=M_n+O_M(n^{-\half})=M_n^{(1)}-M_n^{(2)}+O_M(n^{-\half})$ and 
  $(1,-1)G_\infty^{[2]}(1,-1)^T=\rbr{3\widehat c - 2^{2H+1}\tilde c} \widehat{G}_\infty
  =G_\infty$,
  we have the second convergence.
\end{proof}
\begin{comment}  
\begin{itembox}[l]{Memo}
  \myred ここのCLTの記述(Lemmaの立て方とか)はもう一度読み直してこれでいいのか考えたほうがいい気がする．

  main resultに書いてあった文章案:
  The first estimate implies that $G_\infty$ plays the role of 
the asymptotic variance.
From the above estimates,
we have the weak convergence of $M_n$ to the mixed normal distribution 
$\caln(0,G_\infty)$
by the argument of {\myblue NNP \cite{nourdin2016quantitative}}.
(We can also prove the stable convergence. 
See also Remark 3.4 of {\cite{nualart2019asymptotic}}.)
\end{itembox}
\end{comment}

\begin{remark}
  The detail of the proof of (i) is omitted, 
  since it is essentially contained in 
  the proof of the asymptotic expansion (Proposition \ref{230726.1831}).
  However the two estimates 
  (Lemma \ref{230721.1235} (i) and \eqref{eq:230926.2252}) 
  constitute the key conditions when proving a mixed central limit theorem
  on the Wiener space. 
  % We can show the stable convergence as well.
  (See \cite{nourdin2016quantitative}.)
\end{remark}

\begin{todelete}
\begin{itembox}[l]{Memo}
  \begin{itemize}
    \item 
  $G_\infty =
\rbr{3\hat c - 2^{2H+1}\tilde c}
\rbr{\iv}^{-2} \int^1_0 \rbr{V^{[1]}_t}^4 dt
= 
\rbr{3\hat c - 2^{2H+1}\tilde c} \widehat{G}_\infty.$

\item
Consider $Z_n=M_n+r_nN_n$, $X_\infty$ and $G_\infty=A_\infty^T A_\infty$.
(この$G_\infty$は一般論のものなので，$d\times d$random matrix．)
These satisfy some conditions.
\begin{lemma}\label{230504.1201}
  Assume that
  $\norm{\abr{DM_n, u_n}-G_\infty}_{L^1}\to0,$
  $r_n\norm{\abr{DN_n, u_n}}_{L^1}\to0,$
  $\norm{\abr{DX_\infty, u_n}}_{L^1}\to0,$
  $\norm{\abr{DG_\infty, u_n}}_{L^1}\to0,$ and 
  $r_n\norm{N_n}_{L^1}\to0$ 
  elementwisely.
  In other words, 
  $\norm{\abr{DM_n^{(i)}, u_n^{(j)}}-G_\infty^{(i,j)}}_{L^1}\to0$,
  $r_n\norm{\abr{DN_n^{(i)}, u_n^{(j)}}}_{L^1}\to0,$
  $\norm{\abr{DX_\infty^{(\ell)}, u_n^{(j)}}}_{L^1}\to0,$
  $\norm{\abr{DG_\infty^{(i,j)}, u_n^{(k)}}}_{L^1}\to0,$ and 
  $r_n\norm{N_n^{(i)}}_{L^1}\to0$ 
  for $i,j,k\in\cbr{1,...d}$ and $\ell\in\cbr{1,..,d_1}$.

  Then $(Z_n, X_\infty)$ weakly converges to $(\zeta A_\infty, X_\infty)$, 
  where $\zeta$ is a $d$-dimensional standard normal random variable independent of $\calf$.
\end{lemma}
\end{itemize}
\end{itembox}
\end{todelete}

\subsection{Asymptotic expansion of $Z_n^\circ$}\label{sec:231001.1030}
By Proposition \ref{230721.1408}, 
$Z_n=2\log2 \sqrt{n} (\hat H_n-H)$ 
converges weakly to the mixed normal distribution $\caln(0,G_\infty)$
with $G_\infty$ defined at \eqref{eq:230926.1711}.
To obtain the asymptotic expansion of the distribution of $Z_n$,
we first investigate that of $Z_n^\circ$ defined by
\begin{align}
  Z_n^\circ = \delta(u_n) + r_n N_n,
  \label{230726.1348}
\end{align}
where the definitions of 
$u_n$ and $N_n$ are given at 
\eqref{eq:230926.1446} and \eqref{230725.1532} respectively,
by applying the theory of asymptotic expansion of \cite{nualart2019asymptotic}.
Specifically speaking, we will confirm Condition {\bf[D]} for 
$u_n, N_n$ and $G_\infty$
and prove the asymptotic expansion of $Z_n^\circ$ 
by Theorem \ref{thm:230927.1318}.

\subsubsection{Functionals related to asymptotic expansion}
Here we deal with functionals appearing in the context of 
asymptotic expansion of Skorohod integrals,
that is  $D_{u_n}^2M_n$, $N_n$ and $D_{u_n}N_n$.
These functionals appear in both (ii) and (iii)
of the Condition {\bf[D]}.

\subsubsection*{(i) Decomposition and estimates about $D_{u_n}^2M_n$}
From \eqref{eq:230926.2226}, \eqref{eq:230926.2225} and 
Lemma \ref{230721.1104}, we have 
\begin{align*}
  \hspsm &\babr{D\bar M_n^{(i_1)},u_n^{(i_2)}} =
  2\, \cali_n^{(i_1,i_2)}+ \calj_n^{(i_1;i_2)} = 
  2\, \brbr{\cali_n^{(i_1,i_2)(1)} + {\cali_n^{(i_1,i_2)(2)}}}
  +O_M(n^{-H})
\end{align*}
for $i_1,i_2\in\cbr{1,2}$, where
\begin{align*}
  \cali_n^{(i_1,i_2)(1)} &= 
  (i_1i_2)^{(2H-1)}%2^{(2H-1)(i_1+i_2-2)}
  n^{4H-1}
  \sum_{j\in[i_1 n-1]\times[i_2 n-1]} 
  \rbr{\iv}^{-2} \brbr{V^{[1]}_{t^{i_1 n}_{j_1}}}^2 \brbr{V^{[1]}_{t^{i_2 n}_{j_2}}}^2 
  \babr{\diffker^{i_1 n}_{j_1},\diffker^{i_2n}_{j_2}}^2
  \\
  \cali_n^{(i_1,i_2)(2)} &= 
  (i_1i_2)^{(2H-1)}%2^{(2H-1)(i_1+i_2-2)}
  n^{4H-1}
  \sum_{j\in[i_1 n-1]\times[i_2 n-1]} 
  \rbr{\iv}^{-2} \brbr{V^{[1]}_{t^{i_1 n}_{j_1}}}^2 \brbr{V^{[1]}_{t^{i_2 n}_{j_2}}}^2 
  I_2\brbr{\diffker^{i_1 n}_{j_1}\otimes\diffker^{i_2 n}_{j_2}} 
  \babr{\diffker^{i_1 n}_{j_1},\diffker^{i_2n}_{j_2}}
  \\ 
  u_n^{(i)} &= 
  i^{(2H-1)}%2^{(2H-1)(i-1)} 
  n^{2H-\half} 
  \sum_{j\in[in-1]} \rbr{\iv}^{-1} \brbr{V^{[1]}_{t^{in}_j}}^2 
  I_1\brbr{\diffker^{in}_j} \diffker^{in}_j.
\end{align*}
\begin{lemma}\label{230724.1710}
  For $i\in\cbr{1,2}^{\cbr{1,2,3}}$,
  the functional
  $D_{u_n^{(i_3)}}D_{u_n^{(i_2)}}\bar M_n^{(i_1)}$
  decomposes as 
  \begin{align*}
    D_{u_n^{(i_3)}}D_{u_n^{(i_2)}}\bar M_n^{(i_1)}
    &=
    2\times
    (i_1i_2i_3)^{(2H-1)}%2^{(2H-1)(i_1+i_2+i_3-3)}
    n^{(6H-\frac32)}\cali_n^{(3;i)(1)}
    +O_M(n^{-H})
  \end{align*}
  with 
  \begin{align}
    \cali_n^{(3;i)(1)}=
    2\,\sum_{j\in\bbJ_n([3],i)} 
    \rbr{\iv}^{-3} 
    \brbr{V^{[1]}_{t^{i_3n}_{j_3}}}^2 \brbr{V^{[1]}_{t^{i_1 n}_{j_1}}}^2 \brbr{V^{[1]}_{t^{i_2 n}_{j_2}}}^2 
    \abr{\diffker^{i_1 n}_{j_1},\diffker^{i_2n}_{j_2}}
    \abr{\diffker^{i_2 n}_{j_2},\diffker^{i_3n}_{j_3}}
    \abr{\diffker^{i_1 n}_{j_1},\diffker^{i_3n}_{j_3}},
  \end{align}
  and
  $n^{6H-\frac32}\cali_n^{(3;i)(1)}=O_M(n^{-\half})$.
  In particular, 
  $D_{u_n^{(i_3)}}D_{u_n^{(i_2)}}\bar M_n^{(i_1)}=O_M(n^{-\half})$.
\end{lemma}
\begin{proof}
By Lemma \ref{230721.1101},
\begin{align}\label{230724.1657}
  D_{u_n^{(i_3)}}D_{u_n^{(i_2)}}\bar M_n^{(i_1)}
  &=
  % D_{u_n^{(i_3)}}\rbr{2\, \cali_n^{(i_1,i_2)(1)} + 2\,\cali_n^{(i_1,i_2)(2)} + \calj_n^{(i_1;i_2)}}
  % \\
  % &=
  2\,D_{u_n^{(i_3)}}\rbr{\cali_n^{(i_1,i_2)(1)} + \cali_n^{(i_1,i_2)(2)}}
  +O_M(n^{-H})
\end{align}
Hence we deal with 
$D_{u_n^{(i_3)}}\,\cali_n^{(i_1,i_2)(1)}$ and 
$D_{u_n^{(i_3)}}\,\cali_n^{(i_1,i_2)(2)}$.
As for the first functional, 
\begin{align*}
  &
  % (i_1i_2)^{-(2H-1)}%2^{-(2H-1)(i_1+i_2-2)}
  % n^{-(4H-1)}
  % (i_3)^{-(2H-1)}%2^{-(2H-1)(i_3-1)}
  % n^{-(2H-\half)} 
  % D_{u_n^{(i_3)}}\,\cali_n^{(i_1,i_2)(1)}
  % \\&=
  %(i_1i_2)^{(2H-1)}%2^{(2H-1)(i_1+i_2-2)}
  %(i_1i_2i_3)^{(2H-1)}%2^{(2H-1)(i_1+i_2+i_3-3)}
  (i_1i_2i_3)^{-(2H-1)}%2^{-(2H-1)(i_1+i_2+i_3-3)}
  n^{-4H+\frac32}
  %n^{-6H+\frac32}
  \abr{D\cali_n^{(i_1,i_2)(1)}, u_n^{(i_3)}}
  \\&=
  n^{2H}
  \abr{D\rbr{
  \sum_{j\in[i_1 n-1]\times[i_2 n-1]} 
  \rbr{\iv}^{-2} \brbr{V^{[1]}_{t^{i_1 n}_{j_1}}}^2 \brbr{V^{[1]}_{t^{i_2 n}_{j_2}}}^2 
  \abr{\diffker^{i_1 n}_{j_1},\diffker^{i_2n}_{j_2}}^2},
  \sum_{j\in[i_3n-1]} \rbr{\iv}^{-1} \brbr{V^{[1]}_{t^{i_3n}_j}}^2 
  I_1\rbr{\diffker^{i_3n}_j} \diffker^{i_3n}_j.
  }
  % \\&=
  % %\sum_{j\in[i_1 n-1]\times[i_2 n-1]} \sum_{j\in[i_3n-1]} 
  % \sum_{j\in\bbJ_n([3],i)} 
  % \abr{D\rbr{
  % \rbr{\iv}^{-2} \brbr{V^{[1]}_{t^{i_1 n}_{j_1}}}^2 \brbr{V^{[1]}_{t^{i_2 n}_{j_2}}}^2 
  % \abr{\diffker^{i_1 n}_{j_1},\diffker^{i_2n}_{j_2}}^2},
  % \rbr{\iv}^{-1} \brbr{V^{[1]}_{t^{i_3n}_{j_3}}}^2 
  % I_1\rbr{\diffker^{i_3n}_{j_3}} \diffker^{i_3n}_{j_3}.
  % }
  \\&=
  %\sum_{j\in[i_1 n-1]\times[i_2 n-1]\times[i_3n-1]} 
  %n^{-2H}
  \sum_{j\in\bbJ_n([3],i)} 
  n^{2H}
  \abr{D\rbr{
    \rbr{\iv}^{-2} \brbr{V^{[1]}_{t^{i_1 n}_{j_1}}}^2 \brbr{V^{[1]}_{t^{i_2 n}_{j_2}}}^2},
    \diffker^{i_3n}_{j_3}}
  \rbr{\iv}^{-1} \brbr{V^{[1]}_{t^{i_3n}_{j_3}}}^2 
  I_1\rbr{\diffker^{i_3n}_{j_3}} 
  \abr{\diffker^{i_1 n}_{j_1},\diffker^{i_2n}_{j_2}}^2.
\end{align*}
We used $\bbJ_n([3],i)=[i_1 n-1]\times[i_2 n-1]\times[i_3n-1]$.
The above functional 
corresponds to the weighted graph \eqref{fig:230721.1720},
whose exponent is $(1+(-\half-H)) + (1-4H)=\frac32-5H$.
Hence, 
\begin{align}
  D_{u_n^{(i_3)}}\,\cali_n^{(i_1,i_2)(1)}= 
  \abr{D\cali_n^{(i_1,i_2)(1)}, u_n^{(i_3)}} 
  = O_M(n^{(\frac32-5H)+(4H-\frac32)})
  = O_M(n^{-H})
  \label{230628.1149}
\end{align}

% \noindent{\myblue
% expo($\cali_n^{(i_1,i_2)(1)}$を$A_{n,j}$に入れた場合): 
% $(2H-\half) + (-2H) + (1+(-\half-H))=-H$
% \begin{align*}
%   \abr{D\cali_n^{(i_1,i_2)(1)}, u_n^{(i_3)}}&=
%   \abr{D\cali_n^{(i_1,i_2)(1)}, 
%   2^{(2H-1)(i_3-1)} n^{2H-\half} 
%   \sum_{j\in[i_3n-1]} \rbr{\iv}^{-1} \brbr{V^{[1]}_{t^{i_3n}_j}}^2 
%   I_1\rbr{\diffker^{i_3n}_j} \diffker^{i_3n}_j}
%   \\&=
%   2^{(2H-1)(i_3-1)} n^{2H-\half} 
%   \sum_{j\in[i_3n-1]} \rbr{\iv}^{-1} \brbr{V^{[1]}_{t^{i_3n}_j}}^2 
%   \abr{D\cali_n^{(i_1,i_2)(1)},\diffker^{i_3n}_j}
%   I_1\rbr{\diffker^{i_3n}_j} 
% \end{align*}}

Concerning $D_{u_n^{(i_3)}}\,\cali_n^{(i_1,i_2)(2)}$, 
\begin{align}
  & 
  % (i_1i_2)^{-(2H-1)}%2^{-(2H-1)(i_1+i_2-2)}
  % n^{-(4H-1)}
  % (i_3)^{-(2H-1)}%2^{-(2H-1)(i_3-1)}
  % n^{-(2H-\half)}
  % D_{u_n^{(i_3)}}\,\cali_n^{(i_1,i_2)(2)}
  % \\&=
  (i_1i_2i_3)^{-(2H-1)}%2^{-(2H-1)(i_1+i_2+i_3-3)}
  n^{-(6H-\frac32)}
  \abr{D\cali_n^{(i_1,i_2)(2)},u_n^{(i_3)}}
  % \\&=
  % \Bigg\langle D\rbr{
  % \sum_{j\in[i_1 n-1]\times[i_2 n-1]} 
  % \rbr{\iv}^{-2} \brbr{V^{[1]}_{t^{i_1 n}_{j_1}}}^2 \brbr{V^{[1]}_{t^{i_2 n}_{j_2}}}^2 
  % I_2\rbr{\diffker^{i_1 n}_{j_1}\otimes\diffker^{i_2 n}_{j_2}} 
  % \abr{\diffker^{i_1 n}_{j_1},\diffker^{i_2n}_{j_2}}},
  % \\&\qquad
  % \sum_{j\in[i_3n-1]} \rbr{\iv}^{-1} \brbr{V^{[1]}_{t^{i_3n}_j}}^2 
  % I_1\rbr{\diffker^{i_3n}_j} \diffker^{i_3n}_j
  % \Bigg\rangle
  \nn\\&=
  \sum_{j\in\bbJ_n([3],i)} 
  %\sum_{j\in[i_1 n-1]\times[i_2 n-1]\times[i_3n-1]} 
  \rbr{\iv}^{-1} \brbr{V^{[1]}_{t^{i_3n}_{j_3}}}^2 
  \abr{D\rbr{
  \rbr{\iv}^{-2} \brbr{V^{[1]}_{t^{i_1 n}_{j_1}}}^2 \brbr{V^{[1]}_{t^{i_2 n}_{j_2}}}^2 
  I_2\rbr{\diffker^{i_1 n}_{j_1}\otimes\diffker^{i_2 n}_{j_2}}},
  \diffker^{i_3n}_{j_3}}
  %\\&\hspace{50pt}\times
  I_1\rbr{\diffker^{i_3n}_{j_3}} 
  \abr{\diffker^{i_1 n}_{j_1},\diffker^{i_2n}_{j_2}}
  \nn\\
  % &=
  % \sum_{j\in\bbJ_n([3],i)} 
  % %\sum_{j\in[i_1 n-1]\times[i_2 n-1]\times[i_3n-1]} 
  % \rbr{\iv}^{-1} \brbr{V^{[1]}_{t^{i_3n}_{j_3}}}^2 
  % \abr{D\rbr{
  % \rbr{\iv}^{-2} \brbr{V^{[1]}_{t^{i_1 n}_{j_1}}}^2 \brbr{V^{[1]}_{t^{i_2 n}_{j_2}}}^2},
  % \diffker^{i_3n}_{j_3}}
  % %\\&\hspace{50pt}\times
  % I_2\rbr{\diffker^{i_1 n}_{j_1}\otimes\diffker^{i_2 n}_{j_2}}
  % I_1\rbr{\diffker^{i_3n}_{j_3}} 
  % \abr{\diffker^{i_1 n}_{j_1},\diffker^{i_2n}_{j_2}}
  % \\&\quad+
  % \sum_{j\in\bbJ_n([3],i)} 
  % %\sum_{j\in[i_1 n-1]\times[i_2 n-1]\times[i_3n-1]} 
  % \rbr{\iv}^{-3} 
  % \brbr{V^{[1]}_{t^{i_3n}_{j_3}}}^2 \brbr{V^{[1]}_{t^{i_1 n}_{j_1}}}^2 \brbr{V^{[1]}_{t^{i_2 n}_{j_2}}}^2 
  % I_1\rbr{\diffker^{i_2n}_{j_2}}
  % I_1\rbr{\diffker^{i_3n}_{j_3}} 
  % \abr{\diffker^{i_1 n}_{j_1},\diffker^{i_3n}_{j_3}}
  % \abr{\diffker^{i_1 n}_{j_1},\diffker^{i_2n}_{j_2}}
  % \\&\quad+
  % \sum_{j\in\bbJ_n([3],i)} 
  % %\sum_{j\in[i_1 n-1]\times[i_2 n-1]\times[i_3n-1]} 
  % \rbr{\iv}^{-3} 
  % \brbr{V^{[1]}_{t^{i_3n}_{j_3}}}^2 \brbr{V^{[1]}_{t^{i_1 n}_{j_1}}}^2 \brbr{V^{[1]}_{t^{i_2 n}_{j_2}}}^2 
  % I_1\rbr{\diffker^{i_1 n}_{j_1}}
  % I_1\rbr{\diffker^{i_3n}_{j_3}} 
  % \abr{\diffker^{i_2 n}_{j_2},\diffker^{i_3n}_{j_3}}
  % \abr{\diffker^{i_1 n}_{j_1},\diffker^{i_2n}_{j_2}}
  % \\
  &=
  n^{-2H}
  \sum_{j\in\bbJ_n([3],i)} 
  %\sum_{j\in[i_1 n-1]\times[i_2 n-1]\times[i_3n-1]} 
  \rbr{\iv}^{-1} \brbr{V^{[1]}_{t^{i_3n}_{j_3}}}^2 
  n^{2H}
  \Babr{D\Brbr{
  \rbr{\iv}^{-2} \brbr{V^{[1]}_{t^{i_1 n}_{j_1}}}^2 \brbr{V^{[1]}_{t^{i_2 n}_{j_2}}}^2},
  \diffker^{i_3n}_{j_3}}
  %\\&\hspace{50pt}\times
  I_2\rbr{\diffker^{i_1 n}_{j_1}\otimes\diffker^{i_2 n}_{j_2}}
  I_1\rbr{\diffker^{i_3n}_{j_3}} 
  \abr{\diffker^{i_1 n}_{j_1},\diffker^{i_2n}_{j_2}}
  \nn\\&\quad+
  \sum_{j\in\bbJ_n([3],i)} 
  %\sum_{j\in[i_1 n-1]\times[i_2 n-1]\times[i_3n-1]} 
  \rbr{\iv}^{-3} 
  \brbr{V^{[1]}_{t^{i_3n}_{j_3}}}^2 \brbr{V^{[1]}_{t^{i_1 n}_{j_1}}}^2 \brbr{V^{[1]}_{t^{i_2 n}_{j_2}}}^2 
  I_2\rbr{\diffker^{i_2n}_{j_2}\otimes\diffker^{i_3n}_{j_3}}
  \abr{\diffker^{i_1 n}_{j_1},\diffker^{i_3n}_{j_3}}
  \abr{\diffker^{i_1 n}_{j_1},\diffker^{i_2n}_{j_2}}
  \nn\\&\quad+
  \sum_{j\in\bbJ_n([3],i)} 
  %\sum_{j\in[i_1 n-1]\times[i_2 n-1]\times[i_3n-1]} 
  \rbr{\iv}^{-3} 
  \brbr{V^{[1]}_{t^{i_3n}_{j_3}}}^2 \brbr{V^{[1]}_{t^{i_1 n}_{j_1}}}^2 \brbr{V^{[1]}_{t^{i_2 n}_{j_2}}}^2 
  I_2\rbr{\diffker^{i_1 n}_{j_1}\otimes\diffker^{i_3n}_{j_3}} 
  \abr{\diffker^{i_2 n}_{j_2},\diffker^{i_3n}_{j_3}}
  \abr{\diffker^{i_1 n}_{j_1},\diffker^{i_2n}_{j_2}}
  \nn\\&\quad+
  2\,\sum_{j\in\bbJ_n([3],i)} 
  %2\,\sum_{j\in[i_1 n-1]\times[i_2 n-1]\times[i_3n-1]} 
  \rbr{\iv}^{-3} 
  \brbr{V^{[1]}_{t^{i_3n}_{j_3}}}^2 \brbr{V^{[1]}_{t^{i_1 n}_{j_1}}}^2 \brbr{V^{[1]}_{t^{i_2 n}_{j_2}}}^2 
  \abr{\diffker^{i_1 n}_{j_1},\diffker^{i_2n}_{j_2}}
  \abr{\diffker^{i_2 n}_{j_2},\diffker^{i_3n}_{j_3}}
  \abr{\diffker^{i_1 n}_{j_1},\diffker^{i_3n}_{j_3}} 
  \nn\\&=:
  n^{-2H}\cali_n^{(3;i)(4)}
  +\cali_n^{(3;i)(2)}
  +\cali_n^{(3;i)(3)}
  +\cali_n^{(3;i)(1)}
  \label{eq:231005.1831}
\end{align}
% \begin{itemize}
%   \item exponent ($\cali_n^{(3;i)(1)}$): $1-6H$
%   \item exponent ($\cali_n^{(3;i)(2)}$): $((1-4H)+(-\half-H-H))=\half-6H$
%   \item exponent ($\cali_n^{(3;i)(3)}$): $((1-4H)+(-\half-H-H))=\half-6H$
%   \item exponent ($\cali_n^{(3;i)(4)}$): $((1-2H)+(-\half-H-H)) + (1+(-\half-H))=1-5H$
% \end{itemize}
The functionals
$\cali_n^{(3;i)(1)}$, $\cali_n^{(3;i)(2)}$, $\cali_n^{(3;i)(3)}$ and $\cali_n^{(3;i)(4)}$
correspond to the weighted graphs 
\eqref{fig:230721.1723},
\eqref{fig:230721.1722},
\eqref{fig:230721.1722} and
\eqref{fig:230721.1721}, 
whose exponents are 
$1-6H$, $\half-6H$, $\half-6H$ and $1-5H$,
respectively.
Hence, we have 
$n^{(6H-\frac32)}\cali_n^{(3;i)(1)}=O_M(n^{-\half})$,\; 
$n^{(6H-\frac32)}\cali_n^{(3;i)(2)}=O_M(n^{-1})$,\; 
$n^{(6H-\frac32)}\cali_n^{(3;i)(3)}=O_M(n^{-1})$ and 
$n^{(4H-\frac32)}\cali_n^{(3;i)(4)}
%n^{(6H-\frac32)}n^{-2H}\cali_n^{(3;i)(4)}
=O_M(n^{-\half-H})$,
and we obtain 
\begin{align}
\babr{D\cali_n^{(i_1,i_2)(2)},u_n^{(i_3)}}=
%(i_1i_2i_3)^{-(2H-1)}%2^{-(2H-1)(i_1+i_2+i_3-3)}
(i_1i_2i_3)^{(2H-1)}%2^{(2H-1)(i_1+i_2+i_3-3)}
n^{(6H-\frac32)}\cali_n^{(3;i)(1)}
+O_M(n^{-1}).
\label{230628.1150}
\end{align}

From \eqref{230724.1657}, \eqref{230628.1149} and \eqref{230628.1150},
we obtain 
\begin{align*}
  D_{u_n^{(i_3)}}D_{u_n^{(i_2)}}\bar M_n^{(i_1)}
  &=
  % 2\,D_{u_n^{(i_3)}}\,\cali_n^{(i_1,i_2)(1)}
  % +2\,D_{u_n^{(i_3)}}\cali_n^{(i_1,i_2)(2)}
  % +O_M(n^{-H})
  % \\&=
  2\times
  (i_1i_2i_3)^{(2H-1)}%2^{(2H-1)(i_1+i_2+i_3-3)}
  n^{(6H-\frac32)}\cali_n^{(3;i)(1)}
  +O_M(n^{-H}).
\end{align*}
\end{proof}
Hence we obtain the following estimates about 
$(D_{u_n})^2 M_n$ and $(D_{u_n})^3 M_n$.
\begin{lemma}\label{230721.1259}
  \item[(i)]
  The functional $D_{u_n}D_{u_n} M_n$ decomposes as
  \begin{align}
    &D_{u_n}D_{u_n} M_n
    \nn\\&=
    2^2\times n^{(6H-\frac32)}
    \sum_{j\in[n-1]^3}%\bbJ_n([3],i)} 
    \rbr{\iv}^{-3} 
    \brbr{V^{[1]}_{t^{n}_{j_1}}}^2 \brbr{V^{[1]}_{t^{n}_{j_2}}}^2 \brbr{V^{[1]}_{t^{n}_{j_3}}}^2 
    \abr{\diffker^{n}_{j_1},\diffker^{n}_{j_2}}
    \abr{\diffker^{n}_{j_2},\diffker^{n}_{j_3}}
    \abr{\diffker^{n}_{j_1},\diffker^{n}_{j_3}} 
    \nn\\
    &\quad-3\times
    2^2\times2^{(2H-1)}\times n^{(6H-\frac32)}
    \sum_{j\in[n-1]^2\times[2n-1]} 
    \rbr{\iv}^{-3} 
    \brbr{V^{[1]}_{t^{n}_{j_1}}}^2 \brbr{V^{[1]}_{t^{n}_{j_2}}}^2 \brbr{V^{[1]}_{t^{2n}_{j_3}}}^2 
    \abr{\diffker^{n}_{j_1},\diffker^{n}_{j_2}}
    \abr{\diffker^{n}_{j_2},\diffker^{2n}_{j_3}}
    \abr{\diffker^{n}_{j_1},\diffker^{2n}_{j_3}} 
    \nn\\
    &\quad+3\times
    2^2\times2^{2(2H-1)}\times n^{(6H-\frac32)}
    \sum_{j\in[n-1]\times[2n-1]} 
    \rbr{\iv}^{-3} 
    \brbr{V^{[1]}_{t^{n}_{j_1}}}^2 \brbr{V^{[1]}_{t^{2n}_{j_2}}}^2 \brbr{V^{[1]}_{t^{2n}_{j_3}}}^2 
    \abr{\diffker^{n}_{j_1},\diffker^{2n}_{j_2}}
    \abr{\diffker^{2n}_{j_2},\diffker^{2n}_{j_3}}
    \abr{\diffker^{n}_{j_1},\diffker^{2n}_{j_3}} 
    \nn\\
    &\quad-
    2^2\times2^{3(2H-1)}\times n^{(6H-\frac32)}
    \sum_{j\in[2n-1]^3} 
    \rbr{\iv}^{-3} 
    \brbr{V^{[1]}_{t^{2n}_{j_1}}}^2 \brbr{V^{[1]}_{t^{2n}_{j_2}}}^2 \brbr{V^{[1]}_{t^{2n}_{j_3}}}^2 
    \abr{\diffker^{2n}_{j_1},\diffker^{2n}_{j_2}}
    \abr{\diffker^{2n}_{j_2},\diffker^{2n}_{j_3}}
    \abr{\diffker^{2n}_{j_1},\diffker^{2n}_{j_3}} 
    \nn\\&\quad+O_M(n^{-H}),
    \label{230627.1337}
  \end{align}
  and its order is estimated as $D_{u_n}D_{u_n} M_n=O_M(n^{-\half})$.
  \item[(ii)]
  $(D_{u_n})^3 M_n=O_M(n^{-H})$.
\end{lemma}
\begin{proof}
  By Lemma \ref{230724.1710},
  we have the following decomposition of $D_{u_n}D_{u_n} M_n$:
  \begin{align}
    D_{u_n}D_{u_n} M_n&=
    D_{u_n}D_{u_n} \bar M_n - D_{u_n}D_{u_n} \negTerm_n^{(3)}
    =
    D_{u_n}D_{u_n} \bar M_n + O_M(n^{-H})
    \nn\\&=
    \sum_{i\in[2]^3} (-1)^{(i_1+i_2+i_3-3)}\;
    D_{u_n^{(i_3)}}D_{u_n^{(i_2)}} (\bar M^{(i_1)}_n)
    +O_M(n^{-H})
    \nn% \label{230725.1118}
    \\&=
    \sum_{i\in[2]^3} (-1)^{(i_1+i_2+i_3-3)}\;
    2\times
    (i_1i_2i_3)^{(2H-1)}%2^{(2H-1)(i_1+i_2+i_3-3)}
    n^{(6H-\frac32)}\cali_n^{(3;i)(1)}
    +O_M(n^{-H}).
    \label{230725.1119}
    % \\&=
    % 2^2\times n^{(6H-\frac32)}
    % \sum_{j\in[n-1]^3}%\bbJ_n([3],i)} 
    % \rbr{\iv}^{-3} 
    % \brbr{V^{[1]}_{t^{n}_{j_1}}}^2 \brbr{V^{[1]}_{t^{n}_{j_2}}}^2 \brbr{V^{[1]}_{t^{n}_{j_3}}}^2 
    % \abr{\diffker^{n}_{j_1},\diffker^{n}_{j_2}}
    % \abr{\diffker^{n}_{j_2},\diffker^{n}_{j_3}}
    % \abr{\diffker^{n}_{j_1},\diffker^{n}_{j_3}} 
    % \nn\\
    % &\quad-3\times
    % 2^2\times2^{(2H-1)}\times n^{(6H-\frac32)}
    % \sum_{j\in[n-1]^2\times[2n-1]} 
    % \rbr{\iv}^{-3} 
    % \brbr{V^{[1]}_{t^{n}_{j_1}}}^2 \brbr{V^{[1]}_{t^{n}_{j_2}}}^2 \brbr{V^{[1]}_{t^{2n}_{j_3}}}^2 
    % \abr{\diffker^{n}_{j_1},\diffker^{n}_{j_2}}
    % \abr{\diffker^{n}_{j_2},\diffker^{2n}_{j_3}}
    % \abr{\diffker^{n}_{j_1},\diffker^{2n}_{j_3}} 
    % \nn\\
    % &\quad+3\times
    % 2^2\times2^{2(2H-1)}\times n^{(6H-\frac32)}
    % \sum_{j\in[n-1]\times[2n-1]} 
    % \rbr{\iv}^{-3} 
    % \brbr{V^{[1]}_{t^{n}_{j_1}}}^2 \brbr{V^{[1]}_{t^{2n}_{j_2}}}^2 \brbr{V^{[1]}_{t^{2n}_{j_3}}}^2 
    % \abr{\diffker^{n}_{j_1},\diffker^{2n}_{j_2}}
    % \abr{\diffker^{2n}_{j_2},\diffker^{2n}_{j_3}}
    % \abr{\diffker^{n}_{j_1},\diffker^{2n}_{j_3}} 
    % \nn\\
    % &\quad-
    % 2^2\times2^{3(2H-1)}\times n^{(6H-\frac32)}
    % \sum_{j\in[2n-1]^3} 
    % \rbr{\iv}^{-3} 
    % \brbr{V^{[1]}_{t^{2n}_{j_1}}}^2 \brbr{V^{[1]}_{t^{2n}_{j_2}}}^2 \brbr{V^{[1]}_{t^{2n}_{j_3}}}^2 
    % \abr{\diffker^{2n}_{j_1},\diffker^{2n}_{j_2}}
    % \abr{\diffker^{2n}_{j_2},\diffker^{2n}_{j_3}}
    % \abr{\diffker^{2n}_{j_1},\diffker^{2n}_{j_3}} 
    % \nn\\&\quad+O_M(n^{-H})
  \end{align}
  \item [(i)] 
  This is obvious from \eqref{230725.1119} and Lemma \ref{230724.1710}.
  \item [(ii)] 
  By \eqref{230725.1119} and $u_n=O_M(1)$ from Lemma \ref{230721.1101}, 
  we only need to ensure 
  $D_{u_n}(n^{(6H-\frac32)}\cali_n^{(3;i)(1)})=O_M(n^{-H})$.
  The functional $\cali_n^{(3;i)(1)}$ corresponds to the weighted graph 
  \eqref{fig:230721.1723} for any $i\in[2]^3$,
  and the weighted graph \eqref{fig:230721.1723} has the only component
  $C$ such that $\barq(C)=0$.
  Hence the case (a) of Proposition \ref{221223.1757} applies to $\cali_n^{(3;i)(1)}$,
  and by Proposition \ref{230725.1133}, we have 
  $D_{u_n^{(i')}}\cali_n^{(3;i)(1)}=O_M(n^{(1-6H)-\half})=O_M(n^{\half-6H})$ for $i'=1,2$.
  Thus, we obtain 
  $D_{u_n}(n^{(6H-\frac32)}\cali_n^{(3;i)(1)})
  =O_M(n^{-1})$.
\end{proof}
\begin{remark}
  Although we only prove 
  $(D_{u_n})^3 M_n=O_M(n^{-H})$ here,
  since this estimate is sufficient for the argument 
  to prove the asymptotic expansion,
  we can prove 
  $D_{u_n}^3 M_n=O_M(n^{-1})$
  by examining the functionals represented by $O_M(n^{-H})$ 
  in \eqref{230627.1337}
  and using Proposition \ref{221223.1757}.
\end{remark}
  
\subsubsection*{(ii) Decomposition and estimates about
$N_n$ and $D_{u_n}N_n$}
Recall we have defined $N_n$ by
\begin{align}
  % r_nN_n &= 
  % \sqrt{n}\rbr{\iv}^{-1} \rbr{\convDNeg_n - \convDNeg_{2n}}
  % + \underbrace{\Gamma_n^{(1)}+\Gamma_n^{(2)}+\Gamma_n^{(3)}}_{=O_M(n^{(-H)\vee(-\frac32+H)})}
  % \\
  N_n &= 
  n\rbr{\iv}^{-1} \rbr{\convDNeg_n - \convDNeg_{2n}}
  + O_M(n^{(\half-H)\vee(-1+H)})
  \label{230725.1225}
\end{align}
and $n\rbr{\iv}^{-1} \rbr{\convDNeg_n - \convDNeg_{2n}}$ decomposes as
\begin{align}
  n\rbr{\iv}^{-1} \rbr{\convDNeg_n - \convDNeg_{2n}}&=
  n\rbr{\iv}^{-1} \rbr{
    \convDiff^{(2,1)}_n+\convDiff^{(2,2)}_n + \convDiff^{(1,1)}_n + \convDiff^{(1,2)}_n 
    - \rbr{\convDiff^{(2,1)}_{2n}+\convDiff^{(2,2)}_{2n} + 
    \convDiff^{(1,1)}_{2n} + \convDiff^{(1,2)}_{2n}}}.
    \label{230725.2120}
\end{align}
The summands are written as
\begin{align}
  \rbr{\iv}^{-1} \convDiff^{(2,1)}_n &=
  2\;\rbr{\iv}^{-1} 
  n^{2H-1} \sum_{j=1}^{n-1} V^{[1]}_\tj V^{[(1;1),1]}_\tjm
  I_3\rbr{\diffker^n_j \otimes (\bbone^n_j)^{\otimes2}} 
  =:n^{2H-1}\caliN{1}_n
  \nn\\
  \rbr{\iv}^{-1} \convDiff^{(2,2)}_n &=
  n^{2H-1} \rbr{\iv}^{-1}
    \sum_{j=1}^{n-1} V^{[1]}_\tj {V^{[(1;1),1]}_\tj} 
    I_3\brbr{\diffker^n_j\otimes{(\bbone^n_\jp)}^{\otimes2}}
  \nn\\&\qquad
    -
  n^{2H-1} \rbr{\iv}^{-1} 
    \sum_{j=1}^{n-1} V^{[1]}_\tj {V^{[(1;1),1]}_\tj} 
    I_3\brbr{\diffker^n_j\otimes{(\bbone^n_j)}^{\otimes2}}
  \nn\\&=:
  n^{2H-1}\caliN{2}_n + n^{2H-1}\caliN{3}_n
  \nn\\
  \rbr{\iv}^{-1} \convDiff^{(1,1)}_n &=
  c_{2,H}\rbr{\iv}^{-1} n^{-1}
  \sum_{j=1}^{n-1} a'_\tj V^{[1]}_\tj B(\bbh^n_j)
  =:
  n^{-1}\caliN{4}_n
  \nn\\
  \rbr{\iv}^{-1} \convDiff^{(1,2)}_n &=
  -\frac1{2} \rbr{\iv}^{-1}   c_{2,H} n^{-1}
  \Brbr{\brbr{V^{[1]}_0 }^2 + \brbr{V^{[1]}_1 }^2} 
  =:n^{-1}\AN
  \label{230726.1845}
\end{align}
with
$\bbh^{(1) n}_{j}(\cdot) = \bbone^{2n}_{2j}(\cdot)\; n(\cdot-\ttjm)$,\;
$\bbh^{(2) n}_{j}(\cdot) = \bbone^{2n}_{2j-1}(\cdot)\; n(\ttjm-\cdot)$ and
$\bbh^{n}_{j} = \bbh^{(1) n}_{j} - \bbh^{(2) n}_{j+1}$.
The functionals
$\caliN{1}_n$, $\caliN{2}_n$ and $\caliN{3}_n$ 
correspond to the weighted graph \eqref{fig:230721.1712},
and $\caliN{4}_n$
corresponds to the weighted graph \eqref{fig:230721.1713}.
The exponents of the weighted graphs \eqref{fig:230721.1712} and \eqref{fig:230721.1713}
are $1+(-1-2H)=-2H$ and $1+(-1)=0$, respectively.

\begin{lemma}\label{230721.1242}
  \item[(i)] $N_n=O_M(1)$
  \item[(ii)] $D_{u_n}N_n = O_M(n^{(\half-H)\vee(-1+H)})$
  \item[(iii)] $D_{u_n}^2N_n = O_M(n^{(\half-H)\vee(-1+H)})$
\end{lemma}
\begin{proof}
  (i)
  From the above arguments, we have
\begin{align}
  n\rbr{\iv}^{-1} \rbr{\convDNeg_n - \convDNeg_{2n}}&=
  % n\rbr{\iv}^{-1} \rbr{
  %   \convDiff^{(2,1)}_n+\convDiff^{(2,2)}_n + \convDiff^{(1,1)}_n + \convDiff^{(1,2)}_n 
  %   - \rbr{\convDiff^{(2,1)}_{2n}+\convDiff^{(2,2)}_{2n} + 
  %   \convDiff^{(1,1)}_{2n} + \convDiff^{(1,2)}_{2n}}}
  % \\&=
  n^{2H}(\caliN{1}_n + \caliN{2}_n + \caliN{3}_n)
  +\caliN{4}_n + \AN
  \nn\\&\qquad-\rbr{
  2^{2H-1}n^{2H}(\caliN{1}_{2n} + \caliN{2}_{2n} + \caliN{3}_{2n})
  +2^{-1}\caliN{4}_{2n} + 2^{-1}\AN}
  \nn\\&=
  \sum_{i=1}^2 (-1)^{(i-1)} \rbr{
    i^{(2H-1)}%2^{(2H-1)(i-1)} 
    n^{2H}\sum_{k=1}^3\caliN{k}_{in} 
    + i^{-1}%2^{-(i-1)}
    \caliN{4}_{in}}
  + 2^{-1}\AN
  \label{230725.1221}\\&=O_M(n^0)
  =O_M(1)
  \nn
\end{align}
and $N_n=O_M(1)$ by \eqref{230725.1225}.

\item[(ii)] The functional $D_{u_n}N_n$ is written as
\begin{align*}
  D_{u_n}N_n &= 
  % D_{u_n}\rbr{n\rbr{\iv}^{-1} \rbr{\convDNeg_n - \convDNeg_{2n}}}
  % + O_M(n^{(\half-H)\vee(-1+H)})
  % \\&=
  % \sum_{i=1,2}(-1)^{i-1}D_{u_n^{(i)}} \rbr{\sum_{i=1}^2 (-1)^{(i-1)} \rbr{
  %   i^{(2H-1)}%2^{(2H-1)(i-1)} 
  %   n^{2H}\sum_{k=1}^3\caliN{k}_{in} 
  %   + i^{-1}%2^{-(i-1)}
  %   \caliN{4}_{in}}}
  % \\&\qquad+
  % 2^{-1}\sum_{i=1,2}(-1)^{i-1}D_{u_n^{(i)}}\AN
  % \\&\qquad+ O_M(n^{(\half-H)\vee(-1+H)})
  % \\&=
  \sum_{i_1,i_2=1}^2 (-1)^{i_1+i_2}  
  \rbr{
    i_1^{(2H-1)}%2^{(2H-1)(i_1-1)} 
  n^{2H}\sum_{k=1}^3 D_{u_n^{(i_2)}}\caliN{k}_{i_1n} 
  + i_1^{-1}%2^{-(i_1-1)}
  D_{u_n^{(i_2)}}\caliN{4}_{i_1n}}
  +%\\&\qquad+
  2^{-1}\sum_{i=1,2}(-1)^{i-1}D_{u_n^{(i)}}\AN
  \\&\qquad+ O_M(n^{(\half-H)\vee(-1+H)}),
\end{align*}
where 
\begin{align}
  D_{u_n^{(i)}}\AN&=
  \abr{D\AN,
  i^{(2H-1)}%2^{(2H-1)(i-1)}
  n^{2H-\half} 
  \sum_{j\in[in-1]} \rbr{\iv}^{-1} \brbr{V^{[1]}_{t^{in}_j}}^2 
  I_1\rbr{\diffker^{in}_j} \diffker^{in}_j}
  \nn\\&=
  i^{(2H-1)}%2^{(2H-1)(i-1)}
  n^{-\half} 
  \sum_{j\in[in-1]} \rbr{\iv}^{-1} \brbr{V^{[1]}_{t^{in}_j}}^2 
  n^{2H} \abr{D\AN,\diffker^{in}_j}
  I_1\rbr{\diffker^{in}_j}.
  \label{eq:231005.1843}
\end{align}
The functionals
$\caliN{1}_{i_1n}$, $\caliN{2}_{i_1n}$ and $\caliN{3}_{i_1n}$ 
correspond to the weighted graph \eqref{fig:230721.1712},
which has the only component $C$ such that $\barq(C)=3$.
Hence the case (a) of Proposition \ref{221223.1757} applies, 
and by Proposition \ref{230725.1133}, we have 
$D_{u_n^{(i_2)}}\caliN{k}_{i_1n}
=O_M(n^{-2H-\half})$ for $i_1,i_2\in\cbr{1,2}$ and $k\in[3]$.
By a similar argument, we also have 
$D_{u_n^{(i_2)}}\caliN{4}_{i_1n}
=O_M(n^{-\half})$ for $i_1,i_2\in\cbr{1,2}$.
As for $D_{u_n^{(i)}}\AN$, 
the rescaled functional $n^{\half}D_{u_n^{(i)}}\AN$ 
corresponds to the weighted graph \eqref{fig:230721.1714},
whose exponent is $(1+(-\half-H))=\half-H$.
Thus $D_{u_n^{(i)}}\AN=O_M(n^{-H})$.
Hence, 
\begin{align*}
  D_{u_n}N_n &= 
  n^{2H}\times O_M(n^{-2H-\half})+
  O_M(n^{-\half})+
  O_M(n^{-H})+
  O_M(n^{(\half-H)\vee(-1+H)})
  =
  O_M(n^{(\half-H)\vee(-1+H)})
\end{align*}
\item[(iii)]
This follows immediately 
from (ii) and $u_n=O_M(1)$ (Lemma \ref{230721.1101}). 
\end{proof}

\subsubsection{Condition {\bf [D]}(iii) (Random symbols)}
\label{sec:230726.1315}
Consider the random symbols
$\mS_{0,n}^{(2,0)},\mS_n^{(3,0)},\mS_n^{(1,0)}$ and $\mS_{1,n}^{(2,0)}$
for $u_n=u_n^{(1)}-u_n^{(2)}$, $N_n$ and $G_\infty$ defined at 
\eqref{eq:230926.1446}, \eqref{230725.1532} and \eqref{230725.1532}, respectively.
For every
$\mT_n=\mS_{0,n}^{(2,0)},\mS_n^{(3,0)},\mS_n^{(1,0)}$ and $\mS_{1,n}^{(2,0)}$,
we identify the limit random symbol 
$\mT=\mS_0^{(2,0)},\mS^{(3,0)},\mS^{(1,0)}$ and $\mS_1^{(2,0)}$,
which appear in the asymptotic expansion formula, 
and verify the conditions in {\bf [D]} (iii) for each pair $(\mT_n,\mT)$,
that is,
we will prove 
\begin{align}
  E \sbr{\Psi (\sfz)\mT_n(\tti\sfz)} = E[\Psi (\sfz)\bar\mT_n(\tti\sfz)]
  \tand
  \bar\mT_n\to\mT\text{ in }L^p
  \label{eq:230927.1443}
\end{align}
with some auxiliary random symbols $\bar\mT_n$.

\subsubsection*{(i) The random symbol $\mS^{(2,0)}_{0,n}=\half\qtan_n$ and its limit}
The functional $\abr{DM_n, u_n}$ decomposes as
\begin{align*}
  \abr{DM_n, u_n} &=
  % \abr{D\bar M_n, u_n} - 
  % \underbrace{\babr{D\negTerm_n^{(3)}, u_n}}_{O_M(n^{-H})}
  % \\&=
  % \sum_{i_1,i_2=1}^2(-1)^{i_1+i_2}
  % \brbr{2\,\cali_n^{(i_1,i_2)(1)}+2\,\cali_n^{(i_1,i_2)(2)}+\calj_n^{(i_1;i_2)}}
  % +O_M(n^{-H})
  % \nn\\&=
  \sum_{i_1,i_2=1}^2(-1)^{i_1+i_2}
  2\times\brbr{\cali_n^{(i_1,i_2)(1)}+\cali_n^{(i_1,i_2)(2)}}
  +O_M(n^{-H})
  \nn\\&=
  \sum_{i_1,i_2=1}^2 (-1)^{i_1+i_2}
  2\times\cali_n^{(i_1,i_2)(2)}
  +G_\infty
  +O_M(n^{(-H)\vee(-\frac34)}),
\end{align*}
where we used 
$\cali_n^{(2,2)(1)}=\half\cali_{2n}^{(1,1)(1)}$ and Lemma \ref{230724.1350}.
Notice that 
$\cali_n^{(i_1,i_2)(2)}=O_M(n^{-\half})$ (Lemma \ref{230721.1102}).
The random symbol $\qtan_n$ is written as
\begin{align*}
  \qtan_n [\sfi\sfz]^2&=
  r_n^{-1}(\abr{DM_n, u_n} - G_\infty) [\sfi\sfz]^2
  \\&=
  2\;n^{\half}\cali_n^{(1,1)(2)}[\sfi\sfz]^2
  - 4\;n^{\half}\cali_n^{(1,2)(2)}[\sfi\sfz]^2
  +2\; n^{\half}\cali_{n}^{(2,2)(2)}[\sfi\sfz]^2
  %+ n^{\half}\cali_{2n}^{(1,1)(2)}
  +O_M(n^{(\half-H)\vee(-\frac14)}).
\end{align*}

We want to know the weak limit of $\qtan_n$ in the sense of \eqref{eq:230927.1443}.
By the IBP formula, we decompose
the expectation of the product of 
$\Psi(\sfz)=\exp\rbr{2^{-1} G_\infty[\tti\sfz]^2}$ and 
the random symbol $n^{\half}\cali_n^{(1,1)(2)} [\sfi\sfz]^2$ 
as follows:
\begin{align}
  &E\sbr{\Psi(\sfz) n^{\half}\cali_n^{(1,1)(2)}[\sfi\sfz]^2}%-(4H-1)=-4H+1
  % \\&=
  % n^{4H-\half}
  % E\sbr{
  % \Psi(\sfz)
  % \rbr{\iv}^{-2}
  % \sum_{j\in[n-1]^2}
  % \rbr{V^{[1]}_{t_\jon} }^2  \rbr{V^{[1]}_{t_\jtw} }^2 
  % I_2\rbr{\diffker^n_\jon \otimes \diffker^n_\jtw}
  % \abr{\diffker^n_\jon, \diffker^n_\jtw}
  % }
  % \\&=
  % n^{4H-\half}
  % \sum_{j\in[n-1]^2}
  % E\sbr{
  %   \Psi(\sfz)
  % \rbr{\iv}^{-2}
  % \rbr{V^{[1]}_{t_\jon} }^2  \rbr{V^{[1]}_{t_\jtw} }^2 
  % I_2\rbr{\diffker^n_\jon \otimes \diffker^n_\jtw}
  % \abr{\diffker^n_\jon, \diffker^n_\jtw}
  % }
  % \\&=
  % n^{4H-\half}
  % \sum_{j\in[n-1]^2}
  % E\sbr{
  %   \abr{D\rbr{\Psi(\sfz)\rbr{\iv}^{-2}
  %   \rbr{V^{[1]}_{t_\jon} }^2  \rbr{V^{[1]}_{t_\jtw} }^2},
  %   \diffker^n_\jtw}
  % I_1\rbr{\diffker^n_\jon}
  % \abr{\diffker^n_\jon, \diffker^n_\jtw}
  % }
  % \\&=
  \nn\\&=
  n^{2H-\half}
  \sum_{j\in[n-1]^2}
  E\sbr{
    n^{2H}
    \abr{D\rbr{\Psi(\sfz)\rbr{\iv}^{-2}
    \rbr{V^{[1]}_{t_\jon} }^2  \rbr{V^{[1]}_{t_\jtw} }^2},
    \diffker^n_\jtw}
  I_1\rbr{\diffker^n_\jon}
  \abr{\diffker^n_\jon, \diffker^n_\jtw}
  }[\sfi\sfz]^2.
  % \\&=
  % n^{2H-\half}
  % \sum_{j\in[n-1]^2}
  % E\sbr{
  %   n^{2H}\Psi(\sfz)
  %   \abr{D\rbr{\rbr{\iv}^{-2}
  %   \rbr{V^{[1]}_{t_\jon} }^2  \rbr{V^{[1]}_{t_\jtw} }^2},
  %   \diffker^n_\jtw}
  % I_1\rbr{\diffker^n_\jon}
  % \abr{\diffker^n_\jon, \diffker^n_\jtw}
  % }[\sfi\sfz]^2
  % \\&\hspace{10pt}+
  % n^{2H-\half}
  % \sum_{j\in[n-1]^2}
  % E\sbr{
  %   n^{2H}
  %   \abr{D\Psi(\sfz),
  %   \diffker^n_\jtw}
  %   \rbr{\iv}^{-2}
  %   \rbr{V^{[1]}_{t_\jon} }^2  \rbr{V^{[1]}_{t_\jtw} }^2
  % I_1\rbr{\diffker^n_\jon}
  % \abr{\diffker^n_\jon, \diffker^n_\jtw}
  % }[\sfi\sfz]^2
  \nn\\&=
  n^{2H-\half}
  \sum_{j\in[n-1]^2}
  E\sbr{
    \Psi(\sfz)
    n^{2H}
    \abr{D\rbr{\rbr{\iv}^{-2}
    \rbr{V^{[1]}_{t_\jon} }^2  \rbr{V^{[1]}_{t_\jtw} }^2},
    \diffker^n_\jtw}
  I_1\rbr{\diffker^n_\jon}
  \abr{\diffker^n_\jon, \diffker^n_\jtw}
  }[\sfi\sfz]^2
  \nn\\&\hspace{10pt}+
  \half
  n^{2H-\half}
  \sum_{j\in[n-1]^2}
  E\sbr{
    \Psi(\sfz)
    n^{2H}
    \abr{ DG_\infty,
    \diffker^n_\jtw}
    \rbr{\iv}^{-2}
    \rbr{V^{[1]}_{t_\jon} }^2  \rbr{V^{[1]}_{t_\jtw} }^2
  I_1\rbr{\diffker^n_\jon}
  \abr{\diffker^n_\jon, \diffker^n_\jtw}
  }[\sfi\sfz]^4.
  \nn\\&=
  E\sbr{\Psi(\sfz)\Brbr{\cali_n^{(\qtan)(1,1)(1)}[\sfi\sfz]^2+
  \cali_n^{(\qtan)(1,1)(2)}[\sfi\sfz]^4}},
  \label{eq:231005.1838}
\end{align}
% {\myred
% $\Psi(\sfz)=\exp(-\half G_\infty\sfz^2)
% =\exp(\half G_\infty[\sfi\sfz]^2)$, 
% $D\Psi(\sfz)
% =\Psi(\sfz)\half DG_\infty[\sfi\sfz]^2$}
where we define
\begin{align*}
  \cali_n^{(\qtan)(1,1)(1)}&=
  n^{2H-\half}
  \sum_{j\in[n-1]^2}
    n^{2H}
    \abr{D\rbr{\rbr{\iv}^{-2}
    \rbr{V^{[1]}_{t_\jon} }^2  \rbr{V^{[1]}_{t_\jtw} }^2},
    \diffker^n_\jtw}
  I_1\rbr{\diffker^n_\jon}
  \abr{\diffker^n_\jon, \diffker^n_\jtw}
  \\
  \cali_n^{(\qtan)(1,1)(2)}&=
  \half
  n^{2H-\half}
  \sum_{j\in[n-1]^2}
    n^{2H}
    \abr{ DG_\infty,
    \diffker^n_\jtw}
    \rbr{\iv}^{-2}
    \rbr{V^{[1]}_{t_\jon} }^2  \rbr{V^{[1]}_{t_\jtw} }^2
  I_1\rbr{\diffker^n_\jon}
  \abr{\diffker^n_\jon, \diffker^n_\jtw}.
\end{align*}
Both the rescaled functionals
$n^{\half-2H}\cali_n^{(\qtan)(1,1)(k)}$ ($k=1,2$)
correspond to the weighted graph \eqref{fig:230721.1724},
whose exponent is $\half-3H$.
Hence $\cali_n^{(\qtan)(1,1)(k)}=O_M(n^{-H})$, and 
the weak limit of 
the random symbol $n^{\half}\cali_n^{(1,1)(2)} [\sfi\sfz]^2$ 
in the sense of \eqref{eq:230927.1443}
is zero.
Similar arguments work for $n^{\half} \cali_n^{(1,2)(2)}$ and 
$n^{\half} \cali_n^{(2,2)(2)}$.
Setting $\mS^{(2,0)}_0(\tti\sfz)=0$, we have the following lemma:
\begin{lemma}\label{230725.1941}
  The conditions for 
  % $(\mS^{(3,0)}_{n},\mS^{(3,0)})$ 
  $(\mS^{(2,0)}_{0,n}, \mS^{(2,0)}_0)$
  in {\bf [D]}(iii) hold.
\end{lemma}

\subsubsection*{(ii) The random symbol $\mS^{(3,0)}_{n}=\frac13\qtor_n$ and its limit}
We define the random symbol $\mS^{(3,0)}$ by
\begin{align*}
  \mS^{(3,0)}(\tti\sfz)=
  c_\qtor \;\rbr{\iv}^{-3} \int_0^1 (a_t)^3 dt\; [\sfi\sfz]^3,
\end{align*}
where $a_t=(V^{[1]}(X_t))^2$ and
the constant $c_\qtor$ is defined at \eqref{230726.1339}.
The following lemma confirm that the conditions in {\bf [D]}(iii) 
for the pair $(\mS^{(3,0)}_{n},\mS^{(3,0)})$ hold.
\begin{lemma}\label{230725.1942}
  For some $p>1$, 
  the coefficients of the random symbol $\mS^{(3,0)}$ are in $L^p$, and 
  there exists a polynomial random symbol $\bar\mS^{(3,0)}_{n}$ such that
  the coefficients of $\bar\mS^{(3,0)}_{n}$ belongs to $L^p$,
  the equation
  \begin{align*}
    E\bsbr{\Psi (\sfz)\mS^{(3,0)}_{n}(\tti\sfz)} = 
    E\bsbr{\Psi (\sfz)\bar\mS^{(3,0)}_{n}(\tti\sfz)}
  \end{align*}
  holds for $\sfz\in\bbR$, and
  the convergence $\bar\mS^{(3,0)}_{n}\to\mS^{(3,0)}$ in $L^p$ holds.
\end{lemma}

\begin{proof}
  By the definition of $\mS^{(3,0)}_{n}$ and 
  Lemma \ref{lem:230925.1805}, 
  we have 
  \begin{align*}
    \mS^{(3,0)}_{n}(\sfi\sfz)&=
    \frac13\qtor_n[\sfi\sfz]^3=
    \frac13n^{1/2}D_{u_n}D_{u_n}M_n[\sfi\sfz]^3
    \\&= 
    c_\qtor \; \rbr{\iv}^{-3} \int_0^1 (a_t)^3 dt\;
    [\sfi\sfz]^3
    +o_M(1)=
    \mS^{(3,0)}(\tti\sfz)
    +o_M(1).
  \end{align*}
  Setting 
  $\bar\mS^{(3,0)}_{n}=\mS^{(3,0)}_{n}$,
  the assertion of this lemma stands.
\end{proof}

% \begin{lemma}\label{230906.2127}
%   {\myblue Let $u_n$ and $M_n$ be as \eqref{??} and \eqref{??}.}
%   Then the following decomposition of $n^{1/2}D_{u_n}D_{u_n} (M_n)$ holds:
%   \begin{align*}
%     &n^{1/2}D_{u_n}D_{u_n}M_n= 
%     3\,c_\qtor \;%\times
%     \rbr{\iv}^{-3}
%     \int_0^1 (a_t)^3 dt
%     +o_M(1),
%   \end{align*}
%   where the constant $c_\qtor$ is defined by 
%   \begin{align}
%     c_\qtor&=
%     {\sum_{i_1,i_2\in\bbZ}
%     \widehat\rho(i_1)
%     \widehat\rho(i_1-i_2)
%     \widehat\rho(i_2)}
%     -
%     2^{(2H+1)}
%     {\sum_{i_1,i_2\in\bbZ}
%     \tilde\rho(i_2)
%     \tilde\rho(2i_1-i_2)
%     \widehat\rho(i_1)}
%     +%\\&\hspace{120pt}+
%     2^{2H}
%     {\sum_{i_1,i_2\in\bbZ}
%     \tilde\rho(i_1)
%     \tilde\rho(i_2)
%     \widehat\rho(i_1-i_2)}
%     %\label{230726.1339}
%   \end{align}
%   and the definitions of $\widehat\rho$ and $\tilde\rho$
%   are given at \eqref{eq:230925.1740} and \eqref{eq:230925.1741}.
% \end{lemma}

\subsubsection*{(iii) The random symbols
$\mS^{(1,0)}_{n}$ and $\mS^{(2,0)}_{1,n}$ and their limits}

Recall that we have defined $N_n$ at \eqref{230725.1225}, which decomposes as \eqref{230725.2120}.
We can show that the contribution from 
$n\rbr{\iv}^{-1}\convDiff^{(2,1)}_{in}$ and
$n\rbr{\iv}^{-1}\convDiff^{(2,2)}_{in}$ ($i=1,2$)
to the weak limit of the random symbol 
$\mS^{(1,0)}_{n}(\tti\sfz)=N_n[\sfi\sfz]$ 
is negligible by an argument similar to the one for $\mS^{(2,0)}_{0,n}=\half\qtan_n$.

As for $n\rbr{\iv}^{-1} \convDiff^{(1,1)}_{in}$, we have
\begin{align*}
  E\sbr{\Psi(\sfz) n\rbr{\iv}^{-1} \convDiff^{(1,1)}_n[\sfi\sfz]}
  %\\
  &=
  % E\sbr{\Psi(\sfz) c_{2,H}\rbr{\iv}^{-1}\sum_{j=1}^{n-1} a'_\tj V^{[1]}_\tj B(\bbh^n_j)}[\sfi\sfz]
  % \\&=
  % c_{2,H}\sum_{j=1}^{n-1}
  % E\sbr{\Psi(\sfz) \rbr{\iv}^{-1} a'_\tj V^{[1]}_\tj B(\bbh^n_j)}[\sfi\sfz]
  % \\&=
  % c_{2,H}\sum_{j=1}^{n-1}
  % E\sbr{\abr{D\brbr{\Psi(\sfz) \rbr{\iv}^{-1} a'_\tj V^{[1]}_\tj}, \bbh^n_j}}[\sfi\sfz]
  % \\&=
  % c_{2,H}\sum_{j=1}^{n-1}
  % E\sbr{\abr{D\Psi(\sfz), \bbh^n_j}\rbr{\iv}^{-1} a'_\tj V^{[1]}_\tj}[\sfi\sfz]
  % \\&\quad+
  % c_{2,H}\sum_{j=1}^{n-1}
  % E\sbr{\Psi(\sfz)\babr{D\brbr{ \rbr{\iv}^{-1} a'_\tj V^{[1]}_\tj}, \bbh^n_j}}[\sfi\sfz]
  % \\&=
  \half c_{2,H}\sum_{j=1}^{n-1}
  E\sbr{\Psi(\sfz)\abr{DG_\infty, \bbh^n_j}\rbr{\iv}^{-1} a'_\tj V^{[1]}_\tj}[\sfi\sfz]^3
  \\&\quad+
  c_{2,H}\sum_{j=1}^{n-1}
  E\sbr{\Psi(\sfz)\babr{D\brbr{ \rbr{\iv}^{-1} a'_\tj V^{[1]}_\tj}, \bbh^n_j}}[\sfi\sfz].
\end{align*}
By a similar argument of Lemma \ref{lemma:230616.1754} (i), we can show 
\begin{align*}
  \sup_{n\in\ntwo}\sup_{j\in[n-1]}\norm{n^{2H}\abr{DG_\infty, \bbh^n_j}}_p<\infty, \tand
  \sup_{n\in\ntwo}\sup_{j\in[n-1]}
  \norm{n^{2H}\babr{D\brbr{ \rbr{\iv}^{-1} a'_\tj V^{[1]}_\tj}, \bbh^n_j}}_p<\infty
\end{align*}
for any $p>1$.
Hence the order of the functionals
\begin{align*}
  \half c_{2,H}\sum_{j=1}^{n-1}
  \abr{DG_\infty, \bbh^n_j}\rbr{\iv}^{-1} a'_\tj V^{[1]}_\tj, 
  \tand
  c_{2,H}\sum_{j=1}^{n-1}
  \babr{D\brbr{ \rbr{\iv}^{-1} a'_\tj V^{[1]}_\tj}, \bbh^n_j}
\end{align*}
is of $O_{L^p}(n^{1-2H})$ for $p>1$.
Hence $n\rbr{\iv}^{-1} \convDiff^{(1,1)}_{in}$ has no contribution to the limit.

The contribution from $n\rbr{\iv}^{-1} \rbr{\convDiff^{(1,2)}_n - \convDiff^{(1,2)}_{2n}}$ is obviously 
$\half\AN=
-\frac1{4} \rbr{\iv}^{-1}   c_{2,H} \Brbr{\brbr{V^{[1]}_0 }^2 + \brbr{V^{[1]}_1 }^2}$.
Therefore, we set 
$\mS^{(1,0)}(\tti\sfz)=\half\AN[\sfi\sfz]$, 
and we have the following lemma.
\begin{lemma}\label{230725.2158}
  Condition {\bf [D]}(iii) (b) holds for $(\mS^{(1,0)}_{n},\mS^{(1,0)})$.
\end{lemma}

For $\mS^{(2,0)}_{1,n}(\tti\sfz)=D_{u_n}N_n[\sfi\sfz]^2$,
we have 
$D_{u_n}N_n = O_{L^p}(n^{(\half-H)\vee(-1+H)})= o_{L^p}(1)$ for $p>1$
by Lemma \ref{230721.1242} (ii).
Hence, setting
$\mS^{(2,0)}_{1}(\tti\sfz)=0$, we have the following lemma.
\begin{lemma}\label{230725.2159}
  Condition {\bf [D]}(iii) (b) holds for $(\mS^{(2,0)}_{1,n},\mS^{(2,0)}_{1})$.
\end{lemma}

\subsubsection{Condition {\bf [D]}(iv) (Nondegeneracy of $M_n$)}
  We want to verify the condition [D](iv)(b).
  The following estimate gives the nondegeneracy of $M_n$.
  \begin{lemma}\label{230725.1635}
    For any $p\geq1$, the probability $P\sbr{\Delta_{M_n}<G_\infty}$ is bounded as
    $P\sbr{\Delta_{M_n}<G_\infty}=O(n^{-\frac{p}{2}})$.
  \end{lemma}
  \begin{proof}
  The (one-dimensional) Malliavin matrix 
  $\Delta_{M_n}=\abr{DM_n, DM_n}_\calh$ decomposes as
  \begin{align*}
    \Delta_{M_n}&=%\abr{DM_n, DM_n}_\calh= 
    \abr{D\bar M_n - D\negTerm_n^{(3)}, D\bar M_n - D\negTerm_n^{(3)}}
    =%\\&=
    \abr{D\bar M_n, D\bar M_n} +O_M(n^{-H})
  \end{align*}
  Since we have $D\bar M_n = 2 u_n + v_n$,
  \begin{align*}
    \abr{D\bar M_n, D\bar M_n}&= 
    \abr{2 u_n + v_n, 2 u_n + v_n} =
    4\abr{u_n, u_n} + 4\abr{u_n, v_n} + \abr{v_n, v_n}.
    %\\&=
  \end{align*}
  By Lemma \ref{230503.1530}
  ($\abr{u_n, u_n} = \half G_\infty + O_M(n^{-\half})$) and
  Lemma \ref{230721.1104} ($\abr{u_n,v_n}=O_M(n^{-H})$),
  we have 
  \begin{align*}
    \Delta_{M_n} &= 
    2G_\infty + \norm{v_n}_\calh^2 + O_M(n^{-\half}).
  \end{align*}
  We write $\calr_n = \Delta_{M_n} - \rbr{2G_\infty + \norm{v_n}_\calh^2}.$
  Notice that $\Delta_{M_n}\geq 2G_\infty + \calr_n$ since $\norm{v_n}_\calh^2\geq0$.
  
  Since we have $\norm{G_\infty^{-1}}_p<\infty$ for any $p>1$
  by Assumption \ref{ass:230927.1617} (ii),
  we obtain
  \begin{align*}
    P[\Delta_{M_n}<G_\infty] &\leq P[2G_\infty + \calr_n<G_\infty] =
    P[\calr_n<-G_\infty] 
    \\&\leq 
    P[\abs{\calr_n}>G_\infty] =
    E\sbr{\bbone_{\frac{\abs{\calr_n}}{G_\infty}>1}} \leq 
    E\sbr{\rbr{\frac{\abs{\calr_n}}{G_\infty}}^p} \leq 
    \rbr{\norm{\calr_n}_{2p}\norm{G_\infty^{-1}}_{2p}}^p 
    = O(n^{-\frac{p}{2}})
  \end{align*}
  for any $p\geq1$.
  \end{proof}

%   \begin{lemma}
%   {\myblue nondegeneracyの主張を書く}
%   There exist a sequence of positive random variables $s_n\in\bbD^{\ell-2,\infty}$ and $\kappa>0$ such that 
%   \begin{itemize}
%     \setlength{\parskip}{0cm} \setlength{\itemsep}{0pt} 
%     \item $\sup_{n\in\bbN}\rbr{\norm{s_n^{-1}}_p + \norm{s_n}_{\ell-2,p}}<\infty$ for any $p>1$.
%     \item $P\sbr{\Delta_{M_n}<s_n}=O(r_n^{1+\kappa})$
%   \end{itemize}
% \end{lemma}

\subsubsection{Proof of the asymptotic expansion of $Z_n^\circ$}
Based on the argument in Section \ref{sec:230726.1315}, we define the random symbol $\mS$ by 
\begin{align}
  \mS(\sfi\sfz) &= 
  \rbr{\mS^{(2,0)}_{0}+\mS^{(3,0)}+\mS^{(1,0)}+\mS^{(2,0)}_{1}}(\sfi\sfz)
  =%\nn\\&=
  c_\qtor \;\rbr{\iv}^{-3} \int_0^1 (a_t)^3 dt[\sfi\sfz]^3
  +\half\AN[\sfi\sfz], %\mS^{(1,0)}(\tti\sfz)
  % -\frac1{4} \rbr{\iv}^{-1}   c_{2,H} \Brbr{\brbr{V^{[1]}_0 }^2 + \brbr{V^{[1]}_1 }^2} [\sfi\sfz],
  \label{230726.1342}
\end{align}
where $c_\qtor$ is defined at \eqref{230726.1339}.
We write 
$\mS_n = 1+n^{-1/2}\mS$ and
$\cale(M,\gamma)=\cbr{f:\bbR\to\bbR\mid \abs{f(z)}\leq M(1+\abs{z})^\gamma}$ for $M,\gamma>0$.
\begin{proposition}\label{230726.1831}
  Let the functional $Z_n^\circ$ and random symbol $\mS$ be as in 
  \eqref{230726.1348} and \eqref{230726.1342}, respectively.
  Then, for any $M,\gamma>0$,
   the following asymptotic expansion formula for $Z_n^\circ$ holds as $n\to\infty$:
  \begin{align*}
    \sup_{f\in\cale(M,\gamma)}
    \abs{E\sbr{f(Z_n^\circ)} - \int_{z\in\bbR} f(z)p_n^{Z^\circ}(z)dz}
    =o(n^{-1/2}),
  \end{align*}
  where
  $p_n^{Z^\circ}(z) = E\sbr{\mS_n(\partial_z)^*\phi(z;0,G_\infty)}$.
\end{proposition}
\begin{proof}
Condition {\bf [D]}(i) is the assumption of the regularity of functionals involved, 
and we omit the detail.
We will verify Condition {\bf [D]}(ii). 
By setting $\kappa=(-1+2H)\wedge(2-2H)$($>0$ for $H\in(\half,1)$),
the estimate \eqref{220215.1241} follows from Lemma \ref{230721.1101},
\eqref{220215.1242} from Lemma \ref{230721.1235},
\eqref{220215.1243} and \eqref{220215.1244} from Lemma \ref{230721.1236},
and
\eqref{220215.1246} and \eqref{220215.1247} from Lemma \ref{230721.1242}.
We can prove the estimate \eqref{220215.1245} 
by Lemmas \ref{230721.1259} and \ref{230721.1236},
since 
$D_{u_n}^2G^{(2)}_n = D_{u_n}^3 M_n - D_{u_n}^2G_\infty$.
% {\myblue 
% $D_{u_n}^2N_n = O_M(n^{(\half-H)\vee(-1+H)})$.
% $n^{(\half-H)\vee(-1+H)}=
% n^{-(-\half+H)\wedge(1-H)}=
% n^{-\half2(-\half+H)\wedge(1-H)}$.
% $2(-\half+H)\wedge(1-H)=(-1+2H)\wedge(2-2H)$.
% }

For the pairs of random symbols 
$(\mS^{(2,0)}_{0,n},\mS^{(2,0)}_{0})$, 
$(\mS^{(3,0)}_n,\mS^{(3,0)})$, 
$(\mS^{(1,0)}_n,\mS^{(1,0)})$ and 
$(\mS^{(2,0)}_{1,n},\mS^{(2,0)}_{1})$,
where %the random symbols
$\mS^{(2,0)}_{0}$, $\mS^{(3,0)}$, $\mS^{(1,0)}$ and $\mS^{(2,0)}_{1}$ are 
defined in Section \ref{sec:230726.1315},
we can justify Condition {\bf [D]}(iii).
Indeed, (a) is obvious, and (b) has been verified at
Lemmas \ref{230725.1941}, \ref{230725.1942}, \ref{230725.2158} and \ref{230725.2159}.

Finally we check Condition {\bf [D]}(iv). Assumption \ref{ass:230927.1617} (ii) ensures (a).
For (b), we set $s_n=G_\infty$ for any $n\in\ntwo$ and $\kappa>0$ arbitrarily.
We can verify the condition (b) by Lemma \ref{230725.1635},
Assumption \ref{ass:230927.1617} (ii) and the fact $G_\infty\in\bbD^\infty$.
% {\myred 
% There exists a sequence of positive random variables $s_n\in\bbD^{\ell-2,\infty}$ and $\kappa>0$ such that 
% \begin{itemize}
%   \setlength{\parskip}{0cm} \setlength{\itemsep}{0pt} 
%   \item $\sup_{n\in\bbN}\rbr{\norm{s_n^{-1}}_p + \norm{s_n}_{\ell-2,p}}<\infty$ for any $p>1$.
%   \item $P\sbr{\Delta_{M_n}<s_n}=O(r_n^{1+\kappa})$
% \end{itemize}}
\end{proof}

The adjoint action of random symbols on $\phi(z;0,G_\infty)$
is defined at \eqref{eq:230927.1623}.
In this case, $p_n^{Z^\circ}$ is explicitly written as 
\begin{align*}
  &p_n^{Z^\circ}(z) =%&= 
  E\sbr{(1+n^{-1/2}\mS)(\partial_z)^*\phi(z;0,G_\infty)}
  %E\sbr{\mS_n(\partial_z)^*\phi(z;0,G_\infty)}
  % \\&=
  % E\sbr{\phi(z;0,G_\infty)} + n^{-\half}E\sbr{\mS(\partial_z)^*\phi(z;0,G_\infty)}
  \\&=
  E\sbr{\phi(z;0,G_\infty)} 
  + n^{-\half}\rbr{
  E\sbr{{\frac{c_{\qtor}}{\rbr{\iv}^{3}} \rbr{\int_0^1 (a_t)^3 dt}}
  \rbr{\frac{z^3}{G_\infty^3}-\frac{3z}{G_\infty^2}} \phi(z;0,G_\infty)}+
  E\sbr{{\half\AN}\frac{z}{G_\infty}\phi(z;0,G_\infty)}},
\end{align*}
where $\AN$ is defined at \eqref{230726.1845}.

\subsection{Asymptotic expansion of $Z_n$ and the error of the estimator}
% \subsection{Asymptotic expansion of $Z_n$ and $\sqrt{n}(\hat H_n-H)$}
\label{sec:231001.1649}
% {\myblue memoでは$X_n,X_\infty$と書いていたものを$Z_n^\circ,Z_\infty^\circ$で置き換えた．}
By Lemma \ref{lem:230926.1026}, $Z_n$ is written as 
% the expansion \eqref{230423.1826},
\begin{align}
  Z_n&=
  \sqrt{n} 
  \cbr{{\frac{\convDiff_n}\iv -\half\rbr{\frac{\convDiff_n}\iv}^2}
  - \rbr{\frac{\convDiff_{2n}}\iv -\half\rbr{\frac{\convDiff_{2n}}\iv}^2}
  }
  + \negTerm_n^{(1)},
  \label{230425.1507}
\end{align}
with %$\negTerm_n^{(1)}$ is defined at {\myblue\eqref{??}sto.exp.の証明内に出てくる定義を引用} and 
$\negTerm_n^{(1)}=O_M\rbr{n^{-1}}$.
We set  
$X^{(1)}_n=\sqrt{n} \frac{\convDiff_n}\iv$ and
$X^{(2)}_n=\sqrt{n} \frac{\convDiff_{2n}}\iv$.
Since we have defined 
\begin{align*}
  Y_n
  %=-\half \rbr{\brbr{X^{(1)}_n}^2 - \brbr{X^{(2)}_n}^2}
  =-n \cbr{\half\rbr{\frac{\convDiff_n}\iv}^2 - \half\rbr{\frac{\convDiff_{2n}}\iv}^2}
\end{align*}
and 
$Z_n^\circ = \delta(u_n) + r_n N_n=Z_n-r_nY_n$,
we have
\begin{align*}
  Y_n
  &=-\half \rbr{\brbr{X^{(1)}_n}^2 - \brbr{X^{(2)}_n}^2}
  %=n \cbr{-\half\rbr{\frac{\convDiff_n}\iv}^2 + \half\rbr{\frac{\convDiff_{2n}}\iv}^2}
  \\
  Z_n^\circ&
  =\sqrt{n} \rbr{\frac{\convDiff_n}\iv - \frac{\convDiff_{2n}}\iv}+\negTerm_n^{(1)}
  =X^{(1)}_n-X^{(2)}_n+\negTerm_n^{(1)}
\end{align*}
By Lemmas \ref{lem:230720.1911} and \ref{lem:230925.1301}, and
the proof of Proposition \ref{230720.2014},
we can see
\begin{align*}
  X^{(i)}_n
  %=\sqrt{n} \frac{\convDiff_{in}}\iv
  &=\sqrt{n} \frac{\convDiff_{in}^{(0)}}\iv + O_M(n^{-\half})
  \tand 
  \sqrt{n} \frac{\convDiff_{in}^{(0)}}\iv 
  =\bar M_n^{(i)} 
  =M_n^{(i)} + O_M(n^{-H}). 
\end{align*}
% $X^{(i)}_n
% %=\sqrt{n} \frac{\convDiff_{in}}\iv
% =\sqrt{n} \frac{\convDiff_{in}^{(0)}}\iv + O_M(n^{-\half})
% =\bar M_n^{(i)} + O_M(n^{-\half})
% =M_n^{(i)} + O_M(n^{-\half})$
By the mixed CLT of $(M_n^{(1)},M_n^{(2)})$ (Proposition \ref{230721.1408}),
we have the following weak convergence.
\begin{align*}
  (X^{(1)}_n,X^{(2)}_n)&\overset{d}{\to}
  \zeta A_\infty =: (X^{(1)}_\infty,X^{(2)}_\infty),
  % =   
  % \rbr{\widehat G_\infty}^{\half}
  % \rbr{\sqrt{2\hat c}\:\zeta_1, 
  % \frac{2^{2H}\tilde c}{\sqrt{2\hat c}}\zeta_1 +
  % \sqrt{\hat c -\frac{\rbr{2^{2H}\tilde c}^2}{2\hat c}}\zeta_2}
\end{align*}
where 
$\zeta=(\zeta_1,\zeta_2)$ is a $2$-dimensional 
standard normal random variable independent of $\calf$, and
the $2\times 2$-random matrix $A_\infty$ is defined at \eqref{230726.1711}.
As a result, $(Z_n^\circ,Y_n)$ converges in distribution as follows:
\begin{align*}
  (Z_n^\circ, Y_n) &
  \overset{d}{\to}
  \Brbr{X^{(1)}_\infty-X^{(2)}_\infty,
  -\half \brbr{\nrbr{X^{(1)}_\infty}^2 - \nrbr{X^{(2)}_\infty}^2}}
  =: (Z_\infty^\circ,Y_\infty).
\end{align*}
Notice that $Z_\infty^\circ\sim\caln(0,G_\infty)$.

Hence we will apply Theorem \ref{thm:230927.1632} % S.Y.perturbationの定理を引用する
to $Z^\circ$ and $Y$ in place of $\bbX$ and $\bbY$.
We set $\xi_n$ as follows: 
\begin{align*}
  \xi_n &= 
  \frac{3G_\infty}{2G_\infty + 4\Delta_{Z_n^\circ}}%4\Delta_{X_n}}
\end{align*}
Then the following lemma holds:
\begin{lemma}
  \begin{itemize}
    \setlength{\parskip}{0cm} \setlength{\itemsep}{0pt} 
    \item [(i)] $\sup_n\norm{\xi_n}_{k,p} <\infty$ for $k\geq0$ and $p>1$.
    \item [(ii)] $P\sbr{\abs{\xi_n}\geq\half} = O(r_n^\alpha)$ 
    with some $\alpha>1$.
    \item [(iii)] $\sup_n E\sbr{\bbone_\cbr{\abs{\xi_n}<1} \Delta_{Z_n^\circ}^{-p}} < \infty$ 
    for any $p>1$.
  \end{itemize}
\end{lemma}
\begin{proof}
  (i) This is obvious from 
  $Z_n^\circ=O_M(1)$ and $G_\infty^{-1}\in L^{\infty-}$.
  % To show $\sup_n\norm{F_n/G_n}_{k,p}<\infty$ for general functionals $F_n$ and $G_n$, 
  % it suffices to prove that \newline
  % $\sup_n \norm{F_n}_{k, (2k+2)p}<\infty$, 
  % $\sup_n \norm{G_n}_{k, (2k+2)p}<\infty$ and 
  % $\sup_n \norm{G_n^{-1}}_{(2k+2)p}<\infty$.

  % In this case, the first one is obvious. 
  % The second condition is also easily verified by $Z_n^\circ=O_M(1)$.
  % Since $2G_\infty+4\Delta_{Z_n^\circ}\geq2G_\infty\geq0$ and 
  % $G_\infty^{-1}\in L^{\infty-}$, 
  % we have 
  % $\bnorm{\rbr{2G_\infty+4\Delta_{Z_n^\circ}}^{-1}}_{(2k+2)p}\leq 
  % \bnorm{\rbr{2G_\infty}^{-1}}_{(2k+2)p}<\infty$.
  % This means the third condition.
  
  \noindent(ii)
  Since $\Delta_{Z_n^\circ} = 2G_\infty + \norm{v_n}_{\calh}^2 + \decayDelX_n$
  with $\decayDelX_n = O_M(n^{-\half})$,
  we have, for $L\geq1$,
  \begin{align*}
    P\sbr{\abs{\xi_n}\geq\half} &= 
    P\sbr{G_\infty\geq\Delta_{Z_n^\circ}} \leq 
    P\sbr{G_\infty\geq2G_\infty + \norm{v_n}_{\calh}^2 + \decayDelX_n}\leq 
    P\sbr{-G_\infty\geq \decayDelX_n} \leq
    P\sbr{G_\infty\leq \abs{\decayDelX_n}}
    \\&\leq 
    E\sbr{\bbone_\cbr{{\decayDelX_n}/{G_\infty}\geq1}}\leq
    E\sbr{\rbr{{\decayDelX_n}/{G_\infty}}^L} = 
    \norm{\decayDelX_n}_{2L}^L \norm{G_\infty^{-1}}_{2L}^L 
    = O(n^{-\frac{L}{2}})
  \end{align*}
  and we can take $L$ large enough to prove (ii).

  \noindent(iii)
  Since $G_\infty^{-1}\in L^{\infty-}$ and 
  $%\begin{align*}
    \abs{\xi_n}<1 %\Leftrightarrow \xi_n<1 
    %\Leftrightarrow 3G_\infty < 2G_\infty + 4 \Delta_{Z_n^\circ}
    \Leftrightarrow (0\leq)G_\infty < 4 \Delta_{Z_n^\circ}
    \Leftrightarrow 4 G_\infty^{-1} > \Delta_{Z_n^\circ}^{-1},
  $ %\end{align*}
  we have 
  \begin{align*}
    E\sbr{\bbone_\cbr{\abs{\xi_n}<1}\Delta_{Z_n^\circ}^{-p}} \leq
    E\sbr{\bbone_\cbr{\abs{\xi_n}<1}4^p G_\infty^{-p}} \leq
    4^p E\sbr{G_\infty^{-p}} < \infty
  \end{align*}
  for any $p>1$ and $n\in\ntwo$.
\end{proof}

Thus, we obtain the asymptotic expansion formula of $Z_n$ from Theorem \ref{thm:230927.1632}.
Recall that 
$\cale(M,\gamma)=\{f:\bbR\to\bbR\mid \abs{f(z)}\leq M(1+\abs{z})^\gamma\}$ for $M,\gamma>0$.
\begin{proposition}\label{prop:231005.1738}
  % {\myblue S.Yのperturbation methodより$Z_n$の漸近展開が得られる．}
  Let the functional $Z_n$ be as in \eqref{eq:230926.1554}.
  Then, for any $M,\gamma>0$, it holds that
  %  the following asymptotic expansion formula for $Z_n$ holds
  \begin{align*}
    \sup_{f\in\cale(M,\gamma)}
    \abs{E\sbr{f(Z_n)} - \int_{z\in\bbR} f(z)p_n^{Z}(z)dz}
    =o(n^{-\half})
  \end{align*}
  as $n\to\infty$,
  where
  $p_n^{Z}(z) = p_n^{Z^\circ}(z) + n^{-\half} g_\infty(z)$ with 
  $g_\infty(z)
  =-\partial_z \brbr{E\sbr{Y_\infty\mid Z_\infty^\circ = z} p^{Z_\infty^\circ}(z)}$
  and 
  $p_n^{Z^\circ}(z)$ is given at Proposition \ref{230726.1831}.
\end{proposition}
The modification term $g_\infty(z)$ can be written explicitly by 
\begin{align*}
  g_\infty(z)
  &=
  \frac{1}{2}
  %\frac{\Sigma_{11}-\Sigma_{22}}{\Sigma_{11}+\Sigma_{22}-2\Sigma_{12}}
  \frac{\hat c}{3\hat c-2\tilde c 2^{2H}}
  E\sbr{\rbr{2z-\frac{z^3}{G_\infty}}\phi(z;0,G_\infty)}.
\end{align*}
% \begin{align*}
%   g_\infty(z)
%   &=
%   \frac{1}{2}
%   \frac{\Sigma_{11}-\Sigma_{22}}
%   {\Sigma_{11}+\Sigma_{22}-2\Sigma_{12}}
%   E\sbr{\rbr{2z-\frac{z^3}{G_\infty}}\phi(z;0,G_\infty)}
%   \text{ with }
%   \begin{pmatrix}
%     \Sigma_{11} & \Sigma_{12} \\
%     \Sigma_{12} & \Sigma_{22} \\
%   \end{pmatrix}
%   =
%   \begin{pmatrix}
%     2\hat c & \tilde c 2^{2H} \\
%     \tilde c 2^{2H} & \hat c \\
%   \end{pmatrix}.
% \end{align*}
The main result of this paper follows immediately from Proposition \ref{prop:231005.1738}.
\begin{theorem}\label{thm:231001.1041}
  % $\sqrt{n}(\hat H_n-H)$の漸近展開．
  Suppose that SDE \eqref{eq:230925.1452} satisfies Assumption \ref{ass:230927.1617},
  and the estimator $\hat H_n$ is defined as in \eqref{eq:231001.1046}.
  Let $M>0$ and $\gamma>1$.
  The following asymptotic expansion formula for $\sqrt{n} (\hat H_n - H)$ 
  holds as $n\to\infty$
   \begin{align*}
    \sup_{f\in\cale(M,\gamma)}\abs{E[f(\sqrt{n} (\hat H_n - H))] - \int_{z\in\bbR} f(z) p_n(z) dz} = {o}(n^{-\half})
  \end{align*}
  with $p_n(z)=2\log2\times p_n^Z(2\log2\,z)$.
\end{theorem}
\begin{proof}
  Recall $Z_n=(2\log2)\sqrt{n}(\hat H_n-H)\psi_n$.
  Write $Z_n^*=\sqrt{n}(\hat H_n-H)\psi_n$ and $f^*(x)=f\rbr{(2\log2)^{-1}x}$.
  Since
  $\int_{z\in\bbR} f^*(z) p_n^Z(z) dz=\int_{z\in\bbR} f\rbr{z}p_n(z) dz$,
  we have 
  \begin{align*}
    E\sbr{f(Z_n^*)}=E\sbr{f^*(Z_n)}=\int_{z\in\bbR} f\rbr{z}p_n(z) dz + {o}(n^{-\half})
  \end{align*}
  uniformly in $f$.

  Then, by the estimate \eqref{230704.1630} about $\psi_n$, (see Lemma \ref{230509.2330},)
  \begin{align*}
    \abs{E[f(\sqrt{n} (\hat H_n - H))]-E\sbr{f(Z_n^*)}}
    % &=%\\&=
    % \abs{E\sbr{f(\sqrt{n}(\hat H_n-H))-f(\sqrt{n}(\hat H_n - H)\psi_n)}}
    &\leq%\\&\leq
    E\sbr{\abs{f(\sqrt{n}(\hat H_n-H))-f(\sqrt{n}(\hat H_n - H)\psi_n)}}
    \\&=
    E\sbr{\babs{f(\sqrt{n}(\hat H_n-H))-f(\sqrt{n}(\hat H_n - H)\psi_n)}\; \bbone_\cbr{\psi_n<1}}
    \\&\leq
    E\sbr{\babs{f(\sqrt{n}(\hat H_n - H))}\: \bbone_\cbr{\psi_n<1}}+
    E\sbr{\babs{f(\sqrt{n}(\hat H_n-H)\psi_n)}\: \bbone_\cbr{\psi_n<1}}
    \\&\leq
    2\times M(1+n^{\half})^\gamma\, P\sbr{\psi_n<1}
    \\&= 
    O(n^{\frac{\gamma}{2}+(-L)})
  \end{align*}
  for any $L>1$ and uniformly in $f\in\cale(M,\gamma)$.
  By taking $L$ large enough, we obtain the result.
\end{proof}

% \newpage
\subsection{Weighted graphs}\label{230925.1250}
We collect all the weighted graphs (and their exponents) we encounter in the proof of
the asymptotic expansion formula including the lemmas' in Section \ref{sec:231002.2426}.
The rescaling factors are omitted from the functionals in the caption below the graphs.
The functionals without reference numbers appear 
in the proof of Lemma \ref{lem:230925.1301} (i).

% %\newcommand{\graphcycle}{
%   \begin{tikzpicture}
%     %% vertex labels
%     \node[shape=circle,draw=black,inner sep=1pt] (1) at (120:1){0};
%     \node[shape=circle,draw=black,inner sep=1pt] (2) at (180:1) {0};
%     \node[shape=circle,draw=black,inner sep=1pt] (3) at (0:1) {0};
%     \node[shape=circle,draw=black,inner sep=1pt] (4) at (60:1) {0};
%     %%% edges
%     \draw[thick] (1) to [out=210,in=90,looseness=0.5] (2);
%     \draw[dashed,thick] (2) to [out=270,in=270,looseness=1.6] (3);
%     \draw[thick] (3) to [out=90,in=330,looseness=0.5] (4);
%     \draw[thick] (4) to [out=150,in=30,looseness=0.5] (1);
  
%     \node  at (150:0.8) {1};
%     \node  at (30:0.8) {1};
%     \node  at (90:0.8) {1};
%     \node  at (125:1.45) {$v^G_{1}$};
%     \node  at (170:1.4) {$v^G_{2}$};
%     \node  at (10:1.7) {$v^G_{I(G)-1}$};
%     \node  at (55:1.5) {$v^G_{I(G)}$};
%    \end{tikzpicture}
% %}

%\captionsetup[subfigure]{labelformat=empty}

\begin{figure}[H]
  \begin{tabular}{ccc}
    \begin{minipage}[t]{0.30\textwidth}
      \centering
      \begin{tikzpicture}
        %% vertex labels
        \node[shape=circle,draw=black,inner sep=1pt] (1) at (0,0) [label=below:$v_{0}$]{0};
        \node[shape=circle,draw=black,inner sep=1pt] (2) at (0,0.6) []{2};
        \draw[double, line width=1pt] (1) -- (2);
      \end{tikzpicture}
      \caption{
      exponent: $\half-2H$\\
      $\convDiff^{(0)}_n$ at \eqref{230925.1600}.\hspace{40pt}
      % exponent: $1+(-\half-2H)\\=\half-2H$
      }
      \label{fig:230721.1711}
    \end{minipage} 
    &
    %%%%
    \begin{minipage}[t]{0.30\textwidth}
      \centering
      \begin{tikzpicture}
        %% vertex labels
        \node[shape=circle,draw=black,inner sep=1pt] (1) at (0,0) [label=below:$v_{0}$]{2};
        \node[shape=circle,draw=black,inner sep=1pt] (2) at (0,0.6) []{1};
        \draw[double, line width=1pt] (1) -- (2);
      \end{tikzpicture}
      \caption{
        exponent: $-2H$\\
        $\check\convDiff^{(2,1,1)}_n$, %\eqref{230703.1611}\\
        $\check\convDiff^{(2,4,1)}_n$,
        $\check\convDiff^{(2,4,3)}_n$,\\ %\eqref{230703.1631}, \eqref{230703.1633}\\
        % ($\caliN{1}_n$, $\caliN{2}_n$, $\caliN{3}_n$)
        % \begin{comment}
        %   {\myblue$\convDiff^{(2,1)}_n, \convDiff^{(2,2)}_n$に残るorder}\\
        % \end{comment}
        }
        % exponent: $1+(-1-2H)=-2H$}
      \label{fig:230721.1712}
    \end{minipage} 
    &
    %%%%
    \begin{minipage}[t]{0.30\textwidth}
      \centering
      \begin{tikzpicture}
        %% vertex labels
        \node[shape=circle,draw=black,inner sep=1pt] (1) at (0,0) [label=below:$v_{0}$]{1};
        \node[shape=circle,draw=black,inner sep=1pt] (2) at (0,0.6) []{0};
        \draw[double, line width=1pt] (1) -- (2);
      \end{tikzpicture}
      \caption{
        exponent: $0$\\  
        $\convDiff^{(1,1)}_n$ at \eqref{230804.1411},\\% ($\caliN{4}_n$),\\ %\\%\eqref{230703.1604}, 
        $\check\convDiff^{(2,1,2)}_n$, %\eqref{230703.1612}\\
        $\check\convDiff^{(2,4,2)}_n$,
        $\check\convDiff^{(2,4,4)}_n$. %\eqref{230703.1632}, \eqref{230703.1634}\\
      }
      % exponent: $1-1=0$}
      \label{fig:230721.1713}
    \end{minipage} 
  \end{tabular}
\end{figure}
\begin{figure}[H]
  \begin{tabular}{ccc}
    \begin{minipage}[t]{0.30\textwidth}
      \centering
      \begin{tikzpicture}
        %% vertex labels
        \node[shape=circle,draw=black,inner sep=1pt] (1) at (0,0) [label=below:$v_{0}$]{0};
        \node[shape=circle,draw=black,inner sep=1pt] (2) at (0,0.6) []{1};
        \draw[double, line width=1pt] (1) -- (2);
      \end{tikzpicture}
      \caption{
        exponent: $\half-H$\\
        $\negTerm_n^{(3,1)}$, %:\eqref{230703.1701},
        $\negTerm_n^{(3,2)}$ at \eqref{eq:231005.1844},\\ %:\eqref{230703.1702},\\
        {$\babr{DG_\infty,u_n^{(i)}}$} at \eqref{eq:231005.1842}, \\
        % {$n^{\half}\babr{DG_\infty,u_n^{(i)}}$}, \\
        {$D_{u_n^{(i)}}\AN$} at \eqref{eq:231005.1843}, \\
        % {$n^{\half}D_{u_n^{(i)}}\AN$} %\\
        $\check\convDiff^{(2,1,3)}_n$. %\eqref{230703.1613}\\
      }
      % exponent: $1+(-\half-H)=\half-H$}
      \label{fig:230721.1714}
    \end{minipage} 
    &
    %%%%
    \begin{minipage}[t]{0.30\textwidth}
      \centering
      \begin{tikzpicture}
        %% vertex labels
        \node[shape=circle,draw=black,inner sep=1pt] (1) at (0,0) [label=below:$v_{0}$]{1};
        \node[shape=circle,draw=black,inner sep=1pt] (2) at (0,0.6) []{1};
        \draw[double, line width=1pt] (1) -- (2);
      \end{tikzpicture}
      \caption{
        exponent: $-H$\\  
      $\check\convDiff^{(2,2,1)}_n$, %\eqref{230703.1614}\\
      $\check\convDiff^{(2,3,1)}_n$,  
      $\check\convDiff^{(2,5,1)}_n$, \\$\check\convDiff^{(2,5,2)}_n$, % \eqref{230703.1641}\\
      $\check\convDiff^{(2,6,1)}_n$, $\check\convDiff^{(2,6,2)}_n$. %\eqref{230703.1646}\\ 
      }
      % exponent: $1+(-1-H)=-H$}
    \end{minipage} 
    &
    %%%%
    \begin{minipage}[t]{0.30\textwidth}
      \centering
      \begin{tikzpicture}
        %% vertex labels
        \node[shape=circle,draw=black,inner sep=1pt] (1) at (0,0) [label=below:$v_{0}$]{0};
        \node[shape=circle,draw=black,inner sep=1pt] (2) at (0,0.6) []{0};
        \draw[double, line width=1pt] (1) -- (2);
      \end{tikzpicture}
      \caption{
        exponent: $1$\\
        $\check\convDiff^{(2,2,2)}_n$, %\eqref{230703.1615}\\ 
        $\check\convDiff^{(2,3,2)}_n$, %の2つ目\\ 
        $\check\convDiff^{(2,5,3)}_n$, %の2つ目, \eqref{230703.1642}\\
        $\check\convDiff^{(2,6,3)}_n$. %の2つ目, \eqref{230703.1647}\\ 
        }
      \end{minipage} 
    % \\
  \end{tabular}
\end{figure}
\begin{figure}[H]
  \begin{tabular}{ccc}
    \begin{minipage}[t]{0.30\textwidth}
      \centering
      \begin{tikzpicture}
        %% vertex labels
        \node[shape=circle,draw=black,inner sep=1pt] (1) at (-.8,0) [label=below:$v_{0}$]{0};
        \node[shape=circle,draw=black,inner sep=1pt] (2) at (-.8,0.6) []{0};
        \draw[double, line width=1pt] (1) -- (2);
  
        \node[shape=circle,draw=black,inner sep=1pt] (3) at (0.8,0) [label=below:$v_{1}$]{0};
        \node[shape=circle,draw=black,inner sep=1pt] (4) at (0.8,0.6) []{0};
        \draw[double, line width=1pt] (3) -- (4);
  
        \draw[thick] (2) -- node[midway, above] {2} (4);
  
      \end{tikzpicture}
      \caption{
        exponent: $1-4H$\\
        $\cali^{(k_0;i;0)}_n$ at \eqref{230703.1713},\\
      $\cali^{(k_0;i;2)}_n$ at \eqref{230703.1728},\\
      $\cali_n^{(i_1,i_2)(1)}$ at \eqref{230703.1731}.
      }
      % exponent: $1-4H+0=1-4H$}
      \label{fig:230721.1717}
    \end{minipage} 
    &
    %%%%
    \begin{minipage}[t]{0.30\textwidth}
      \centering
      \begin{tikzpicture}
        %% vertex labels
        \node[shape=circle,draw=black,inner sep=1pt] (1) at (-.8,0) [label=below:$v_{0}$]{0};
        \node[shape=circle,draw=black,inner sep=1pt] (2) at (-.8,0.6) []{1};
        \draw[double, line width=1pt] (1) -- (2);
  
        \node[shape=circle,draw=black,inner sep=1pt] (3) at (0.8,0) [label=below:$v_{1}$]{0};
        \node[shape=circle,draw=black,inner sep=1pt] (4) at (0.8,0.6) []{1};
        \draw[double, line width=1pt] (3) -- (4);
  
        \draw[thick] (2) -- node[midway, above] {1} (4);
      \end{tikzpicture}
      \caption{
    exponent: $\half-4H$\\  
      $\cali^{(k_0;i;1)}_n$ at \eqref{230703.1713},\\
    $\cali_n^{(i_1,i_2)(2)}$ at \eqref{230703.1732}.
    }
    % exponent: $(1-2H)+(-\half-2H)\\=\half-4H$}
    \label{fig:230721.1718}
    \end{minipage} 
    &
    %%%%
    \begin{minipage}[t]{0.30\textwidth}
      \centering
      \begin{tikzpicture}
        %% vertex labels
        \node[shape=circle,draw=black,inner sep=1pt] (1) at (-.6,0) [label=below:$v_{0}$]{0};
        \node[shape=circle,draw=black,inner sep=1pt] (2) at (-.6,0.6) []{2};
        \draw[double, line width=1pt] (1) -- (2);
  
        \node[shape=circle,draw=black,inner sep=1pt] (3) at (0.6,0) [label=below:$v_{1}$]{0};
        \node[shape=circle,draw=black,inner sep=1pt] (4) at (0.6,0.6) []{1};
        \draw[double, line width=1pt] (3) -- (4);
  
        %\draw[thick] (2) -- node[midway, above] {1} (4);
      \end{tikzpicture}
      \caption{
        exponent: $1-3H$\\  
      $\calj_n^{(i_1;i_2)}$ at \eqref{eq:231005.1823}.%{230703.1733}
    }
    % exponent: $(1+(-\half-2H))\\+(1+(-\half-H))=1-3H$}
    \label{fig:230721.1719}
    \end{minipage} 
    % \vspsm
    % \\
  \end{tabular}
\end{figure}
\begin{figure}[H]
  \begin{tabular}{ccc}
    \begin{minipage}[t]{0.30\textwidth}
      \centering
      \begin{tikzpicture}
        %% vertex labels
        \node[shape=circle,draw=black,inner sep=1pt] (1) at (-1.4,0) [label=below:$v_{1}$]{0};
        \node[shape=circle,draw=black,inner sep=1pt] (2) at (-1.4,0.6) []{0};
        \draw[double, line width=1pt] (1) -- (2);
  
        \node[shape=circle,draw=black,inner sep=1pt] (3) at (0,0) [label=below:$v_{2}$]{0};
        \node[shape=circle,draw=black,inner sep=1pt] (4) at (0,0.6) []{0};
        \draw[double, line width=1pt] (3) -- (4);
  
        \node[shape=circle,draw=black,inner sep=1pt] (5) at (1.4,0) [label=below:$v_{3}$]{0};
        \node[shape=circle,draw=black,inner sep=1pt] (6) at (1.4,0.6) []{1};
        \draw[double, line width=1pt] (5) -- (6);
  
        \draw[thick] (2) -- node[midway, above] {2} (4);
  
      \end{tikzpicture}
      \caption{
        exponent: $\frac32-5H$\\  
      $D_{u_n^{(i_3)}}\,\cali_n^{(i_1,i_2)(1)}$ at \eqref{230628.1149}.
      % $n^{-4H+\frac32}D_{u_n^{(i_3)}}\,\cali_n^{(i_1,i_2)(1)}$\\
      }
      % exponent: $(1-4H)+(1+(-\half-H))\\=\frac32-5H$}
      \label{fig:230721.1720}
    \end{minipage} 
    &
    %%%%
    \begin{minipage}[t]{0.30\textwidth}
      \centering
      \begin{tikzpicture}
        %% vertex labels
        \node[shape=circle,draw=black,inner sep=1pt] (1) at (-1.4,0) [label=below:$v_{1}$]{0};
        \node[shape=circle,draw=black,inner sep=1pt] (2) at (-1.4,0.6) []{1};
        \draw[double, line width=1pt] (1) -- (2);
  
        \node[shape=circle,draw=black,inner sep=1pt] (3) at (0,0) [label=below:$v_{2}$]{0};
        \node[shape=circle,draw=black,inner sep=1pt] (4) at (0,0.6) []{1};
        \draw[double, line width=1pt] (3) -- (4);
  
        \node[shape=circle,draw=black,inner sep=1pt] (5) at (1.4,0) [label=below:$v_{3}$]{0};
        \node[shape=circle,draw=black,inner sep=1pt] (6) at (1.4,0.6) []{1};
        \draw[double, line width=1pt] (5) -- (6);
  
        \draw[thick] (2) -- node[midway, above] {1} (4);
  
      \end{tikzpicture}
      \caption{
        exponent: $1-5H$\\  
      $\cali_n^{(3;i)(4)}$ at \eqref{eq:231005.1831}.
      }
      % exponent: $((1-2H)+(-\half-2H))\\+(1+(-\half-H))=1-5H$}
      \label{fig:230721.1721}
    \end{minipage} 
    &
    %%%%
    \begin{minipage}[t]{0.30\textwidth}
      \centering
      \begin{tikzpicture}
        %% vertex labels
        \node[shape=circle,draw=black,inner sep=1pt] (1) at (-1.4,0) [label=below:$v_{1}$]{0};
        \node[shape=circle,draw=black,inner sep=1pt] (2) at (-1.4,0.6) []{1};
        \draw[double, line width=1pt] (1) -- (2);
  
        \node[shape=circle,draw=black,inner sep=1pt] (3) at (0,0) [label=below:$v_{2}$]{0};
        \node[shape=circle,draw=black,inner sep=1pt] (4) at (0,0.6) []{0};
        \draw[double, line width=1pt] (3) -- (4);
  
        \node[shape=circle,draw=black,inner sep=1pt] (5) at (1.4,0) [label=below:$v_{3}$]{0};
        \node[shape=circle,draw=black,inner sep=1pt] (6) at (1.4,0.6) []{1};
        \draw[double, line width=1pt] (5) -- (6);
  
        \draw[thick] (2) -- node[midway, above] {1} (4);
        \draw[thick] (4) -- node[midway, above] {1} (6);
  
      \end{tikzpicture}
    \caption{
      exponent: $\half-6H$\\
      $\cali_n^{(3;i)(2)}$, $\cali_n^{(3;i)(3)}$
      at \eqref{eq:231005.1831}.
    }
    % exponent: $(1-4H)+(-\half-2H)\\=\half-6H$}
    \label{fig:230721.1722}
    \end{minipage} 
    % \vspsm
    % \\
  \end{tabular}
\end{figure}
\begin{figure}[H]
  \begin{tabular}{ccc}
    \begin{minipage}[t]{0.30\textwidth}
      \centering
      \begin{tikzpicture}
        %% vertex labels
        \node[shape=circle,draw=black,inner sep=1pt] (1) at (-1.5,0) [label=below:$v_{1}$]{0};
        \node[shape=circle,draw=black,inner sep=1pt] (2) at (-1.5,0.6) []{0};
        \draw[double, line width=1pt] (1) -- (2);
  
        \node[shape=circle,draw=black,inner sep=1pt] (3) at (0,-1.6) [label=below:$v_{2}$]{0};
        \node[shape=circle,draw=black,inner sep=1pt] (4) at (0,-1) []{0};
        \draw[double, line width=1pt] (3) -- (4);
  
        \node[shape=circle,draw=black,inner sep=1pt] (5) at (1.5,0) [label=below:$v_{3}$]{0};
        \node[shape=circle,draw=black,inner sep=1pt] (6) at (1.5,0.6) []{0};
        \draw[double, line width=1pt] (5) -- (6);
  
        \draw[thick] (2) -- node[midway, above] {1} (4);
        \draw[thick] (4) -- node[midway, above] {1} (6);
        \draw[thick] (6) -- node[midway, above] {1} (2);
  
      \end{tikzpicture}
    \caption{
      exponent: $1-6H$\\
      $\cali_n^{(3;i)(1)}$ at \eqref{eq:231005.1831}.
    }
    \label{fig:230721.1723}
    \end{minipage} 
    &
    %%%%
    \begin{minipage}[t]{0.30\textwidth}
      \centering
      \begin{tikzpicture}
        %% vertex labels
        \node[shape=circle,draw=black,inner sep=1pt] (1) at (-.8,0) [label=below:$v_{0}$]{0};
        \node[shape=circle,draw=black,inner sep=1pt] (2) at (-.8,0.6) []{1};
        \draw[double, line width=1pt] (1) -- (2);
  
        \node[shape=circle,draw=black,inner sep=1pt] (3) at (0.8,0) [label=below:$v_{1}$]{0};
        \node[shape=circle,draw=black,inner sep=1pt] (4) at (0.8,0.6) []{0};
        \draw[double, line width=1pt] (3) -- (4);
  
        \draw[thick] (2) -- node[midway, above] {1} (4);
      \end{tikzpicture}
      \caption{
        exponent: $\half-3H$\\
        $\cali_n^{(\qtan)(1,1)(k)}$ at \eqref{eq:231005.1838}.
        % $r_n^{-1}\cali_n^{(i_1,i_2)(2)}$のWeak limitに関連
    }
    % exponent: $(1-2H)+(-\half-H)\\=\half-3H$}
    \label{fig:230721.1724}
    \end{minipage} 
    &
    %%%%
    % \begin{minipage}[t]{0.30\textwidth}
    %   \centering
    %   \begin{tikzpicture}
    %     %% vertex labels
    %     \node[shape=circle,draw=black,inner sep=1pt] (1) at (0,0) [label=below:$v_{0}$]{2};
    %     \node[shape=circle,draw=black,inner sep=1pt] (2) at (0,0.6) []{0};
    %     \draw[double, line width=1pt] (1) -- (2);
    %   \end{tikzpicture}
    % \caption{
    %   exponent: $-H$\\  
    % $n\,\rbr{\iv}^{-1} \convDiff^{(2,1)}_n$, 
    %   $n\,\rbr{\iv}^{-1} \convDiff^{(2,2)}_n$のweak limit関連\\
    %   }
    %   % exponent: $1+(-1-H)=-H$}
    % \label{fig:230721.1725}
    % \end{minipage} 

  \end{tabular}
\end{figure}

%% file: subfiles/6-sec_technical_lemmas.tex
\section{Technical lemmas}
\label{sec:231002.2426}
In this section, we collect the lemmas used in Section \ref{sec:231002.2423}
and give the proofs of them.
In Section \ref{sec:231004.1533}, we give some basic estimates related to SDE driven 
by fractional Brownian motion with $H>1/2$.
Lemmas related to decomposing $Z_n$ to obtain its stochastic expansion
are considered in Section \ref{230804.1348}, 
while in Section \ref{sec:231005.2102}
we deal with convergence of functionals related to random symbols 
appearing in Section \ref{sec:231002.2423}.

\subsection{Estimates about stochastic differential equations driven by fBm with $H>1/2$}
\label{sec:231004.1533}

We consider estimates about SDE \eqref{eq:230925.1452} 
under Assumption \ref{ass:230927.1617} (i) on the coefficients
$V^{[1]}$ and $V^{[2]}$.
% We have the following lemma.

\begin{lemma}\label{230924.1800}
  Suppose that $X_t$ is the solution of SDE \eqref{eq:230925.1452}
  under Assumption \ref{ass:230927.1617} (i).
  Let $T>0$ and
  $f\in C^\infty(\bbR)$ with bounded derivatives of any order.
  \item[(i)]
  Let $r_0<r_1\in[0,T]$.
  Writing $Y^{(0)}_t=t$ and $Y^{(1)}_t=B_t$,
  consider the $k$-th $(k\in\bbN)$ iterated integral 
  \begin{align*}
    F_{r_0,r_1}^{(i)}&=
    \int^{r_1}_{r_0}
    \int^{t^{(1)}}_{r_0}\cdots\rbr{
    \int^{t^{(k-1)}}_{r_0}
    f(X_{t^{(k)}})dY^{(i_k)}_{t^{(k)}}}
    \cdots dY^{(i_2)}_{t^{(2)}}dY^{(i_1)}_{t^{(1)}}
  \end{align*}
  for $(i_j)_{j=1}^k\in\cbr{0,1}^k$.
  Then for any $N\in\bbZ_{\geq0}$ and $p\geq1$,
  the following estimate holds with any $\epsilon>0$:
  \begin{align*}
    \norm{F_{r_0,r_1}^{(i)}}_{N,p}&\leq
    C \abs{r_0-r_1}^{k_0+k_1H-\epsilon\bbone_\cbr{k_1>0}},
  \end{align*}
  where 
  $k_\ell=\sharp\cbr{1\leq j\leq k \mid i_j=\ell}$ for $\ell=0,1$, and
  the constant $C$ is independent of $r_0$ and $r_1$.

  \item[(ii)] %\ref{230807.1543}
  Let $r_0<r_1\in[0,T]$, $N\in\bbZ_{\geq0}$ and $p>1$.
  For any $\epsilon>0$,
  it holds that 
  \begin{align*}
    \norm{f(X_{r_1})-f(X_{r_0})}_{N,p} \leq C \abs{r_1-r_0}^{H-\epsilon}
  \end{align*}
  with the constant $C$ independent of $r_0$ and $r_1$.

  \item[(iii)]
  The difference of the integral of $f(X_t)$ %with respect to $dt$ 
  and its Riemann sum is estimated as
  \begin{align}\label{230924.1812}
    \int^T_0 f(X_t)dt - 
    \frac{T}{n}\sum_{j=0}^{n-1} f(X_{t_{j}})
    =O_M(n^{-1}),
  \end{align}
  where $t_j=jT/n$.
\end{lemma}

\begin{remark}
  The proof of (i) and (ii) of the following lemma are based on a similar arguments 
in Lemma 5.5 of \cite{2024Yamagishi-AsymptoticEO}.
Although in that paper the authors studied a functionals
related to fBm with $H\in(\frac12,\frac34)$,
the proof there does not depend on the assumption $H\in(\frac12,\frac34)$ and
works in the case of $H\in[\frac34,1]$ without change.
This is because the estimates used there rely on the arguments developed 
in Section 3 and Lemma A.1 in \cite{hu2016rate}, where only $H\in(\half,1)$ is assumed.

\end{remark}

\begin{proof}[Proof of Lemma \ref{230924.1800}]
  The estimate (i) is a generalization of Lemma 5.5 (ii) of \cite{2024Yamagishi-AsymptoticEO}, and 
  (ii) is the same as Lemma 5.5 (i) of \cite{2024Yamagishi-AsymptoticEO}.

  We only give the proof of (iii).
  The difference decomposes as 
  \begin{align*}
    \int^T_0 f(X_t)dt - 
    \frac{T}{n}\sum_{j=0}^{n-1} f(X_{t_{j}})
    % \\&=
    % \sum_{j=0}^{n-1}\cbr{
    % \int^{t_{j+1}}_{t_j} f(X_t)dt - 
    % \frac{T}{n} f(X_{t_{j}})}
    &=
    \sum_{j=0}^{n-1}
    \int^{t_{j+1}}_{t_j} \rbr{f(X_t) - f(X_{t_{j}})}dt
    \\&=
    \sum_{j=0}^{n-1}\cbr{
    \int^{t_{j+1}}_{t_j} 
    \int^t_{t_{j}} f'(X_{t'}) V^{[1]}_{t'} dB_{t'}dt
    +
    \int^{t_{j+1}}_{t_j} 
    \int^t_{t_{j}} f'(X_{t'}) V^{[2]}_{t'} dt'dt}
    \\&=
    S_{1,n}+S_{2,n}+S_{3,n}
  \end{align*}
  with 
  \begin{align*}
    S_{1,n}&=
    \sum_{j=0}^{n-1}\cbr{
    \int^{t_{j+1}}_{t_j} 
    \int^t_{t_{j}} \rbr{f'(X_{t'}) V^{[1]}_{t'} - f'(X_{t_j}) V^{[1]}_{t_j}}
    dB_{t'}dt}
    \\
    S_{2,n}&=
    \sum_{j=0}^{n-1}
    f'(X_{t_j}) V^{[1]}_{t_j}
    \int^{t_{j+1}}_{t_j} \int^t_{t_{j}}  dB_{t'}dt
    \\
    S_{3,n}&=
    \sum_{j=0}^{n-1}
    \int^{t_{j+1}}_{t_j} 
    \int^t_{t_{j}} f'(X_{t'}) V^{[2]}_{t'} dt'dt.
  \end{align*}
  For $S_{1,n}$, decomposing the difference 
  $f'(X_{t'}) V^{[1]}_{t'} - f'(X_{t_j}) V^{[1]}_{t_j}$ 
  by the change of variables, and using the estimate (i),
  we have 
  $S_{1,n}=O_M(n^{1-2H-1+\epsilon})=O_M(n^{-2H+\epsilon})$ with any $\epsilon>0$.
  Similarly, we have 
  $S_{3,n}=O_M(n^{-1})$.

  Since we can write 
  \begin{align*}
    S_{2,n}&=
    \sum_{j=0}^{n-1}
    f'(X_{t_j}) V^{[1]}_{t_j}
    \int^{t_{j+1}}_{t_j} \int^t_{t_{j}}  dB_{t'}dt
    =
    n^{-1}
    \sum_{j=0}^{n-1}
    f'(X_{t_j}) V^{[1]}_{t_j}
    I_1(f_j)
  \end{align*}
  with $f_j(\cdot)=n(t_{j+1}-\cdot)\bbone_{[t_{j},t_{j+1}]}$,
  by arguments of exponent and 
  % theory of exponent introduced in Section \ref{sec:231002.2156} 
  Proposition \ref{230725.1133},
  we have $S_{2,n}=O_M(n^{-1})$.
  Hence we obtain \eqref{230924.1812}.
\end{proof}

\subsection{Related to the stochastic expansion of $Z_n$}
\label{230804.1348}
\subsubsection{Decomposition of $\convDiff_n$}
\label{230926.1330}

Recall that we have defined 
$\convDiff_n = \rv_n-\iv$ at \eqref{230804.1403}, and 
we have 
$\convDiff_n=\sum_{k=0}^3 \convDiff^{(k)}_n$ with 
$\convDiff^{(k)}_n$ ($k=0,...,3$)
defined at 
\eqref{230925.1600},
\eqref{230925.1601},
\eqref{230925.1602} and
\eqref{230925.1603}.
We decompose 
$\convDiff^{(1)}_n$ at Lemma \ref{lem:230720.1911}, and
$\convDiff^{(2)}_n$ and $\convDiff^{(3)}_n$ at Lemma \ref{lem:230925.1301}.

Recall the definition of $\convDiff^{(1)}_n$:
\begin{align*}
  \convDiff^{(1)}_n &=
  c_{2,H}\;n^{-1} \sum_{j=1}^{n-1} \rbr{V^{[1]}_{t^n_j}}^2
  - c_{2,H}\int^1_0 \rbr{V^{[1]}_t}^2 dt,
\end{align*}
where the constant $c_{2,H}$ is defined at \eqref{230925.1607}.
Recall also that
we write 
$F_n=\hat O_M(n^\alpha)$,
when $F_n=O_M(n^{\alpha+\epsilon})$ for any $\epsilon>0$, 
\begin{lemma}\label{lem:230720.1911}
  The functional $\convDiff^{(1)}_n$ decomposes as
  \begin{align*}    
    \convDiff^{(1)}_n =
    \convDiff^{(1,1)}_n + \convDiff^{(1,2)}_n
    + \hat O_M(n^{-2H}),
  \end{align*}
  where 
  $\convDiff^{(1,1)}_n$ and $\convDiff^{(1,2)}_n$ are given at 
  \eqref{eq:230809.1743} and \eqref{eq:230809.1745} below, 
  respectively, 
  and $\convDiff^{(1,1)}_n,\convDiff^{(1,2)}_n=O_M(n^{-1})$.
  In particular, it holds that $\convDiff^{(1)}_n=O_M(n^{-1})$.
\end{lemma}
\begin{proof}
  We expand $c_{2,H}^{-1}\convDiff^{(1)}_n$ as
  \begin{align*}
    c_{2,H}^{-1}\convDiff^{(1)}_n &=
    n^{-1} \sum_{j=1}^{n-1} \brbr{V^{[1]}_{t^n_j}}^2
    - \int^1_0 \brbr{V^{[1]}_t}^2 dt
    =
    \sum_{k=1}^3 \check\convDiff^{(1,k)}_n
  \end{align*}
  with
  \begin{align*}
    \check\convDiff^{(1,1)}_n&=
    n^{-1} \sum_{j=1}^{n-1} \brbr{V^{[1]}_{t^n_j} }^2
    - \int^{t^{2n}_{2n-1}}_{t^{2n}_1} \brbr{V^{[1]}_t}^2 dt
    \\ 
    \check\convDiff^{(1,2)}_n&=
    -\frac1{2n}\rbr{\brbr{V^{[1]}_0 }^2 + \brbr{V^{[1]}_1 }^2}
    \\
    \check\convDiff^{(1,3)}_n&=
    \frac1{2n}\rbr{\brbr{V^{[1]}_0 }^2 + \brbr{V^{[1]}_1 }^2}
      - \int^{t^{2n}_1}_0 \brbr{V^{[1]}_t}^2 dt
      - \int^1_{t^{2n}_{2n-1}} \brbr{V^{[1]}_t}^2 dt.
  \end{align*}
  By Lemma \ref{230924.1800} (i),
  we can show 
  $\check\convDiff^{(1,3)}_n=O_M(n^{-1-\beta})$ with any $\beta\in(0,H)$.
  
  Writing 
  $a(x):=\rbr{V^{[1]}(x)}^2$,
  $a_t:= \brbr{V^{[1]}_t}^2$ and
  $a^{[1]}_t:=a'(X_t)$, %=\frac{da}{dx}\big|_{x=X_t}$,
  we decompose $\check\convDiff^{(1,1)}_n$ as 
  \begin{align*}
    \check\convDiff^{(1,1)}_n&=
    % n^{-1} \sum_{j=1}^{n-1} \brbr{V^{[1]}_\tj }^2
    % - \int^{t^{2n}_{2n-1}}_{t^{2n}_1} \brbr{V^{[1]}_t}^2 dt
    % =
    n^{-1} \sum_{j=1}^{n-1} a_{t^n_j}
    - \int^{t^{2n}_{2n-1}}_{t^{2n}_1} a_t dt
    \\&=
    \sum_{j=1}^{n-1} \rbr{
    \int^\ttj_\ttjm (a_\ttj - a_t) dt
    +
    % \sum_{j=1}^{n-1} 
    \int^\ttjp_\ttj (a_\ttj - a_t) dt}
    %=-
    %\sum_{j=1}^n \int^\tj_\tjm (a_t - a_\tjm) dt
    \\&=
    \sum_{j=1}^{n-1}\Bigg\{ 
      a^{[1]}_\ttj V^{[1]}_\ttj 
    \int^\ttj_\ttjm \int^\ttj_t dB_s dt
    +%\\&\hspsm+
    \underbrace{
    \int^\ttj_\ttjm \int^\ttj_t 
    \rbr{a^{[1]}_s V^{[1]}_s - a^{[1]}_\ttj V^{[1]}_\ttj} dB_s dt }
    _{=\hat O_M(n^{-2H-1})}
    \\&\hspace{40pt}+
    \underbrace{a^{[1]}_\ttj V^{[2]}_\ttj 
    \int^\ttj_\ttjm \int^\ttj_t ds dt}_{\clubsuit}
    +
    \underbrace{
    \int^\ttj_\ttjm \int^\ttj_t 
    \rbr{a^{[1]}_s V^{[2]}_s - a^{[1]}_\ttj V^{[2]}_\ttj} ds dt }
    _{=\hat O_M(n^{-2-H})}\Bigg\}
    \\&\hspsm
    -\sum_{j=1}^{n-1} \Bigg\{
    a^{[1]}_\ttj V^{[1]}_\ttj 
    \int^\ttjp_\ttj \int^t_\ttj dB_s dt
    +%\\&\hspsm+
    \underbrace{\int^\ttjp_\ttj \int^t_\ttj
    \rbr{a^{[1]}_s V^{[1]}_s - a^{[1]}_\ttj V^{[1]}_\ttj} dB_s dt }
    _{=\hat O_M(n^{-2H-1})}
    \\&\hspace{50pt}
    +\underbrace{a^{[1]}_\ttj V^{[2]}_\ttj 
    \int^\ttjp_\ttj \int^t_\ttj ds dt}_{\clubsuit}
    +%\\&\hspsm+
    \underbrace{\int^\ttjp_\ttj \int^t_\ttj
    \rbr{a^{[1]}_s V^{[2]}_s - a^{[1]}_\ttj V^{[2]}_\ttj} ds dt }
    _{=\hat O_M(n^{-2-H})}\Bigg\}
    \\&=
    \sum_{j=1}^{n-1} a^{[1]}_\ttj V^{[1]}_\ttj \cbr{
    \int^\ttj_\ttjm \int^\ttj_t dB_s dt
    - 
    \int^\ttjp_\ttj \int^t_\ttj dB_s dt}
    +\hat O_M(n^{-2H}),
    % \\&=
    % \underbrace{
    % \sum_{j=1}^{n-1} a^{[1]}_\ttj V^{[1]}_\ttj 
    % \rbr{\int^\ttj_\ttjm \int^s_\ttjm dt dB_s 
    % -
    % \int^\ttjp_\ttj \int^\ttjp_s dt dB_s}}
    % _{{\mygreen=:c_{2,H}^{-1}\convDiff^{(1,1)}_n} }
    % +
    % O_M(n^{-2H})
    \\&=:
    c_{2,H}^{-1}\convDiff^{(1,1)}_n+\hat O_M(n^{-2H})
  \end{align*}
  where we have used the change of variables of Young integral and
  Lemma \ref{230924.1800} (i).
  We can write
  % $\check\convDiff^{(1,1)}_n=
  % c_{2,H}^{-1}\convDiff^{(1,1)}_n+O_M(n^{-2H})$ with the following definition of $\convDiff^{(1,1)}_n$
  \begin{align*}
    \convDiff^{(1,1)}_n &= 
    c_{2,H}\;
    \sum_{j=1}^{n-1} a^{[1]}_{t^n_j} V^{[1]}_{t^n_j}
    \rbr{\int^\ttj_\ttjm (s-\ttjm) dB_s 
    -\int^\ttjp_\ttj (\ttjp-s) dB_s}.
  % \\&=
  % c_{2,H}\;
  % \sum_{j=1}^{n-1} a^{[1]}_{t^n_j} V^{[1]}_{t^n_j}
  % \rbr{\int^\ttj_\ttjm \int^s_\ttjm dt dB_s 
  % -\int^\ttjp_\ttj \int^\ttjp_s dt dB_s}.
  \end{align*}

  Define 
  $\bbh^{(1) n}_{j}$, $\bbh^{(2) n}_{j}$ and $\bbh^{n}_{j}$  for $j\in[n-1]$ by
  \begin{align*}
    \bbh^{(1) n}_{j} = \rbr{\bbone^{2n}_{2j}(s)\, n(s-\ttjm)}_{s\in[0,1]}, \qquad
    \bbh^{(2) n}_{j} = \rbr{\bbone^{2n}_{2j-1}(s)\, n(\ttjm-s)}_{s\in[0,1]} \tand
    \bbh^{n}_{j} = \bbh^{(1) n}_{j} - \bbh^{(2) n}_{j+1}.
  \end{align*}
  Then we have
  \begin{align}\label{eq:230809.1743}
    \convDiff^{(1,1)}_n&=
    c_{2,H}\,n^{-1}
    \cbr{\sum_{j=1}^{n-1} a^{[1]}_{t^n_j} V^{[1]}_{t^n_j}
    B(\bbh^{(1) n}_{j})
    -
    % c_{2,H}\,n^{-1}
    \sum_{j=1}^{n-1} a^{[1]}_{t^n_j} V^{[1]}_{t^n_j}
    B(\bbh^{(2) n}_{j+1})}
    =%\\&=
    c_{2,H}\,n^{-1}
    \sum_{j=1}^{n-1} a^{[1]}_{t^n_j} V^{[1]}_{t^n_j}
    B(\bbh^n_j)
  \end{align}
  and by an argument of the exponent and Proposition \ref{230725.1133},
  we can show 
  $\convDiff^{(1,1)}_n=O_M(n^{-1})$.
  % {\myred [ToDelete??] Remark: by extending the theory of exponent itself, we can prove 
  % $\convDiff^{(1,1)}_n=O_M(n^{-\half-H})$.
  % }

  Define 
  \begin{align}\label{eq:230809.1745}
    \convDiff^{(1,2)}_n=
    c_{2,H}\check\convDiff^{(1,2)}_n=
    -\frac{c_{2,H}}{2n}\rbr{\brbr{V^{[1]}_0 }^2 + \brbr{V^{[1]}_1 }^2},
  \end{align}
  and 
  it is obvious that 
  $\convDiff^{(1,2)}_n=O_M(n^{-1})$.
  Thus, we obtain
  \begin{align*}
    \convDiff^{(1)}_n &=
    \sum_{k=1}^3 c_{2,H}\check\convDiff^{(1,k)}_n
    =
    % \convDiff^{(1,1)}_n+\hat O_M(n^{-2H})
    % +\convDiff^{(1,2)}_n
    % +\hat O_M(n^{-1-H})
    % \\&=
    \convDiff^{(1,1)}_n
    +\convDiff^{(1,2)}_n
    +\hat O_M(n^{-2H})
    =O_M(n^{-1}).
  \end{align*}
\end{proof}

The functionals $\convDiff^{(2)}_n$ and $\convDiff^{(3)}_n$
include the factor $\resincj$,
which appeared in the decomposition of $\secDiff{n}{j} X$, and
was defined at \eqref{230925.1518}.
We first decompose this term $\resincj$.
Recall its definition: 
\begin{align}
  \resincj &=
  \rbr{V^{[1]}_\tj - V^{[1]}_\tjm} \Delta^n_{j}B
  + \rbr{V^{[2]}_\tj - V^{[2]}_\tjm} n^{-1}
  +\resint
  \nn\\
  \resint &= 
  \int^\tjp_\tj \rbr{V^{[1]}_t - V^{[1]}_\tj} dB_t + 
  \int^\tjp_\tj \rbr{V^{[2]}_t - V^{[2]}_\tj} dt
  \nn\\&\quad
  -\rbr{
    \int^\tj_\tjm \rbr{V^{[1]}_t - V^{[1]}_\tjm} dB_t + 
    \int^\tj_\tjm \rbr{V^{[2]}_t - V^{[2]}_\tjm} dt
  }.
  \label{230809.1624}
\end{align}
We define $\kerinc^n_j$ and $\kerdec^n_j$ by
\begin{align}\label{eq:230925.1628}
  \kerinc^n_j = \rbr{\bbone^n_j(t)\: n\:(t-t^n_{j-1})}_{t\in[0,1]},\tand
  \kerdec^n_j = (\bbone^n_j(t)\: n\:(t^n_{j}-t))_{t\in[0,1]}.
\end{align}
\begin{lemma}[Decomposition of $\resincj$]\label{230809.1737}
  Let $\resincj$ be as above for $j\in[n-1]$.
  It holds that
  \begin{align*}
    \resincj &=
    V^{[(1;1),1]}_\tjm B(\bbone^n_j)^2
    + V^{[(1;1),2]}_\tjm n^{-1} B(\bbone^n_j)
    + V^{[(2;1),1]}_\tjm B(\bbone^n_j) n^{-1}
    \\&\hspsm
    + \half {V^{[(1;1),1]}_\tj} \rbr{B(\bbone^n_\jp)^2 - B(\bbone^n_j)^2} 
    + V^{[(1;1),2]}_\tj \rbr{B(\kerinc^n_\jp)-B(\kerinc^n_j)} n^{-1}
    \\&\hspsm
    + V^{[(2;1),1]}_\tj \rbr{B(\kerdec^n_\jp)-B(\kerdec^n_j)} n^{-1}
    +O_M(n^{-2})+\hat O_M(n^{-3H})
  % +O_M(n^{(-2)\vee(-3H)})
  \end{align*}
  and $\resincj=O_M(n^{-2H})$
  uniformly in $j$ as $n\to\infty.$
\end{lemma}
\begin{proof}
  Decompose the first stochastic integral of \eqref{230809.1624} as follows:
  \begin{align*}
    \int^\tjp_\tj \rbr{V^{[1]}_t - V^{[1]}_\tj} dB_t
    &=
    % \int^\tjp_\tj 
    % \int^t_\tj \rbr{V^{[1;1]}_s V^{[1]}_s dB_s + V^{[1;1]}_s V^{[2]}_s ds}
    % dB_t
    % \\&=
    \int^\tjp_\tj \int^t_\tj V^{[(1;1),1]}_s dB_s dB_t +
    \int^\tjp_\tj \int^t_\tj  V^{[(1;1),2]}_s ds dB_t
    \\&=
    {V^{[(1;1),1]}_\tj}~\half B(\bbone^n_\jp)^2 +
    \int^\tjp_\tj \int^t_\tj \rbr{V^{[(1;1),1]}_s - V^{[(1;1),1]}_\tj} dB_s dB_t
    \\
    &\hspsm+%&=
    V^{[(1;1),2]}_\tj \int^\tjp_\tj (t-\tj) dB_t
    +
    \int^\tjp_\tj \int^t_\tj  \rbr{V^{[(1;1),2]}_s-V^{[(1;1),2]}_\tj} ds dB_t
    \\
    &=
    \half {V^{[(1;1),1]}_\tj} B(\bbone^n_\jp)^2+
    V^{[(1;1),2]}_\tj 
    B(\kerinc^n_\jp) n^{-1}
    + \hat O_M(n^{-3H}),
  \end{align*}
  where $\kerinc^n_j$ is defined at \eqref{eq:230925.1628}.
  We used Lemma \ref{230924.1800} (i) and 
  the estimate $\hat O_M(n^{-3H})$ holds uniformly in $j$.
  % \redb{ iterated integralに関して補題をreferする．\koko}
  Hence we can write the difference of the first and third integrals of \eqref{230809.1624} 
  as below:
  \begin{align*}
    &\int^\tjp_\tj \rbr{V^{[1]}_t - V^{[1]}_\tj} dB_t - 
    \int^\tj_\tjm \rbr{V^{[1]}_t - V^{[1]}_\tjm} dB_t
    \\&=
    \half {V^{[(1;1),1]}_\tj} B(\bbone^n_\jp)^2 -
    \half {V^{[(1;1),1]}_\tjm} B(\bbone^n_j)^2
    + V^{[(1;1),2]}_\tj B(\kerinc^n_\jp) n^{-1}
    - V^{[(1;1),2]}_\tjm B(\kerinc^n_j) n^{-1}
    + \hat O_M(n^{-3H})
    \\&=
    \half {V^{[(1;1),1]}_\tj} \rbr{B(\bbone^n_\jp)^2 - B(\bbone^n_j)^2}+
    \half \rbr{V^{[(1;1),1]}_\tj - V^{[(1;1),1]}_\tjm} B(\bbone^n_j)^2 
    \\&\hspsm
    + V^{[(1;1),2]}_\tj \rbr{B(\kerinc^n_\jp)-B(\kerinc^n_j)} n^{-1}
    + \rbr{V^{[(1;1),2]}_\tj - V^{[(1;1),2]}_\tjm} B(\kerinc^n_j) n^{-1}
    + \hat O_M(n^{-3H}).
    \\&=
    \half {V^{[(1;1),1]}_\tj} \rbr{B(\bbone^n_\jp)^2 - B(\bbone^n_j)^2}
    + V^{[(1;1),2]}_\tj \rbr{B(\kerinc^n_\jp)-B(\kerinc^n_j)} n^{-1}
    + \hat O_M(n^{-3H}).
  \end{align*}
  
  The second integral of \eqref{230809.1624} decomposes as
  \begin{align*}
    \int^\tjp_\tj \rbr{V^{[2]}_t - V^{[2]}_\tj} dt
    &=
    % \int^\tjp_\tj \rbr{
    % \int^t_\tj \rbr{V^{[2;1]}_s V^{[1]}_s dB_s + V^{[2;1]}_s V^{[2]}_s ds}
    % } dt
    % \\&=
    \int^\tjp_\tj \int^t_\tj V^{[(2;1),1]}_s dB_s dt
    +
    \int^\tjp_\tj \int^t_\tj  V^{[(2;1),2]}_s dsdt
    \\&=
    V^{[(2;1),1]}_\tj \int^\tjp_\tj \int^t_\tj dB_s dt
    + \int^\tjp_\tj \int^t_\tj \rbr{V^{[(2;1),1]}_s-V^{[(2;1),1]}_\tj} dB_s dt
    \\&\quad+
    \int^\tjp_\tj \int^t_\tj  V^{[(2;1),2]}_s dsdt
    \\&=
    % V^{[(2;1),1]}_\tj B(\kerdec^n_\jp) n^{-1}
    % +O_M(n^{-1-2H})+O_M(n^{-2})
    % =
    V^{[(2;1),1]}_\tj B(\kerdec^n_\jp) n^{-1}
    +O_M(n^{-2}),
  \end{align*}
  where $\kerdec^n_j$ is defined at \eqref{eq:230925.1628},
  % with $\kerdec^n_j = (\bbone^n_j(t)\: n\:(\tj-t))_{t\in[0,1]}$,
  and the difference of the second and fourth integrals of \eqref{230809.1624} is
  \begin{align*}
    &\int^\tjp_\tj \rbr{V^{[2]}_t - V^{[2]}_\tj} dt
    - \int^\tj_\tjm \rbr{V^{[2]}_t - V^{[2]}_\tjm} dt
    \\&=
    V^{[(2;1),1]}_\tj B(\kerdec^n_\jp) n^{-1} - 
    V^{[(2;1),1]}_\tjm B(\kerdec^n_j) n^{-1}
    +O_M(n^{-2})
    % +O_M(n^{-2\vee(-1-2H)})
    \\&=
    V^{[(2;1),1]}_\tj \rbr{B(\kerdec^n_\jp)-B(\kerdec^n_j)} n^{-1} +
    \rbr{V^{[(2;1),1]}_\tj-V^{[(2;1),1]}_\tjm} B(\kerdec^n_j) n^{-1}
    +O_M(n^{-2})
    % +O_M(n^{-2\vee(-1-2H)})
    \\&=
    V^{[(2;1),1]}_\tj \rbr{B(\kerdec^n_\jp)-B(\kerdec^n_j)} n^{-1} 
    +O_M(n^{-2})
  \end{align*}
  Therefore, $\resint$ is written as 
  \begin{align*}
    \resint &= 
    \half {V^{[(1;1),1]}_\tj} \rbr{B(\bbone^n_\jp)^2 - B(\bbone^n_j)^2} 
    + V^{[(1;1),2]}_\tj \rbr{B(\kerinc^n_\jp)-B(\kerinc^n_j)} n^{-1}
    \\&\hspsm
    + V^{[(2;1),1]}_\tj \rbr{B(\kerdec^n_\jp)-B(\kerdec^n_j)} n^{-1}
    +O_M(n^{-2})
    +\hat O_M(n^{-3H}).
  \end{align*}

Since 
\begin{align*}
  \rbr{V^{[i]}_\tj - V^{[i]}_\tjm}
  &=
  \int^\tj_\tjm \brbr{V^{[(i;1),1]}_t dB_t + V^{[(i;1),2]}_t dt}
  \\&=
  V^{[(i;1),1]}_\tjm B(\bbone^n_j)
  + \int^\tj_\tjm \rbr{V^{[(i;1),1]}_t - V^{[(i;1),1]}_\tjm} dB_t 
  + V^{[(i;1),2]}_\tjm n^{-1}
  + \int^\tj_\tjm \rbr{V^{[(i;1),2]}_t - V^{[(i;1),2]}_\tjm} dt
  \\&=
  V^{[(i;1),1]}_\tjm B(\bbone^n_j)
  + V^{[(i;1),2]}_\tjm n^{-1}
  +\hat O_M(n^{-2H})
\end{align*}
for $i=1,2$, we obtain 
\begin{align*}
  \resincj &=
  \rbr{V^{[1]}_\tj - V^{[1]}_\tjm} \Delta^n_{j}B
  + \rbr{V^{[2]}_\tj - V^{[2]}_\tjm} n^{-1}
  +\resint
  \\&=
  V^{[(1;1),1]}_\tjm B(\bbone^n_j)^2
  + V^{[(1;1),2]}_\tjm n^{-1} B(\bbone^n_j)
  + V^{[(2;1),1]}_\tjm B(\bbone^n_j) n^{-1}
  \\&\hspsm
  + \half {V^{[(1;1),1]}_\tj} \rbr{B(\bbone^n_\jp)^2 - B(\bbone^n_j)^2} 
  + V^{[(1;1),2]}_\tj \rbr{B(\kerinc^n_\jp)-B(\kerinc^n_j)} n^{-1}
  \\&\hspsm
  + V^{[(2;1),1]}_\tj \rbr{B(\kerdec^n_\jp)-B(\kerdec^n_j)} n^{-1}
  +O_M(n^{-2})+\hat O_M(n^{-3H}).
  % +O_M(n^{(-2)\vee(-3H)})
\end{align*}
Notice that $\resincj=O_M(n^{-2H})$ uniformly in $j$. 
\end{proof}

Before we proceed to the decompositions of 
$\convDiff^{(2)}_n$ and $\convDiff^{(3)}_n$, 
we calculate some inner products used 
in the proof of Lemma \ref{lem:230925.1301}.
Recall that we have defined 
$\diffker^n_j = \bbone^n_\jp - \bbone^n_j$
with
$\bbone^n_j=\bbone_{[t^n_{j-1},t^n_j]}$, and
$c_{2,H}=4-2^{2H}
=\abr{\diffker^1_0, \diffker^1_0}
=\abr{\bbone^1_1-\bbone^1_0, \bbone^1_1-\bbone^1_0}$,
where the inner product $\abr{\cdot,\cdot}$ should be read as the extended one.
(See \eqref{def:231005.2143}.)
Note that we read 
$\bbone^1_0=\bbone_{[-1,0]}$ and 
$\bbone^1_1=\bbone_{[0,1]}$.
\begin{lemma}\label{230809.1430}
  The following relations hold:
  % {\myred for $n\in\bbN$ and $j\in\bbZ$}:
  \begin{itemize}
\item [(i)]
$\abr{\diffker^1_0, \bbone^1_1}=\half c_{2,H}$ and
$\abr{\diffker^1_0, \bbone^1_0}=-\half c_{2,H}$.

\item [(ii)]
$\abr{\bbone^n_j,\bbone^n_j}=
% \abr{\bbone^n_0,\bbone^n_0}=
n^{-2H}\abr{\bbone^1_0,\bbone^1_0}=
n^{-2H}$ .

\item [(iii)]
$\babr{\diffker^n_j,\kerinc^n_\jp-\kerinc^n_j}=\half c_{2,H}\,n^{-2H}$
and 
% \item [(iv)]
$\babr{\diffker^n_j,\kerdec^n_\jp-\kerdec^n_j}=\half c_{2,H}\,n^{-2H}$.
\end{itemize}
\end{lemma}

\begin{proof}
  % $-\half \abr{\diffker^1_0, \diffker^1_0}=\abr{\diffker^1_0, \bbone^1_0}$
  (i) % \item[(i)]
  Since the following equality holds
  \begin{align*}
    \abr{\diffker^1_0, \bbone^1_1}=
    \abr{\bbone^1_1-\bbone^1_0, \bbone^1_1}=
    \abr{\bbone^1_0-\bbone^1_1, \bbone^1_0}=
    -\abr{\bbone^1_1-\bbone^1_0, \bbone^1_0}=
    -\abr{\diffker^1_0, \bbone^1_0},
  \end{align*}
  we have
  $\abr{\diffker^1_0, \bbone^1_1}=
  \half(\abr{\diffker^1_0, \bbone^1_1}-\abr{\diffker^1_0, \bbone^1_0})=
  \half(\abr{\diffker^1_0, \diffker^1_0})=
  \half c_{2,H}$.

  \item [(ii)] By self-similarity, it is obvious.
  
  \item [(iii)] By self-similarity, we have
  $\babr{\diffker^n_j,\kerinc^n_\jp-\kerinc^n_j}=
  n^{-2H}
  \babr{\diffker^1_0,\kerinc^1_1-\kerinc^1_0}$,
  and the following calculation and (i) prove (iii):
  \begin{align*}
    \babr{\diffker^1_0,\kerinc^1_1-\kerinc^1_0}&=
    % \babr{\bbone^1_1-\bbone^1_0,\kerinc^1_1-\kerinc^1_0}
    % \\&=
    \babr{\bbone^1_1,\kerinc^1_1}-
    \babr{\bbone^1_1,\kerinc^1_0}-
    \babr{\bbone^1_0,\kerinc^1_1}+
    \babr{\bbone^1_0,\kerinc^1_0}
    \\&=
    \babr{\bbone^1_1,\kerinc^1_1}-
    \babr{\bbone^1_1,\kerinc^1_0}-
    \babr{\bbone^1_1,\kerdec^1_0}+
    \babr{\bbone^1_1,\kerdec^1_1}
    \\&=
    \babr{\bbone^1_1,\kerinc^1_1+\kerdec^1_1}-
    \babr{\bbone^1_1,\kerinc^1_0+\kerdec^1_0}
    \\&=
    \babr{\bbone^1_1,\bbone^1_1}-
    \babr{\bbone^1_1,\bbone^1_0}
    \\&=
    \babr{\bbone^1_1,\diffker^1_0}.
  \end{align*}
\end{proof}

Recall the definitions of $\convDiff^{(k)}_n$ for $k=0,2,3$:
\begin{align*}
  \convDiff^{(0)}_n &=
  n^{2H-1} \sum_{j=1}^{n-1} 
  \rbr{V^{[1]}_\tj }^2 I_2\rbr{(\diffker^n_j)^{\otimes2}},
  \qquad
  % &%\\
  % \convDiff^{(1)}_n &=
  % c_{2,H}\;n^{-1} \sum_{j=1}^{n-1} \rbr{V^{[1]}_\tj }^2
  % - c_{2,H}\int^1_0 \rbr{V^{[1]}_t}^2 dt,
  % \\
  \convDiff^{(2)}_n =
  2 n^{2H-1} \sum_{j=1}^{n-1} V^{[1]}_\tj 
  B\rbr{\diffker^n_j} \resinc{j},
  \qquad
  \convDiff^{(3)}_n =
  n^{2H-1} \sum_{j=1}^{n-1} \resinc{j}^2.
\end{align*}
\begin{lemma}\label{lem:230925.1301}
  (i) % \item [(i)]
  The following decomposition of $\convDiff^{(2)}_n$ holds:
  $\convDiff^{(2)}_n =
  \convDiff^{(2,1)}_n+\convDiff^{(2,2)}_n
  + O_M(n^{(-\half-H)\vee(H-2)})$,
  where we write 
  \begin{align*}
    \convDiff^{(2,1)}_n &=
    2\;n^{2H-1} \sum_{j=1}^{n-1} V^{[1]}_\tj V^{[(1;1),1]}_\tjm
    I_3\rbr{\diffker^n_j \otimes (\bbone^n_j)^{\otimes2}}
    \\
    \convDiff^{(2,2)}_n &=
    n^{2H-1} \sum_{j=1}^{n-1} V^{[1]}_\tj {V^{[(1;1),1]}_\tj} 
    I_3\brbr{\diffker^n_j\otimes{\bbone^n_\jp}^{\otimes2}}
    -
    n^{2H-1} \sum_{j=1}^{n-1} V^{[1]}_\tj {V^{[(1;1),1]}_\tj} 
    I_3\brbr{\diffker^n_j\otimes{\bbone^n_j}^{\otimes2}}
  \end{align*}  
  and 
  $\convDiff^{(2,1)}_n,\convDiff^{(2,2)}_n=O_M(n^{-1})$.
  In particular, it holds that $\convDiff^{(2)}_n=O_M(n^{-1})$.

  \item[(ii)]
  $\convDiff^{(3)}_n=n^{2H-1} \sum_{j=1}^{n-1} \resinc{j}^2=O_M(n^{-2H})$.

  \item[(iii)]
  The orders of $\convDiff^{(0)}_n$ and $\convDiff_n$
  are estimated as 
  $\convDiff^{(0)}_n=O_M(n^{-\half})$ and
  $\convDiff_n=O_M(n^{-\half})$.
\end{lemma}
\begin{proof}
  (i)
  The functional $\convDiff^{(2)}_n$ decomposes as 
  \begin{align*}
    \half\convDiff^{(2)}_n &=
    \sum_{k=1}^6 \check\convDiff^{(2,k)}_n 
    + O_M(n^{-2+H}) + \hat O_M(n^{-2H})
    % + O_M(n^{(-2H)\vee(-2+H)})
  \end{align*}
  with 
  \begin{align*}
    \check\convDiff^{(2,1)}_n&=
    n^{2H-1} \sum_{j=1}^{n-1} V^{[1]}_\tj B\rbr{\diffker^n_j} 
    V^{[(1;1),1]}_\tjm B(\bbone^n_j)^2
    %\label{221126.1311}
    \\
    \check\convDiff^{(2,2)}_n&=
    n^{2H-1} \sum_{j=1}^{n-1} V^{[1]}_\tj B\rbr{\diffker^n_j} 
    V^{[(1;1),2]}_\tjm n^{-1} B(\bbone^n_j)
    %\label{221126.1312}
    \\
    \check\convDiff^{(2,3)}_n&=
    n^{2H-1} \sum_{j=1}^{n-1} V^{[1]}_\tj B\rbr{\diffker^n_j} 
    V^{[(2;1),1]}_\tjm B(\bbone^n_j) n^{-1}
    %\label{221126.1313}
    \\
    \check\convDiff^{(2,4)}_n&=
    n^{2H-1} \sum_{j=1}^{n-1} V^{[1]}_\tj B\rbr{\diffker^n_j} 
    \half {V^{[(1;1),1]}_\tj} \rbr{B(\bbone^n_\jp)^2 - B(\bbone^n_j)^2} 
    %\label{221126.1314}
    \\
    \check\convDiff^{(2,5)}_n&=
    n^{2H-1} \sum_{j=1}^{n-1} V^{[1]}_\tj B\rbr{\diffker^n_j} 
    V^{[(1;1),2]}_\tj \rbr{B(\kerinc^n_\jp)-B(\kerinc^n_j)} n^{-1}
    %\label{221126.1315}
    \\
    \check\convDiff^{(2,6)}_n&=
    n^{2H-1} \sum_{j=1}^{n-1} V^{[1]}_\tj B\rbr{\diffker^n_j} 
    V^{[(2;1),1]}_\tj \rbr{B(\kerdec^n_\jp)-B(\kerdec^n_j)} n^{-1}
    %\label{221126.1316}
  \end{align*}

  By the product formula, we have
\begin{align*}
  B\rbr{\diffker^n_j} B(\bbone^n_j)^2 &=
  I_1\rbr{\diffker^n_j}
  \rbr{I_2\rbr{(\bbone^n_j)^{\otimes2}} + \abr{\bbone^n_j,\bbone^n_j}}
  =%\\&=
  I_3\rbr{\diffker^n_j \otimes (\bbone^n_j)^{\otimes2}} 
  + 2 \abr{\diffker^n_j, \bbone^n_j} I_1\rbr{\bbone^n_j} 
  + I_1\rbr{\diffker^n_j} \abr{\bbone^n_j,\bbone^n_j},
\end{align*}
and $\check\convDiff^{(2,1)}_n$ expands as
$%\begin{align*}
  \check\convDiff^{(2,1)}_n=
  \sum_{k=1}^3\check\convDiff^{(2,1,k)}_n
$ %\end{align*}
with
\begin{align*}
  \check\convDiff^{(2,1,1)}_n&=
  n^{2H-1} \sum_{j=1}^{n-1} V^{[1]}_\tj V^{[(1;1),1]}_\tjm
  I_3\rbr{\diffker^n_j \otimes (\bbone^n_j)^{\otimes2}} 
  \\ 
  \check\convDiff^{(2,1,2)}_n&=
  2 n^{2H-1} \sum_{j=1}^{n-1} V^{[1]}_\tj V^{[(1;1),1]}_\tjm
  \abr{\diffker^n_j, \bbone^n_j} I_1\rbr{\bbone^n_j} 
  =
  -c_{2,H}\,n^{-1} 
  \sum_{j=1}^{n-1} V^{[1]}_\tj V^{[(1;1),1]}_\tjm
  I_1\rbr{\bbone^n_j} 
  \\
  \check\convDiff^{(2,1,3)}_n&=
  n^{2H-1} \sum_{j=1}^{n-1} V^{[1]}_\tj V^{[(1;1),1]}_\tjm
  I_1\rbr{\diffker^n_j} \abr{\bbone^n_j,\bbone^n_j}
  =
  n^{-1} \sum_{j=1}^{n-1} V^{[1]}_\tj V^{[(1;1),1]}_\tjm
  I_1\rbr{\diffker^n_j} 
\end{align*}
% where 
% $\abr{\diffker^1_0, \bbone^1_0}=-\half c_{2,H}$,
% $c_{[1,2]}=\abr{\diffker^1_0, \bbone^1_0}$,
% % $c_{[1,2]}=\abr{\bbone_{[0,1]}-\bbone_{[-1,0]}, \bbone_{[-1,0]}}$ and
since 
$\abr{\diffker^n_j, \bbone^n_j}=
n^{-2H}\abr{\diffker^1_0, \bbone^1_0}=
-\half c_{2,H}n^{-2H}$
and
$\abr{\bbone^n_j,\bbone^n_j}=
% \abr{\bbone^n_0,\bbone^n_0}=
n^{-2H}\abr{\bbone^1_0,\bbone^1_0}=
n^{-2H}$.
The rescaled functionals 
$n^{1-2H}\check\convDiff^{(2,1,1)}_n$,
$n^{1}\check\convDiff^{(2,1,2)}_n$ and
$n^{1}\check\convDiff^{(2,1,3)}_n$ 
correspond to the weighted graphs
\eqref{fig:230721.1712},
\eqref{fig:230721.1713} and 
\eqref{fig:230721.1714},
and whose exponents are
$-2H$, $0$ and $\half-H$,
respectively.
Hence, we obtain the estimates
\begin{align*}
  \check\convDiff^{(2,1,1)}_n=O_M(n^{-1}),\quad
  \check\convDiff^{(2,1,2)}_n=O_M(n^{-1})\tand
  \check\convDiff^{(2,1,3)}_n=O_M(n^{-\half-H}).
\end{align*}

Since 
\begin{align*}
  B(\bbone^n_\jp)^2 - B(\bbone^n_j)^2
  =
  I_2\nrbr{{\bbone^n_\jp}^{\otimes2}} 
  - I_2\nrbr{{\bbone^n_j}^{\otimes2}},
\end{align*}
by the product formula,
the functional $\check\convDiff^{(2,4)}_n$ decomposes as
$\check\convDiff^{(2,4)}_n=
\sum_{k=1}^4\check\convDiff^{(2,4,k)}_n$
with
\begin{align*}
  \check\convDiff^{(2,4,1)}_n&=
  \half 
  n^{2H-1} \sum_{j=1}^{n-1} V^{[1]}_\tj {V^{[(1;1),1]}_\tj} 
  I_3\nrbr{\diffker^n_j\otimes{\bbone^n_\jp}^{\otimes2}}
  \\
  \check\convDiff^{(2,4,2)}_n&=
  \half 
  n^{2H-1} \sum_{j=1}^{n-1} V^{[1]}_\tj {V^{[(1;1),1]}_\tj} 
  2I_1\rbr{{\bbone^n_\jp}} \abr{\diffker^n_j,\bbone^n_\jp}
  =
  \half c_{2,H}
  n^{-1} \sum_{j=1}^{n-1} V^{[1]}_\tj {V^{[(1;1),1]}_\tj} 
  I_1\rbr{{\bbone^n_\jp}}  
  \\
  \check\convDiff^{(2,4,3)}_n&=
  -\half 
  n^{2H-1} \sum_{j=1}^{n-1} V^{[1]}_\tj {V^{[(1;1),1]}_\tj} 
  I_3\nrbr{\diffker^n_j\otimes{\bbone^n_j}^{\otimes2}}
  \\
  \check\convDiff^{(2,4,4)}_n&=
  -\half 
  n^{2H-1} \sum_{j=1}^{n-1} V^{[1]}_\tj {V^{[(1;1),1]}_\tj} 
  2I_1\rbr{{\bbone^n_j}} \abr{\diffker^n_j,\bbone^n_j}
  =
  \half c_{2,H}
  n^{-1} \sum_{j=1}^{n-1} V^{[1]}_\tj {V^{[(1;1),1]}_\tj} 
  I_1\rbr{{\bbone^n_j}}
\end{align*}
since
$\abr{\diffker^n_j,\bbone^n_{j+1}}=
% \abr{\diffker^n_0,\bbone^n_{1}}=
n^{-2H}\abr{\diffker^1_0, \bbone^1_1}=
\half c_{2,H}n^{-2H}$
and 
$\abr{\diffker^n_j,\bbone^n_j}=
n^{-2H}\abr{\diffker^1_0, \bbone^1_0}=
-\half c_{2,H}n^{-2H}$.
The rescaled functionals 
$n^{1-2H}\check\convDiff^{(2,4,1)}_n$,
$n^{1}\check\convDiff^{(2,4,2)}_n$,
$n^{1-2H}\check\convDiff^{(2,4,3)}_n$and
$n^{1}\check\convDiff^{(2,4,4)}_n$
correspond to
\eqref{fig:230721.1712},
\eqref{fig:230721.1713},
\eqref{fig:230721.1712} and
\eqref{fig:230721.1713}.
Thus we have 
$\check\convDiff^{(2,4,k)}_n=O_M(n^{-1})$ for $k\in[4]$.
Notice that 
\begin{align*}
  \check\convDiff^{(2,4,2)}_n+\check\convDiff^{(2,4,4)}_n&=
  \half c_{2,H}
  n^{-1} \sum_{j=1}^{n-1} V^{[1]}_\tj {V^{[(1;1),1]}_\tj} 
  I_1\rbr{{\bbone^n_\jp}+{\bbone^n_j}}  
  =% \\&=
  % \half c_{2,H}
  % n^{-1} \sum_{j=1}^{n-1} V^{[1]}_\tj {V^{[(1;1),1]}_\tj} 
  % I_1\rbr{\diffker^n_j+2{\bbone^n_j}}  
  % \\&=
  c_{2,H}
  n^{-1} \sum_{j=1}^{n-1} V^{[1]}_\tj {V^{[(1;1),1]}_\tj} 
  I_1\rbr{\half\diffker^n_j+{\bbone^n_j}}  
\end{align*}
Hence we have the following estimate:
\begin{align*}
  &
  \check\convDiff^{(2,1,2)}_n+
  \check\convDiff^{(2,4,2)}_n+
  \check\convDiff^{(2,4,4)}_n
  % \\&=
  % -c_{2,H}\,n^{-1} 
  % \sum_{j=1}^{n-1} V^{[1]}_\tj V^{[(1;1),1]}_\tjm
  % I_1\rbr{\bbone^n_j} 
  % +
  % \half c_{2,H}
  % n^{-1} \sum_{j=1}^{n-1} V^{[1]}_\tj {V^{[(1;1),1]}_\tj} 
  % I_1\rbr{{\bbone^n_\jp}}
  % +
  % \half c_{2,H}
  % n^{-1} \sum_{j=1}^{n-1} V^{[1]}_\tj {V^{[(1;1),1]}_\tj} 
  % I_1\rbr{{\bbone^n_j}}
  \\&=
  c_{2,H}\,n^{-1}\cbr{
  -
  \sum_{j=1}^{n-1} V^{[1]}_\tj V^{[(1;1),1]}_\tjm
  I_1\rbr{\bbone^n_j} 
  +
  \sum_{j=1}^{n-1} V^{[1]}_\tj {V^{[(1;1),1]}_\tj} 
  I_1\rbr{\half\diffker^n_j+{\bbone^n_j}}
  }
  \\&=
  c_{2,H}\,n^{-1}\cbr{
  \sum_{j=1}^{n-1} V^{[1]}_\tj \rbr{V^{[(1;1),1]}_\tj - V^{[(1;1),1]}_\tjm} 
  I_1\rbr{\bbone^n_j}
  +
  \half
  \sum_{j=1}^{n-1} V^{[1]}_\tj {V^{[(1;1),1]}_\tj} 
  I_1\rbr{\diffker^n_j}
  }
  \\&=
  \hat O_M(n^{-1+1-H-H}) + O_M(n^{-1+(\half-H)})
  =O_M(n^{-\half-H}),
\end{align*}
where we used Lemma \ref{230924.1800} (ii),
and we obtain 
\begin{align*}
  \check\convDiff^{(2,1)}_n+\check\convDiff^{(2,4)}_n=
  \check\convDiff^{(2,1,1)}_n+
  \check\convDiff^{(2,4,1)}_n+
  \check\convDiff^{(2,4,3)}_n+
  O_M(n^{-\half-H}).
\end{align*}

As for $\check\convDiff^{(2,2)}_n$,
we have 
$\check\convDiff^{(2,2)}_n=
\sum_{k=1}^2\check\convDiff^{(2,2,k)}_n$
with
\begin{align*}
  \check\convDiff^{(2,2,1)}_n&=
  n^{2H-2} \sum_{j=1}^{n-1} V^{[1]}_\tj V^{[(1;1),2]}_\tjm  
  I_2(\diffker^n_j \otimes \bbone^n_j) 
  \tand%\\
  \check\convDiff^{(2,2,2)}_n=%&=
  % n^{2H-2} \sum_{j=1}^{n-1} V^{[1]}_\tj V^{[(1;1),2]}_\tjm  
  % \abr{\diffker^n_j, \bbone^n_j}
  % =
  -\half c_{2,H}
  n^{-2} \sum_{j=1}^{n-1} V^{[1]}_\tj V^{[(1;1),2]}_\tjm.
\end{align*}
For $\check\convDiff^{(2,5)}_n$,
we have 
$\check\convDiff^{(2,5)}_n=
\sum_{k=1}^3\check\convDiff^{(2,5,k)}_n$
with
\begin{align*}
  \check\convDiff^{(2,5,1)}_n=&
  n^{2H-2} \sum_{j=1}^{n-1} V^{[1]}_\tj V^{[(1;1),2]}_\tj 
  I_2\brbr{\diffker^n_j\otimes \kerinc^n_\jp},
  \qquad&%\\
  \check\convDiff^{(2,5,2)}_n=&
  - n^{2H-2} \sum_{j=1}^{n-1} V^{[1]}_\tj V^{[(1;1),2]}_\tj 
  I_2\brbr{\diffker^n_j\otimes \kerinc^n_j}
  \\
  \check\convDiff^{(2,5,3)}_n=&
  % n^{2H-2} \sum_{j=1}^{n-1} V^{[1]}_\tj V^{[(1;1),2]}_\tj 
  % \babr{\diffker^n_j,\kerinc^n_\jp-\kerinc^n_j}
  % =
  \half c_{2,H}\,
  n^{-2} \sum_{j=1}^{n-1} V^{[1]}_\tj V^{[(1;1),2]}_\tj. 
\end{align*}
Here we used 
$\babr{\diffker^n_j,\kerinc^n_\jp-\kerinc^n_j}=\half c_{2,H}\,n^{-2H}$
from Lemma \ref{230809.1430}.
By an argument of exponent, we have 
$\check\convDiff^{(2,2,1)}_n,
\check\convDiff^{(2,5,1)}_n, 
\check\convDiff^{(2,5,2)}_n=O_M(n^{H-2})$.
By Lemma \ref{230924.1800} (ii), we have
\begin{align*}
  \check\convDiff^{(2,2,2)}_n+
  \check\convDiff^{(2,5,3)}_n=
  % -\half c_{2,H}
  % n^{-2} \sum_{j=1}^{n-1} V^{[1]}_\tj V^{[(1;1),2]}_\tjm
  % +
  % \half c_{2,H}\,
  % n^{-2} \sum_{j=1}^{n-1} V^{[1]}_\tj V^{[(1;1),2]}_\tj 
  % =
  \half c_{2,H} n^{-2}
  \sum_{j=1}^{n-1} V^{[1]}_\tj 
  \brbr{V^{[(1;1),2]}_\tj - V^{[(1;1),2]}_\tjm}
  % =\hat O_M(n^{-2+1-H})
  =\hat O_M(n^{-1-H}),
\end{align*}
and we obtain 
$\check\convDiff^{(2,2)}_n+\check\convDiff^{(2,5)}_n
=O_M(n^{H-2}) + \hat O_M(n^{-1-H})
=O_M(n^{H-2})$.

Similarly, we have 
$\check\convDiff^{(2,3)}_n=
\sum_{k=1}^2\check\convDiff^{(2,3,k)}_n$
with
\begin{align*}
  \check\convDiff^{(2,3,1)}_n=
  n^{2H-2} \sum_{j=1}^{n-1} V^{[1]}_\tj V^{[(2;1),1]}_\tjm  
  I_2(\diffker^n_j \otimes \bbone^n_j)
  \tand
  \check\convDiff^{(2,3,2)}_n=
  -\half c_{2,H}
  n^{-2} \sum_{j=1}^{n-1} V^{[1]}_\tj V^{[(2;1),1]}_\tjm,
\end{align*}
and 
$\check\convDiff^{(2,6)}_n=
\sum_{k=1}^3\check\convDiff^{(2,6,k)}_n$
with
\begin{align}
  \check\convDiff^{(2,6,1)}_n&=
  n^{2H-2} \sum_{j=1}^{n-1} V^{[1]}_\tj V^{[(2;1),1]}_\tj 
  I_2\brbr{\diffker^n_j\otimes \kerdec^n_\jp},
  \qquad&%\nn\\
  \check\convDiff^{(2,6,2)}_n&=
  -
  n^{2H-2} \sum_{j=1}^{n-1} V^{[1]}_\tj V^{[(2;1),1]}_\tj 
  I_2\brbr{\diffker^n_j\otimes \kerdec^n_j}
  \nn\\
  \check\convDiff^{(2,6,3)}_n&=
  % n^{2H-2} \sum_{j=1}^{n-1} V^{[1]}_\tj V^{[(2;1),1]}_\tj 
  % \babr{\diffker^n_j,\kerdec^n_\jp-\kerdec^n_j}
  % =
  \half c_{2,H}
  n^{-2} \sum_{j=1}^{n-1} V^{[1]}_\tj V^{[(2;1),1]}_\tj,
\end{align}
where we used 
$\abr{\diffker^n_j,\kerdec^n_\jp-\kerdec^n_j}=\half c_{2,H}\,n^{-2H}$
from Lemma \ref{230809.1430}.
Thus, by the same argument as 
$\check\convDiff^{(2,2)}_n+\check\convDiff^{(2,5)}_n$,
we obtain 
$\check\convDiff^{(2,3)}_n+\check\convDiff^{(2,6)}_n
=O_M(n^{H-2}) + \hat O_M(n^{-1-H})
=O_M(n^{H-2})$.

Therefore, we obtain 
\begin{align*}
  \convDiff^{(2)}_n &=
  2\sum_{k=1}^6 \check\convDiff^{(2,k)}_n 
  + O_M(n^{-2+H}) + \hat O_M(n^{-2H})
  % + O_M(n^{(-2H)\vee(-2+H)})
  % \\&=
  % \check\convDiff^{(2,1,1)}_n+
  % \check\convDiff^{(2,4,1)}_n+
  % \check\convDiff^{(2,4,3)}_n
  % +\underbrace{O_M(n^{-\half-H})}_{
  %   \check\convDiff^{(2,1,3)}_n
  % }
  % +\underbrace{O_M(n^{-\half-H})}_{
  %   \check\convDiff^{(2,1,2)}_n+
  %   \check\convDiff^{(2,4,2)}_n+
  %   \check\convDiff^{(2,4,4)}_n
  % }
  % +O_M(n^{H-2}) + \hat O_M(n^{-1-H})
  % + O_M(n^{(-2H)\vee(-2+H)})
  =% \\&=
  \convDiff^{(2,1)}_n+\convDiff^{(2,2)}_n
  +O_M(n^{(-\half-H)\vee(H-2)}),
  % +O_M(n^{-\half-H})+O_M(n^{H-2})
\end{align*}
where we define 
$\convDiff^{(2,1)}_n=2\check\convDiff^{(2,1,1)}_n$
and 
$\convDiff^{(2,2)}_n=
2(\check\convDiff^{(2,4,1)}_n+\check\convDiff^{(2,4,3)}_n)$.

\item[(ii)] 
It is straightforward since we have $\resincj=O_M(n^{-2H})$ 
from Lemma \ref{230809.1737}.

\item [(iii)]
The functional $\convDiff_n$ decomposes as 
$\convDiff_n=\sum_{k=0}^3 \convDiff^{(k)}_n$
with
\begin{align*}
  \convDiff^{(0)}_n &=
  n^{2H-1} \sum_{j=1}^{n-1} 
  \rbr{V^{[1]}_\tj }^2 I_2\rbr{(\diffker^n_j)^{\otimes2}}.
\end{align*}
We have $\convDiff^{(0)}_n=O_M(n^{-\half})$,
since the rescaled functional $n^{-2H+1}\convDiff^{(0)}_n$ 
corresponds to the weighted graph \eqref{fig:230721.1711} 
in Section \ref{230925.1250},
whose exponent is $\half-2H$.
By Lemma \ref{lem:230720.1911}, and (i) and (ii) of this lemma,
we obtain 
$\convDiff_n=O_M(n^{-\half})$.
\end{proof}

\subsubsection{About the truncation functional $\psi_n$}
Recall that we have defined the truncation functional $\psi_n$ at \eqref{eq:230926.1352} as
\begin{align*}
  \psi_n = 
  \psi\rbr{ \frac{R_n}{\eta_0\iv}}\psi\rbr{\frac{R_{2n}}{\eta_0\iv}}.
\end{align*}
with some smooth function $\psi:\bbR\to[0,1]$ satisfying 
$\psi(x)=1$ on $\cbr{x\mid\abs{x}\leq\half}$ and 
$\psi(x)=0$ on $\cbr{x\mid\abs{x}\geq1}$,
and 
the constant $\eta_0\in\bbR_{>0}$
is taken small enough so that $\hat H_n\in(0,1)$ when $\psi_n>0$.
Note that $\convDiff_n$ is defined at \eqref{230804.1403} 
and by Lemma \ref{lem:230925.1301} (iii) we have
$\convDiff_n=O_M(n^{-\half})$.

We also write 
% $\psi_n^{(1)} = 
% \psi\rbr{\eta_0^{-1} \frac{\convDiff_n}{\iv}}$ and
$\psi_n^{(i)} = 
\psi\rbr{\eta_0^{-1} \frac{\convDiff_{in}}{\iv}}$ ($i=1,2$)
for notational convenience.
The following estimate is basic to Lemma \ref{230510.1711}. 
\begin{lemma}\label{230509.2330}
  For any $L>1$, as $n\to\infty$,
  \begin{align*}
    P\sbr{\eta_0^{-1} \abs{\frac{\convDiff_n}{\iv}}\geq\half} = O(n^{-L}).
  \end{align*}
\end{lemma}
\begin{proof}
  For any $L>1$,
  \begin{align*}
    P\sbr{\eta_0^{-1} \abs{\frac{\convDiff_n}{\iv}}\geq\half} &%=
    % P\sbr{2\eta_0^{-1} \abs{\frac{\convDiff_n}{\iv}}\geq1} = 
    % E\sbr{\bbone_\cbr{2\eta_0^{-1} \abs{\frac{\convDiff_n}{\iv}}\geq1}}
    \leq%\\&\leq
    E\sbr{\rbr{2\eta_0^{-1} \abs{\frac{\convDiff_n}{\iv}}}^L}
    \leq
    \rbr{2\eta_0^{-1}}^L \;
    \norm{\convDiff_n}_{2L}^L\; \norm{\bbV_\infty^{-1}}_{2L}^L
    =O(n^{-\frac{L}{2}}),
  \end{align*}
  where we used 
  $\convDiff_n=O_M(n^{-\half})$ (Lemma \ref{lem:230925.1301} (iii)) and
  $\bbV_\infty^{-1}\in L^{\infty-}$.
\end{proof}

The truncation functional $\psi_n$ satisfies the properties in the next lemma.
\begin{lemma}\label{230510.1711}
  The following estimates on $\psi_n$ hold:
  \item[(i)]
  $P\sbr{\psi_n<1}=O(n^{-L})$
  for any $L>0$.

  \item[(ii)] 
  % \begin{align}%\label{230510.1200}
  $\norm{\psi_n-1}_{k,p} = O(n^{-L})$
  % \end{align}
  for any $k\in\bbZ_{\geq0}$, $p>1$ and $L>0$.
\end{lemma}
\begin{proof}
  We write $\psiKer_n = \eta_0^{-1}R_n/\iv$ for notational convenience.
  Notice that $\psiKer_n=O_M(n^{-\half})$.
  Then we have $\psi_n^{(i)}=\psi(\psiKer_{in})$ for $i=1,2$.
  
  We can prove the estimate (i) by
  \begin{align}
    P\sbr{\psi_n<1} = 
    P\bsbr{\psi_n^{(1)}\psi_n^{(2)}<1} %=
    %P\sbr{\psi_n^{(1)}\neq1 \text{ or } \psi_n^{(2)}\neq1}
    &\leq
    P\bsbr{\psi_n^{(1)}<1} + P\bsbr{\psi_n^{(2)}<1}
    \nn\\&\leq
    P\Bsbr{\abs{\psiKer_n}\geq\half}
    %P\sbr{\abs{\eta_0^{-1}\frac{\convDiff_n}{\iv}}\geq\half}
    +P\Bsbr{\abs{\psiKer_{2n}}\geq\half}
    %P\sbr{\abs{\eta_0^{-1}\frac{\convDiff_{2n}}{\iv}}\geq\half}
    =O(n^{-L}),
    \label{230704.1630}
  \end{align}
  where the last estimate follows from Lemma \ref{230509.2330}. 

  \item[(ii)]
  The estimate $\norm{\psi_n-1}_{p} = O(n^{-L})$ for any $L>0$ 
  follows from 
  \begin{align*}
    \norm{\psi_n-1}_{p} &\leq
    E\sbr{\bbone_\cbr{\psi_n<1}}^{1/p} =
    P\sbr{\psi_n<1}^{1/p}
  \end{align*}
  and (i).

  For $k\geq1$, the Malliavin derivative $D^k \psi_n^{(1)}$ is decomposed as 
  \begin{align*}
    D^k \psi_n^{(1)}&= 
    \sum_{\ell, k_1,...,k_\ell} 
    c_{\ell,k_{.}} 
    \psi^{(\ell)}\rbr{\psiKer_n}\;
    %\psi^{(\ell)}\rbr{\eta_0^{-1} \frac{\convDiff_n}{\iv}}\;
    D^{k_1}\psiKer_n
    %D^{k_1}\rbr{\eta_0^{-1} \frac{\convDiff_n}{\iv}} 
    \tilde\otimes \cdots \tilde\otimes
    D^{k_\ell}\psiKer_n,
    %D^{k_\ell}\rbr{\eta_0^{-1} \frac{\convDiff_n}{\iv}}
  \end{align*}
  where $\ell=1,...,k$ and $k=k_1+\cdots+k_\ell$ with $k_i\geq1$.
  By Lemma \ref{230509.2330}, for any $p>1$ and $L>0$ we have
  \begin{align*}
    \norm{\psi^{(\ell)}\rbr{\psiKer_n}}_p &= 
    %\norm{\psi^{(\ell)}\rbr{\eta_0^{-1} \frac{\convDiff_n}{\iv}}}_p &= 
    E\sbr{\abs{\psi^{(\ell)}\rbr{\psiKer_n}}^p}^{\frac1p} 
    %E\sbr{\abs{\psi^{(\ell)}\rbr{\eta_0^{-1} \frac{\convDiff_n}{\iv}}}^p}^{\frac1p} 
    \\&\leq
    c_\ell
    E\sbr{\rbr{\bbone_\cbr{\abs{\psiKer_n}\geq\half}}^p}^{\frac1p} 
    %E\sbr{\bbone_\cbr{\abs{\eta_0^{-1} \frac{\convDiff_n}{\iv}}\geq\half}}^{\frac1p} 
    = 
    c_\ell
    P\sbr{\abs{\psiKer_n}\geq\half}^{\frac1p} 
    = O(n^{-\frac{L}{p}})
  \end{align*}
  where $c_\ell=\sup_{x\in\bbR} \abs{\psi^{(\ell)}(x)}<\infty$. 
  Hence 
  \begin{align*}
    \norm{\bnorm{D^k \psi_n^{(1)}}_{\calh^{\otimes k}}}_p& \leq
    \sum_{\ell, k_1,...,k_\ell} 
    c_{\ell,k_{.}} 
    \norm{\psi^{(\ell)}\rbr{\psiKer_n}}_{(\ell+1)p}\;
    \Bnorm{\norm{D^{k_1}\psiKer_n}_{\calh^{\otimes k_1}}}_{(\ell+1)p}
    \times \cdots \times
    \Bnorm{\norm{D^{k_\ell}\psiKer_n}_{\calh^{\otimes k_\ell}}}_{(\ell+1)p}
    \\&\leq
    \sum_{\ell, k_1,...,k_\ell} 
    O(n^{-\frac{L}{(\ell+1)p} - \half\ell})
    =O(n^{-L'})
  \end{align*}
  for any $L'>0$ and $p>1$. 
\end{proof}

\subsubsection{Decomposition of $Z_n$}
The main part of 
$\hat H_n =
0\vee\rbr{\half + \frac1{2\log2} \log\frac{\secVar_n}{\secVar_{2n}}}\wedge1$
is written as the difference of two $\log$'s:
\begin{align}\label{eq:230804.2030}
  \log\frac{\secVar_n}{\secVar_{2n}} &= 
  \log\rbr{2^{2H-1} \frac{\rv_n}{\rv_{2n}}}
  =
  \rbr{H-\half}2\log2 +
  \log\rbr{\frac{\convDiff_n}{\iv}+1}-\log\rbr{\frac{\convDiff_{2n}}{\iv}+1}.
\end{align}
By Taylor expansion, 
$\log\rbr{\frac{\convDiff_n}{\iv}+1}$ decomposes as
\begin{align}
  \label{230423.1826}
  \log\rbr{\frac{\convDiff_n}\iv+1}
  % =
  % {\frac{\convDiff_n}\iv -\half\rbr{\frac{\convDiff_n}\iv}^2}
  % + \reslogone\rbr{\frac{\convDiff_n}\iv+1},
  =
  {\frac{\convDiff_n}\iv -\half\rbr{\frac{\convDiff_n}\iv}^2}
  + r\rbr{\frac{\convDiff_n}\iv},
\end{align}
where $r(x)$ is defined by 
% where $R^{(3)}_{\log,1}(x)$ is defined by 
\begin{align*}
  r(x) &= 
  \int^{x+1}_1 \int^{t_1}_1 \int^{t_2}_1
  \log^{(3)}(t_3) dt_3 dt_2 dt_1
  = 
  \int^{x+1}_1 \int^{t_1}_1 \int^{t_2}_1
  2(t_3)^{-3} dt_3 dt_2 dt_1
\end{align*}
for $x>-1$, where 
$\log^{(3)}$ is the third derivative of $\log$.
We denote by $r^{(k)}(x)$ the $k$-th derivative of $r(x)$.
First let us prove the next elementary lemma.
\begin{lemma}
  Let $\eta\in(0,1)$.
  For $\abs{x}<\eta$, the following estimates hold:
  \begin{align*}
    % \abs{r(x)}&\leq 2(1-\eta)^{-3}\;\frac1{3!}\abs{x}^3
    % \\ 
    % \abs{r'(x)}&\leq 2(1-\eta)^{-3}\;\frac1{2!}\abs{x}^2
    % \\
    % \abs{r''(x)}&\leq 2(1-\eta)^{-3}\;\abs{x}
    \babs{r^{(k)}(x)}&\leq 2(1-\eta)^{-3}\;\frac1{(3-k)!}\abs{x}^{(3-k)}
    && \text{for } k=0,1,2.
    \\
    \babs{r^{(k)}(x)}&\leq (k-1)! (1-\eta)^{-k}
    && \text{for } k\geq3.
  \end{align*}
\end{lemma}
\begin{proof}
  If $\eta>x\geq0$, then
  \begin{align*}
    \abs{r(x)} &= 
    \int^{x+1}_1 \int^{t_1}_1 \int^{t_2}_1
    2(t_3)^{-3} dt_3 dt_2 dt_1
    \leq%\\&\leq
    2\int^{x+1}_1 \int^{t_1}_1 \int^{t_2}_1 dt_3 dt_2 dt_1
    =%\\&= 
    2\; \frac1{3!}\; x^3.\hspace{70pt}
  \end{align*}
  For $-\eta<x<0$, 
  \begin{align*}
    \abs{r(x)} &= 
    \int^{1}_{x+1} \int^{1}_{t_1} \int^1_{t_2}
    2(t_3)^{-3} dt_3 dt_2 dt_1
    \leq%\\&\leq
    2(1-\eta)^{-3}
    \int^{1}_{x+1} \int^{1}_{t_1} \int^1_{t_2} dt_3 dt_2 dt_1
    =%\\&= 
    2(1-\eta)^{-3}\; \frac1{3!}\; \abs{x}^3.
  \end{align*}
  Hence in the both case, we obtain 
  \begin{align*}
    \abs{r(x)} &\leq
    2(1-\eta)^{-3}\; \frac1{3!}\; \abs{x}^3.
  \end{align*}
  Similar argument stands for $r^{(k)}(x)$ ($k\geq1$).
\end{proof}
The following lemma shows that 
the residual term of the Taylor expansion \eqref{230423.1826} 
% of $\log\rbr{{\convDiff_n}/\iv+1}$
is of $O_M(n^{-\frac32})$ under the truncation by $\psi_n$.
\begin{lemma}\label{230510.1721}
  As $n\to\infty$, the following estimates hold:
  \begin{align*}
    r\rbr{\frac{\convDiff_n}{\iv}} \psi_n 
    &= O_M(n^{-\frac32})
    \tand%\\
    r\rbr{\frac{\convDiff_{2n}}{\iv}} \psi_n 
    = O_M(n^{-\frac32}).
  \end{align*}
\end{lemma}

\begin{proof}
  % \begin{align*}
  %   \psi_n &=
  %   \psi\rbr{\eta_0^{-1} \frac{\convDiff_n}{\iv}} \psi\rbr{\eta_0^{-1} \frac{\convDiff_{2n}}{\iv}}
  %   %\\&= 
  % \end{align*}
  Writing $\widehat\psi(x) = \psi(\eta_0^{-1}x)$,
  we can have
    $\psi_n^{(1)} =
    \psi\rbr{\eta_0^{-1} \frac{\convDiff_n}{\iv}}
    =\widehat\psi\rbr{\frac{\convDiff_n}{\iv}}$,
    $\psi_n^{(2)} = 
    \psi\rbr{\eta_0^{-1} \frac{\convDiff_{2n}}{\iv}}
    =\widehat\psi\rbr{\frac{\convDiff_{2n}}{\iv}}$
  and
  $\psi_n = \psi_n^{(1)} \psi_n^{(2)}$. 
  For notational convenience, we write
  $\phi(x) = r(x) \widehat\psi(x)%= \reslogone(x+1) \widehat\psi(x)
  $ for $x>-1$.
  Since the function $\hat\psi$ satisfies $\hat\psi(x)=0$ for $\abs{x}\geq\eta_0$,
  we extend the domain of $\phi$ to $x\leq-1$ by $\phi(x)=0$.
  We can prove $\phi=r\hat\psi\in C^\infty_{\text{cpt}}(\bbR)$.
  Then we can write 
  \begin{align*}
    r\rbr{\frac{\convDiff_n}{\iv}} \psi_n 
    &= r\rbr{\frac{\convDiff_n}{\iv}} \psi_n^{(1)} \psi_n^{(2)}
    = r\rbr{\frac{\convDiff_n}{\iv} } \widehat\psi\rbr{\frac{\convDiff_n}{\iv}} \psi_n^{(2)}
    =%\\&= 
    \phi\rbr{\frac{\convDiff_n}{\iv}} \psi_n^{(2)}.
  \end{align*}

  We are going to prove 
  $\norm{\phi\rbr{\frac{\convDiff_n}{\iv}}}_{k,p}=O_M(n^{-\frac32})$
  for $k\geq0$ and $p>1$.
  For $\eta_1>0$ satisfying $(0<)\eta_0<\eta_1<1$,
  we have 
  \begin{align*}  
    \abs{\phi\rbr{\frac{\convDiff_n}\iv}}
    =\abs{r\rbr{\frac{\convDiff_n}\iv}}\abs{\hat\psi\rbr{\frac{\convDiff_n}\iv}}
    \leq\underbrace{\sup_x\nabs{\hat\psi(x)}}_{=1}\;
    2(1-\eta_1)^{-3}\; \frac1{3!}\; \abs{\frac{\convDiff_n}\iv}^3\;
    \bbone_{\cbr{\babs{\frac{\convDiff_n}\iv}<\eta_0}}
    \simleq %2(1-\eta_1)^{-3}\; \frac1{3!}
    \abs{\frac{\convDiff_n}\iv}^3,
  \end{align*}
  since $\babs{\frac{\convDiff_n}\iv}<\eta_0$
  if ${\hat\psi\rbr{\frac{\convDiff_n}\iv}}>0$.
  We obtain 
  $\norm{\phi\rbr{\frac{\convDiff_n}\iv}}_p 
  \simleq 
  \norm{\abs{\frac{\convDiff_n}\iv}^3}_p
  % \leq
  % \norm{\convDiff_n}_{6p}^3 \norm{\rbr{\iv}^{-1}}_{6p}^3
  =O(n^{-\frac32})$. 

  We write $F_n = \frac{\convDiff_n}\iv$ for notational convenience.
  For $k\geq1$, the Malliavin derivative $D^k \phi(F_n)$ is decomposed as 
  \begin{align*}
    D^k \phi(F_n)&= 
    \sum_{\ell, k_1,...,k_\ell} 
    c_{\ell,k_{.}} 
    \phi^{(\ell)}\rbr{\psiKer_n}\;
    D^{k_1}\psiKer_n
    \tilde\otimes \cdots \tilde\otimes
    D^{k_\ell}\psiKer_n
    \\&= 
    \sum_{\ell, k_1,...,k_\ell} 
    c_{\ell,k_{.}} 
    \rbr{\sum_{i=0}^\ell c_i\; r^{(i)}(F_n) \hat\psi^{(\ell-i)}(F_n)}\;
    D^{k_1}\psiKer_n
    \tilde\otimes \cdots \tilde\otimes
    D^{k_\ell}\psiKer_n,
  \end{align*}
  where $\ell=1,...,k$ and $k=k_1+\cdots+k_\ell$ with $k_i\geq1$.
  Hence 
  \begin{align*}
    \norm{\norm{D^k \phi(F_n)}_{\calh^{\otimes k}}}_p&\simleq
    \sum_{\ell, k_1,...,k_\ell} \sum_{i=0}^\ell  
    \norm{\babs{r^{(i)}(F_n)}\;\babs{\hat\psi^{(\ell-i)}(F_n)}\;
    \norm{D^{k_1}\psiKer_n}_{\calh^{\otimes k_1}}
    \times \cdots \times
    \norm{D^{k_\ell}\psiKer_n}_{\calh^{\otimes k_\ell}}}_p
    \\&\simleq
    \sum_{\ell, k_1,...,k_\ell} 
    \norm{\babs{r(F_n)}\;\babs{\hat\psi^{(\ell-i)}(F_n)}\;
    \norm{D^{k_1}\psiKer_n}_{\calh^{\otimes k_1}}
    \times \cdots \times
    \norm{D^{k_\ell}\psiKer_n}_{\calh^{\otimes k_\ell}}}_p
    \\&\quad+
    \sum_{\ell, k_1,...,k_\ell}   
    \norm{\babs{r^{(1)}(F_n)}\;\babs{\hat\psi^{(\ell-i)}(F_n)}\;
    \norm{D^{k_1}\psiKer_n}_{\calh^{\otimes k_1}}
    \times \cdots \times
    \norm{D^{k_\ell}\psiKer_n}_{\calh^{\otimes k_\ell}}}_p
    \\&\quad+
    \sum_{\ell\geq2, k_1,...,k_\ell}
    \norm{\babs{r^{(2)}(F_n)}\;\babs{\hat\psi^{(\ell-i)}(F_n)}\;
    \norm{D^{k_1}\psiKer_n}_{\calh^{\otimes k_1}}
    \times \cdots \times
    \norm{D^{k_\ell}\psiKer_n}_{\calh^{\otimes k_\ell}}}_p
    \\&\quad+
    \sum_{\ell\geq3, k_1,...,k_\ell} \sum_{i=3}^\ell  
    \norm{\babs{r^{(i)}(F_n)}\;\babs{\hat\psi^{(\ell-i)}(F_n)}\;
    \norm{D^{k_1}\psiKer_n}_{\calh^{\otimes k_1}}
    \times \cdots \times
    \norm{D^{k_\ell}\psiKer_n}_{\calh^{\otimes k_\ell}}}_p
    \\&\simleq
    \sum_{\ell(\geq1), k_1,...,k_\ell}
    \bnorm{F_n}_{(\ell+3)p}^3
    \Bnorm{\norm{D^{k_1}\psiKer_n}_{\calh^{\otimes k_1}}}_{(\ell+3)p}
    \times \cdots \times
    \Bnorm{\norm{D^{k_\ell}\psiKer_n}_{\calh^{\otimes k_\ell}}}_{(\ell+3)p}
    \\&\quad+
    \sum_{\ell(\geq1), k_1,...,k_\ell}   
    \bnorm{F_n}_{(\ell+2)p}^2
    \Bnorm{\norm{D^{k_1}\psiKer_n}_{\calh^{\otimes k_1}}}_{(\ell+2)p}
    \times \cdots \times
    \Bnorm{\norm{D^{k_\ell}\psiKer_n}_{\calh^{\otimes k_\ell}}}_{(\ell+2)p}
    \\&\quad+
    \sum_{\ell\geq2, k_1,...,k_\ell}
    \bnorm{F_n}_{(\ell+1)p}
    \Bnorm{\norm{D^{k_1}\psiKer_n}_{\calh^{\otimes k_1}}}_{(\ell+1)p}
    \times \cdots \times
    \Bnorm{\norm{D^{k_\ell}\psiKer_n}_{\calh^{\otimes k_\ell}}}_{(\ell+1)p}
    \\&\quad+
    \sum_{\ell\geq3, k_1,...,k_\ell} \sum_{i=3}^\ell  
    \Bnorm{\norm{D^{k_1}\psiKer_n}_{\calh^{\otimes k_1}}}_{\ell p}
    \times \cdots \times
    \Bnorm{\norm{D^{k_\ell}\psiKer_n}_{\calh^{\otimes k_\ell}}}_{\ell p}
    \\&=
    \sum_{\ell(\geq1), k_\cdot} O(n^{n^{-\frac{\ell+3}2}}) + 
    \sum_{\ell(\geq1), k_\cdot} O(n^{n^{-\frac{\ell+2}2}}) + 
    \sum_{\ell\geq2, k_\cdot} O(n^{n^{-\frac{\ell+1}2}}) + 
    \sum_{\ell\geq3, k_\cdot} O(n^{n^{-\frac{\ell}2}})
    \\&=
    O(n^{-\frac32})
  \end{align*}
  Hence we obtain $\phi(F_n)=\phi(\convDiff_n/\iv)=O_M(n^{-\frac32})$

  By the proof of Lemma \ref{230510.1711} (ii), 
   we can see $\psi_n^{(2)}=O_M(1)$ and 
  \begin{align*}
    r\rbr{\frac{\convDiff_n}{\iv}} \psi_n 
    &= 
    \phi\rbr{\frac{\convDiff_n}{\iv}} \psi_n^{(2)}
    = O_M(n^{-\frac32}).
  \end{align*}
  We can also prove 
  \begin{align*}
    r\rbr{\frac{\convDiff_{2n}}{\iv}} \psi_n
    &= 
    \phi\rbr{\frac{\convDiff_{2n}}{\iv}} \psi_n^{(1)}
    = O_M(n^{-\frac32}).
  \end{align*}
\end{proof}

Hence we obtain the decomposition of $Z_n$, 
which is defined at \eqref{eq:230926.1554}.
\begin{lemma}\label{lem:230926.1026}
  The functional $Z_n$ decomposes as
  \begin{align*}
    Z_n &=
    \sqrt{n} 
    \cbr{\frac{\convDiff_n}\iv - \frac{\convDiff_{2n}}\iv}
    -  \sqrt{n} 
    \cbr{\half\rbr{\frac{\convDiff_n}\iv}^2 -\half\rbr{\frac{\convDiff_{2n}}\iv}^2}
    + \negTerm_n^{(1)}
  \end{align*}
  with $\negTerm_n^{(1)}=O_M\rbr{n^{-1}}$.
\end{lemma}
\begin{proof}
For notational convenience, we set 
$\widetilde H_n = 
\rbr{\half + \frac1{2\log2} \log\frac{\secVar_n}{\secVar_{2n}}} \psi_n$.
When $\psi_n>0$, we have 
$\half + \frac1{2\log2} \log\frac{\secVar_n}{\secVar_{2n}}\in(0,1)$, and hence
$\hat H_n = \half + \frac1{2\log2} \log\frac{\secVar_n}{\secVar_{2n}}$.
Thus we have 
$\widetilde H_n = \hat H_n\psi_n$ whether $\psi_n>0$ or $\psi_n=0$.
Then, by \eqref{eq:230804.2030} and \eqref{230423.1826}, we can write  
\begin{align*}
  Z_n &= 2\log2\,\sqrt{n} (\hat H_n - H) \psi_n 
  = 2\log2\,\sqrt{n} (\widetilde H_n - H \psi_n)
  \\&= 
  2\log2\,\sqrt{n} \cbr{
    \rbr{H 
    +\frac1{2\log2}
    \Brbr{
      %\log\frac{\rv_n}{\iv} - \log\frac{\rv_{2n}}{\iv}
      \log\Brbr{\frac{\convDiff_n}{\iv}+1}-\log\Brbr{\frac{\convDiff_{2n}}{\iv}+1}
      }
    } \psi_n - H \psi_n
  }
  % \\&= 
  % \sqrt{n}\Brbr{
  %     \log\Brbr{\frac{\convDiff_n}{\iv}+1}-\log\Brbr{\frac{\convDiff_{2n}}{\iv}+1}}
  %   \psi_n
  \\&=
  \sqrt{n}
  \cbr{
    \frac{\convDiff_n}\iv -\half\rbr{\frac{\convDiff_n}\iv}^2
      + r\rbr{\frac{\convDiff_n}\iv}
    -\rbr{
      \frac{\convDiff_{2n}}\iv -\half\rbr{\frac{\convDiff_{2n}}\iv}^2
      + r\rbr{\frac{\convDiff_{2n}}\iv}}
  }\psi_n
  \\&=
  \sqrt{n} 
  \cbr{\frac{\convDiff_n}\iv - \frac{\convDiff_{2n}}\iv}
  -  \sqrt{n} 
  \cbr{\half\rbr{\frac{\convDiff_n}\iv}^2 -\half\rbr{\frac{\convDiff_{2n}}\iv}^2}
  + \negTerm_n^{(1)},
\end{align*}
where 
\begin{align*}
  \negTerm_n^{(1)}=
  \sqrt{n}
  \rbr{r\rbr{\frac{\convDiff_n}\iv}-r\rbr{\frac{\convDiff_{2n}}\iv}}
  \psi_n
  +
  \sqrt{n}
  \cbr{\frac{\convDiff_n}\iv -\half\rbr{\frac{\convDiff_n}\iv}^2
  -\rbr{\frac{\convDiff_{2n}}\iv -\half\rbr{\frac{\convDiff_{2n}}\iv}^2}}(\psi_n-1).
\end{align*}
Then, by Lemmas \ref{230510.1721}, \ref{lem:230925.1301} (iii)
and \ref{230510.1711} (ii),
% $\convDiff_n=O_M(n^{-\half})$,
we obtain $\negTerm_n^{(1)} = O_M(n^{-1})$.
\end{proof}

\subsection{Limits of functionals}\label{sec:231005.2102}

Here we will identify the limits of 
$\cali_n^{(i_1,i_2)(1)}$ ($i_1,i_2\in\cbr{1,2}$)
and 
$n^{1/2}D_{u_n}D_{u_n} M_n$.
The former functional appears as the principal term of
$\abr{DM_n,u_n}$, which plays an imporatnt role in
the coefficient of the random symbol ``quasi-tangent'', 
and its limit is a part of the asymptotic variance $G_\infty$ of $Z_n$.
The latter functional is the coefficient of the random symbol 
``quasi-torsion'' and its limit appears in the asymptotic expansion formula.

We define $\widehat\rho$ and $\tilde\rho$ by
\begin{align}
  \widehat\rho(j)
  &=
  \abr{\bbone_{[0,1]}-\bbone_{[-1,0]},
  \bbone_{[j,j+1]}-\bbone_{[j-1,j]}}
  \label{eq:230925.1740}
  \\
  \tilde\rho(j)
  &=
  2^{-2H}\abr{
    \bbone_{[0,1]}-\bbone_{[-1,0]},
    \bbone_{[j,j+2]}-\bbone_{[j-2,j]}}
  \label{eq:230925.1741}
\end{align}
for $j\in\bbZ$, and set
\begin{align}\label{eq:230925.1914}
  \widehat c = \sum_{j\in\bbZ} \widehat\rho(j)^2, \tand
  \tilde c = \sum_{j\in\bbZ} \tilde\rho(j)^2.
\end{align}
Here the inner product in \eqref{eq:230925.1740} and \eqref{eq:230925.1741} 
should be read as the extended one.
(See \eqref{def:231005.2143}.)
It holds that 
\begin{align*}
  n^{2H}\abr{\diffker^n_\jon, \diffker^{n}_\jtw}=
  \widehat\rho(\jon-\jtw)
  \tand
  n^{2H}\abr{\diffker^n_\jon, \diffker^{2n}_\jtw}=
  \tilde\rho(2\jon-\jtw).
\end{align*}

The following lemma is from
Lemma 3.4 and 3.6 of \cite{mishura2023asymptotic}.
By this lemma, we have $\hat c, \tilde c<\infty$.
\begin{lemma}\label{230807.1545}
  (i) It holds that 
  $\widehat\rho(k) = O(\abs{k}^{2H-4})$ as $\abs{k}\to\infty$.
  In other words,
  % $\tilde\rho(k) \simleq(\abs{k}\vee1)^{2H-4}$;
  there exists $C>0$ such that
  $\abs{\widehat\rho(k)} \leq C(\abs{k}\vee1)^{2H-4}$ for $k\in\bbZ$.

  \item[(ii)]
  $\tilde\rho(k) = O(\abs{k}^{2H-4})$ as $\abs{k}\to\infty$.
  
\end{lemma}

\subsubsection{Limits of $\cali_n^{(i_1,i_2)(1)}$} % ($i_1,i_2\in\cbr{1,2}$)}
Recall the definitions of 
$\cali_n^{(i_1,i_2)(1)}$ ($i_1,i_2\in\cbr{1,2}$)
and $\widehat{G}_\infty$
defined at \eqref{230703.1731} and \eqref{230728.2030}:
\begin{align*}
  \cali_n^{(i_1,i_2)(1)} &= 
  (i_1\,i_2)^{(2H-1)} n^{4H-1}
  \sum_{j\in[i_1 n-1]\times[i_2 n-1]} 
  \rbr{\iv}^{-2} \brbr{V^{[1]}_{t^{i_1 n}_{j_1}}}^2 \brbr{V^{[1]}_{t^{i_2 n}_{j_2}}}^2 
  \abr{\diffker^{i_1 n}_{j_1},\diffker^{i_2n}_{j_2}}^2
  \\
  \widehat{G}_\infty &= \rbr{\iv}^{-2} \int^1_0 \brbr{V^{[1]}_t}^4 dt.
\end{align*}
% \begin{comment}
%   {\myblue $a(x) = (V^{[1]}(x))^2$
%   (hence $a_t = \brbr{V^{[1]}_{t} }^2 $)を使ったほうが記述が軽くなるはず．\koko}
% \end{comment}

\begin{lemma}[Lemma \ref{230724.1350}]\label{lem:230926.1700}
  For any $\beta\in(\half,H)$,
  the functionals $\cali_n^{(1,1)(1)}$ and $\cali_n^{(1,2)(1)}$
  are written as:
  \begin{align*}
    \cali_n^{(1,1)(1)} 
    &= 
    % \hat c \;
    % \rbr{\iv}^{-2}
    % \int^1_0 \rbr{V^{[1]}_t}^4 dt
    % +O_M(n^{(-H)\vee(-\frac34)})
    % = 
    \widehat c \; \widehat{G}_\infty
    +O_M(n^{(-\beta)\vee(-\frac34)})
    % +O_M(n^{(-H)\vee(-\frac34)})
    \tand%\\
    \cali_n^{(1,2)(1)}%=\cali_n^{(2,1)(1)} 
    =%&=
    % 2^{2H-1} \tilde c\;
    % \rbr{\iv}^{-2}
    % \int^1_0 \rbr{V^{[1]}_t}^4 dt
    % +O_M(n^{(-H)\vee(-\frac34)})
    % =
    2^{2H-1} \tilde c\; \widehat{G}_\infty
    +O_M(n^{(-\beta)\vee(-\frac34)}),
    % +O_M(n^{(-H)\vee(-\frac34)}),
  \end{align*}
  as $n$ tends to $\infty$.
\end{lemma}
\begin{proof}
  By the definition of $\tilde\rho$, 
  the functional $\cali_n^{(1,2)(1)}$ is 
  written as 
  \begin{align*}
    \cali_n^{(1,2)(1)} &=
    2^{2H-1} n^{-1} \rbr{\iv}^{-2}
    \sum_{\jon\in[n-1]} \sum_{\jtw\in[2n-1]}
    \brbr{V^{[1]}_{t^n_\jon} }^2 \brbr{V^{[1]}_{t^{2n}_\jtw} }^2 
    \tilde\rho(2\jon-\jtw)^2
  \end{align*}
  The main factor of $\cali_n^{(1,2)(1)}-2^{2H-1} \tilde c\; \widehat{G}_\infty$, 
  namely
  $2^{-(2H-1)}\rbr{\iv}^{2}
  \Brbr{\cali_n^{(1,2)(1)}-2^{2H-1} \tilde c\; \widehat{G}_\infty}$,
  decomposes as 
  \begin{align*}
    2^{-(2H-1)}\rbr{\iv}^{2}
    \Brbr{\cali_n^{(1,2)(1)}-2^{2H-1} \tilde c\; \widehat{G}_\infty}
    &=%\\&=
    2^{-(2H-1)}\rbr{\iv}^{2}\cali_n^{(1,2)(1)} 
    - \tilde c \int^1_0 \brbr{V^{[1]}_t}^4 dt
    \\&=
    \Delta_1+\Delta_2+\Delta_3
  \end{align*}
  with 
  \begin{align*}
    \Delta_1&=
     n^{-1} 
    \sum_{\jon\in[n-1]} \sum_{\jtw\in[2n-1]}
    \brbr{V^{[1]}_{t^n_\jon} }^2 \brbr{V^{[1]}_{t^{2n}_\jtw} }^2 
    \tilde\rho(2\jon-\jtw)^2
    - 
    n^{-1} 
    \sum_{\jon\in[n-1]} \sum_{\jtw\in[2n-1]}
    \brbr{V^{[1]}_{t^{n}_{\jon}} }^2 \brbr{V^{[1]}_{t^{2n}_{2\jon}} }^2 
    \tilde\rho(2\jon-\jtw)^2
    \\
    \Delta_2&=
    n^{-1} 
    \sum_{\jon\in[n-1]} \sum_{\jtw\in[2n-1]}
    \brbr{V^{[1]}_{t^{n}_{\jon}} }^4
    \tilde\rho(2\jon-\jtw)^2
    -
    n^{-1} 
    \sum_{\jon\in[n-1]} 
    \brbr{V^{[1]}_{t^{n}_{\jon}} }^4
    \sum_{\jtw\in\bbZ}
    \tilde\rho(2\jon-\jtw)^2
    \\
    \Delta_3&=
    \tilde c\,
    n^{-1} 
    \sum_{\jon\in[n-1]} 
    \brbr{V^{[1]}_{t^{n}_{\jon}} }^4
    -
    \tilde c \int^1_0
    \brbr{V^{[1]}_{t} }^4dt.
  \end{align*}
  Here recall that we have defined 
  $\tilde c = \sum_{j\in\bbZ} \tilde\rho(j)^2$.
  
  For $k\in\bbZ_{\geq0}$ and $p>1$,
  the norm of $\Delta_1$ is bounded as
  \begin{align}
    \norm{\Delta_1}_{k,p}&=
    n^{-1} 
    \norm{\sum_{\jon\in[n-1]} \sum_{\jtw\in[2n-1]}
    \brbr{V^{[1]}_{t^n_\jon} }^2 
    \Brbr{\brbr{V^{[1]}_{t^{2n}_\jtw} }^2 - \brbr{V^{[1]}_{t^{2n}_{2\jon}} }^2}
    \tilde\rho(2\jon-\jtw)^2}_{k,p}
    % \\&\leq 
    % n^{-1} 
    % \sum_{\jon\in[n-1]} \sum_{\jtw\in[2n-1]} 
    % \norm{
    % \brbr{V^{[1]}_{t^n_\jon} }^2 
    % \Brbr{\brbr{V^{[1]}_{t^{2n}_\jtw} }^2 - \brbr{V^{[1]}_{t^{2n}_{2\jon}} }^2}}_{k,p}
    % \tilde\rho(2\jon-\jtw)^2
    \nn\\&\simleq 
    n^{-1} 
    \sum_{\jon\in[n-1]} \sum_{\jtw\in[2n-1]} 
    \norm{\brbr{V^{[1]}_{t^n_\jon} }^2}_{k,2p}
    \norm{\brbr{V^{[1]}_{t^{2n}_\jtw}}^2 - \brbr{V^{[1]}_{t^{2n}_{2\jon}}}^2}_{k,2p}
    \tilde\rho(2\jon-\jtw)^2
  % \end{align*}
  % \begin{align*}
  %   \norm{\Delta_1}_{k,p}&\simleq 
  %   n^{-1} 
  %   \sum_{\jon\in[n-1]} \sum_{\jtw\in[2n-1]} 
  %   \norm{\brbr{V^{[1]}_{t^n_\jon} }^2}_{k,2p}
  %   \norm{\brbr{V^{[1]}_{t^{2n}_\jtw}}^2 - \brbr{V^{[1]}_{t^{2n}_{2\jon}}}^2}_{k,2p}
  %   \tilde\rho(2\jon-\jtw)^2
    \nn\\&\simleq 
    n^{-1} 
    \sum_{\jon\in[n-1]} \sum_{\jtw\in[2n-1]} 
    % \norm{\brbr{V^{[1]}_{t^n_\jon} }^2}_{k,2p}
    \abs{2\jon-\jtw}^\beta \;n^{-\beta}\;
    (\abs{2\jon-\jtw}\vee1)^{4H-8}
    \label{230807.1550}
    \\&\leq 
    % n^{-1-\beta}
    % \sum_{\jon\in[n-1]} \sum_{\jtw\in[2n-1]} 
    % % \norm{\brbr{V^{[1]}_{t^n_\jon} }^2}_{k,2p}
    % (\abs{2\jon-\jtw}\vee1)^{4H+\beta-8}
    % \nn\\&\leq 
    n^{-1-\beta} (n-1) \sum_{j\in\bbZ} (\abs{j}\vee1) ^{4H+\beta-8}
    =O(n^{-\beta}).
    \nn
  \end{align}
  Here $\beta$ is taken arbitrarily from $(\half,H)$.
  At \eqref{230807.1550},
  we used Lemma \ref{230807.1545} (i) and 
   an estimate  
  \begin{align*}
    \bnorm{\brbr{V^{[1]}_{t^{2n}_\jtw}}^2 - \brbr{V^{[1]}_{t^{2n}_{2\jon}}}^2}_{k,2p}
    \simleq
    \abs{2\jon-\jtw}^\beta \;n^{-\beta}
  \end{align*}
  from Lemma \ref{230924.1800} (ii). %\ref{230807.1543}.

  Concerning with $\Delta_2$, 
  we have the following estimate for $\epsilon\in(0,1)$:
  \begin{align*}
    &
    n^{-1} 
    \sum_{\jon\in[n-1]} ~\sum_{k\leq2\jon-2n}
    \tilde\rho(k)^2
    +
    n^{-1} 
    \sum_{\jon\in[n-1]} ~\sum_{k\geq2\jon}
    \tilde\rho(k)^2
    \\&=
    n^{-1} 
    \sum_{1\leq\jon<n-n^\epsilon} ~\sum_{k\leq2\jon-2n}
    \tilde\rho(k)^2
    +
    n^{-1} 
    \sum_{n-n^\epsilon\leq\jon\leq n-1} ~\sum_{k\leq2\jon-2n}
    \tilde\rho(k)^2
    \\&\hspsm+
    n^{-1} 
    \sum_{1\leq\jon\leq n^\epsilon} ~\sum_{k\geq2\jon}
    \tilde\rho(k)^2
    +
    n^{-1} 
    \sum_{n^\epsilon<\jon\leq n-1} ~\sum_{k\geq2\jon}
    \tilde\rho(k)^2
    \\&=:
    \sum_{i=1}^4 \Delta_2^{(i)}
  \end{align*}
  We have 
  $\Delta_2^{(1)},\Delta_2^{(4)} = O(n^{-1}n^{1}n^{(4H-7)\epsilon})
  = O(n^{(4H-7)\epsilon})$, and 
  $\Delta_2^{(2)},\Delta_2^{(3)} = O(n^{-1}n^{\epsilon})O(1) 
  = O(n^{-1+\epsilon})$.
  Hence, by choosing 
  $\epsilon = \frac1{4}$,
  % $\epsilon = \frac1{8-4H}$
  we obtain 
  $\sum_{i=1}^4 \Delta_2^{(i)} = O(n^{-\frac34})$ and
  \begin{align*}
    \norm{\Delta_2}_{k,p}&=
    % \norm{n^{-1} 
    % \sum_{\jon\in[n-1]} 
    % \brbr{V^{[1]}_{t^{2n}_{2\jon}} }^4
    % \sum_{\jtw\in[2n-1]}
    % \tilde\rho(2\jon-\jtw)^2
    % -
    % n^{-1} 
    % \sum_{\jon\in[n-1]} 
    % \brbr{V^{[1]}_{t^{2n}_{2\jon}} }^4
    % \sum_{\jtw\in\bbZ}
    % \tilde\rho(2\jon-\jtw)^2}_{k,p}
    % \\&=
    \norm{n^{-1} 
    \sum_{\jon\in[n-1]} 
    \brbr{V^{[1]}_{t^{n}_{\jon}} }^4
    \rbr{\sum_{\jtw\leq0}
    \tilde\rho(2\jon-\jtw)^2
    +
    \sum_{\jtw\geq2n}
    \tilde\rho(2\jon-\jtw)^2}}_{k,p}
    = O(n^{-\frac34}).
  \end{align*}

  From Lemma \ref{230924.1800}(iii), %\ref{230807.1750}, 
  $\Delta_3=O_M(n^{-1})$.
  Therefore, 
  \begin{align*}
    {\cali_n^{(1,2)(1)}-2^{2H-1} \tilde c\; \widehat{G}_\infty}
    &=
    2^{(2H-1)}\rbr{\iv}^{-2}
    (\Delta_1+\Delta_2+\Delta_3)
    =O_M(n^{(-\beta)\vee(-\frac34)})
  \end{align*}
  for any $\beta\in(\half,H)$.

  For $\cali_n^{(1,1)(1)}$, by the definition of $\widehat\rho$,
  we can write 
  \begin{align*}
    \cali_n^{(1,1)(1)} &=
    n^{-1} \rbr{\iv}^{-2}
    \sum_{\jon\in[n-1]} \sum_{\jtw\in[n-1]}
    \brbr{V^{[1]}_{t^n_\jon} }^2 \brbr{V^{[1]}_{t^{n}_\jtw} }^2 
    \widehat\rho(\jon-\jtw)^2.
  \end{align*}
  Since $\widehat\rho(j)$ decays as $\abs{j}\to\infty$ 
  at the same order as $\tilde\rho$ (Lemma \ref{230807.1545}),
  an argument parallel to that of $\cali_n^{(1,2)(1)}$ shows 
  \begin{align*}
    \cali_n^{(1,1)(1)} -
    \hat c \; \widehat{G}_\infty
    =O_M(n^{(-\beta)\vee(-\frac34)}).
  \end{align*}
  for any $\beta\in(\half,H)$.
\end{proof}

\subsubsection{Limits of $n^{\half}D_{u_n}D_{u_n} M_n$
(The coefficient of $\qtor$)}
% \begin{itembox}[l]{Memo}
%   $n^{2H}\abr{\diffker^n_\jon, \diffker^{n}_\jtw}=
%   \widehat\rho(\jon-\jtw)$,
%   $n^{2H}\abr{\diffker^n_\jon, \diffker^{2n}_\jtw}=
%   \tilde\rho(2\jon-\jtw)$
% \end{itembox}
The next estimate is used in the proof of Lemma \ref{lem:230925.1805}.
We denote  $a_t = \brbr{V^{[1]}_{t} }^2 $ for brevity.
\begin{lemma}\label{230109.1833}
  The following estimates hold as $n$ tends to $\infty$:
  \begin{align*}
    &n^{-1} 
    \sum_{j\in[n-1]^3}
    (a_{t_\jon})^3\;
    \widehat\rho(\jon-\jtw)
    \widehat\rho(\jtw-\jth)
    \widehat\rho(\jth-\jon)
    \\&\hspace{150pt}=
    % =%\\&=
    n^{-1} 
    \sum_{j_1\in[n-1]}
    (a_{t_\jon})^3\;
    \sum_{i_1,i_2\in\bbZ}
    \widehat\rho(i_1)
    \widehat\rho(i_1-i_2)
    \widehat\rho(i_2)
    +o_M(1)
    \\
    &n^{-1} \sum_{j\in[n-1]^2\times[2n-1]}
    \brbr{a_{t_{j_1}^n}}^3 
    \tilde\rho(2j_2-j_3)
    \tilde\rho(2j_1-j_3)
    \widehat\rho(j_1-j_2)
    \\&\hspace{150pt}=
    n^{-1} \sum_{j_1\in[n-1]}
    \brbr{a_{t_{j_1}^n}}^3 
    \sum_{i_1,i_2\in\bbZ}
    \tilde\rho(i_2)
    \tilde\rho(2i_1-i_2)
    \widehat\rho(i_1)
    +o_M(1)
    \\
    &n^{-1} \sum_{j\in[n-1]\times[2n-1]^2}
    \brbr{a_{t_{j_1}^n}}^3 
    \tilde\rho(j_2-2j_1)
    \tilde\rho(j_3-2j_1)
    \widehat\rho(j_2-j_3)
    \\&\hspace{150pt}=
    n^{-1} \sum_{j_1\in[n-1]}
    \brbr{a_{t_{j_1}^n}}^3 
    \sum_{i_1,i_2\in\bbZ}
    \tilde\rho(i_1)
    \tilde\rho(i_2)
    \widehat\rho(i_1-i_2)
    +o_M(1)
  \end{align*}
\end{lemma}
\begin{proof}
  We only consider the second convergence above. 
  (The others can be proved by similar arguments.)
  Since the LHS is written as 
  \begin{align*}
    % &n^{-1} \sum_{j\in[n-1]^2\times[2n-1]}
    % \nrbr{a_{t_{j_1}^n}}^3 
    % \tilde\rho(2j_2-j_3)
    % \tilde\rho(2j_1-j_3)
    % \widehat\rho(j_1-j_2)
    &%\\&=
    n^{-1} 
    \sum_{j_1\in[n-1]}
    \brbr{a_{t_{j_1}^n}}^3 
    \sum_{(j_2,j_3)\in[n-1]\times[2n-1]}
    \tilde\rho(2j_2-j_3)
    \tilde\rho(2j_1-j_3)
    \widehat\rho(j_1-j_2)
    % \\&=
    % n^{-1} 
    % \sum_{j_1\in[n-1]}
    % \brbr{a_{t_{j_1}^n}}^3 
    % \sum_{j_2\in[n-1]}
    % \sum_{j_3\in[2n-1]}
    % \widehat\rho(j_2-j_1)
    % \tilde\rho(2j_2-j_3)
    % \tilde\rho(j_3-2j_1)
    \\&=
    n^{-1} 
    \sum_{j_1\in[n-1]}
    \brbr{a_{t_{j_1}^n}}^3 
    \sum_{\substack{i_1\in [n-1]-j_1\\i_2\in[2n-1]-2j_1}}
    % \sum_{i_1\in [n-1]-j_1}
    % \sum_{i_2\in[2n-1]-2j_1}
    \tilde\rho(2i_1-i_2)
    \tilde\rho(i_2)
    \widehat\rho(i_1),
  \end{align*}
  the error of convergence is written as 
  \begin{align*}
    n^{-1} 
    \sum_{j_1\in[n-1]}
    \brbr{a_{t_{j_1}^n}}^3 
    \bbrbr{
    \sum_{\substack{i_1\in [n-1]-j_1\\i_2\in[2n-1]-2j_1}}
    \widehat\rho(i_1)
    \tilde\rho(2i_1-i_2)
    \tilde\rho(i_2)
    -
    \sum_{i_1,i_2\in\bbZ}
    \widehat\rho(i_1)
    \tilde\rho(2i_1-i_2)
    \tilde\rho(i_2)
    }=:\Delta_n.
  \end{align*}
  We will consider the estimate of the following difference
  for $n^\epsilon<j_1<n-n^\epsilon$
  with $\epsilon\in(0,1)$:
  \begin{align*}
    &\bbabs{
    \sum_{i_1,i_2\in\bbZ}
    \widehat\rho(i_1)
    \tilde\rho(2i_1-i_2)
    \tilde\rho(i_2)
    -
    \sum_{\substack{i_1\in [n-1]-j_1\\i_2\in[2n-1]-2j_1}}
    \widehat\rho(i_1)
    \tilde\rho(2i_1-i_2)
    \tilde\rho(i_2)}
    \\&=
    \bbabs{
      \sum_{\substack{i_1\not\in [n-1]-j_1\\\tor\\i_2\not\in[2n-1]-2j_1}}
      \widehat\rho(i_1)
      \tilde\rho(2i_1-i_2)
      \tilde\rho(i_2)}
    \\&\leq
    \bbabs{
      \sum_{\substack{\abs{i_1}>n^\epsilon\\\tor\\\abs{i_2}>2 n^\epsilon}}
      \widehat\rho(i_1)
      \tilde\rho(2i_1-i_2)
      \tilde\rho(i_2)}
    \\&\leq
    \bbabs{
      \sum_{\substack{\abs{i_1}>n^\epsilon\\\tandsm\\\abs{i_2}>2 n^\epsilon}}
      \widehat\rho(i_1)
      \tilde\rho(2i_1-i_2)
      \tilde\rho(i_2)}
    +
    \bbabs{
      \sum_{\substack{\abs{i_1}\leq n^\epsilon\\\tandsm\\\abs{i_2}>2 n^\epsilon}}
      \widehat\rho(i_1)
      \tilde\rho(2i_1-i_2)
      \tilde\rho(i_2)}
    +
    \bbabs{
      \sum_{\substack{\abs{i_1}>n^\epsilon\\\tandsm\\\abs{i_2}\leq2 n^\epsilon}}
      \widehat\rho(i_1)
      \tilde\rho(2i_1-i_2)
      \tilde\rho(i_2)}
  \end{align*}
  We used that 
  for $n^\epsilon<j_1<n-n^\epsilon$,
  $i_1\not\in [n-1]-j_1\Rightarrow \abs{i_1}>n^\epsilon$ and 
  $i_2\not\in [2n-1]-2j_1 \Rightarrow \abs{i_2}>2n^\epsilon$.
  The second term above is estimated as
  \begin{align*}
    &\bbabs{\sum_{\substack{\abs{i_1}\leq n^\epsilon%\\\tandsm
      \\\abs{i_2}>2 n^\epsilon}}
      \widehat\rho(i_1)
      \tilde\rho(2i_1-i_2)
      \tilde\rho(i_2)}
    \leq 
    \sum_{\substack{\abs{i_1}\leq n^\epsilon%\\\tandsm
    \\\abs{i_2}>2 n^\epsilon}}
    %   \substack{n^\epsilon\geq\abs{i_1}%\\\tandsm
    % \\2 n^\epsilon<\abs{i_2}}
    \abs{\widehat\rho(i_1)}\,
    \abs{\tilde\rho(2i_1-i_2)}\,
    \abs{\tilde\rho(i_2)}
    \\&\leq
    \sup_{i_1,i_2}
    \abs{\tilde\rho(2i_1-i_2)}
    \sum_{\abs{i_1}\leq n^\epsilon}
    \abs{\widehat\rho(i_1)}
    \sum_{\abs{i_2}>2 n^\epsilon}
    \abs{\tilde\rho(i_2)}
    \\&=O(1)O(1)O(n^{\epsilon(2H-3)})
    =O(n^{\epsilon(2H-3)}),
  \end{align*}
  since 
  $\sup_{i\in\bbZ}\abs{\tilde\rho(i)}<\infty$,
  $\sum_{i\in\bbZ}\abs{\widehat\rho(i)}<\infty$ and 
  $\sum_{\abs{i}>k}\abs{\tilde\rho(i)}=O(k^{2H-3})$ as $k\to\infty$.
  Similar arguments show that 
  \begin{align*}
    \bbabs{
      \sum_{\substack{\abs{i_1}>n^\epsilon\\\tandsm\\\abs{i_2}>2 n^\epsilon}}
      \widehat\rho(i_1)
      \tilde\rho(2i_1-i_2)
      \tilde\rho(i_2)}
      =O(n^{2\epsilon(2H-3)}),\tand
    \bbabs{
      \sum_{\substack{\abs{i_1}>n^\epsilon\\\tandsm\\\abs{i_2}\leq2 n^\epsilon}}
      \widehat\rho(i_1)
      \tilde\rho(2i_1-i_2)
      \tilde\rho(i_2)}
    =O(n^{\epsilon(2H-3)}).
  \end{align*}

  Hence the order of the norm of the difference is bounded as
  \begin{align*}
    % &\Bigg|\Bigg|
    % n^{-1} \sum_{j\in[n-1]^2\times[2n-1]}
    % \brbr{a_{t_{j_1}^n}}^3 
    % \tilde\rho(2j_2-j_3)
    % \tilde\rho(2j_1-j_3)
    % \widehat\rho(j_1-j_2)
    % \\&\hspace{150pt}-
    % n^{-1} \sum_{j_1\in[n-1]}
    % \brbr{a_{t_{j_1}^n}}^3 
    % \sum_{i_1,i_2\in\bbZ}
    % \tilde\rho(i_2)
    % \tilde\rho(2i_1-i_2)
    % \widehat\rho(i_1)
    % \Bigg|\Bigg|_{k,p}
    % \\&=
    % \Bigg|\Bigg|
    % n^{-1} \sum_{j_1\in[n-1]}
    % \brbr{a_{t_{j_1}^n}}^3 
    % \bbrbr{
    % \sum_{\substack{j_2\in[n-1]\\j_3\in[2n-1]}}
    % % \sum_{(j_2,j_3)\in[n-1]\times[2n-1]}
    % \tilde\rho(2j_2-j_3)
    % \tilde\rho(2j_1-j_3)
    % \widehat\rho(j_1-j_2)
    % % \\&\hspace{50pt}-
    % -
    % \sum_{i_1,i_2\in\bbZ}
    % \tilde\rho(i_2)
    % \tilde\rho(2i_1-i_2)
    % \widehat\rho(i_1)}
    % \Bigg|\Bigg|_{k,p}
    % &% \\&
    % \Bigg|\Bigg|
    % n^{-1} 
    % \sum_{j_1\in[n-1]}
    % \brbr{a_{t_{j_1}^n}}^3 
    % \bbrbr{
    % \sum_{\substack{i_1\in [n-1]-j_1\\i_2\in[2n-1]-2j_1}}
    % \widehat\rho(i_1)
    % \tilde\rho(2i_1-i_2)
    % \tilde\rho(i_2)
    % -
    % \sum_{i_1,i_2\in\bbZ}
    % \tilde\rho(i_2)
    % \tilde\rho(2i_1-i_2)
    % \widehat\rho(i_1)}
    % \Bigg|\Bigg|_{k,p}
    % \\&\leq
    \norm{\Delta_n}_{k,p}&\leq
    n^{-1} \sum_{j_1\in[n-1]}
    \norm{\brbr{a_{t_{j_1}^n}}^3}_{k,p}
    \bbabs{
    \sum_{\substack{i_1\in[n-1]-j_1\\i_2\in[2n-1]-2j_1}}
    % \sum_{\substack{j_2\in[n-1]\\j_3\in[2n-1]}}
    \widehat\rho(i_1)
    \tilde\rho(2i_1-i_2)
    \tilde\rho(i_2)
    % \tilde\rho(2j_2-j_3)
    % \tilde\rho(2j_1-j_3)
    % \widehat\rho(j_1-j_2)
    % \\&\hspace{50pt}-
    -
    \sum_{i_1,i_2\in\bbZ}
    \widehat\rho(i_1)
    \tilde\rho(2i_1-i_2)
    \tilde\rho(i_2)}
    \\&\simleq
    n^{-1} \sum_{j_1\leq n^\epsilon\torsm n-n^\epsilon\leq j_1}
    \bbabs{
    \sum_{\substack{i_1\in[n-1]-j_1\\i_2\in[2n-1]-2j_1}}
    \widehat\rho(i_1)
    \tilde\rho(2i_1-i_2)
    \tilde\rho(i_2)
    -
    \sum_{i_1,i_2\in\bbZ}
    \widehat\rho(i_1)
    \tilde\rho(2i_1-i_2)
    \tilde\rho(i_2)}
    \\&\quad+
    n^{-1} \sum_{n^\epsilon<j_1<n-n^\epsilon}
    \bbabs{
    \sum_{\substack{i_1\in[n-1]-j_1\\i_2\in[2n-1]-2j_1}}
    \widehat\rho(i_1)
    \tilde\rho(2i_1-i_2)
    \tilde\rho(i_2)
    -
    \sum_{i_1,i_2\in\bbZ}
    \widehat\rho(i_1)
    \tilde\rho(2i_1-i_2)
    \tilde\rho(i_2)}
    \\&=
    O(n^{-1+\epsilon}) + O(n^{-1+1+\epsilon(2H-3)})=o(1),
  \end{align*}
  since $\epsilon\in(0,1)$.
  Here we used the following estimate 
  % We used Lemma \ref{230808.1804}
  for the case %where
  $j_1\leq n^\epsilon$ or $n-n^\epsilon\leq j_1$:
  \begin{align*}
    \sum_{i_1,i_2\in\bbZ}
    \abs{\widehat\rho(i_1)}\,
    \abs{\tilde\rho(2i_1-i_2)}\,
    \abs{\tilde\rho(i_2)}
    \leq
    \sup_{i_1,i_2}\abs{\tilde\rho(2i_1-i_2)}\,
    \sum_{i_1\in\bbZ}
    \abs{\widehat\rho(i_1)}\,
    \sum_{i_2\in\bbZ}
    \abs{\tilde\rho(i_2)}
    <\infty.
  \end{align*}

\end{proof}

% \begin{lemma}\label{230808.1804}
%   \begin{align*}
%     \sum_{i_1,i_2\in\bbZ}
%     \abs{\widehat\rho(i_1)}\,
%     \abs{\tilde\rho(2i_1-i_2)}\,
%     \abs{\tilde\rho(i_2)}<\infty
%   \end{align*}
% \end{lemma}
% \begin{proof}
%   \begin{align*}
%     \sum_{i_1,i_2\in\bbZ}
%     \abs{\widehat\rho(i_1)}\,
%     \abs{\tilde\rho(2i_1-i_2)}\,
%     \abs{\tilde\rho(i_2)}
%     \leq
%     \sup_{i_1,i_2}\abs{\tilde\rho(2i_1-i_2)}\,
%     \sum_{i_1\in\bbZ}
%     \abs{\widehat\rho(i_1)}\,
%     \sum_{i_2\in\bbZ}
%     \abs{\tilde\rho(i_2)}
%     <\infty.
%   \end{align*}
% \end{proof}

Recall that we have defined
$u_n^{(i)}$ ($i=1,2$) at \eqref{eq:230926.1446}, %{230725.1531},
and 
$u_n=u_n^{(1)}-u_n^{(2)}$ and $M_n=\delta(u_n)$.
\begin{lemma}%[Lemma \ref{230906.2127}]
  \label{lem:230925.1805}
  The following decomposition of $n^{1/2}D_{u_n}D_{u_n} M_n$ holds:
  \begin{align}
    &n^{1/2}D_{u_n}D_{u_n}M_n = 
    3\,c_\qtor \; \rbr{\iv}^{-3} \int_0^1 (a_t)^3 dt
    +o_M(1),
    \label{230807.1953}
  \end{align}
  where the constant $c_\qtor$ is defined by 
  \begin{align}
    c_\qtor&=
    {\sum_{i_1,i_2\in\bbZ}
    \widehat\rho(i_1)
    \widehat\rho(i_1-i_2)
    \widehat\rho(i_2)}
    -
    2^{(2H+1)}
    {\sum_{i_1,i_2\in\bbZ}
    \tilde\rho(i_2)
    \tilde\rho(2i_1-i_2)
    \widehat\rho(i_1)}
    +%\\&\hspace{120pt}+
    2^{2H}
    {\sum_{i_1,i_2\in\bbZ}
    \tilde\rho(i_1)
    \tilde\rho(i_2)
    \widehat\rho(i_1-i_2)}.
    \label{230726.1339}
  \end{align}

\end{lemma}
\begin{proof}
  By Lemma \ref{230721.1259},
  the functional $n^\half D_{u_n}D_{u_n} M_n$ decomposes as
  \begin{align*}
    &n^{1/2}D_{u_n}D_{u_n} M_n
    \nn\\&=
    2^2\times n^{6H-1}
    \sum_{j\in[n-1]^3}%\bbJ_n([3],i)} 
    \rbr{\iv}^{-3} 
    \brbr{V^{[1]}_{t^{n}_{j_1}}}^2 \brbr{V^{[1]}_{t^{n}_{j_2}}}^2 \brbr{V^{[1]}_{t^{n}_{j_3}}}^2 
    \abr{\diffker^{n}_{j_1},\diffker^{n}_{j_2}}
    \abr{\diffker^{n}_{j_2},\diffker^{n}_{j_3}}
    \abr{\diffker^{n}_{j_1},\diffker^{n}_{j_3}} 
    \nn\\
    &\quad-3\times
    2^2\times2^{(2H-1)}\times n^{6H-1}
    \sum_{j\in[n-1]^2\times[2n-1]} 
    \rbr{\iv}^{-3} 
    \brbr{V^{[1]}_{t^{n}_{j_1}}}^2 \brbr{V^{[1]}_{t^{n}_{j_2}}}^2 \brbr{V^{[1]}_{t^{2n}_{j_3}}}^2 
    \abr{\diffker^{n}_{j_1},\diffker^{n}_{j_2}}
    \abr{\diffker^{n}_{j_2},\diffker^{2n}_{j_3}}
    \abr{\diffker^{n}_{j_1},\diffker^{2n}_{j_3}} 
    \nn\\
    &\quad+3\times
    2^2\times2^{2(2H-1)}\times n^{6H-1}
    \sum_{j\in[n-1]\times[2n-1]} 
    \rbr{\iv}^{-3} 
    \brbr{V^{[1]}_{t^{n}_{j_1}}}^2 \brbr{V^{[1]}_{t^{2n}_{j_2}}}^2 \brbr{V^{[1]}_{t^{2n}_{j_3}}}^2 
    \abr{\diffker^{n}_{j_1},\diffker^{2n}_{j_2}}
    \abr{\diffker^{2n}_{j_2},\diffker^{2n}_{j_3}}
    \abr{\diffker^{n}_{j_1},\diffker^{2n}_{j_3}} 
    \nn\\
    &\quad-
    2^2\times2^{3(2H-1)}\times n^{6H-1}
    \sum_{j\in[2n-1]^3} 
    \rbr{\iv}^{-3} 
    \brbr{V^{[1]}_{t^{2n}_{j_1}}}^2 \brbr{V^{[1]}_{t^{2n}_{j_2}}}^2 \brbr{V^{[1]}_{t^{2n}_{j_3}}}^2 
    \abr{\diffker^{2n}_{j_1},\diffker^{2n}_{j_2}}
    \abr{\diffker^{2n}_{j_2},\diffker^{2n}_{j_3}}
    \abr{\diffker^{2n}_{j_1},\diffker^{2n}_{j_3}} 
    \nn\\&\quad+O_M(n^{1/2-H}),
  \\&=
  2^2 \brbr{
  \: \cali^{(3.0)}_n
  - 3\times2^{(2H-1)} \times \cali^{(3.1)}_n
  + 3\times2^{(2H-1)2}\times \cali^{(3.2)}_n
  - 2^{(2H-1)3}\times\cali^{(3.3)}_n}
  +O_M(n^{1/2-H}),
\end{align*}
where we define 
$a(x) = (V^{[1]}(x))^2$
(hence $a_t = \brbr{V^{[1]}_{t} }^2 $)
and 
\begin{align*}
  \cali^{(3.0)}_n
  &=
  % n^{6H-1}
  %   \sum_{j\in[n-1]^3}%\bbJ_n([3],i)} 
  %   \rbr{\iv}^{-3} 
  %   \brbr{V^{[1]}_{t^{n}_{j_1}}}^2 \brbr{V^{[1]}_{t^{n}_{j_2}}}^2 \brbr{V^{[1]}_{t^{n}_{j_3}}}^2 
  %   \abr{\diffker^{n}_{j_1},\diffker^{n}_{j_2}}
  %   \abr{\diffker^{n}_{j_2},\diffker^{n}_{j_3}}
  %   \abr{\diffker^{n}_{j_1},\diffker^{n}_{j_3}} 
  % \\&=
  \rbr{\iv}^{-3}
  n^{-1} 
  \sum_{j\in[n-1]^3}
  a_{t^{n}_{j_1}}  a_{t^{n}_{j_2}} a_{t^{n}_{j_3}}
  \widehat\rho(\jon-\jtw)
  \widehat\rho(\jtw-\jth)
  \widehat\rho(\jth-\jon)
  \\
  \cali^{(3.1)}_n
  &=
  % n^{6H-1} 
  % \sum_{j\in[n-1]^2\times[2n-1]}
  % \rbr{\iv}^{-3}
  % \rbr{V^{[1]}_{t_\jon} }^2  \rbr{V^{[1]}_{t_\jtw} }^2 
  % \rbr{V^{[1]}_{t_{j_3}} }^2
  % \abr{\diffker^{n}_\jtw, \diffker^{2n}_\jth}  
  % \abr{\diffker^{n}_\jon, \diffker^{2n}_\jth}  
  % \abr{\diffker^{n}_\jon, \diffker^{n}_\jtw}
  % \\&=
  \rbr{\iv}^{-3}
  n^{-1} 
  \sum_{j\in[n-1]^2\times[2n-1]}
  a_{t^{n}_{j_1}}  a_{t^{n}_{j_2}} a_{t^{2n}_{j_3}}\;
  \widehat\rho(\jon-\jtw)
  \tilde\rho(2\jtw-\jth)
  \tilde\rho(2\jon-\jth)
  \\
  \cali^{(3.2)}_n
  &=
  % n^{6H-1} 
  % \sum_{j\in[n-1]\times[2n-1]^2}
  % \rbr{\iv}^{-3}
  % \rbr{V^{[1]}_{t_\jon} }^2  \rbr{V^{[1]}_{t_\jtw} }^2 
  % \rbr{V^{[1]}_{t_{j_3}} }^2
  % \abr{\diffker^{n}_\jon, \diffker^{2n}_\jtw}
  % \abr{\diffker^{2n}_\jtw, \diffker^{2n}_\jth}  
  % \abr{\diffker^{n}_\jon, \diffker^{2n}_\jth}  
  % \\&=
  2^{-2H}
  \rbr{\iv}^{-3}
  n^{-1} 
  \sum_{j\in[n-1]\times[2n-1]^2}
  a_{t^{n}_{j_1}}  a_{t^{2n}_{j_2}} a_{t^{2n}_{j_3}}\;
  \tilde\rho(2\jon-\jtw)
  \widehat\rho(\jtw-\jth)
  \tilde\rho(2\jon-\jth)
  \\
  \cali^{(3.3)}_n
  &=
  % n^{6H-1} 
  % \sum_{j\in[2n-1]^3}
  % \rbr{\iv}^{-3}
  % \rbr{V^{[1]}_{t_\jon} }^2  \rbr{V^{[1]}_{t_\jtw} }^2 
  % \rbr{V^{[1]}_{t_{j_3}} }^2
  % \abr{\diffker^{2n}_\jon, \diffker^{2n}_\jtw}
  % \abr{\diffker^{2n}_\jtw, \diffker^{2n}_\jth}  
  % \abr{\diffker^{2n}_\jon, \diffker^{2n}_\jth}  
  % \\&=
  2^{-6H}
  \rbr{\iv}^{-3}
  n^{-1} 
  \sum_{j\in[2n-1]^3}
  a_{t^{2n}_{j_1}}  a_{t^{2n}_{j_2}} a_{t^{2n}_{j_3}}\;
  \widehat\rho(\jon-\jtw)
  \widehat\rho(\jtw-\jth)
  \widehat\rho(\jth-\jon).
\end{align*}
The functional $\cali^{(3.1)}_n$ decomposes as 
\begin{align}
  \cali^{(3.1)}_n
  &=
  % \rbr{\iv}^{-3}
  % n^{-1} 
  % \sum_{j\in[n-1]^2\times[2n-1]}
  % a_{t^n_\jon}  a_{t^n_\jtw} a_{t^{2n}_{j_3}}
  % \tilde\rho(2\jtw-\jth)
  % \tilde\rho(2\jon-\jth)
  % \widehat\rho(\jon-\jtw)
  % \nn\\&=
  \rbr{\iv}^{-3}
  n^{-1} 
  \sum_{j\in[n-1]^2\times[2n-1]}
  a_{t^n_\jon}  \brbr{a_{t^n_\jtw}-a_{t^n_\jon}} a_{t^{2n}_{j_3}}\;
  \tilde\rho(2\jtw-\jth)
  \tilde\rho(2\jon-\jth)
  \widehat\rho(\jon-\jtw)
  \label{230512.1251}
  \\&\hspsm+
  \rbr{\iv}^{-3}
  n^{-1} 
  \sum_{j\in[n-1]^2\times[2n-1]}
  \brbr{a_{t^n_\jon}}^2 \brbr{a_{t^{2n}_{j_3}}-a_{t^n_\jon}}\; 
  \tilde\rho(2\jtw-\jth)
  \tilde\rho(2\jon-\jth)
  \widehat\rho(\jon-\jtw)
  \label{230512.1252}
  \\&\hspsm+
  \rbr{\iv}^{-3}
  n^{-1} 
  \sum_{j\in[n-1]^2\times[2n-1]}
  \brbr{a_{t^n_\jon}}^3 % a_{t^n_\jtw} a_{t^{2n}_{j_3}}
  \tilde\rho(2\jtw-\jth)
  \tilde\rho(2\jon-\jth)
  \widehat\rho(\jon-\jtw)
  \nn
\end{align}
The order of \eqref{230512.1251} is bounded by the following sum
with any $\beta\in(\half,H)$:
\begin{align*}
  %\rbr{\iv}^{-3}
  &n^{-1} 
  \sum_{j\in[n-1]^2\times[2n-1]}
  % \rbr{\frac{\abs{\jtw-\jon}}{n}}^\beta\;
  \abs{\jtw-\jon}^\beta n^{-\beta}\;
  \abs{\tilde\rho(2\jtw-\jth)}
  \abs{\tilde\rho(2\jon-\jth)}
  \abs{\widehat\rho(\jon-\jtw)}
  % \\&\leq
  % n^{-1-\beta} 
  % \sum_{j\in[n-1]^2\times[2n-1]}
  % \sup_{\jon,\jtw\in[n-1]}\rbr{\abs{\jtw-\jon}^\beta\;\widehat\rho(\jon-\jtw)}
  % \tilde\rho(2\jtw-\jth)
  % \tilde\rho(2\jon-\jth)
  \\&\leq
  n^{-1-\beta} 
  \sup_{\jon,\jtw\in[n-1]}\rbr{\abs{\jtw-\jon}^\beta\;\abs{\widehat\rho(\jon-\jtw)}}
  \sum_{j\in[n-1]^2\times[2n-1]}
  \abs{\tilde\rho(2\jtw-\jth)}
  \abs{\tilde\rho(2\jon-\jth)}
  \\&\leq
  n^{-1-\beta} 
  \sup_{\jon,\jtw\in[n-1]}\rbr{\abs{\jtw-\jon}^\beta\;\abs{\widehat\rho(\jon-\jtw)}}
  \\&\hspace{50pt}\times
  \sup_{\jth\in[2n-1]}\sum_{\jtw\in[n-1]}\abs{\tilde\rho(2\jtw-\jth)}
  \sup_{\jth\in[2n-1]}\sum_{\jon\in[n-1]}\abs{\tilde\rho(2\jon-\jth)}
  \sum_{\jth\in[2n-1]}1
\end{align*}
Here we used Lemma \ref{230924.1800} (ii). % \ref{230807.1543}.
Since
$\hat \rho(k), \tilde \rho(k) = O(\abs{k}^{2H-4})$ as $\abs{k}\to\infty$
by Lemma \ref{230807.1545},
we have 
\begin{align*}
  &\sup_{\jon,\jtw\in[n-1]}\rbr{\abs{\jtw-\jon}^\beta\;\abs{\widehat\rho(\jon-\jtw)}}
  <\infty
  \\
  &\sup_{\jth\in[2n-1]}\sum_{\jtw\in[n-1]}\abs{\tilde\rho(2\jtw-\jth)}
  % \rbr{=\sup_{\jth\in[2n-1]}\sum_{\jon\in[n-1]}\tilde\rho(2\jon-\jth)}
  <\sum_{j\in\bbZ} \tilde\rho(j)<\infty.
\end{align*}
Hence the term \eqref{230512.1251} is estimated as $O_M(n^{-\beta})$.
A similar calculation shows the term \eqref{230512.1252}
is of $O_M(n^{-\beta})$ as well.
By Lemma \ref{230109.1833} and \ref{230924.1800} (iii), % \ref{230807.1750}},
\begin{align*}
  \cali^{(3.1)}_n
  &=
  \rbr{\iv}^{-3} n^{-1} 
  \sum_{j\in[n-1]^2\times[2n-1]}
  \brbr{a_{t^n_\jon}}^3 % a_{t^n_\jtw} a_{t^{2n}_{j_3}}
  \tilde\rho(2\jtw-\jth)
  \tilde\rho(2\jon-\jth)
  \widehat\rho(\jon-\jtw)
  + O_M(n^{-\beta})
  \\&=
  \rbr{\iv}^{-3}
  n^{-1} \sum_{j_1\in[n-1]}
  \brbr{a_{t_{j_1}^n}}^3 
  \sum_{i_1,i_2\in\bbZ}
  \tilde\rho(i_2)
  \tilde\rho(2i_1-i_2)
  \widehat\rho(i_1)
  +o_M(1)
  \\&=
  \rbr{\iv}^{-3}
  \rbr{\sum_{i_1,i_2\in\bbZ}
  \tilde\rho(i_2)
  \tilde\rho(2i_1-i_2)
  \widehat\rho(i_1)}
  \int_0^1 (a_t)^3 dt
  +o_M(1)
\end{align*}

Similar arguments prove 
\begin{align*}
  \cali^{(3.0)}_n
  &=
  \rbr{\iv}^{-3}\;
  \rbr{\sum_{i_1,i_2\in\bbZ}
  \widehat\rho(i_1)
  \widehat\rho(i_1-i_2)
  \widehat\rho(i_2)}
  \int_0^1 (a_t)^3 dt
  +
  o_M(1).
  \\
  \cali^{(3.2)}_n
  &=
  2^{-2H}
  \rbr{\iv}^{-3}
  \rbr{\sum_{i_1,i_2\in\bbZ}
  \tilde\rho(i_1)
  \tilde\rho(i_2)
  \widehat\rho(i_1-i_2)}
  \int_0^1 (a_t)^3 dt
  +o_M(1)
  \\
  \cali^{(3.3)}_n&=
  2^{1-6H}
  \rbr{\iv}^{-3}
  \rbr{\sum_{i_1,i_2\in\bbZ}
  \widehat\rho(i_1)
  \widehat\rho(i_2)
  \widehat\rho(i_1-i_2)}
  \int_0^1 (a_t)^3 dt
  +o_M(1)
\end{align*}
and \eqref{230807.1953}.
\end{proof}

%% file: subfiles/7-sec_numerical_analysis.tex
\section{Numerical analysis}
\label{sec:240701.2110}
\newcommand{\NAimages}{images}
\newcommand{\imageFile}[3]{#1/HurstFSDE_centered_H#3_#1_n#2_sd911_MC1e+05MCAE10000.png}

Here we give results of numerical experiments.
We consider the following two SDEs for this experiment
\begin{align}
  dX_t &= X_t dt + (2+\sin(X_t)) dB_t
  \label{eq:231203.1401}
  \\
  dX_t &= \sin(X_t) dt + (2+\cos(X_t)) dB_t.
  \label{eq:231203.1402}
\end{align}
with $X_0=1$ for both.
These SDEs are considered as examples for numerical simulation in \cite{kubilius2012rate}.

To see the fit of the curves of approximate densities,
we generate the values of the estimator without the cutoffs at $0$ and $1$,
that is 
\begin{align*}
  \hat H_n' &=
  \half + \frac1{2\log2} \log\frac{\secVar_n}{\secVar_{2n}},
\end{align*}
and the number of paths to make the histograms is $10^5$ for all plots.
The dashed curves are obtained by the (mixed) normal approximation, and
the solid ones by the asymptotic expansion.
These curves are calculated based on $10^4$-times Monte Carlo simulations.

The following four plots are in the case of SDE \eqref{eq:231203.1401},
$H=0.55$ and $n=16,32,64,128$.
% \imageFile{d1s1}{32}
\begin{figure}[H]
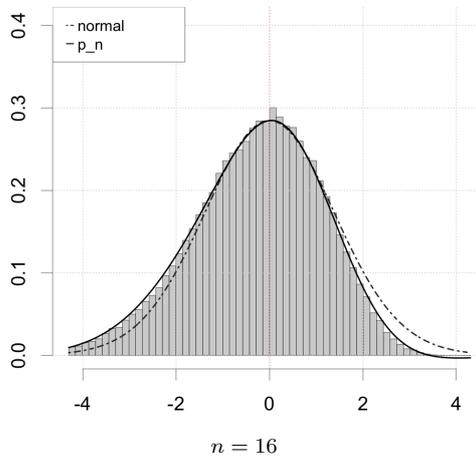
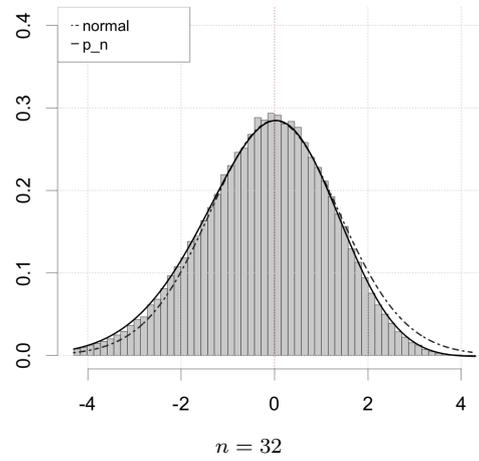

    \centering
  \begin{subfigure}[t]{0.45\textwidth}
    \centering
    \includegraphics[trim=110 140 120 190,clip,width=0.8\linewidth]{\NAimages/\imageFile{d1s1}{16}{55}}
    \caption*{$n=16$}
		% \caption*{$\sqrt{n}(\wh H_n'-H)$}
		% \includegraphics[clip,width=\linewidth]{\NAimages/\imageFile{d1s1}{16}{55}}
		% \includegraphics%[trim=110 140 120 190,clip,width=0.8\linewidth]
  \end{subfigure}
  \begin{subfigure}[t]{0.45\textwidth}
    \centering
		% \caption*{$\sqrt{n}(\wh H_n'-H)$}
		% \includegraphics[clip,width=\linewidth]{\NAimages/\imageFile{d1s1}{32}{55}}
    \includegraphics[trim=110 140 120 190,clip,width=0.8\linewidth]{\NAimages/\imageFile{d1s1}{32}{55}}
    \caption*{$n=32$}
  \end{subfigure}

  \vspsm
  % \begin{subfigure}[t]{0.45\textwidth}
  %   \centering
	% 	\caption*{$\sqrt{n}(\wh H_n^{(2)'}-H)$}
% 	\caption*{$H=0.55$, $n=16,32$}
% \end{figure}
% %
% \begin{figure}[H]
%   \centering
\begin{subfigure}[t]{0.45\textwidth}
  \centering
  \includegraphics[trim=110 140 120 190,clip,width=0.8\linewidth]{\NAimages/\imageFile{d1s1}{64}{55}}
  \caption*{$n=64$}
\end{subfigure}
\begin{subfigure}[t]{0.45\textwidth}
  \centering
  \includegraphics[trim=110 140 120 190,clip,width=0.8\linewidth]{\NAimages/\imageFile{d1s1}{128}{55}}
  \caption*{$n=128$}
\end{subfigure}
% \caption*{$H=0.55$, $n=16,32,64,128$}
\caption*{The plot of $\sqrt{n}(\wh H_n'-H)$ for SDE \eqref{eq:231203.1401} with $H=0.55$.}

\end{figure}

The first two plots are for 
SDE \eqref{eq:231203.1401}, $H=0.85$ and $n=16,32$;
the second for 
SDE \eqref{eq:231203.1402}, $H=0.55$ and $n=16,32$;
the last for 
SDE \eqref{eq:231203.1402}, $H=0.85$ and $n=16,32$;
\begin{figure}[H]
  \centering
\begin{subfigure}[t]{0.45\textwidth}
  \centering
  % \caption*{$\sqrt{n}(\wh H_n'-H)$}
  % \includegraphics[clip,width=\linewidth]{\NAimages/\imageFile{d1s1}{16}{85}}
  \includegraphics[trim=110 140 120 190,clip,width=0.8\linewidth]{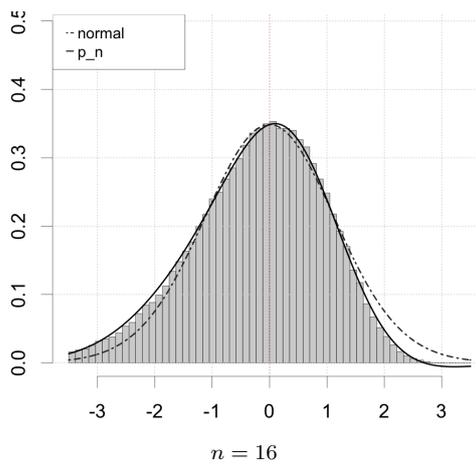}
  \caption*{$n=16$}
\end{subfigure}
\begin{subfigure}[t]{0.45\textwidth}
  \centering
  % \caption*{$\sqrt{n}(\wh H_n'-H)$}
  % \includegraphics[clip,width=\linewidth]{\NAimages/\imageFile{d1s1}{32}{85}}
  \includegraphics[trim=110 140 120 190,clip,width=0.8\linewidth]{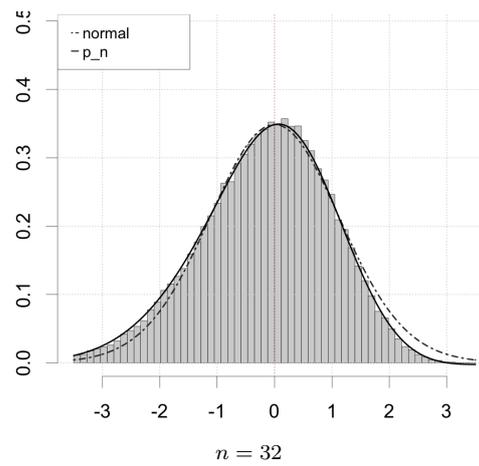}
  \caption*{$n=32$}
\end{subfigure}
\caption*{The plot of $\sqrt{n}(\wh H_n'-H)$ for SDE \eqref{eq:231203.1401} with $H=0.85$.}
% \caption*{$H=0.85$, $n=16,32$}
\end{figure}

\begin{figure}[H]
  \centering
\begin{subfigure}[t]{0.45\textwidth}
  \centering
  % \caption*{$\sqrt{n}(\wh H_n'-H)$}
  % \includegraphics[clip,width=\linewidth]{\NAimages/\imageFile{d2s2}{16}{55}}
  \includegraphics[trim=110 140 120 190,clip,width=0.8\linewidth]{\NAimages/\imageFile{d2s2}{16}{55}}
  \caption*{$n=16$}
\end{subfigure}
\begin{subfigure}[t]{0.45\textwidth}
  \centering
  % \caption*{$\sqrt{n}(\wh H_n'-H)$}
  % \includegraphics[clip,width=\linewidth]{\NAimages/\imageFile{d2s2}{32}{55}}
  \includegraphics[trim=110 140 120 190,clip,width=0.8\linewidth]{\NAimages/\imageFile{d2s2}{32}{55}}
  \caption*{$n=32$}
\end{subfigure}
\caption*{The plot of $\sqrt{n}(\wh H_n'-H)$ for SDE \eqref{eq:231203.1402} with $H=0.55$.}
% \caption*{$H=0.55$, $n=16,32$}
\end{figure}

\begin{figure}[H]
  \centering
\begin{subfigure}[t]{0.45\textwidth}
  \centering
  % \caption*{$\sqrt{n}(\wh H_n'-H)$}
  % \includegraphics[clip,width=\linewidth]{\NAimages/\imageFile{d2s2}{16}{95}}
  \includegraphics[trim=110 140 120 190,clip,width=0.8\linewidth]{\NAimages/\imageFile{d2s2}{16}{85}}
  \caption*{$n=16$}
\end{subfigure}
\begin{subfigure}[t]{0.45\textwidth}
  \centering
  % \caption*{$\sqrt{n}(\wh H_n'-H)$}
  % \includegraphics[clip,width=\linewidth]{\NAimages/\imageFile{d2s2}{32}{95}}
  \includegraphics[trim=110 140 120 190,clip,width=0.8\linewidth]{\NAimages/\imageFile{d2s2}{32}{85}}
  \caption*{$n=32$}
\end{subfigure}
\caption*{The plot of $\sqrt{n}(\wh H_n'-H)$ for SDE \eqref{eq:231203.1402} with $H=0.85$.}
	% \caption*{$H=0.95$, $n=16,32$}
\end{figure}

% 比較の方向性
% - d1s1,d2s2
% - nについて
% - Hについて